\documentclass[aos]{imsart}

\RequirePackage{amsthm,amsmath,amsfonts,amssymb}
\RequirePackage{natbib}
\RequirePackage{graphicx,color}
\RequirePackage{url}
\RequirePackage{xcolor}
\RequirePackage[colorlinks,citecolor=blue!70!black,linkcolor=blue,urlcolor=blue]{hyperref}
\RequirePackage{subcaption}   
\RequirePackage{enumitem}
\RequirePackage{ulem}
\RequirePackage{comment}
\RequirePackage{pgfplots}
\usepackage[ruled,linesnumbered]{algorithm2e} 
\RequirePackage{tikz}
\usetikzlibrary{shapes.geometric}

\usepackage{url}
\usepackage{titlesec}
\usepackage{easybmat}
\usepackage{float}
\usepackage{array}
\newcolumntype{M}[1]{>{\centering\arraybackslash}p{#1}}
\usepackage{multirow,bigdelim}
\newcommand*\hexbrace[2]{%
  \underset{#2}{\underbrace{\rule{#1}{0pt}}}}

\DeclareFontEncoding{LS1}{}{}
\DeclareFontSubstitution{LS1}{stix}{m}{n}
\DeclareSymbolFont{symbols2}{LS1}{stixfrak} {m} {n}
\DeclareMathSymbol{\operp}{\mathbin}{symbols2}{"A8}

\newcommand{\blue}[1]{\textcolor{blue}{#1}}

\newcommand{\RNum}[1]{\uppercase\expandafter{\romannumeral #1\relax}}

\startlocaldefs
\newtheorem{theorem}{Theorem}

\newtheorem{condition}{Condition}

\newtheorem{definition}{Definition}

\newtheorem{lemma}[theorem]{Lemma}
\newtheorem{notation}[theorem]{Notation}

\newtheorem{proposition}[theorem]{Proposition}
\newtheorem{remark}{Remark}

\newtheorem{fact}{Fact}[section]

\newcommand{\EXPT}{\mathbb{E}}

\newcommand{\mytrans}{\top}
\newcommand{\myCt}{\mathcal{C}_t}

\newcommand{\myBt}{\mathcal{B}_t}
\newcommand{\myAt}{\mathcal{A}_t}
\newcommand{\myUt}{\mathcal{V }_t}

\newcommand{\myDone}{\mathcal{D}_1}
\newcommand{\myXt}{\mathcal{X}_t}

\newcommand{\mM}{\mathcal{M}_k}
\newcommand{\bmMp}{\bar{\mathcal{M}}^{\perp}_{k}}
\newcommand{\DZ}{\mathcal{D}_Z}
\newcommand{\UZ}{\mathcal{U}_Z}

\newcommand{\UZb}{\mathcal{U}_{\check Z}}
\newcommand{\ZQs}{Z^{\star}\check Q}
\newcommand{\myDtri}{\mathcal{D}_3}

\newcommand{\initerrprob}{1 - O((nT)^{-s}) }
\newcommand{\Seff}{S_{eff} }
\newcommand{\Heff}{H_{eff} }
\newcommand{\Ieff}{I_{eff} }
\newcommand{\rankU}{\mathrm{r}_{n,k} }
\newcommand{\Ginit}{\hat G_{0}}

\DeclareMathOperator*{\argmin}{arg\,min}
\DeclareMathOperator*{\argmax}{arg\,max}

\pgfplotsset{compat=1.18}

\endlocaldefs

\begin{document}

\begin{frontmatter}
\title{Semiparametric Modeling and Analysis for Longitudinal Network Data}

\runtitle{Semiparametric Analysis for Longitudinal Network Data}

\begin{aug}
\author[A]{\fnms{Yinqiu}~\snm{He}\ead[label=e1]{yinqiu.he@wisc.edu}},
\author[B]{\fnms{Jiajin}~\snm{Sun}\ead[label=e2]{jsun5@fsu.edu}},
\author[C]{\fnms{Yuang}~\snm{Tian}\ead[label=e3]{yatian20@fudan.edu.cn}},
\author[D]{\fnms{Zhiliang}~\snm{Ying}\ead[label=e4]{zying@stat.columbia.edu}},
\and
\author[E]{\fnms{Yang}~\snm{Feng}\ead[label=e5]{yang.feng@nyu.edu}}
\thanks{He, Sun, and Tian contribute equally to this work. Corresponding Author: Yang Feng.}
\address[A]{Department of Statistics,
University of Wisconsin-Madison\printead[presep={,\ }]{e1}}

\address[B]{Department of Statistics, 
Florida State University\printead[presep={,\ }]{e2}}

\address[C]{Shanghai Center for Mathematical Sciences, Fudan University\printead[presep={,\ }]{e3}}

\address[D]{Department of Statistics, 
Columbia University\printead[presep={,\ }]{e4}}

\address[E]{Department of Biostatistics, School of Global Public Health, New York University\printead[presep={,\ }]{e5}}

\end{aug}

\begin{abstract}
We introduce a semiparametric latent space model for analyzing longitudinal network data. The model consists of a static latent space component and a time-varying node-specific baseline component. We develop a semiparametric efficient score equation for the latent space parameter by adjusting for the baseline nuisance component. 
Estimation is accomplished through a one-step update estimator and an appropriately penalized maximum likelihood estimator. We derive oracle error bounds for the two estimators and address identifiability concerns from a quotient manifold perspective. Our approach is demonstrated using the New York Citi Bike Dataset.
\end{abstract}




\end{frontmatter}

\section{Introduction}

Recent years have seen increased availability of time-varying interaction/network data \citep{butts20084,linderman2014discovering,holme2015modern}, with
examples such as email exchange history of coworkers \citep{klimt2004enron} and transport records among bike-sharing stations \citep{bikedata}. This type of data not only includes counts of pairwise interaction events but also the timestamps of these interactions. 

Non-network count data have been studied extensively in the survival analysis literature. \cite{andersen1982cox} proposed a Poisson-type intensity-based regression model. For extensions to handle non-Poisson count data, see \cite{pepe1993some}, \cite{cook2007statistical}, \cite{lin2000semiparametric}, \cite{sun2013statistical}.

For cross-sectional single network data, an important development is the latent space model, in which each node is represented by a latent vector and the relationship between two nodes is quantified through their inner product \citep{hoff2002latent}.
The latent vectors may be treated as random \citep{hoff2003random,hoff2005bilinear,bickel2013asymptotic} or fixed \citep{athreya2018survey,ma2020universal}. 

Motivated by various scientific applications, multiple-network models and longitudinal network models have also been developed; see, for example, \cite{zhang2020mixed} for application in neuroscience and \cite{butts20084} for application in social science. Similar to the single-network models, latent space modeling in the context of multiple-network settings has been examined from different perspectives \citep{sewell2015latent,nielsen2018multiple,matias2014modeling}.
One line of research views the latent vectors as random \citep{hoff2011hierarchical,salter2017latent}, and another line considers fixed latent vectors \citep{jones2020multilayer,levin2017central}. The latter often leads to frequentist approaches, in which it is of common interest to study 
the estimation of shared structure across multiple networks, and in particular, how the information accumulation across multiple networks improves the fitting quality  \citep{arroyo2021inference,zheng2022limit,macdonald2020latent,zhang2020flexible}. 
\nocite{yang2011detecting,xu2014dynamic,matias2017statistical,matias2018semiparametric,zhang2020modularity,chen2022global}


In this paper, we develop a semiparametric modeling framework for longitudinal networks with interaction event counts. 
The model comprises a static latent space component, which accounts for inherent node interactions, and a baseline component that accommodates heterogeneity across different time points and nodes. 
We introduce two semiparametric estimation methods. The first method forms a generalized semiparametric one-step estimator by solving a local linear approximation of the efficient score equation.
We demonstrate that the estimator automatically eliminates the identifiability issue through an interpretation in the quotient manifold induced by orthogonal transformation equivalence classes.
The second method optimizes a convex surrogate loss function based on the log-likelihood, relaxing the non-convex rank constraint of the latent space component. We establish that both methods achieve, up to a logarithmic factor, the {oracle} estimation error rates for the latent space component. 

The rest of the paper is organized as follows. 
Section \ref{sec::2model} introduces the semiparametric latent space model. Section \ref{sec:newton} introduces the one-step estimator, provides its connection to the quotient manifold theory, and establishes the error bound for the estimator. Section \ref{sec:penalizedmethod} introduces the penalized maximum likelihood estimator and establishes its error bound. Simulation results are presented in Section \ref{sec::simulations}. Section \ref{sec::bikedata} applies our framework to analyze the New York Citi Bike Dataset \citep{bikedata}. Section \ref{sec::discussion} concludes the main body of the paper with some discussions. All the technical derivations are deferred to the Supplementary Material.

\section{Semiparametric Poisson Latent Space Model}
\label{sec::2model}
\subsection{Notation and Model Specification}
\label{subsec::1.2model}
We consider longitudinal pairwise interaction counts of $n$ subjects (nodes) over $T$ discrete time points.  
Specifically, for a time point $t\in \{1,\ldots, T\}$ and nodes $i, j\in \{1,\ldots, n\}$, $A_{t,ij}$ denotes the number of $i$-$j$ interactions  at the  time point $t$.
We propose a Poisson-based latent space model  
\begin{align}
A_{t,ij} =&~ A_{t,ji}\ \ {\sim}\ \ \text{Poisson}\{\EXPT(A_{t,ij}\mid z,\alpha) \}, \quad \text{independently with }\label{model_DLSM} \\
\EXPT(A_{t,ij}\mid z,\alpha)=&~\exp(\alpha_{it} + \alpha_{jt} + 
\langle z_{i}, z_j\rangle ), \notag 
\end{align} 
which naturally adopts the {exponential} link function $\exp(\cdot)$ to model the event counts. 
For any two nodes $i$ and $j$, their interaction effect is modeled through the inner product of two corresponding latent vectors $\langle z_i, z_j \rangle = z_i^{\top}z_j$, 
similarly to the inner product model of a single network \citep{ma2020universal}. 
The latent vectors $z_i$'s do not change with respect to the time point $t$ and represent the shared latent structures across $T$ heterogeneous networks. 
For example, $z_i$'s can encode the time-invariant geographic information in multiple transportation networks. At a given time point $t$, when $\alpha_{it}$ increases and all the other parameters are fixed, edges connecting the node $i$ tend to have higher numbers of counts at the time point $t$, indicating higher baseline activity levels. 
Therefore, $\alpha_{it}$'s model the degree heterogeneity across different nodes $i\in \{1,\ldots, n\}$ and time points $t\in \{1,\ldots, T\}$, and are called baseline degree heterogeneity parameters of nodes and time.  
In the hourly bike-sharing networks, $\alpha_{it}$'s can represent distinct baseline activity levels across different stations and hours. 

Model specification \eqref{model_DLSM} may be expressed in vector-matrix notation as 
\begin{align} \label{eq:modelmatrixform}
    \EXPT(\mathbf A_t \mid  Z, \alpha) = \exp(\alpha_t {1}_n^{\top} + {1}_n \alpha_t ^{\top} + Z Z^{\top}), 
\end{align}
where $\alpha_t=(\alpha_{1t},\ldots, \alpha_{nt})^{\top} \in \mathbb{R}^{n\times 1}$, $Z=(z_1,\ldots, z_n)^{\top} \in \mathbb{R}^{n\times k}$, $1_n = (1,\ldots, 1)^{\top} \in \mathbb{R}^{n\times 1}$, and $\exp(\cdot)$ is the elementwise exponential operation.  
Throughout this paper, we consider the asymptotic regime in which the number of nodes $n$ and the number of time periods $T$ increase to infinity while the dimension of the latent space $k$ is fixed. Thus, $\alpha_t {1}_n^{\top} + {1}_n \alpha_t ^{\top} + Z Z^{\top}$ is a low-rank matrix.
To ensure identifiability,
we assume that column means of $Z$ are zero, i.e., $
1_n^{\top} Z/n = 0$.  
This centering assumption is analogous to the classical two-way analysis of variance (ANOVA) modeling with interaction \citep{scheffe1999analysis}. 
Additionally, since $ZZ^{\top}=ZQQ^{\top}Z^{\top}$ for any 
$Q \in \mathcal{O}(k)$, 
where  $\mathcal{O}(k) = \{ Q\in \mathbb{R}^{k\times k} : Q Q^{\top} =\mathrm{I}_k \}$,
$Z$ is identifiable up to a common orthogonal transformation of its rows. 

Let $\lambda_{ij,t}=\exp(\alpha_{it}+\alpha_{jt})$. Then
$\EXPT(A_{t,ij}) = \lambda_{ij,t}\times \exp(\langle z_i, z_j \rangle )$.
This form resembles the multiplicative counting process modeling in the analysis of recurrent event times \citep{cook2007statistical,sun2013statistical}. In particular,  the $\lambda_{ij,t}$ corresponds to the baseline rate/intensity function changing over time while the $z_i$ corresponds to the time-invariant parameters of interest. Note that the number of baseline (nuisance) parameters is $nT$, 
and the size of the shared latent positions $Z$ or number of parameters of interest is $nk$; 
 the former increases with both $n$ and $T$, whereas the latter does not increase with $T$. 
Thus, we may be able to estimate $Z$ more accurately than the nuisance parameter. 

Some notations used in this paper are summarized as follows.
For a matrix $X=[x_{ij}]\in \mathbb{R}^{n\times n}$,  $\mathrm{tr}(X)=\sum_{i=1}^n x_{ii}$ stands for its trace and $\sigma_{\min}(X)$ represents its minimum eigenvalue. 
For $X, Y\in \mathbb{R}^{n\times m}$,  $\langle X, Y\rangle = \mathrm{tr}(X^{\top}  Y) $.
If $m\geqslant n$, for any matrix $X$ with the singular value decomposition $X=\sum_{i=1}^n s_i u_i v_i^{\top}$, we let $\|X\|_{*}= \sum_{i=1}^n s_i$,  $\|X\|_{\mathrm{F}}= \sqrt{\sum_{i=1}^n s_i^2}$, and $\|X\|_{\mathrm{op}}= \max_{i=1, \ldots, n} s_i$ stand for the nuclear norm, the Frobenius norm, and the operator norm of the matrix, respectively. For a vector $x\in \mathbb{R}^{n}$, $\|x\|_2=\sqrt{  x^{\top}x }$. 

For two sequences of real numbers $\{f_n\}$ and $\{h_n\}$, $f_n\lesssim h_n $ and $f_n = O(h_n) $ mean that $|f_n|\leqslant c_1|h_n|$ for a constant $c_1>0$; $f_n\asymp h_n$ means $c_2h_n\leqslant  f_n\leqslant c_1  h_n $ for some constants $c_1,c_2>0$; $f_n=o(h_n)$ and $f_n\ll h_n$ mean $\lim_{n\to \infty} f_n/h_n=0$; $f_n\gg h_n$ means $\lim_{n\to \infty}h_n/f_n=0$. 
For a sequence of random variables $X_n$ and a sequence of real numbers $f_n$, write $X_n=O_p(f_n)$ if for any $\epsilon>0$, there exists finite $M>0$ such that $\sup_n \Pr(|X_n/f_n|>M)<\epsilon$ ($X_n/f_n$ is stochastically bounded); write $X_n=o_p(f_n)$ if for any $\epsilon>0$, $\lim_{n\to \infty}\Pr(|X_n/f_n|>\epsilon) =0$ ($X_n/f_n$  converges to 0 in probability). 

\subsection{Semiparametric Oracle Error Rate and Technical Challenges}\label{sec:oraclerateandchallenges}
To gain insights into the best possible error rate for estimating high-dimensional parameters,  
let's consider a simpler setting of the regular exponential family with a $p_n$-dimensional natural parameter vector $\theta_n$. Suppose we have $n$ independent observations from this family. 
 \cite{portnoy1988asymptotic} showed that the  MLE  $\hat{\theta}_n$ satisfies $\|\hat{\theta}_n-\theta_n\|_2^2=O_p( {p_n / n} )$, where $p_n$ can increase with the sample size $n$. For our problem of estimating $Z$, if the nuisance parameters $\alpha$'s were known,
then each $k$-dimensional latent vector $z_i$ would be measured by $nT$ independent edges: $\{A_{t,ij}: j =1,\ldots, n; t=1,\ldots, T \}$. 
The result of \cite{portnoy1988asymptotic} indicates that the oracle estimation error rate of each $z_i$ would be $ {k}/(nT)$.
When $k$ is fixed, the aggregated estimation error of $n$ latent vectors $Z=[z_1,\ldots, z_n]^{\top}$ is expected to be  of the order of
$O_p(n \times {k}/{(nT)}) = O_p({1}/{T})$,  
which is referred to as the oracle estimation error rate throughout this paper. Since the $\alpha$'s are unknown, the asymptotic theory for the classical semiparametric models \citep{bickel1993efficient} may be modified to yield efficient score equations (projections) so that the same error rate can be achieved. 

To achieve the above oracle bounds, 
several technical challenges remain. Firstly, the baseline parameters $\alpha_{it}$ in  \eqref{model_DLSM} not only characterize the time-specific heterogeneity over $t$ but also represent node-specific heterogeneity over $i$. Due to this two-way heterogeneity, the estimation errors of $\alpha_{it}$'s are intertwined with those of node-specific parameters $z_i$'s in a complicated way. In particular, the partial likelihood \citep{andersen1982cox} cannot eliminate the baseline nuisance parameters.
Secondly, the target parameter matrix $Z$ is only identifiable up to an orthogonal transformation. 
Therefore, given a true matrix $Z^{\star}$ and an estimate $\hat{Z}$, the ordinary Euclidean distance between the two matrices is not a proper measure for the estimation error of $\hat{Z}$. 
Indeed, it is more natural to consider their distance defined up to an orthogonal group transformation $\operatorname{dist}( \hat{Z}, Z^{\star}) = \min_{Q \in \mathcal{O}(k)} \| \hat{Z}-Z^{\star} Q\|_{\mathrm{F}}$. This distance metric indicates that the intrinsic geometry of $Z$ is non-Euclidean. Thirdly, the estimation of $Z$ is further complicated by the high-dimensionality of the parameters and non-linearity of the transformation in \eqref{model_DLSM}. Unlike classical semiparametric problems, the dimension of the parameters of interest also grows with the sample size.
Furthermore, the non-linearity of the exponential link function in \eqref{model_DLSM} 
 makes it difficult for the spectral-based analyses \citep{arroyo2021inference} to deal with the entanglement between $Z$ and $\alpha$ and achieve the semiparametric oracle rates of $Z$. 
 Lastly, the different oracle rates for different parameters also render techniques under the single network setting not easily applicable: in the single network setting $z_i$ and $\alpha_{i}$ have the same oracle rates, so it suffices to jointly control their overall error rates via the natural parameter $\Theta_{ij} = \alpha_i + \alpha_j + \langle z_i, z_j\rangle $ from a generalized linear model perspective. In contrast, in our setting, 
 each $\alpha_{it}$ is measured by $n$ independent edges: $\{A_{t,ij}: j =1,\ldots, n \}$, and the oracle squared error rate of each $\alpha_{it}$ would be $O_p(1/n)$, compared to the $O_p(1/(nT))$ error rate of each $z_i$.
Our goal is to attain the oracle estimation error rate of $Z$, which is a refinement of the overall error rate dominated by $\alpha$'s error. 

To accurately estimate the latent space component $Z$, we next develop two methods: a generalized semiparametric one-step estimator in Section \ref{sec:newton}, and a semiparametric penalized maximum likelihood estimator in Section \ref{sec:penalizedmethod}. We show that both methods achieve, up to a logarithmic factor, the {oracle} estimation error rates for the target parameters $Z$.

\section{Generalized Semiparametric One-Step Estimator} \label{sec:newton}
In this section, we introduce our generalized semiparametric one-step estimator of $Z$ and provide theoretical guarantees. We first introduce some notation. Model \eqref{model_DLSM} leads to the following form for the log-likelihood function
\begin{align}
\label{eq:LZalpha}
L(Z,\alpha) =L(Z_v,\alpha_v)=\sum_{t=1}^T \sum_{1\leqslant i\leqslant j\leqslant n } \left\{A_{t,ij} (\alpha_{it}+\alpha_{jt}+\langle z_i, z_j\rangle ) - \exp(\alpha_{it}+\alpha_{jt}+\langle z_i, z_j\rangle) \right\} 
\end{align}
where, for notational convenience in the differentiation of the likelihood, we use $Z_v$ and $\alpha_v$ to denote vectorizations of $Z$ and $\alpha$, respectively, i.e.,
$Z_v=(z_1^{\top},\ldots, z_n^{\top})^{\top} \in \mathbb{R}^{nk\times 1}$
and ${{\alpha}_v=(\alpha_1^{\top}, \ldots, \alpha_T^{\top})^{\top}} \in \mathbb{R}^{nT\times 1}$.  
Then we let $\dot{L}_Z(Z,\alpha)$ and $\dot{L}_{\alpha}(Z,\alpha)$ 
denote the partial derivatives  of $L(Z,\alpha)$ 
with respect to vectors $Z_v$ and $\alpha_v$, respectively (see the precise formulae in Section \ref{subsec:formula_efficient} of the Supplementary Material). 
Unless otherwise specified, 
such vectorization is applied when considering partial derivatives (for both first and higher orders)  throughout the paper. 
Following the semiparametric literature \citep{tsiatis2006semiparametric},  
the efficient score and the efficient Fisher information matrix for  $Z$ can be expressed as
\begin{align}
       \Seff(Z,\alpha) =&~\dot{L}_Z - \mathbb{E}\big(\dot{L}_{Z}\, \dot{L}_{\alpha}^{\top}\big) \big\{\mathbb{E}\big(\dot{L}_{\alpha}\, \dot{L}_{\alpha}^{\top}\big)\big\}^{-1} \dot{L}_{\alpha} \ \in \  \mathbb{R}^{nk\times 1},\notag\\
    \Ieff(Z, \alpha)=&~\mathbb{E}\big(\Seff(Z,\alpha)\, \Seff^{\top}(Z,\alpha) \big) \ \in \ \mathbb{R}^{nk\times nk},\label{eq:ieffdef}
\end{align}
where, for notational simplicity, $(Z,\alpha)$ is omitted in $\dot{L}_Z$ and $\dot{L}_{\alpha}$, and  $\mathbb{E}(\cdot)$ refers to the expectation  taken when data follows  \eqref{eq:modelmatrixform} with parameters $(Z,\alpha)$. We have derived analytic formulas for the expectation terms in \eqref{eq:ieffdef}, and express $ \Seff(Z,\alpha)$ and $\Ieff(Z,\alpha)$ as close-form functions of $(Z,\alpha)$ and the observed data; see Section \ref{subsec:formula_efficient} of the Supplementary Material.

With the above preparations, we construct our generalized semiparametric one-step estimator as
\begin{align} \label{eq:newtonsolpseudo}
 \hat{Z}_v = \check{Z}_v  + \big\{ \Ieff( \check{Z}, \check{\alpha})\big\}^{+}\Seff( \check{Z}, \check{\alpha}),  
\end{align} 
where $(\check Z,\check \alpha) $ denotes an initial estimate, and
$B^+ $ represents the  Moore–Penrose inverse of a matrix $B$, which is uniquely defined and also named pseudo inverse \citep{ben2003generalized}. By the property of pseudo inverse, 
\eqref{eq:newtonsolpseudo} can be equivalently written as 
\begin{equation}
\label{newtonsolUZ} 
 \hat Z_v = \check Z_v + \check{\mathcal U}\, \left\{\check{\mathcal U}^{\top}\,    \Ieff(\check Z, \check\alpha)\, \, \check{\mathcal U} \right\}^{-1} \,  \check{\mathcal U}^\top \, \Seff(\check Z,\check\alpha),
\end{equation}
where $\check{\mathcal U}$ is a matrix whose columns can be any set of basis of the column space of $\Ieff(\check Z, \check\alpha)$. 

The estimator originates from the one-step estimator for semiparametric problems \citep[Section 25.8]{van2000asymptotic}, which solves the efficient score equation in the presence of nuisance parameters through a linear approximation at an initial estimate.
However, different from the classical one-step estimator, the efficient information matrix is singular in our setting; see Remark \ref{lb:singularity}. As a result, the corresponding linear approximation equations are underdetermined. Taking pseudo inverse corresponds to solving the linear approximation equation with minimum $\ell_2$ norm; the solution lies in the column space of $\check {\mathcal U}$, as indicated in \eqref{newtonsolUZ}.

Owing to the intricate technical challenges outlined in Section \ref{sec:oraclerateandchallenges} and the singularity of $\Ieff(\check{Z},\check{\alpha})$, straightforward conclusions regarding the classical one-step estimator for \eqref{eq:newtonsolpseudo} are elusive. However, we have constructed a detailed semiparametric analysis and demonstrated that our proposed one-step estimator nearly achieves the oracle error rate, given suitable initial estimators. We will present the comprehensive theoretical results in the following Section \ref{sec:formulaofupdateandtheory}. Additionally, in Section \ref{sec:initest}, we will introduce an initial estimator that demonstrates both statistical validity and computational efficiency.

\begin{remark}\label{lb:singularity}
In effect, the singularity of the efficient information matrix  $\Ieff(Z,\alpha)$ is caused by the redundancy and unidentifiability of parameters in  $Z$  \citep{little2010parameter}. 
Particularly, we find that the rank of $\Ieff(Z,\alpha)$ equals  $nk-k(k+1)/2$, which is identical to the number of free parameters in the $n\times k$ matrix $Z$. Below, we provide a high-level  explanation, and please find a formal  justification in 
Lemma \ref{lem_presentingUZmain} in Section \ref{sec:eigenspace}. 
Firstly, recall that 
to avoid the mean shift issue in \eqref{eq:modelmatrixform}, 
we have imposed $k$ linear constraints  $1_n^{\top}Z = 0$.  
Secondly, even with the centering constraints, 
$Z$ can only be identified up to an orthogonal group transformation $\mathcal{O}(k)$.  
In particular, $\mathcal{O}(k)$ is the orthogonal Stiefel manifold, and its dimension is $k(k-1)/2$ \citep[see, e.g.,][Section 3.3.2]{absil2009optimization}.
Intuitively,  $k+k(k-1)/2 = k(k+1)/2$ degrees of freedom are removed due to the constraints and unidentifiability. This leads to  $nk-k(k+1)/2$ free parameters in $Z$.  
\end{remark} 

\subsection{Theory} \label{sec:formulaofupdateandtheory}
Throughout the sequel, we use $(Z^{\star}, \alpha^{\star})$ to denote the true value of $(Z,\alpha)$. In other words,  our observed data follow the model \eqref{model_DLSM}
 with $(Z,\alpha) = (Z^{\star}, \alpha^{\star})$.
Besides, we denote $\Theta_{t,ij} = \alpha_{it}+\alpha_{jt} + \langle z_i, z_j\rangle $. We define the estimation error from the $i$-th row of $\check{Z}$ 
 as $\mathrm{dist}_i(\check{z}_i, z_i^{\star}) = \|\check z_{i} - \check Q^\top z^\star_i \|_2$,
where 
\begin{align}\label{eq:defcheckq}
\check Q = \argmin_{Q\in \mathcal{O}(k)} \|\check Z - Z^\star  Q\|_{\mathrm{F}},
\end{align}
so that $\mathrm{dist}^2(\check{Z}, Z^{\star})=\sum_{i=1}^n \mathrm{dist}_i^2(\check{z}_i, z_i^{\star})$. 
The goal in this subsection is to establish an error bound for the proposed one-step estimator \eqref{eq:newtonsolpseudo} in terms of the distance defined above. To do so, we need the following regularity conditions.

\begin{condition}\label{cond:truvalueregularity}
Assume the true parameters $\alpha^{\star}$'s and $z^{\star}$'s satisfy: 
\begin{itemize}\setlength{\itemsep}{0pt}
    \item[(i)] There exist positive constants $M_{Z,1}$, $M_{\alpha}$, and  $M_{\Theta, 1}$ such that $\left\|z_{i}^{\star}\right\|^2_2\leqslant  M_{Z,1}$, 
 $|\alpha_{it}^{\star}| \leqslant M_{\alpha}$,  and $ \Theta_{t,ij}^{\star} \leqslant -{M_{\Theta, 1}}$ for  $1\leqslant i,j\leqslant n$ and $1\leqslant t \leqslant T$;

\item[(ii)] There exists a positive constant $M_{Z,2}$ such that $\sigma_{\min}[(Z^{\star})^{\top}Z^\star/n]
\geqslant M_{Z,2}$;
 
\item[(iii)] $1_n^{\top}Z^{\star}=0$. 
\end{itemize}
\end{condition}

\begin{condition}
\label{cond_elem_init}
Assume the initial estimates $\check{\alpha}$'s and  $\check{z}$'s satisfy : 
\begin{itemize}\setlength{\itemsep}{0pt}
    \item[(i)] 
There exist constants $M_b >0$ and $\epsilon$ such that
 $|\check\alpha_{it} - \alpha_{it}^\star| + \operatorname{dist}_i(\check{z}_i, z_i^{\star})\leqslant M_bn^{-{1}/{2}}\log^{\epsilon}(nT)$ for  $ 1\leqslant i \leqslant n$ and $1\leqslant t \leqslant T$; 

\item[(ii)] $1_n^{\top} \check{Z}= 0$. 
\end{itemize}
\end{condition}

Condition \ref{cond:truvalueregularity}(i) assumes boundedness of the true parameters.
It also implies $- M_{\Theta, 2}\leqslant  \Theta_{t,ij}^{\star}\leqslant - M_{\Theta, 1}$ with $M_{\Theta, 2}=M_{Z,1}+2M_{\alpha}$. 
Condition \ref{cond:truvalueregularity}(ii) excludes degenerate cases of $Z^{\star}$.  
This, together with $\|z_i^{\star}\|_2^2\leqslant M_{Z,1}$,  indicates  $\|Z^{\star}\|_{\mathrm{op}}\asymp \sqrt{n}$. 
Condition \ref{cond:truvalueregularity}(iii)  
is needed in order to identify the interaction effects  $Z^{\star}( Z^{\star})^{\top}$ in \eqref{eq:modelmatrixform}.  Condition \ref{cond_elem_init}(ii) requires $\check{Z}$ to be a feasible estimator satisfying the centering constraint. With this condition, 
$\hat{Z}$ is guaranteed to be feasible, namely  $1_n^{\top} \hat{Z} = 0$.  
Condition \ref{cond_elem_init}(i) is 
 analogous to the ``$\sqrt n$-consistent initial estimator''  used in the classical one-step estimator  \citep{van2000asymptotic}. 
It assumes that the estimation errors of $\check{\alpha}_{it}$'s and $\check{z}_i$'s are uniformly bounded by $M_b \log^\epsilon(nT) /\sqrt{n}$. In Section \ref{sec:initest}, we will construct an initial estimator satisfying Condition \ref{cond_elem_init} with the error bound of order $ \log^2(nT)/\sqrt{n}$
under high probability. Recall that the classical MLE theory indicates that the optimal estimation errors for individual $\alpha_{it}$ and $z_i$ would be of the orders of $n^{-1/2}$ and $(nT)^{-1/2}$, respectively. 
In this regard, the constructed initial $\check{\alpha}_{it}$ achieves its optimal rate up to a logarithmic factor.
On the other hand, 
Condition \ref{cond_elem_init}(i)
assumes that $\check{z}_i$ achieves the same error rate as $\check{\alpha}_{it}$ but not its own optimal rate $(nT)^{-1/2}$. 
In this connection, 
Condition \ref{cond_elem_init} imposes a mild assumption on the initial $\check{z}_i$'s. 
We also would like to point out that the deterministic upper bounds in Condition \ref{cond_elem_init}(i) may be relaxed to probabilistic upper bounds. 
The current non-probabilistic bounds are used to simplify the presentation. 
Theorem \ref{thm_NRerr} below establishes 
the near optimality of the proposed one-step estimator \eqref{eq:newtonsolpseudo}.

\begin{theorem}
\label{thm_NRerr} 
Assume Conditions
\ref{cond:truvalueregularity}--\ref{cond_elem_init}. 
Let $\hat Z$ be the generalized one-step estimator defined as in \eqref{eq:newtonsolpseudo}.
Let $\varsigma = \max\{\epsilon,1/2\} $. 
For any constant $s>0$,  there exists a constant $C_s > 0$ such that when $n/ \log^{2\varsigma}(T) $ is sufficiently large, 
\begin{align*}
\Pr\left\{
\operatorname{dist}^2(\hat Z, Z^\star)\, >
  \frac{1}{T}\times C_s r_{n,T}
\right\} =O(n^{-s}),
\end{align*}
where $r_{n,T}=\max\left\{1,\,  \frac{T}{n}\right\} \log^{4\varsigma}(nT)$.
\end{theorem} 


Theorem \ref{thm_NRerr} implies that with high probability, 
the estimation error $\operatorname{dist}^2(\hat Z,Z^\star)$ is $O(r_{n,T}/T)$. For ease of understanding, we first ignore the logarithmic term in $r_{n,T}$, i.e.,
set $\varsigma=0$. 
Then, the error order is reduced to
\begin{align}\label{eq:simplifiedboundonestep}
\frac{1}{T} \times \max\left\{1,\  \frac{T}{n}\right\}=   \max\left\{\frac{1}{T} ,\  \frac{1}{n}\right\}.
\end{align}
When $T=1$, i.e., for a single network,  $\eqref{eq:simplifiedboundonestep} = O(1)$ achieves the oracle error rate.
This is similar to \cite{ma2020universal}, which studied a single network; also see Section \ref{sec:detailedcomareliterature} for a detailed discussion. When $1< T \lesssim n$, 
$\eqref{eq:simplifiedboundonestep} = O(1/T)$, indicating that our proposed generalized one-step estimator achieves the oracle estimation error rate in this case.
When $T \gg n $, $\eqref{eq:simplifiedboundonestep} = O(1/n)$, which unfortunately, is not exactly inverse proportional to $T$. 
Intuitively, when $T$ is very large, 
there is too much time heterogeneity,
resulting in too many nuisance parameters $\alpha_{it}$'s over time. 
In turn, this causes technical difficulty in analyzing the target $Z$. 
Nevertheless, we still have $\eqref{eq:simplifiedboundonestep} =O(1/n)$, which is smaller than the optimal error rate 
$O(1)$ for a single network. 
This suggests that given $T>1$ networks, the estimation error of $Z$ can always be improved compared to using only a single network. 
This can be interpreted, in a broad sense, as an ``inverse proportional to $T$'' property. 
We call $O(1/n)$ a sub-oracle rate throughout the paper.

Compared to \eqref{eq:simplifiedboundonestep}, 
the established error bound $r_{n,T}/T$ contains an  extra logarithmic term
$\log^{4\varsigma}(nT)$, 
which comes from the uniform error control over high-dimensional parameters $(\alpha,Z)$ and the initialization error in Condition \ref{cond_elem_init}(i). 
As an example, 
we propose an initial estimator in Section \ref{sec:initest} satisfying  
$\varsigma=\epsilon=2$. 
For a finite $\varsigma$, the error bound $O(r_{n,T}/T)$, similarly to  \eqref{eq:simplifiedboundonestep},  decreases as $T$ increases. 
To sum up, up to logarithmic factors,
$\hat{Z}$ 
achieves the oracle error rate $ O(1/T)$  , when $1\leqslant T \lesssim n$, 
and  $\hat{Z}$ 
achieves the sub-oracle error rate $O(1/n)$,
when $T \gg n$.

\begin{remark}
\label{RMK:Heff_replace_Ieff}
The efficient Fisher information matrix $\Ieff$ in \eqref{newtonsolUZ} can also be replaced by the ``observed efficient information matrix'' defined as $-\Heff $, where
\begin{align*}
\Heff (Z,\alpha) = \ddot{L}_{Z,Z}(Z,\alpha) - \ddot{L}_{Z,\alpha}(Z,\alpha) \ddot{L}^{-1}_{\alpha,\alpha}(Z,\alpha)\ddot{L}_{\alpha,Z}(Z,\alpha)
\end{align*}
 satisfies $\mathbb{E}\{-\Heff(Z,\alpha)\}=\Ieff(Z,\alpha)$ (an analytic expression for $H_{eff}(Z,\alpha)$ is provided in Section \ref{sec:obseffinf} of the Supplementary Material). 
Then, the one-step estimator is given by
  \begin{equation}
\label{newtonsolUZ2} 
 \check Z_v - \check{\mathcal U}\, \left\{\check{\mathcal U}^{\top}\,    \Heff(\check Z, \check\alpha)\, \, \check{\mathcal U} \right\}^{-1} \,  \check{\mathcal U}^\top \, \Seff(\check Z,\check\alpha),
\end{equation}
 and theoretically, we could establish the same rate as that in Theorem \ref{thm_NRerr}.
More details on this are provided in Section \ref{sec:pfRMK} of the Supplementary Material.
\end{remark}

\begin{remark} \label{rmk::pf_thm1}
To establish the near-oracle rate in Theorem \ref{thm_NRerr}, 
the key idea is to 
show that through our semiparametric efficient construction, the estimation error of $\alpha_{it}$'s would not mask that of $z_i$'s.
Intuitively, 
this is achievable due to the fact 
$\mathbb{E}\{\frac{\partial \Seff(Z,\alpha)}{\partial \alpha_{it}}\}= 0$. 
This shows that a small perturbation error of $\alpha_{it}$'s would not change the efficient score of $z_i$'s in terms of the first-order expansion and therefore, its influence on estimating $z_i$'s can be reduced \citep{bickel1993efficient}. 
Nevertheless, 
establishing the 
near-oracle rate of $Z$ is still challenged by the complex  structures of  $\alpha$ and  $Z$ under our model \eqref{model_DLSM}.  


For the nuisance parameters $\alpha$, 
we do not impose any structural assumptions, such as smoothness or sparsity. Therefore,   
$\alpha_{it}$'s can be completely  heterogeneous over all $i\in \{1,\ldots, n\}$ and $t\in \{1,\ldots, T\}$,
and the free dimension of $\alpha$ is $nT$. 
Such two-way heterogeneity and high dimensionality of $\alpha$
make it difficult to separate the estimation error of $\alpha$ from that of $Z$ as discussed in Section \ref{sec:oraclerateandchallenges}. 
Technically, 
we reduce the influence of estimating $\alpha$ 
by carefully investigating the semiparametric efficient score. 
But that influence cannot be  eliminated entirely, 
which causes a sub-oracle rate  when $T$ is very large, as shown in \eqref{eq:simplifiedboundonestep}.  


For the target parameters $Z$, 
the non-Euclidean geometry and identifiability constraints of $Z$ lead to  
the singularity of the efficient information matrix $\Ieff(\check{Z},\check{\alpha})$ as mentioned in Remark \ref{lb:singularity}. 
To address the issue, we explicitly characterize a restricted subspace of  $\mathbb{R}^{nk}$ on which  $\Ieff(\check{Z},\check{\alpha})$  has positive eigenvalues. 
The restricted eigenspace exhibits interesting properties for characterizing the estimation error of $\hat{Z}$ and has fundamental connections with two sources of identifiability constraints of $Z$ in Remark \ref{lb:singularity}. 
Please see details in  Section \ref{sec:eigenspace} of the Supplementary Material. 
  \end{remark}


\begin{remark}\label{rm:tnratio}
To achieve the desired oracle rate, 
our constraint $T\lesssim n$   
has fundamental connections to classical results in semiparametric analysis. 
We first elaborate on typical rates in classical semiparametric analyses and explain the connections with our result. 
Consider 
a canonical semiparametric problem with  
 target parameters $\theta$ and nuisance parameters $\eta$ 
 with their estimates denoted by $\hat{\theta}$ and $\hat{\eta}$, respectively. 
Under 
suitable regularity conditions, 
the contribution from the nuisance parameters to $\|\hat{\theta}-\theta\|_2$ can often  be bounded by $O(\|\hat{\eta}-\eta\|_2^2)$, a {squared} $\ell_2$-error of nuisance parameters; 
see, e.g., Section 25.59 in \cite{van2000asymptotic}. 
When $\theta$ is fixed-dimensional, its oracle estimation error rate is typically quadratic, i.e., $\|\hat{\theta}-\theta\|_2=O_p(n^{-1/2})$.
Then for the contributed squared error $O(\|\hat{\eta} -\eta\|^2_2)$  to be no larger than $\|\hat{\theta}-\theta \|_2$,  
it is common to impose a quarter rate $\|\hat{\eta} -\eta\|_2=O_p(n^{-1/4})$ on estimating the nuisance parameters. 
Other examples include  
Section 7.6 in \cite{bickel1993efficient},
Section 3 in \cite{murphy2000profile}, and 
Eq. (3.8) in \cite{chernozhukiv.ectj.12097}, etc. 
We next show that our constraint  $T\lesssim n$ aligns with 
this classical quarter rate. 
Specifically, as discussed in Section \ref{sec:oraclerateandchallenges} of the main text,
the oracle rates for estimating each target parameter $z_i^\star$ and each nuisance parameter  $\alpha_{it}^\star$ are expected to be 
$O_p((nT)^{-1/2})$ and $O_p(n^{-1/2})$, respectively. 
Following the above ``squared contribution'' idea,
 when 
each squared oracle error rate of $\alpha_{it}^\star$
is no larger than the oracle error rate of $z_i^\star$, the constraint $T\lesssim n$ is required. 
In this sense, our constraint on $T$ aligns with the common quarter rate requirement on nuisance parameters in the classical literature.

Furthermore, we would like to point out 
that similar error rates and constraints on $T$ and $n$ exist in many other related network studies. 
In particular, we have provided a summary of existing estimation error rates in Table \ref{tab:errorrates} below, showing that all those studies require certain constraints on the growth rate of $T$ with respect to $n$ to achieve the oracle rate for estimating $Z$.  
In terms of the error rate and the constraints on $(n, T)$, 
 our result  is consistent with or more relaxed than the existing results in Table \ref{tab:errorrates}.  	
\end{remark}

 \subsection{Related Works} \label{sec:detailedcomareliterature} 
Various models have been proposed for multiple heterogeneous networks. 
We briefly review related works focusing on fixed-effect continuous latent vectors.  

One class of multiple-network models generalizes the random dot product graphs  \citep{young2007random,athreya2018survey,rubin2017statistical}. 
For example, 
 \cite{arroyo2021inference} considered settings with $T$ layers of networks, 
 and in the $t$-th layer network with $1\leqslant t\leqslant T$, 
 each edge $A_{t,ij}=A_{t,ji}\sim \operatorname{Bernoulli}\{\EXPT(A_{t,ij})\}$ independently for $1 \leqslant t\leqslant T$  and $1 \leqslant i<j\leqslant n $, where 
 $ \EXPT(A_{t,ij})= \langle z_i, z_j\rangle_{\Lambda_t}$ with $\langle z_i, z_j \rangle_{\Lambda_t} = z_i^{\top}\Lambda_t z_j.$
Here $z_i$'s denote the common latent positions shared across $T$ networks, and $\Lambda_{t}\in \mathbb{R}^{k\times k}$ characterizes the layer-specific latent information. 
A similar model for directed random graphs has been studied in \cite{zheng2022limit} with non-symmetric latent positions shared across $T$ networks. 
\cite{macdonald2020latent} proposed models based on generalized random dot product graphs 
and established theoretical results when 
$A_{t,ij}\sim \operatorname{Normal} \{\EXPT(A_{t,ij}), \sigma^2\}$ independently,  
where $\sigma$ is a nuisance parameter,  and 
$ \EXPT(A_{t,ij})= \langle z_i, z_j \rangle_{\mathrm{I}_{p_1,q_1}}  + \langle u_{i,t}, u_{j,t} \rangle_{\mathrm{I}_{p_{2,t}, q_{2,t}}}.$
Here   
$z_i$'s denote the shared latent vectors, whereas $u_{i,t}$'s denote the layer-specific latent vectors, 
and $\langle \cdot,\cdot \rangle_{\mathrm{I}_{p,q}}$ defines an indefinite inner product, where   $\mathrm{I}_{p,q}=\mathrm{diag}(1,\ldots, 1, -1,\ldots, -1)$ with $p$ ones followed by $q$ negative ones on its diagonal for $p\geqslant 1$ and $q\geqslant 0$.  

From the perspective of generalized linear models, 
the above models 
can be viewed as linear responses with the identity link function.
 Considering a single network with binary edges, 
 \cite{hoff2002latent} proposed models with the logit link function, 
 and \cite{ma2020universal} proposed and analyzed two estimation methods. 
 Generalizing the idea to multilayer binary networks,
 \cite{zhang2020flexible} proposed the following model: 
for $t=1,\ldots,T$ layers of networks, each edge $A_{t,ij}\sim \operatorname{Bernoulli}\{\EXPT(A_{t,ij})\}$ independently for $1 \leqslant i, j\leqslant n $, where
$\EXPT(A_{t,ij})=\operatorname{logit}^{-1}(\alpha_{it}+\alpha_{jt}+\langle z_i, z_j\rangle_{\Lambda_t})$ with $\operatorname{logit}^{-1}(x)={e^x}/(1+e^x)$. 

Under the above models, 
one common interest is to estimate the  latent vectors $Z=[z_1,\ldots, z_n]^{\top}$ shared across multiple heterogeneous networks. 
Table \ref{tab:errorrates} summarizes the squared estimation error rates for $Z$ in four studies with a focus on the order of the error with respect to $n$ and $T$ only.   
Our results and intuition on the oracle error rate are consistent with some existing results.
But establishing similar results under our model  \eqref{model_DLSM} requires non-trivial technical developments due to various challenges, including two-way heterogeneity, high-dimensionality, and non-linearity of the link function as discussed in Section \ref{sec:oraclerateandchallenges}.   
To the best of our knowledge, for multiple networks $(T>1)$ with \textit{non-linear} link functions in models, 
there is no existing result that is comparable to our derived nearly optimal rate in Theorem \ref{thm_NRerr}. 
\begin{table}[!htbp]
\centering
\caption{Squared Estimation Error Rates of Shared Latent Vectors $z_i$'s in Different Models.}\label{tab:errorrates} 
\centering
\renewcommand{\arraystretch}{2}
\setlength{\tabcolsep}{2pt}
\begin{tabular}{c|c|c|c|c}
\hline
 & \multicolumn{2}{|c|}{\ Models with Identity Link \ }  & \multicolumn{2}{c}{Models with Logit Link} \\\hline
Model & Arroyo et al. (2021) & MacDonald et al. (2022) &Ma et al. (2020) $(T=1)$ & Zhang et al. (2020) \\ \hline
Error Order & $\displaystyle O\biggr(\frac{1}{T} {+\frac{1}{n}}\biggr)$ & $\displaystyle  O\biggr(\frac{n^{\eta}}{T}\biggr)$ \ \big($\eta \geqslant 0, T\ll n^{1/2}$\big) & $\displaystyle O(1)$ & $\displaystyle O\biggr(1+\frac{1}{T}\biggr)$\\[5pt] \hline
\end{tabular}
\end{table}

 
 Besides the aforementioned works, there are other models for analyzing longitudinal or multiple networks that focus on different goals and challenges. 
 For instance, a prevalent subclass within the latent space models is the stochastic block model (SBM) \citep{holland1983stochastic,lee2019review}.
For a single network, 
our proposed latent space model can reduce to  SBM or degree-corrected block models \citep{karrer2011stochastic,zhao2012consistency,jin2015fast} under suitable parametrizations. 
Notably, the baseline component  $\exp(\alpha_{it})$ for each vertex $i$ is  similar to degree-correction parameters in the degree-corrected block model, showing  
flexibility for accommodating degree heterogeneity. 
Typically, SBM assumes that vertices form stochastically equivalent clusters, which  
leads to low-rank expected adjacency matrix and discrete latent structures. 
The intrinsic discreteness of SBMs makes them structurally simple and helpful in identifying clusters. 
Alternatively, the general latent space models 
are advantageous when there is no apparent clustering structure, and the expected adjacency matrix violates low-rankness. 
As a result, the model properties and analysis techniques under those two types of models can be quite different. 
For multiple networks,
various extensions based on block models
have been proposed  \citep{pensky2019spectral, lei2019consistent,zhang2020modularity, lei2023bias,bazzi2020framework,paul2021null,agterberg2022joint,cai2022consensus}. They often assume the expected adjacency matrix is low-rank and emphasize cluster recovery and community detection, which differs from the objectives and challenges in this work. 


Another line of research proposes to impose a low-rank tensor structure for modeling the longitudinal network \citep{lyu2023latent,zhangefficient}.
Such models differ from our model in terms of interpretation and complexity. 
In particular, under fixed ranks, 
 the total number of parameters under \eqref{model_DLSM} is of the order of $O(nT)$, whereas that under a tensor structure is often reduced to the order of $O(n+T)$ or even smaller. 
 Under \eqref{model_DLSM}, $nT$ number of baseline parameters $\alpha_{it}$ are used to flexibly model the degree heterogeneity of different nodes and time points. This, as mentioned above, is similar to the use of degree-correction parameters in degree-corrected block models and does not have a straightforward counterpart in those related low-rank tensor models. 
Due to this difference in the number of parameters, the model complexities are distinct, and thus, the estimation error rates are not directly comparable.

\subsection{Geometric Interpretation on Quotient Manifold}\label{subsec:manifold}
The proposed \eqref{newtonsolUZ} treats the $n\times k$ matrix $Z$ as an $nk$-dimensional vector $Z_v$. 
This ignores the constraint of the space 
$\mathbb{R}_0^{n\times k} = \{ Z\in \mathbb R^{n\times k}: 1_n^\top Z = 0, \text{det}(Z^\top Z)\neq 0 \}$. 
Moreover,  viewing  $Z$ as an element of $\mathbb{R}^{n\times k}_0$ ignores that it is identified up to a rotation under our considered log-likelihood $L(Z,\alpha)$, as $L(Z,\alpha)=L(ZQ, \alpha)$ for any $Q\in \mathcal{O}(k)$. 
Interestingly, in our considered problem, we find that 
\eqref{newtonsolUZ}
implicitly takes the underlying parameter structure into account. 

To demonstrate this, we first show the intrinsic parameter space of $Z$ is a quotient set. 
In particular, 
the ``rotation invariance''  naturally induces an equivalence relation on $\mathbb R_0^{n\times k} $:
$Z_1 \sim Z_2$ if and only if there exists $ Q\in \mathcal{O}(k)$ such that \ $Z_2 = Z_1 Q$. 
Given  the equivalence relation $\sim$ on $\mathbb{R}_0^{n\times k}$ and an element $Z\in \mathbb{R}_0^{n\times k}$, elements in $\mathbb{R}_0^{n\times k}$ that are equivalent to $Z$ form an equivalence class of $Z$, denoted as $[Z]$. 
The set of equivalence classes of all elements in  $\mathbb{R}_0^{n\times k}$ is called the quotient set, denoted as $\mathbb{R}_0^{n\times k}/\sim $. 
Equipped with the equivalence relation $\sim$, 
 the quotient set $\mathbb{R}_0^{n\times k}/\sim $ naturally incorporates the rotation invariance of $Z$ and  thus is  the intrinsic parameter space to examine.
For the simplicity of notation, we let $\mathcal{M}$ denote $\mathbb{R}_0^{n\times k}/\sim$  below.

We next show that \eqref{newtonsolUZ} can be viewed as an approximation to the one-step Newton-Raphson update of the log profile likelihood on the quotient set $\mathcal{M}$. 
Specifically, the log profile likelihood of $Z$ is $pL(Z)=L(Z,\hat{\alpha}(Z))$ with $\hat{\alpha}(Z)=\argmax_{\alpha \in \mathbb{R}^{n \times T}} L(Z,\alpha)$. 
The MLE of $Z$, i.e., the first component of $(\hat{Z}_{\text{MLE}},\hat{\alpha}_{\text{MLE}})$ that maximizes $L(Z,\alpha)$, is also the maximizer of 
$pL(Z)$. 
Therefore, to investigate the MLE of the target parameter $Z$, it suffices to examine $pL(Z)$ \citep{murphy2000profile}. 
On the quotient set $\mathcal{M}$, 
$pL(Z)$  naturally induces a function  $pL_Q: \mathcal{M} \to \mathbb{R}$ with  $pL_Q([Z])= pL(Z)$ for any $Z \in \mathbb{R}_0^{n \times k }$. 
To incorporate the rotation invariance relation, it is natural to search for the maximizer of $pL_Q$ on the quotient set $\mathcal{M}$.  
Fortunately, a Newton-Raphson step for $pL_Q$ can be properly constructed given the nice properties of $\mathcal{M}$ and $pL_Q$ in Lemma \ref{lm:manifoldsmooth}.  
 
\begin{lemma}[\citealp{lee2013smooth}, Chapter 21]\label{lm:manifoldsmooth}
 
 On  $\mathbb R^{n\times k}_0$, 
 consider the canonical atlas $\varphi : \mathbb{R}^{n \times k}_0 \to \mathbb{R}^{(n-1)k}$ given by  $\varphi(Z)=(z_1^{\top}, \ldots, z_{n-1}^{\top})^{\top}$, and the canonical Riemannian metric 
 $g_Z(Z_1 , Z_2) = \operatorname{tr}(Z_1^{\top} Z_2)$ for $Z_1$ and $Z_2$ in the tangent space to $\mathbb R^{n\times k}_0$ at $Z$. 
 (i) Endowed with the differential structure and the Riemannian metric induced by that of $\mathbb{R}_0^{n \times k}$,
the quotient set $\mathcal{M}$ is a smooth  Riemannian quotient manifold of dimension $nk-k(k+1)/2$.   
(ii) $pL_Q: \mathcal{M}\to \mathbb{R}$ is a smooth function.  
\end{lemma}

 In particular, consider 
 the Riemannian manifold $\mathcal{M}$
 equipped with the Riemannian connection and the exponential retraction $R$ \citep[Section 5]{absil2009optimization}. 
For a smooth function $f: \mathcal{M}\to  \mathbb{R}$, 
the Newton-Raphson update at an equivalence class $[Z]\in \mathcal{M}$ is given by $R_{[Z]}(\nu)$, where $\nu$ is a tangent vector 
in the tangent space 
$T_{[Z]}\mathcal{M}$ and 
specified by  
$\operatorname{Hess}f([Z]) [\nu] = -\operatorname{Grad}f([Z]),$
and $ \operatorname{Hess}f([Z])$ and $\operatorname{Grad}f([Z]) $ denote the Hessian and gradient of $f$ at $[Z]$, respectively  \citep{boumal2020introduction}. 
Conceptually, $R_{[Z]}(\nu)$ 
 defines a move in the direction of $\nu$ while staying on the manifold. 
 Interestingly, for the function $pL_Q$, we prove that $R_{[Z]}(\nu)$ has an analytic (matrix) form in Proposition \ref{prop:newtonstepform}.  

\begin{proposition}
\label{prop:newtonstepform}
Under the same setting as Lemma \ref{lm:manifoldsmooth},
equip $\mathcal{M}$ with the Riemannian connection and the exponential retraction $R$.
Consider an initial value $\check{Z}$ and its equivalence class $[\check{Z}]$.
For $pL_Q$ at $[\check{Z}]$, the Newton-Raphson update $R_{[\check{Z}]}(\nu) = [\operatorname{mat}(\check{Z}_v + \bar{\nu}_{\check{Z}})]$, where 
$\check{Z}_v$ is the vectorization of $\check{Z}$, 
$\operatorname{mat}((z_1^{\top},\ldots, z_n^{\top})^{\top}) = (z_1,\ldots, z_n)^\top\in \mathbb{R}^{n\times k}$, 
 and  
\begin{align}\notag
  \bar{\nu}_{\check{Z}} = -\,  \check{\mathcal{U}}\left[\,  \check{\mathcal{U}}^{\top}\   \Heff\left\{ \check{Z}, \hat{\alpha}(\check{Z})\right\} \check{\mathcal{U}}\right]^{-1}\  \check{\mathcal{U}}^{\top}\, \Seff\left\{ \check{Z}, \hat{\alpha}(\check{Z})\right\},
\end{align}
where $\check{\mathcal{U}}$ is the same as that in \eqref{newtonsolUZ}, and $\Heff (Z,\alpha) $ is defined same as in Remark \ref{RMK:Heff_replace_Ieff}. 
\end{proposition}

Note that $\check{Z}_v+\bar{\nu}_{\check{Z}}$ is similar to  \eqref{newtonsolUZ2}  except that 
 $\hat{\alpha}(\check{Z})$ is replaced with $\check{\alpha}$. 
 Given initial values $(\check{Z},\check{\alpha})$ that are close to the true values $(Z^{\star},\alpha^{\star})$, 
we can show that $\hat{\alpha}(\check{Z})$ and $\check{\alpha}$ are close, and thus $\check{Z}_v+\bar{\nu}_{\check{Z}}$ and \eqref{newtonsolUZ2} are close; see rigorous results in Section \ref{sec:addgeores} of the Supplementary Material. 
By Remark \ref{RMK:Heff_replace_Ieff}, 
\eqref{eq:newtonsolpseudo}, 
\eqref{newtonsolUZ}, and \eqref{newtonsolUZ2} can all be viewed as approximations to a one-step Newton-Raphson update on the quotient manifold $\mathcal{M}$,
whose construction 
 implicitly incorporates the underlying geometric structures of the quotient manifold.  
 \begin{figure*}
\quad
\begin{tikzpicture}
\pgfplotsset{width=0.49\textwidth,compat=1.9}
\begin{axis}
[
    title={3D plot of $f(z_{1},z_{2})=-(z_{1}^2+z_{2}^2)$},
    ymin=-5,ymax=5,
    xmin=-5,xmax=5,
]
\addplot3[
    surf,
]
{(-x^2-y^2)};

\addplot3 [black,domain=-5:0,samples y=1,-stealth, line width=1pt, dashed] (-x, x, -2 * x^2);

\end{axis}
\end{tikzpicture}\quad \  
\includegraphics[height=0.4\textwidth]{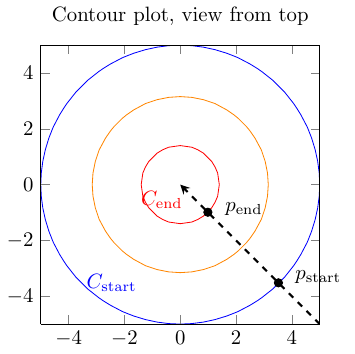} 



\caption{Each circle represents one equivalence class of $(z_1,z_2)$ giving the same value of $\bar{f}(z_1,z_2)$. It suffices to search for maximization along one given direction.} \label{fig:fzillustration}
\end{figure*}


\begin{remark}
We provide an intuitive explanation of   Proposition \ref{prop:newtonstepform}.
Motivated by the fact that the log profile likelihood $pL(Z)$ depends on parameters in $Z$  only through  $ZZ^{\top}$,
we next focus on $ZZ^{\top}$. 
For ease of visualization, 
 we first consider  the  parameter space 
$\mathbb{S}_0^{2\times 2}=\{ Z \in \mathbb{R}^{2\times 2}: 1_2^{\top}Z=0_{1\times 2} \}$. 
For $Z\in \mathbb{S}_0^{2\times 2}$, $ZZ^{\top}=(z_{11}^2+z_{12}^2)(1, -1)^{\top}(1, -1)$ depends on $Z$ only through the  scalar $z_{11}^2+z_{12}^2$. 
This 
suggests that 
 the number of free parameters in $ZZ^{\top}$ is 1, which equals $nk-k(k+1)/2$ as analyzed in Remark \ref{lb:singularity}.  
 We visualize the function $\bar{f}(z_{11},z_{12})=-(z_{11}^2+z_{12}^2) $ in Figure  \ref{fig:fzillustration}. 
The function has a constant value on each circle of $(z_{11},z_{12})$ in the contour plot. 
From this perspective, each circle is an equivalence class giving the same value of $\bar{f}$.  
The set of all 
equivalence classes (circles) yields a quotient set, denoted as  $\mathcal{M}_2$, 
and $\bar{f}$ induces a function $f: \mathcal{M}_2 \to \mathbb{R}$ given by $f([z]) = \bar{f}( z_{11}, z_{12})$, where  $[z] \in \mathcal{M}_2$ represents the equivalence class specified by $(z_{11}, z_{12}) $. 

Consider the problem of maximizing the function $f$ on $\mathcal{M}_2$. 
One natural idea is to search for the maximizer of $f$ across its domain $\mathcal{M}_2$.
In this case, an update step in $\mathcal{M}_2$  yields a move from one circle to another, e.g., $C_{\text{start}}$ to $C_{\text{end}}$ in Figure \ref{fig:fzillustration}. 
Such an abstract update can also be described in the Euclidean space $\mathbb{R}^2$. 
In particular, 
fix one direction that is orthogonal to a tangent line of circles, e.g., the dashed line in Figure \ref{fig:fzillustration},
move along the tangent line from the point $p_{\text{start}}$ to the point $p_{\text{end}}$ in Figure \ref{fig:fzillustration}, 
and map $p_{\text{end}}$ to the corresponding equivalence class $C_{\text{end}}$.  
The update from $p_{\text{start}}$ to  $p_{\text{end}}$ in $\mathbb{R}^2$ can be analytically represented in a matrix form. 
Generalizing this idea, an update step on the quotient manifold $\mathbb{R}_0^{n\times k}/\sim $ can be described through a fixed search space 
in the Euclidean space $\mathbb{R}^{n\times k}$ (analogous to the search direction in Figure \ref{fig:fzillustration}), 
and a common choice is the 
so-called horizontal space \citep{boumal2020introduction}.  
Interestingly, 
we prove that the horizontal space at $\check{Z}$ 
can be 
{explicitly expressed} through $\check{\mathcal{U}}$. 
This enables us to derive an analytic form  
$\bar{\nu}_{\check{Z}}$, which, after properly aligned into a matrix, defines an update in the horizontal space. 
The retraction $R$ maps the post-update  matrix $\operatorname{mat}(\check{Z}_v + \bar{\nu}_{\check{Z}})$ to an element in $\mathbb{R}_0^{n\times k}/\sim $, which gives the targeted one-step estimator in the manifold. 
\end{remark}

\subsection{Initial Estimation} \label{sec:initest}
In this section, we construct an initial estimator $(\check Z,\check\alpha)$ and prove that it satisfies  Condition \ref{cond_elem_init} with $\epsilon=2$ and high probability. The proposed initialization algorithm consists of two stages outlined below, while the implementation details are provided in Section \ref{sec:details} of the Supplementary Material. 
In the following, we define $E_t^{\star}=\mathbb{E}(\mathbf A_t)$,  
$\Theta_{t}^{\star}=\log (E_t^{\star})$ for $t=1,\ldots, T$,
and $G^{\star}=Z^{\star}(Z^{\star})^{\top}$. 
We note that
\begin{align}\label{eq:relationshipstar}
\alpha_t^{\star}=H_n^{-1}\Theta_t^{\star}1_n,\hspace{2em} \text{and} \hspace{2em} G^{\star}= \sum_{t=1}^T(\Theta_t^{\star}-\alpha_t^{\star}1_n^{\top} - 1_n (\alpha_t^{\star})^\top)/T, 
\end{align}
where $H_n= n\mathrm{I}_n + 1_n1_n^{\top}$. 

\subsubsection*{Outline of the Initialization Algorithm} 
\begin{itemize}
\item[(1)] The first stage obtains an ``initial of initial'' estimator $(\mathring{Z},\mathring{\alpha})$ as follows.
A similar ``double-SVD'' idea has been introduced in \cite{zhang2020note}.
\begin{itemize}
\item[(a)] 
Construct an estimator $\mathring{E}_{t}$  for $E_t^{\star}$ by applying the universal singular value thresholding  \citep{chatterjee2015matrix} to each matrix $\mathbf A_t$. Take $\mathring{\Theta}_{t}=\log(\mathring{E}_{t})$ as an estimator for $\Theta_{t}^{\star}$. 
\item[(b)] Let $\mathring{\alpha}_{t}=H_n^{-1} \mathring{\Theta}_{t}1_n$, motivated by \eqref{eq:relationshipstar}. 
Take $\mathring{\alpha} = (\mathring{\alpha}_{1}, \ldots,\mathring{\alpha}_{T})$.
\item[(c)] Let $\mathring{G} =  \sum_{t=1}^T(\mathring{\Theta}_{t}-\mathring{\alpha}_{t}1_n^{\top} - 1_n \mathring{\alpha}_{t}^\top)/T$, motivated by \eqref{eq:relationshipstar}.\\
Let $\mathring{Z}=\mathring{U}_{k}\mathring{D}_{k}^{1/2}$, where $\mathring{U}_{k}\mathring{D}_{k}\mathring{U}_{k}^{\top}$ denote the top-$k$ eigen-decomposition of $\mathring{G}$. 
\\
\end{itemize}
\item[(2)] 
The second stage utilizes the projected gradient descent method \citep{chen2015fast} with $L(Z,\alpha)$ as an objective function to maximize, and $(\mathring Z ,\mathring \alpha )$ from the first stage as initial values.  In particular,  

\begin{itemize}
\item[(a)] 
Update $Z$ and $\alpha$ along their corresponding gradient directions with pre-specified step sizes. 
\item[(b)] 
Project the updated estimates to 
the constraint set $\mathcal{S}_C$ induced by the conditions of identifiability and boundedness of parameters in Condition \ref{cond:truvalueregularity}. Specifically, 
$$\mathcal{S}_C=\{(Z,\alpha): Z \in \mathbb{R}^{n \times k},\, 1_n^{\mytrans} Z=0,\, \max_{ i}\left\|z_i\right\|^2_2 \leqslant M_{Z,1},\, \alpha \in \mathbb{R}^{n\times T},  \, \max_{i,t}|\alpha_{it}| \leqslant M_{\alpha} \}. $$

\item[(c)] 
 Repeat (a) and (b) until convergence. The resulting values, denoted as $(\check{Z}, \check{\alpha})$, are used as the initial estimator in the proposed one-step estimator. 
\end{itemize}

\end{itemize} 

\medskip
In Theorem \ref{thm:initialerror}, we establish an {elementwise} error bound of $(\check Z,\check\alpha)$ obtained through the proposed initialization algorithm above. 
\begin{theorem} \label{thm:initialerror}
    Assume Conditions \ref{cond:truvalueregularity} and \ref{cond:tuning}. 
    For any constant $s>0$,
    there exists a  constant $C_s>0$
    such that when $n /\log^{k+3}(T)$ is sufficiently large, 
    $$
    \Pr\left[\max_{1\leqslant i \leqslant n,\, 1\leqslant t \leqslant T} \left\{ \left|\check\alpha_{it} - \alpha_{it}^\star\right| +
    \mathrm{dist}_i(\check{z}_i, z_i^{\star}) \right\} > 
     \frac{  C_s\log^2(nT) }{\sqrt{n}} \right] = O( (nT)^{-s}).
    $$
\end{theorem}


Theorem \ref{thm:initialerror} implies that with high probability, the initial estimator $(\check Z,\check\alpha )$ satisfy Condition \ref{cond_elem_init} with $\epsilon = 2$.

    
\section{Semiparametric Penalized Maximum Likelihood Estimator} \label{sec:penalizedmethod}

The log-likelihood $L(Z,\alpha)$ is non-convex in $Z$, and also requires specification of the dimension of each $z_i$. On the other hand, when the log-likelihood \eqref{eq:LZalpha} is reparametrized through $G = (G_{ij})_{n \times n}$ with $G_{ij} = \langle z_i,z_j\rangle$, the log-likelihood function 
$$
l(G,\alpha) = \sum_{t=1}^T \sum_{1\leqslant i\leqslant j \leqslant n} \left[A_{t,ij} (\alpha_{it}+\alpha_{jt}+G_{ij}) - \exp(\alpha_{it}+\alpha_{jt}+G_{ij}) \right]
$$
is convex in $\alpha$ and $G$. Nevertheless, $G=ZZ^{\top}$ itself overparameterizes the model, and the dimension of $z_i$ corresponds to the rank of $G$ in this model specification. Subject to the rank constraint, solving MLE of $G$ is a non-convex optimization problem, even though  $l(G,\alpha)$ is convex. 
To overcome this issue, we construct a penalized maximum likelihood estimator of $G$ with a nuclear norm penalization that relaxes the exact rank constraint on $G$. We further show that the proposed estimator can achieve the corresponding almost-oracle error rate. 

We first characterize the parameter space of $G$ and its relationship with $Z$. From the perspective of likelihood functions, matrices $G$ and $Z$ such that $L(Z,\alpha)=l(G,\alpha)$ can be viewed as equivalent parameters. 
Accordingly, any $Z \in \mathbb{R}^{n\times k}_0$ uniquely  defines  $G=ZZ^{\top} \in \mathbb{S}_{0,+}^{n,k}$, 
where $\mathbb{S}_{0,+}^{n,k}$ represents the class of all $n\times n$ positive semidefinite symmetric matrices with rank $k$ and $G1_n=0_n$. Moreover, any matrix $G\in \mathbb{S}_{0,+}^{n,k}$ induces  $G^{1/2}=U_kD_k^{1/2} \in \mathbb{R}_{0}^{n\times k}$, where $U_kD_kU_k^{\top}$ is the top-$k$ eigenvalue components of $G$. Then $G$ uniquely identifies an equivalence class $\{G^{1/2}Q: Q \in \mathcal{O}(k) \} = \{ Z\in \mathbb{R}_0^{n\times k}: ZZ^{\top} = G \}$. 

The penalized maximum likelihood estimator is defined as the solution $(\hat G,\hat\alpha)$ to the following convex optimization problem
\begin{align}\label{eq:penalizedMLEproblem}
    \max_{G,\ \alpha}\ \ &~ \  l(G, \alpha) - \lambda_{n,T}\|G\|_* \\
\text{subject to} \ \ &~  G\in \mathbb{S}_{+}^n, \quad  G1_n=0_n, \quad |G_{ij}| \leqslant {M_{Z,1}}, \notag 
\end{align}
where $\mathbb{S}_{+}^n$ represents the class of $n\times n$ positive semidefinite matrices, 
$M_{Z,1}$ is a constant same as in Condition \ref{cond:truvalueregularity}, and $\lambda_{n,T}$ is a prespecified tuning parameter.

To scale the error of $G$ to the error of $Z$, we present the estimation error of $G$ in terms of $\|\hat G - G^\star\|_{\mathrm{F}}^2 / n$; 
see Remark \ref{rm:gtozoracleerror}.
 The following Theorem \ref{MLE_mainthm} shows that the penalized MLE $\hat{G}$ solved from \eqref{eq:penalizedMLEproblem} almost achieves the oracle rate, similarly to Theorem \ref{thm_NRerr}. 
 \begin{theorem} \label{MLE_mainthm} 
Assume $G^{\star}=Z^{\star}Z^{\star\top}$ with $Z^{\star} \in \mathbb{R}_0^{n\times k}$ satisfying Condition \ref{cond:truvalueregularity}.  
Suppose we choose $\lambda_{n,T} \asymp \sqrt{nT}\log(nT) $. 
Then for any constant $s>0$, there exists a constant $C_s>0$ such that when $n/\log(T)$ is sufficiently large, 
\begin{align*} 
\Pr\left\{
\|\hat{G} - G^{\star}\|_{\mathrm{F}}^2 /n\, >
  \frac{1}{T}\times C_s r_{n,T}'
\right\} =O(n^{-s}),  
\end{align*}
  where $r_{n,T}'=\max\left\{1,\,  \frac{T}{n}\right\} \log^2(nT)$.
\end{theorem}


Theorem \ref{MLE_mainthm} implies that with high probability, 
the estimation error $\|\hat G-G^\star\|_{\mathrm{F}}^2/n$ is $O(r_{n,T}'/T)$. 
Similarly to Theorem \ref{thm_NRerr}, when ignoring the logarithmic term in $r_{n, T}'/T$, the error order reduces to
\eqref{eq:simplifiedboundonestep}. 
Therefore, up to the logarithmic term,  the penalized MLE achieves the oracle estimation error rate $O(1/T)$ when $1\leqslant T \lesssim n$, and it achieves the sub-oracle rate $O(1/n)$
when $T\gg n$. 

\begin{remark}\label{rm:gtozoracleerror}
Intuitively, 
$n^2$ parameters in $G^{\star}$ are redundant and essentially determined by only $nk$ parameters in $Z^{\star}$.
To align with the error rate in Theorem \ref{thm_NRerr}, 
we consider the estimation error of $\hat{G}$  with a scaling factor $1/n$ in 
 Theorem \ref{MLE_mainthm}.   
When the rank $k$ is known, the one-step estimator $\hat{Z}$ induces the estimator $\hat{Z}\hat{Z}^{\top}$ for $G^{\star}$ satisfying   
 $\|\hat{Z}\hat{Z}^{\top} - G^\star\|_{\mathrm{F}}^2 / n  \asymp \operatorname{dist}^2(\hat Z,Z^\star)\asymp 1/T$; see Lemmas \ref{Lemma-28} and \ref{Lemma-29} in the Supplementary Material.
 When $k$ is unknown, 
Theorem \ref{MLE_mainthm} suggests 
the penalized MLE $\hat{G}$ achieves the same error rate as that of $\hat{Z}\hat{Z}^{\top}$.  \end{remark}

\begin{remark}\label{rm:techchanlplm}
The idea of using the nuclear norm penalization as a convex relaxation of the exact rank constraint originates from the low-rank matrix recovery literature 
\citep{candes2010power,candes2012exact,davenport20141}. 
It was also utilized by \cite{ma2020universal} to fit the latent space model for a single network. But different from these studies, we have two components $G$ and $\alpha$ in our model with different but entangling oracle error rates. To establish a sharp error rate for $\hat G$, we develop a novel semiparametric analysis for the penalized MLE through the profile likelihood \citep{murphy2000profile}.

In particular, we note that $\hat{G}$ is also the penalized MLE of  the log profile likelihood  $ pl(G)= l(G,\hat{\alpha}(G))$, where $\hat{\alpha}(G)=\arg\max_{\alpha\in \mathbb{R}^{n\times T}}l(G,\alpha)$. 
Let $\dot{l}_{\alpha}$ and $\ddot{l}_{\alpha\alpha}$ represent the first-order derivatives and the second-order derivatives of $l(G,\alpha)$ with respect to $\alpha_v$, respectively. 
We can obtain 
$\hat{\alpha}(G)-\alpha= (\ddot{l}_{\alpha\alpha})^{-1}\dot{l}_{\alpha}+\text{high-order terms}$,  and $\mathbb{E}( \dot{l}_{\alpha}) = 0$. 
 Intuitively, this suggests that the perturbation error from $\hat{\alpha}(G)$ can be small in terms of the first-order expansions, so its impact on the profile likelihood of $G$ can be reduced. 

Despite the wide use of the profile likelihood method \citep{murphy2000profile,severini1992profile},
our analysis is faced with unique technical challenges.  
First, not only are the nuisance parameters $\alpha$  high-dimensional but also they encode two-way heterogeneity over both nodes and time. 
As a result, to separate the estimation error of $G$ from that of $\alpha$ requires delicate analysis that differs significantly from the existing studies. 
Second, the target parameter $G$ is high-dimensional and intrinsically redundant with a low-rank structure. 
The redundancy of parameters in $G$, as explained in Remark \ref{lb:singularity},  
leads to singularity issues when analyzing the likelihood function and calls for new technical developments.  
Meanwhile, the low-rank structure of $G$ motivates the use of the nuclear norm penalization, which needs to be properly taken care of in our semiparametric analysis. 
\end{remark}

\begin{remark}[Relationship between the two  estimators] \label{rm:relationtwoest}
The two estimators are fundamentally connected as they are both motivated by maximizing the semiparametric profile likelihood. 
In particular,  
the penalized MLE is a convex relaxation of maximizing the log-likelihood function $l(ZZ^{\top},\alpha)$, which is equivalent to maximizing the log profile likelihood $pL(Z)$. 
The one-step estimator solves the maximization of the quadratic function
$q(x)=pL(\check{Z})+(x-\check{Z}_v)^{\top}\Seff(\check{Z},\check{\alpha})-\frac{1}{2}(x-\check{Z}_v)^{\top}\Ieff(\check{Z},\check{\alpha})(x-\check{Z}_v)$. 
In settings with fixed-dimensional target parameters, 
the log profile likelihood  
 is approximated by its local quadratic expansion similar to $q(x)$  \citep{murphy2000profile}.  
 From an alternative point of view, the one-step estimator gives an approximate solution to the efficient score equation, which, in general, is equivalent to maximizing the profile likelihood.
Therefore,  we expect that the two proposed estimators are similar up to the convex relaxation, and they can achieve the oracle error rate as the nuisance $\alpha$ is eliminated through the profile likelihood.

In practice, the penalized MLE  differs from the non-convex one-step estimator since it adopts a convex relaxation and does not require a prespecification of the true latent dimension $k$.  
Nevertheless, the computation of the one-step estimator is more efficient than solving penalized MLE, and it can achieve a better performance when oracle $k$ is specified. Both convex and non-convex estimators are widely examined in the related literature \citep{ma2020universal}, and their theoretical properties are of interest to study.  
\end{remark} 

\begin{remark}\label{rm:fromgtoz}
Given $\hat{G}$ solved from \eqref{eq:penalizedMLEproblem}, 
we can further construct an estimator for $Z^{\star}$. If the true $k$ is given, 
we let $\hat{Z}_k=\hat{U}_k\hat{D}_k^{1/2}$, 
where $\hat{U}_k\hat{D}_k\hat{U}_k^{\top}$ is the top-$k$ eigenvalue components of the penalized MLE $\hat{G}$. 
When $k$ is unknown, we can consistently estimate $k$. 
 Intuitively, since $k$ equals the rank of $G^{\star}$, 
we can estimate $k$ as the number of significant eigenvalues of an estimator $\Ginit$ that approximates $G^{\star}$ sufficiently well. 
For instance, we can set $\Ginit$ to be  
 the penalized MLE $\hat{G}$ in 
\eqref{eq:penalizedMLEproblem} or  the ``double-SVD'' estimator $\mathring{G}$ in Section \ref{sec:initest}, both of which can be obtained without prespecifying $k$. 
We have developed rigorous theoretical guarantees and conducted numerical experiments in Section \ref{sec:estimatekdetails} of the Supplementary Material.
Notably for the one-step estimator in Section \ref{sec:newton}, we can also consistently estimate $k$ when it is unknown.   
\end{remark}

\section{Simulation Studies}
\label{sec::simulations}
We conduct simulation studies to examine how the estimation errors vary with respect to $n$ and $T$.  
To this end, we consider the following two scenarios of $(n, T)$:  
\begin{enumerate}
    \item[(a)] fixing $n=200$ and  varying $T\in \{5, 10, 20, 40, 80\}$;
    \item[(b)] fixing $T=20$ and varying $n \in \{100, 200, 400, 800\}$. 
\end{enumerate}
In each scenario, we allow varying $k\in \{2,4,8\}$.   
Given $(n,T,k)$, 
we generate the shared latent vectors $Z^{\star}$ as follows: 
 independently sample  $w_i$  following a uniform distribution on $\mathbb{B}_2^k = \{x \in \mathbb{R}^k : \|x\|_2 \leqslant 1\}$,  
let $W=[\tilde{w}_1,\ldots, \tilde{w}_n]^{\top}$ where $\tilde{w}_j=w_j-\sum_{l=1}^nw_l/n$,
and set $Z^{\star}=\sqrt{n} W/\|WW^{\top}\|_{\mathrm{F}}^{1/2}$.
In this way, $Z^{\star}$ has zero column means, and $G^{\star}=Z^{\star}Z^{\star\top}$ satisfies  $\|G^{\star}\|_{\mathrm{F}}/n=1$.  
For the baseline heterogeneity parameters $\alpha^{\star}$, we consider two cases
\begin{itemize}
    \item[(I)]   
    $\alpha_{it}^{\star} $ are independently sampled from $ U(-2,0)$; 
    \item[(II)] $ \alpha_{it}^{\star} = t/T + \beta_{it}$ with $ \beta_{it}$ independently sampled from $U(-3,-1)$ for $1 \leqslant i \leqslant \lfloor n/2 \rfloor$, and $ \alpha_{it}^{\star} = - 2t/T + \beta_{it}$ with $ \beta_{it}$ independently sampled from $U(-2,0)$  for $\lfloor n/2 \rfloor < i \leqslant n$. 
\end{itemize}
Intuitively, $\alpha^{\star}$'s in Cases (I) and (II) represent different types of heterogeneity; the former is uniform over $i$, whereas the latter has a two-block structure over $i$. 


Under each model configuration and for each estimator, we estimate the error  $\operatorname{dist}^2(\hat Z, Z^\star)$ over 50 Monte Carlo simulations. 
We compute $\operatorname{dist}^2(\hat Z, Z^\star) = \|\hat{Z}-Z^{\star} VU^{\top}\|_{\mathrm{F}}^2$ where $U\Sigma V^{\top}$ is the singular value decomposition of $\hat{Z}^{\top}Z^{\star}$; see  \cite{schonemann1966generalized}.
For the penalized MLE, we obtain $\hat{Z}$  following Remark \ref{rm:fromgtoz} given true $k$; 
 the corresponding errors of $\hat{G}$ follow  similar patterns and are provided in Section \ref{sec:suppsimultaion} of the Supplementary Material.

\begin{figure}[!htbp] 
	\centering
  \begin{subfigure}{0.32\textwidth}
		\centering
		\includegraphics[width=1\linewidth]{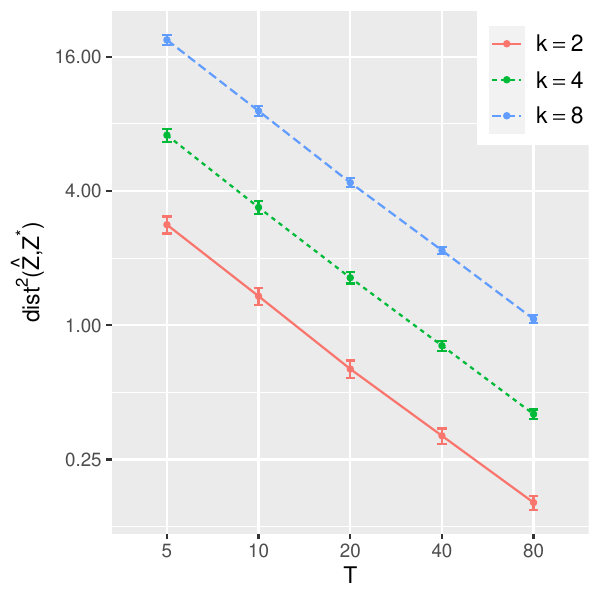}
  \caption{$\operatorname{dist}^2(\hat Z, Z^\star)$ versus $T$}
\end{subfigure}
\begin{subfigure}{0.32\textwidth}
		\centering
		\includegraphics[width=1\linewidth]{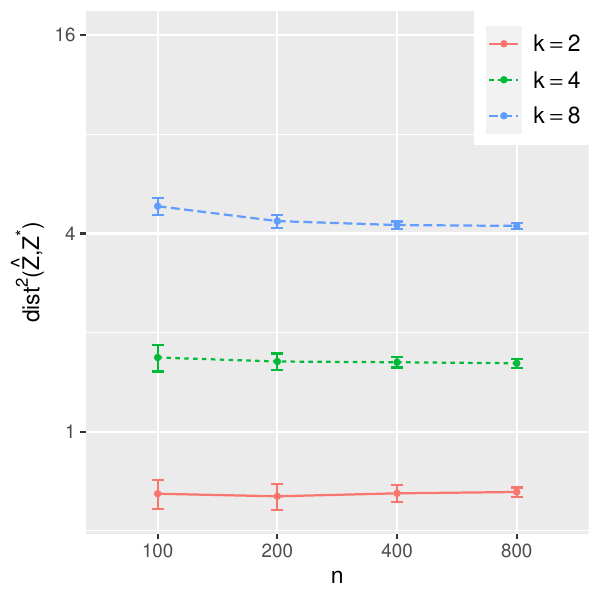}   
\caption{$\operatorname{dist}^2(\hat Z, Z^\star)$ versus $n$}
	\end{subfigure}
\begin{subfigure}{0.32\textwidth}
		\centering
		\includegraphics[width=1\linewidth]{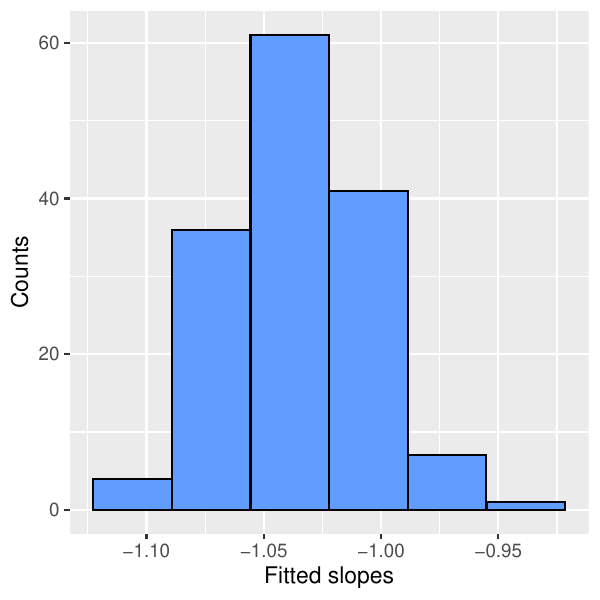}
    \caption{Histogram of fitted slopes}
\end{subfigure}
 \caption{Case (I): Empirical estimation errors of the one-step estimator. Panel (a) presents    $\operatorname{dist}^2(\hat Z, Z^\star)$ (averaged over 50 repetitions) versus $T$ in the scenario (a).  
 Panel (b) presents  $\operatorname{dist}^2(\hat Z, Z^\star)$ (averaged over 50  repetitions)  versus $n$ in the scenario (b). 
In (a) and (b), axes are in the log scale, three lines correspond to results under $k\in \{2,4,8\}$, respectively, and error bars are obtained by $\pm$ the {standard deviation} from 50 repetitions.
 Panel (c) presents the slopes from regressing $\log \operatorname{dist}^2(\hat Z, Z^\star)$ on $\log T$ with fixed $(n,k)\in \{200\}\times \{2,4,8\}$ in the 50 repetitions under the scenario  (a). 
} \label{fig:resultscasea1}
\end{figure}

\begin{figure}[!htbp]
	\centering
\begin{subfigure}{0.32\textwidth}
		\centering
		\includegraphics[width=1\linewidth]{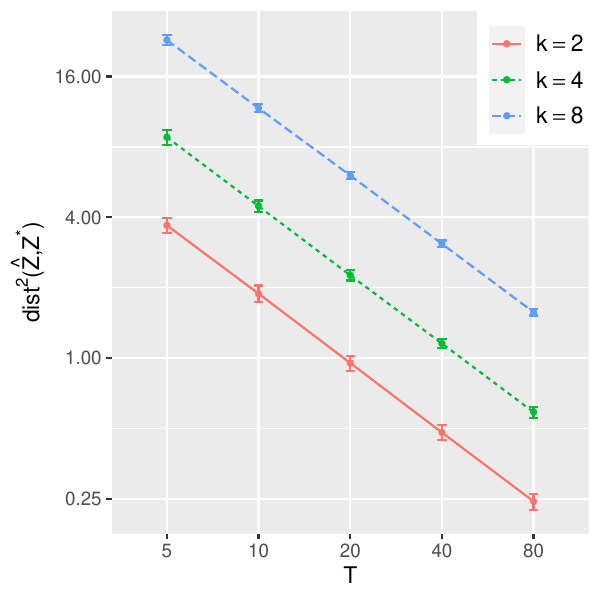}
  \caption{$\operatorname{dist}^2(\hat Z, Z^\star)$ versus $T$}
\end{subfigure}
\begin{subfigure}{0.32\textwidth}
		\centering
		\includegraphics[width=1\linewidth]{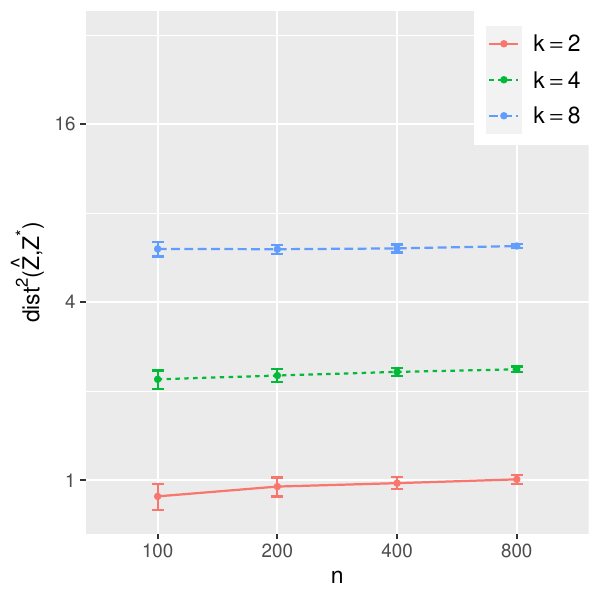}
    \caption{$\operatorname{dist}^2(\hat Z, Z^\star)$ versus $n$}
	\end{subfigure}
 \begin{subfigure}{0.32\textwidth}
		\centering
		\includegraphics[width=1\linewidth]{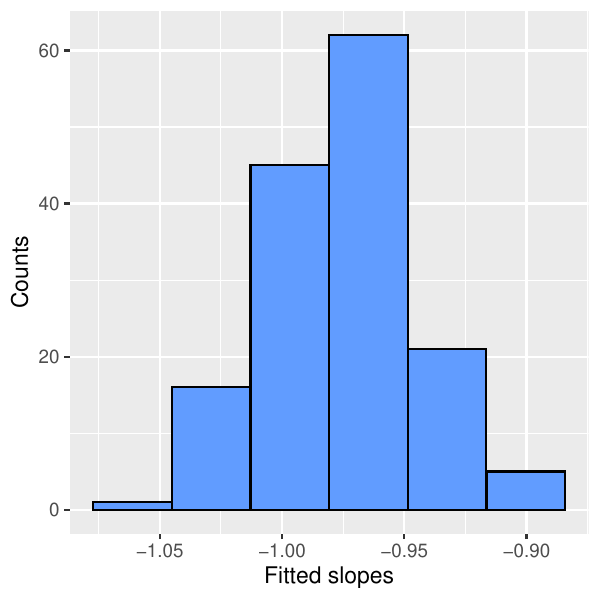}
    \caption{Histogram of fitted slopes} 
\end{subfigure}
 \caption{Case (I): Empirical estimation errors of the penalized MLE. Panels (a)--(c) are presented similarly to Figure \ref{fig:resultscasea1}.}  \label{fig:resultscasea2}
\end{figure}

\begin{figure}[!htbp]
\centering
 \begin{subfigure}{0.32\textwidth}
		\centering
		\includegraphics[width=1\linewidth]{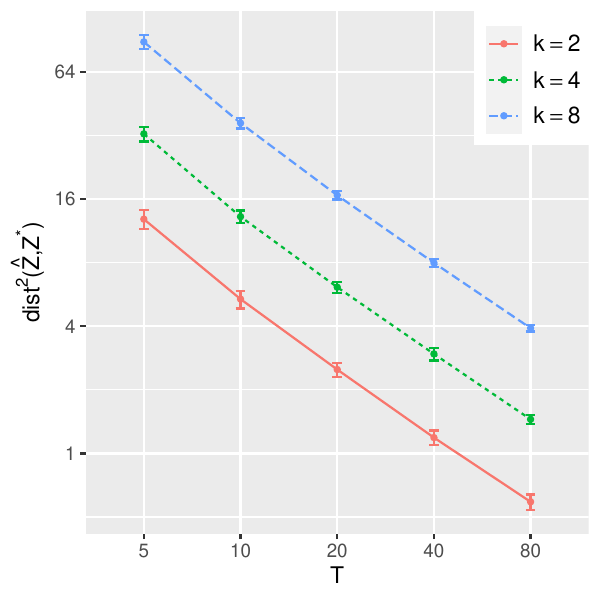}
  \caption{$\operatorname{dist}^2(\hat Z, Z^\star)$ versus $T$}
\end{subfigure}
\begin{subfigure}{0.32\textwidth}
		\centering
		\includegraphics[width=1\linewidth]{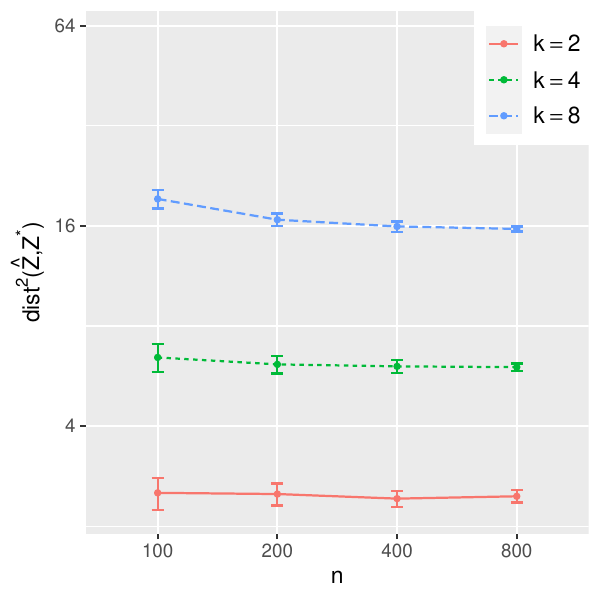}
  \caption{$\operatorname{dist}^2(\hat Z, Z^\star)$ versus $n$}
	\end{subfigure}
\begin{subfigure}{0.32\textwidth}
		\centering
		\includegraphics[width=1\linewidth]{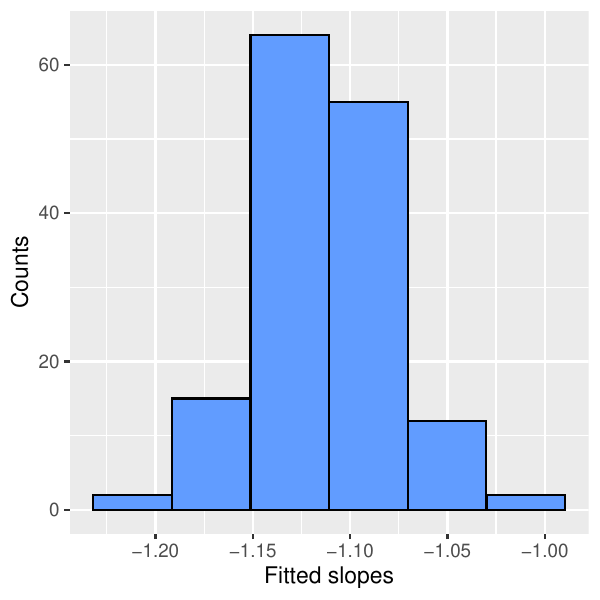}
  \caption{Histogram of fitted slopes}
\end{subfigure}
 \caption{Case (II): Empirical estimation errors of the one-step estimator. Panels (a)--(c) are presented similarly to Figure \ref{fig:resultscasea1}.}  \label{fig:resultscase2nonconvex}
\end{figure}

\begin{figure}[!htbp]
\centering
\begin{subfigure}{0.32\textwidth}
		\centering
		\includegraphics[width=1\linewidth]{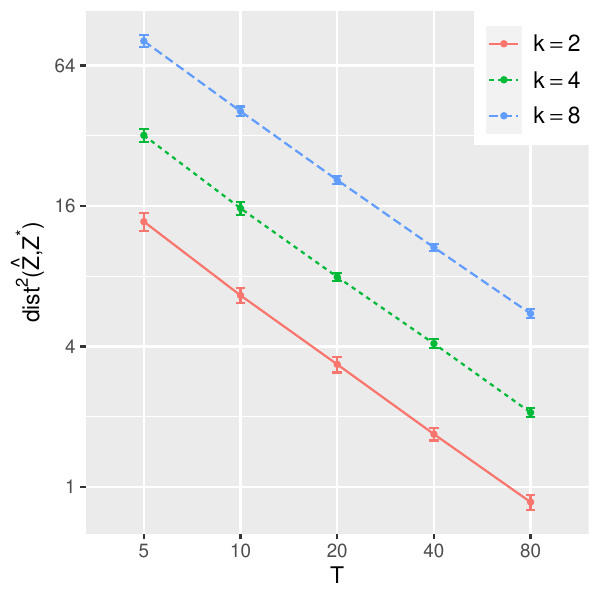}
  \caption{$\operatorname{dist}^2(\hat Z, Z^\star)$ versus $T$}
\end{subfigure}
\begin{subfigure}{0.32\textwidth}
		\centering
		\includegraphics[width=1\linewidth]{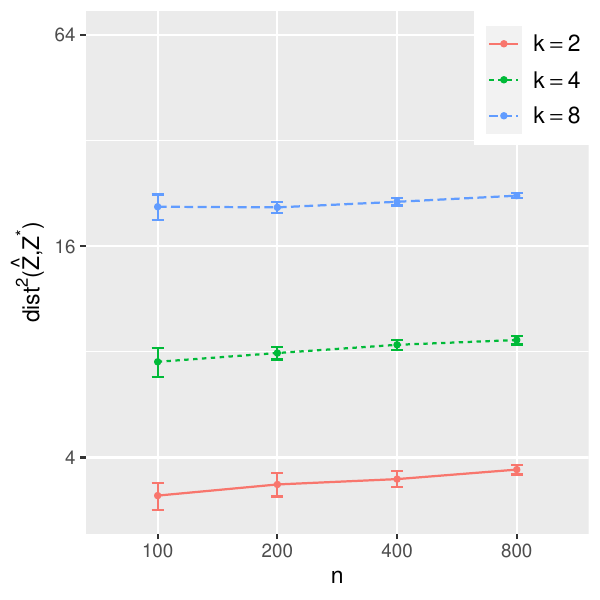}
  \caption{$\operatorname{dist}^2(\hat Z, Z^\star)$ versus $n$}
	\end{subfigure}
 \begin{subfigure}{0.32\textwidth}
		\centering
        \includegraphics[width=1\linewidth]{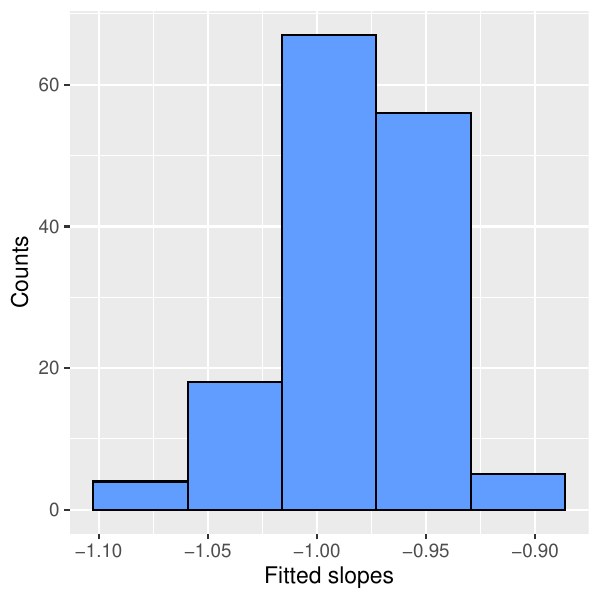}
  \caption{Histogram of fitted slopes}
\end{subfigure}
 \caption{Case (II): Empirical estimation errors of the penalized MLE. Panels (a)--(c) are presented similarly to Figure \ref{fig:resultscasea1}.}  \label{fig:resultscase2convex}
\end{figure}

Under Case (I), we present the empirical estimation errors of the one-step estimator and the penalized MLE in Figures \ref{fig:resultscasea1} and \ref{fig:resultscasea2}, respectively. 
In each figure,
panel (a) suggests that $\operatorname{dist}^2(\hat Z, Z^\star)$ is inverse proportional to $T$ when $n$ is fixed. 
Furthermore, panel (c) shows that all the fitted slopes are close to $-1$. 
In addition, 
panel (b) shows that $\operatorname{dist}^2(\hat Z, Z^\star)$ does not change too much as $n$ increases. 
In conclusion,  
the numerical results are consistent with the oracle theoretical rate  $O(1/T)$ for the estimation error of $Z^{\star}$.   
Comparing  Figures \ref{fig:resultscasea1} and \ref{fig:resultscasea2}, we find that under the same $(n,T,k)$,  the one-step estimator can achieve a slightly smaller estimation error than the penalized MLE. This might be because biases are introduced with the convex relaxation. 
However, the penalized MLE could be more flexible and robust when the true $k$ is unknown in applications. 
Under Case (II), 
 we present the empirical estimation errors $\mathrm{dist}^2(\hat{Z},Z^{\star})$ of the one-step estimator and the penalized MLE in Figures \ref{fig:resultscase2nonconvex} and \ref{fig:resultscase2convex}, respectively.
We can see that the patterns are similar to those in Figures \ref{fig:resultscasea1}--\ref{fig:resultscasea2}. The results suggest the oracle theoretical rate $O(1/T)$ can hold under a variety of types of heterogeneity of  $\alpha^{\star}$.

\begin{remark}\label{rm:choiceoflambda}
When estimating latent vectors with the penalized MLE, 
we set the tuning parameter $\lambda_{n,T}= c_{\lambda}\sqrt{ 
 nT \hat{\mu}}$ with a fixed constant $c_{\lambda}=0.5$ and $\hat{\mu}=\sum_{t,i,j} A_{t,ij} /(n^2T)$ being the average over all entries. In general applications, 
practitioners may also choose $c_{\lambda}$  using the network  cross-validation \citep{chen2018network,li2020network} or based on their data interpretation. Notably, when focusing on 
the aggregation effect over $T$ networks, i.e., how the estimation error changes with respect to $T$,  
we find that this phenomenon is not sensitive to the choice of $c_{\lambda}$. See more discussions in Section \ref{sec:lambdatuning} of the Supplementary Material. 
\end{remark} 

\section{Analysis of New York Citi Bike Dataset}
\label{sec::bikedata}
\begin{figure}[!htbp]
	\centering
	\begin{subfigure}{0.32\linewidth}
		\centering
		\includegraphics[width=1\linewidth]{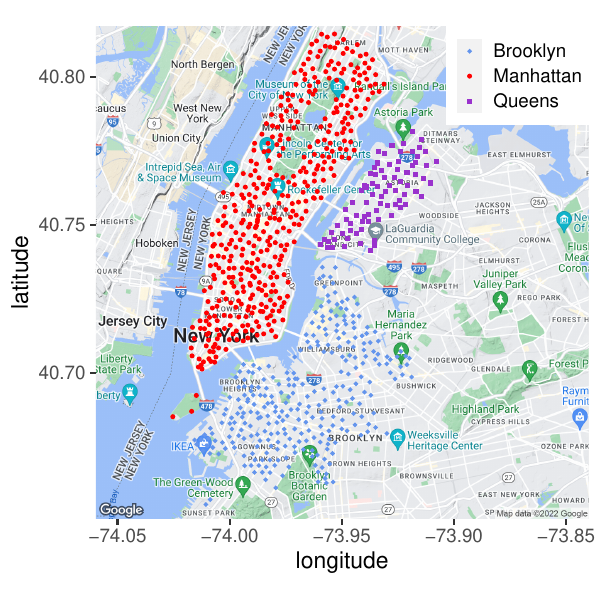}
    \caption{True Station Locations}
	\end{subfigure} \hspace{0.5em}
    \begin{subfigure}{0.32\linewidth}
		\centering
	    \includegraphics[width=1\linewidth]{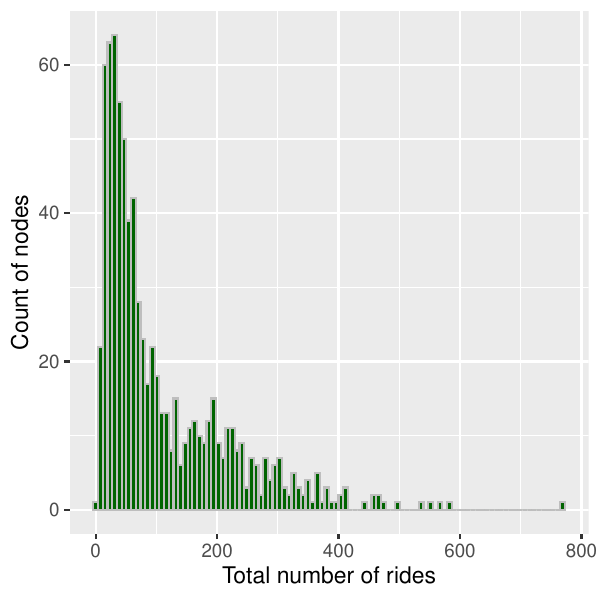}
     \caption{Node-Level Histogram}
	\end{subfigure} \hfill
	\begin{subfigure}{0.32\linewidth}
		\centering
		\includegraphics[width=1\linewidth]{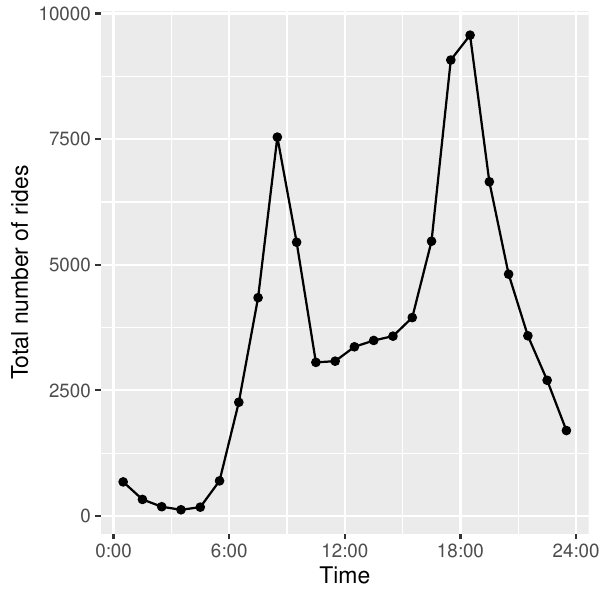}
    \caption{Ride Counts Over Time}
	\end{subfigure}
 \caption{Illustration of Citi Bike Data: Panel (a)  presents true locations of bike stations on Google Maps, colored by three boroughs of  New York City. Panel (b) presents the histogram of the total number of rides by nodes (bike stations). Panel (c) presents the total number of rides over the 24-hour period.} \label{fig:citidata}
\end{figure}

We illustrate the use of the proposed methods by analyzing the New York Citi Bike data \citep{bikedata}. 
The data set contains over $2.3$ million rides between bike stations in New York in August 2019. Each ride is identified by two stations and a time stamp (the hire starting time).  We focus on a weekday, August 1st, 2019, and keep the rides that last between one minute and $3$ hours. The processed data contains $85,854$ rides between $782$ stations over $24$ hours. 
The $782$ bike stations form a common set of nodes across hours.
A pair of stations defines a network edge, and the number of events between them in each hour gives the hourly edge weight. 
We provide exploratory data visualization in 
Figure \ref{fig:citidata}.
Panels (b) and (c) of Figure \ref{fig:citidata} 
show that the number of ride events is very heterogeneous across different bike stations and hours, respectively, motivating the application of the model \eqref{model_DLSM}. 

\begin{figure}[!htbp]

	\centering
 \hspace{0.2em}
\begin{subfigure}{0.225\textwidth}
		\centering
		\includegraphics[width=1\linewidth]{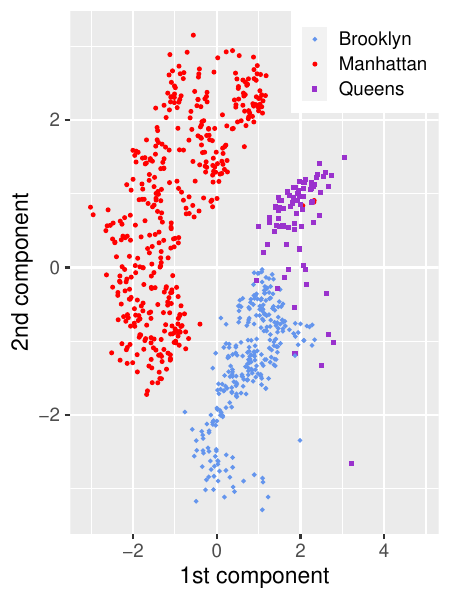}
  \caption{  24-Hour Estimates}\label{fig:24data}
\end{subfigure}
\begin{subfigure}{0.225\textwidth}
		\centering
		\includegraphics[width=1\linewidth]{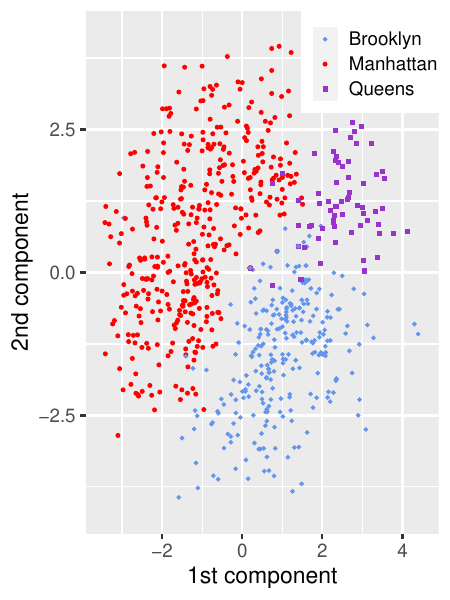}
    \caption{1-Hour Estimates }
	\end{subfigure}
 \begin{subfigure}{0.3\textwidth}
		\centering
		\includegraphics[width=1\linewidth]{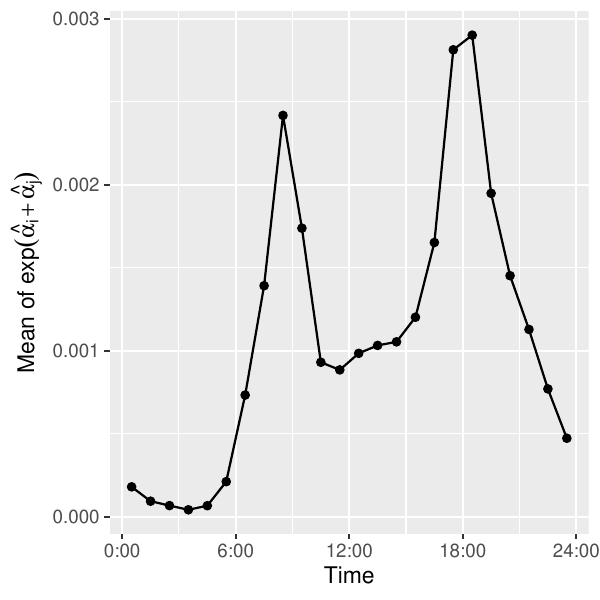}
    \caption{Baseline Estimates}
	\end{subfigure} 
 \begin{subfigure}{0.225\linewidth}
		\centering
		\includegraphics[width=1\textwidth]{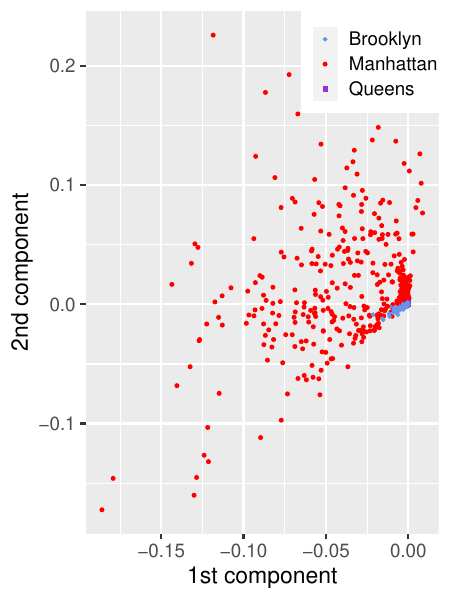}
  \caption{COSIE Estimates} \label{fig:cosie}
	\end{subfigure}
 \caption{Estimation Results: Panel (a) shows the estimated latent positions $\hat{z}_i\in \mathbb{R}^2$ using all the data across the 24-hour period. Each point is colored based on which borough the corresponding bike station is located, and we present the results after a rotation for better visualization.  Panel (b) shows the estimated latent positions $\hat{z}_i\in \mathbb{R}^2$ using one-hour data during 21:00-22:00 on the selected day, where points are similarly colored and rotated to those in the panel (a). Panel (c) presents the mean of estimated baseline levels $\sum_{i,j=1}^n\exp(\hat{\alpha}_{it}+\hat{\alpha}_{jt})/n^2$ for $t=1,\ldots, 24$, i.e., across the 24-hour period under the model \eqref{model_DLSM}. Panel (d) presents the estimated two-dimensional shared latent vectors under COSIE with 24-hour data. } \label{fig:estfig}
\end{figure}

We next fit the model \eqref{model_DLSM} with 2-dimensional latent vectors $z_i$'s. 
The choice of $k=2$ is for the ease of interpretation below and is consistent with the strategy in Remark \ref{rm:fromgtoz}. 
Figure \ref{fig:estfig} (a)  visualizes the estimated positions $\hat{z}_i$'s obtained by the one-step estimator, while the results of the penalized MLE are similar and deferred to the Supplementary Material.

To interpret the results, 
we compare the estimated latent space vectors $\hat{z}_i$'s, visualized in Figure \ref{fig:estfig} (a), with true geographic locations (in latitudes and longitudes) of stations, visualized in Figure  \ref{fig:citidata} (a). 
We can see that, overall, the estimated latent positions can match the true geographic positions. 
Particularly, the stations located in Manhattan are well-separated from those in Brooklyn and Queens.
This could be because 
{the borough of Manhattan is separated from the other two boroughs by the East River, and}
it is less common to use shared bikes to travel between Manhattan and the other two boroughs. 
On the other hand, the estimated station locations in Brooklyn and Queens overlap. This could be because Queens and Brooklyn are not physically divided by a river, and it is easier to travel between these two boroughs by biking. 
In the meantime, some bike stations were estimated to be in the Central Park area. 
This could be because people would ride through Central Park frequently, so two stations that are far away geographically might be viewed as close to each other based on the interactions.  
Besides the latent vectors, we visualize the estimated baseline levels across the 24-hour period in Figure  \ref{fig:estfig} (c). 
The results show a time-varying pattern that is consistent with the true changing pattern of the total number of rides across the 24-hour periods, visualized in Figure  \ref{fig:citidata}  (c). 
In summary, the fitted parameters $\hat{z}_i$'s and $\hat{\alpha}_{it}$'s under the model \eqref{model_DLSM} are reasonable and interpretable.  

In addition, comparing 
panels (a) and (b) in  Figure \ref{fig:estfig},
we can see that the alignment between the estimated latent space vectors and the true geographic locations improves as we utilize observations over more hours.
The results suggest that the proposed methods can effectively extract static latent space information and accommodate time heterogeneity simultaneously.

As a comparison, we fit the data by the common subspace independent edge (COSIE) model in \cite{arroyo2021inference} with two-dimensional shared latent vectors. 
We estimate the model 
by the multiple adjacency spectral embedding (MASE) method in \cite{arroyo2021inference}. We present the estimated latent vectors in Figure \ref{fig:estfig} (d). 
Compared to Figure \ref{fig:estfig} (a),  
Figure \ref{fig:estfig} (d) does not appear to show a direct connection between the estimated latent positions and the geographic positions of bike stations.  
This may be because  \cite{arroyo2021inference} targeted at $\mathbf{A}_t$'s that are adjacency matrices with binary edges.
We present the results by COSIE only to illustrate the differences between different mean models of $\mathbf{A}_t$'s.

\begin{remark}\label{rm:dataundirect}
The analyzed networks are undirected and mainly encode the total usage of bike rides between stations.
The learned embeddings $z_i$'s can be easy to interpret and useful for understanding the usage of shared bikes and for guiding the design of stations. Alternatively, 
one may examine bike rides with directional information, giving directed networks. Then the  model \eqref{model_DLSM} may be generalized as 
 $\exp(\alpha_{it} + \beta_{jt} + \langle z_{i}, w_{j}\rangle)$,
where for each vertex $i$, $\alpha_{it}$ and $\beta_{it}$  characterize the ``leave'' and ``arrival'' baseline activity levels at the time point $t$,
and $z_i$ and $w_i$ are time-invariant latent vectors characterizing the out-degree and in-degree properties, respectively.  
It would be an interesting future research direction to 
extend this study to directed networks.  
\end{remark}



\section{Discussion}
\label{sec::discussion}
In this work, we propose a longitudinal latent space model tailored for recurrent interaction events. We develop two novel semiparametric estimation techniques, i.e., the generalized semiparametric one-step updating and the penalized maximum likelihood estimation, and show that the resulting estimators attain the oracle estimation error rate for the shared latent structure. 
The first approach utilizes the semiparametric efficient score equation to construct a second-order updating estimator. 
We show that the estimator possesses a geometric interpretation on the quotient manifold, which automatically overcomes the non-uniqueness issue due to overparametrization.
The second approach corresponds to a convex relaxation of the low-rank static latent space component.

By separating the (primary) parameters of interest associated with the static latent space from the dynamic nuisance parameters, a strategy commonly found in semiparametrically efficient parameter estimation, we are able to delineate the oracle rates of convergence for the primary and the nuisance parameters according to their dimensions. This strategy also helps us to untangle the static and time-heterogeneous components inherent in the network model and construct the oracle estimators.

There are a few other interesting future works. 
First, the ability to accurately estimate latent structures could enable important downstream analysis such as prediction, hypothesis testing, and change-point detection. For instance, it may be useful to ascertain a change point in the structure of the latent space \citep{ bhattacharjee2020change,enikeeva2021change,padilla2022change}. The achievement of oracle estimation error rates could facilitate the quantification of uncertainty in estimators, which in turn lays a strong foundation for conducting reliable statistical inference.


Second, this paper focuses on the variance in estimation error rates as a function of $n$ and $T$, while treating the latent dimension $k$ and network sparsity level as fixed. Extending the current methodology and theory to the cases when $k$ grows \citep{choi2012stochastic} as well as sparse networks \citep{qin2013regularized,le2017concentration} are important topics.  

 Third, this work aims to unveil the fundamental relationship between the estimation errors and the degree of baseline heterogeneity. We focus on the most challenging scenario where the degree of baseline heterogeneity increases linearly with respect to $n$ and $T$. It is possible to impose additional structures to reduce baseline heterogeneity, such as assuming $\{\alpha_{i1},\ldots, \alpha_{iT}\}$ to be piecewise constants.  
Nevertheless, as $T$ increases, the intrinsic number of parameters would eventually become large to keep up with the increasing data complexity. The developed results would also provide us with techniques for investigating such scenarios. 
Which structural assumptions are appropriate may vary across different applications and require case-by-case analyses in future research.  

 Fourth, the proposed model has the potential for further extensions to capture more complex network structures.  
Currently, 
  heterogeneity across networks is only characterized at the first-order baseline levels $\alpha_{it}$'s, while the second-order interaction terms are modeled by the time-invariant components $z_{i}$'s.  
We find that this model adequately describes the analyzed dataset.
But more generally, it may also be of interest to incorporate time-varying interaction terms, which would further increase the model complexity and pose new theoretical challenges. Moreover, the proposed model adopts the Euclidean inner product to describe the interactions between nodes, which can be limited to capturing homophilic network structures. 
Recently, researchers proposed to use indefinite inner products to capture heterophilic structures  \citep{rubin2017statistical,10.1214/20-AOS1976,macdonald2020latent}. 
{We discuss the feasibility of generalizing  the proposed analysis to heterophilic products  
in Section \ref{sec:genhetero} of the Supplementary Material. 
While developing comprehensive results would require careful consideration of specific model properties, the foundational results and techniques developed under \eqref{model_DLSM} pave the way for studying more general models.} 



Finally, it would be worthwhile to generalize our results to other models, including 
various distributions for weighted edges, continuous time stamps, or additional covariates influencing the network structure \citep{hoff2002latent,vu2011continuous,perry2013point,hoff2015multilinear,kim2018dynamic, Sit2020event,weng2021community, huang2018pairwise}. We believe that the proposed semiparametric analysis framework can function as a valuable building block for establishing sharp estimation error rates under those models.

\begin{appendix}

\section{Details for the Initialization Algorithm}\label{sec:appendixalgo}

 We present the details of the initialization  Algorithm \ref{algor:init} below. 
More explanations on the details can be found in Section \ref{sec:techinitial} of the Supplementary Material.


\subsection*{Stage 1}
Stage 1 utilizes the universal singular value thresholding (USVT) method in \cite{chatterjee2015matrix} to obtain an ``initial of initial'' estimator. 
In particular, for $t=1,\ldots, T$, 
lines 3--4 of Algorithm \ref{algor:init} finds the singular value decomposition of $\mathbf{A}_t$ and retains the 
top singular components with the singular values  $s_{t,l} > \tau_t$. 
Line 5 applies $f_{M_{\Theta,2}}: \mathbb{R}^{n\times n}\to \mathbb{R}^{n\times n}$ to $\widetilde{E}_t$ to ensure that all elements in $\mathring{E}_t$ are in $[e^{-M_{\Theta,2}}, 1]$. 
Specifically, 
given a matrix $X=(X_{ij})_{1\leqslant i,j\leqslant n}\in \mathbb{R}^{n\times n}$, the $(i,j)$-th element of $f_{M_{\Theta,2}}(X)$ is defined as
\begin{align}\label{eq:fmtheta2x}
      (f_{M_{\Theta,2}}(X))_{ij} = \begin{cases}
      X_{ij} & \text{ if } e^{-M_{\Theta,2}}\leqslant X_{ij} \leqslant 1;\\
      e^{-M_{\Theta,2}} & \text{ if }  X_{ij}< e^{-M_{\Theta,2}};\\
      1 & \text{ if } X_{ij} >1.
      \end{cases}
  \end{align}
Line 10 projects $\widetilde{G}$ into $ \mathbb{S}_+^n$ and obtain $\mathring{G}$. 
In particular, 
 $\mathcal{P}_{\mathbb{S}_{+}^n}(\cdot)$ represents the projection mapping onto $\mathbb{S}_{+}^n$, i.e., for $ X\in \mathbb{R}^{n\times n},$ we define
\begin{align}\label{eq:projmap}
\mathcal{P}_{\mathbb{S}_{+}^n}(X) = \argmin_{Y\in \mathbb{S}_{+}^n} \|X-Y\|, 
\end{align}
where $\|\cdot\|$ represents the Euclidean norm.
In line 11, we define 
$U_k=(u_1,\ldots, u_k)\in \mathbb{R}^{n\times k}$ and $D_k=\mathrm{diag}(d_1,\ldots, d_k) \in \mathbb{R}^{k\times k}$,
where $d_1\geqslant \ldots \geqslant d_k$ are the largest $k$ eigenvalues of $\mathring{G}$,
and $u_j\in \mathbb{R}^n$ is the eigenvector of  $\mathring{G}$ corresponding to $d_j$ for $j=1,\ldots, k$.


\vspace{-0.5em}
\subsection*{Stage 2}
Stage 2 sets $(\mathring{Z},\mathring{\alpha})$ obtained in Stage 1 as the initial estimator (line 13) and consists of two iterative sub-stages.

\textit{Stage 2-1 (lines 14--18):}  For $r=1,\ldots, R_1-1$,
let  $(Z^r, \alpha^r)$ denote the updated parameters along gradients after the $r$-th iteration. 
Lines 15--16 of Algorithm \ref{algor:init}  updates the parameters $(Z^r,\alpha^r)$ along the  directions 
\begin{align}
    g_Z(Z^r,\alpha^r) = &~ \sum_{t=1}^{T} (\mathbf{A}_t-\exp(\Theta_t^r)) Z^r\ \in  \mathbb{R}^{n\times k},\label{eq:grandzformulaalo}\\ 
   g_{\alpha_t}(Z^r,\alpha^r) =&~(\mathbf{A}_t-\exp(\Theta_t^r)) 1_{n} \ \in  \mathbb{R}^{n\times 1}\label{eq:grandalphaformulaalo}
\end{align}
 with $ \Theta^{r}_t = \alpha_t^r 1_n^{\top} + 1_n (\alpha_{t}^{r})^{\top} +Z^r (Z^{r})^{\top}$ for $t=1,\ldots, T$. 
Lines 17--18  
projects the updated estimator $(\widetilde{Z}^{r+1},\widetilde{\alpha}^{r+1})$ onto the following  sets 
\begin{align}
        \mathcal{C}_{Z}=&~\{Z \in \mathbb{R}^{n \times k}: 1_n^{\mytrans} Z=0, \max_{1 \leqslant i \leqslant n}\left\|z_i\right\|^2_2 \leqslant M_{Z,1}\}, \label{eq:czsetproj}\\
        \mathcal{C}_{\alpha}=&~\{\alpha \in \mathbb{R}^{n\times T}:\max_{1\leqslant i\leqslant n,\, 1\leqslant t \leqslant T}|\alpha_{it}| \leqslant M_{\alpha}\}, \label{eq:calphasetproj} \\
    \mathcal{C}_{Z}^{\prime}=&~\{Z \in \mathbb{R}^{n \times k}: \max_{1\leqslant i\leqslant n}\left\|z_i\right\|^2_2 \leqslant M_{Z,1}\}. \label{eq:czprimesetproj}
\end{align}
Similarly to \eqref{eq:projmap},  for $\mathcal{C}\subseteq\mathbb{R}^{m\times q}$ and $X \in \mathbb{R}^{m\times q}$,  let $\mathcal{P}_{\mathcal{C}}(X) = \argmin_{Y\in \mathcal{C}}\|Y-X\|$.
The projection mapping  is uniquely defined  
as $\mathcal{C}_{Z}$,  $ \mathcal{C}_{\alpha}$, and $   \mathcal{C}_{Z}^{\prime}$ are closed and convex sets \citep[Section 8.1,][]{boyd2004convex} .

\textit{Stage 2-2 (lines 20--25):} 
The second part of Stage 2 utilizes an  alternating version of the projected gradient descent to obtain the error bound of 
individual rows of $Z$ and individual elements in $\alpha$. 
In particular, 
in lines 21--22, we define 
\begin{align}
    g_{Z,R_1}(Z^r,\alpha^r) = &~ \sum_{t=1}^{T} (\mathbf{A}_t-\exp(\Theta_t^{r,R_1})) Z^{R_1}\ \in  \mathbb{R}^{n\times k},\label{eq:grandzformulaalo2}\\ 
   g_{\alpha_t,R_1}(Z^r,\alpha^r) =&~(\mathbf{A}_t-\exp(\Theta_t^{r,R_1})) 1_{n} \ \in  \mathbb{R}^{n\times 1}\label{eq:grandalphaformulaalo2}
\end{align}
where $ \Theta^{r,R_1}_t = \alpha_t^r 1_n^{\top} + 1_n (\alpha_{t}^{R_1})^{\top} +Z^r (Z^{R_1})^{\top}$ with $\alpha_t^{R_1}$ and $Z^{R_1}=(z_1^{R_1},\ldots, z_n^{R_1})^{\top}$ being estimates from \textit{Stage 2-1}.  Compared to $\Theta_t^r$ in \eqref{eq:grandzformulaalo}--\eqref{eq:grandalphaformulaalo}, $\Theta_t^{r,R_1}$ is defined with part of the parameters fixed at $\alpha_t^{R_1}$ and $Z^{R_1}$.

The final output of Algorithm \ref{algor:init} is  $\check{\alpha}=\alpha^{R_2}$ and $\check{Z}= {JZ^{R_2}}$ with 
  \begin{align}\label{eq:jmatrix}
    J = \mathrm{I}_n - \frac{1}{n} 1_n 1_n^{\top}.
\end{align}
We call $J$ the centering matrix as  $1_n^{\top}J=0$, which gives 
 $1_n^{\top}\check{Z}=0$. 

\begin{condition} \label{cond:tuning}
Assume the tuning parameters in Algorithm \ref{algor:init} satisfy the following: 
\begin{enumerate}
\setlength{\itemsep}{0pt} 
\item[(i)] Thresholds
$\tau_t = \tau n^{1/2}$ for some large enough constant $\tau > 0$.
\item[(ii)]
Step sizes $\eta_{Z} = \eta / (T{\|Z^0\|^2_{\operatorname{op}}})$ and $\eta_{\alpha} = \eta / (2n)$ for some small constant $\eta >0$. 
\item[(iii)]
Numbers of iterations $R_1$ and $R_2$ are sufficiently large.
 \end{enumerate}
\end{condition}

\begin{remark} \label{rm:e3}
 In simulations, we choose parameters following Remark \ref{rm:e2} in the Supplementary Material. 
Condition \ref{cond:tuning} on parameters are imposed for the convenience of proofs. 
It essentially requires that the tuning parameters satisfy certain orders with respect to $n$ and $T$ and can be approximated by the practical choices in Remark \ref{rm:e2}.
\end{remark} 

\normalem 
\begin{algorithm}[!htbp]
	\caption{Initialization algorithm for the semiparametric one-step estimator.}
	\label{algor:init}
	\KwIn{Data: $\mathbf A_1, \ldots, \mathbf A_T \in \mathbb{R}^{n\times n}$; thresholds: $\tau_1,\ldots, \tau_T$; latent space dimension: $k$; step sizes: $\eta_{Z}$, $\eta_{\alpha}$; 
    number of iterations: $R_1, R_2$.}
	\KwOut{$(\check{\alpha},\check{Z}) 
    .$}
	\BlankLine
 \BlankLine
 \textbf{Stage 1:}
 
    \For{$t = 1, \ldots, T$}{
     Find the singular value decomposition of $\mathbf A_{t}=\sum_{l=1}^{n} s_{t,l} u_{t,l}v_{t,l}^{\top}$,  where for $l=1,\ldots, n$, $ s_{t,l} \in \mathbb{R}$, $ u_{t,l}\in \mathbb{R}^n$, and $v_{t,l}\in \mathbb{R}^n$ represent the singular values, the left singular vectors, and the right singular vectors of $\mathbf{A}_t$, respectively. 
    
	Let $\widetilde{E}_{t}=\sum_{\{l:\, s_{t,l} > \tau_{t}\}} s_{t,l} u_{t,l} v_{t,l}^{\top}$. 

  Let $\mathring{E}_{t}=f_{M_{\Theta,2}}(\widetilde{E}_{t})$, where $f_{M_{\Theta,2}}(x)$ is defined in \eqref{eq:fmtheta2x}. 
  

  
   Let $\mathring{\Theta}_{t} = \log(\mathring{E}_{t})$.
	
		Let $\mathring{\alpha}_{t} =\arg \min _{\alpha}\|\mathring{\Theta}_{t}-\alpha 1_{n}^{\top}-1_{n} \alpha^{\top}\|_{\mathrm{F}}^{2} = (n \mathrm{I}_n + 1_n1_n^{\top})^{-1} \mathring{\Theta}_{t} 1_n$.
	}

Let $\widetilde{G} = \sum_{t = 1}^{T} (\mathring{\Theta}_{t}-\mathring{\alpha}_{t}1_{n}^{\top}-1_{n} \mathring{\alpha}_{t}^{\top})/T$. 
 

 	Let $\mathring{G} = \mathcal{P}_{\mathbb{S}_{+}^{n}}(\widetilde{G})$, where  $\mathcal{P}_{\mathbb{S}_+^n}(\cdot)$ is defined in \eqref{eq:projmap}.

	Let $\mathring{Z}=\mathring{U}_{k} \mathring{D}_{k}^{1 / 2}$ where $\mathring{U}_{k} \mathring{D}_{k} \mathring{U}_{k}^{\top}$ is the top-$k$ eigen components of $\mathring{G}$.

     	\BlankLine
      \BlankLine
 \textbf{Stage 2:}

    Let $Z^0 = \mathring{Z}$ \text{ and } $\alpha^0 = \mathring{\alpha} = (\mathring{\alpha}_1, \ldots, \mathring{\alpha}_T)$.

	\For{$r = 0, 1, \ldots, R_1 - 1$}{
			$\widetilde{Z}^{r+1}=Z^{r}+\eta_{Z}    g_Z(Z^r,\alpha^r)$ with $ g_Z(Z^r,\alpha^r)$  in \eqref{eq:grandzformulaalo}.
            
            $\widetilde{\alpha}^{r+1}_{t} = \alpha^{r}_{t} + {\eta_{\alpha}}g_{\alpha_t}(Z^r,\alpha^r)$ with $g_{\alpha_t}(Z^r,\alpha^r)$  in \eqref{eq:grandalphaformulaalo} for $t=1,\ldots, T$.

			 $Z^{r+1}=\mathcal{P}_{\mathcal{C}_{Z}}(\widetilde{Z}^{r+1})$ with $\mathcal{C}_Z$ defined in \eqref{eq:czsetproj}.

            $\alpha^{r+1}=\mathcal{P}_{\mathcal{C}_{\alpha}}(\widetilde{\alpha}^{r+1})$ with $\mathcal{C}_{\alpha}$  in \eqref{eq:calphasetproj}. 
		}

	    \For{$r = R_1, \ldots, R_2 - 1$}{
   			$\widetilde{Z}^{r+1}=Z^{r} +  \eta_{Z}  g_{Z,R_1}(Z^r,\alpha^r)$ with  $  g_{Z,R_1}(Z^r,\alpha^r)$  in \eqref{eq:grandzformulaalo2}. 
   
            $\widetilde{\alpha}^{r+1}_{t}  = \alpha^{r}_{t} +  \eta_{\alpha}  g_{\alpha_t,R_1}(Z^r,\alpha^r) $ with $g_{\alpha_t,R_1}(Z^r,\alpha^r)$  in \eqref{eq:grandalphaformulaalo2} for $ t = 1, \ldots, T$.

            $Z^{r+1}=\mathcal{P}_{\mathcal{C}_{Z}^{\prime}}(\widetilde{Z}^{r+1})$ with  $\mathcal{C}_Z^{\prime}$ defined in \eqref{eq:czprimesetproj}. 
            
            $\alpha^{r+1}=\mathcal{P}_{\mathcal{C}_{\alpha}}(\widetilde{\alpha}^{r+1})$ with $\mathcal{C}_{\alpha}$ defined in \eqref{eq:calphasetproj}. 
		}

\BlankLine
  Let $\check{\alpha}=\alpha^{R_2}$ and $\check{Z}= {JZ^{R_2}},$ where  $J$ is defined in \eqref{eq:jmatrix}. 
\end{algorithm}
\ULforem

\end{appendix}

\begin{acks}[Acknowledgments]
We are grateful to 
the editor Dr. Enno Mammen,  
an associate editor, and two anonymous referees for their helpful comments and suggestions.  
Ying was partially supported by NSF Grant DMS-2015417. 
Feng was partially supported by NIH grant 1R21AG074205-01, NSF Grant DMS-2324489, NYU University Research Challenge Fund, and a grant from NYU School of Global Public Health.
\end{acks}



\begin{supplement}
\stitle{Supplement to ``Semiparametric Modeling and Analysis for Longitudinal Network Data.''} 
\sdescription{Due to space limitation, additional results and proofs are deferred to the supplementary material.  The codes are published at \cite{dlsm}.} 

\end{supplement}

\bibliographystyle{ims}
\bibliography{Citation.bib}

\numberwithin{equation}{section}

\renewcommand{\thetable}{S\arabic{table}}
\renewcommand{\thefigure}{S\arabic{figure}}

\newcommand{\DZt}{\mathcal{D}_{Z\Theta}}
\newcommand{\ZQssub}{{Z}^{\star}_q}
\newcommand{\myDonet}{\mathcal{D}_{\alpha\Theta}}
\newcommand{\bmM}{\bar{\mathcal{M}}_k}

\newcommand{\bdmu}{\boldsymbol{\mu} }

\newcommand{\Dmut}{\mathcal{D}_{\mu_t}}
\newcommand{\Dmutza}{\mathcal{D}_{\mu_t}(Z,\alpha)}




\newpage

\begin{center}
\Large 
    \textbf{Supplementary Material for ``Semiparametric Modeling and Analysis for Longitudinal Network Data.''}
\end{center}



\appendix

\setcounter{section}{+1}


\section{Preliminaries}
\label{sec:technical_prelim}

In this section, we provide preliminaries for developing the theoretical results. 
Section \ref{subsec:formula_efficient} provides analytical expressions of efficient score $\Seff(Z,\alpha)$ and efficient information matrix $ \Ieff(Z, \alpha)$ in \eqref{eq:ieffdef}.

\subsection{Efficient Score and Efficient Information Matrix for \texorpdfstring{$Z$}{}} 
\label{subsec:formula_efficient}

In this section, we derive the explicit formulae of the efficient score and efficient information matrix. 
Recall that the log-likelihood function of $(Z,\alpha)$
is $L(Z,\alpha)=\sum_{t=1}^T L_t(Z,\alpha_t)$, 
where 
\begin{align*}
L_t(Z,\alpha_t) =& ~\sum_{1\leqslant i\leqslant j\leqslant n} \left[A_{t,ij} (\alpha_{it}+\alpha_{jt}+z_i^\top z_j) - \exp(\alpha_{it}+\alpha_{jt}+z_i^\top z_j) \right].
\end{align*}

\subsubsection{Vectorization of variables}
We consider the vectorization of $Z, \alpha, \Theta_t$ defined as
\begin{align}\label{eq:vecdef}
  Z_v=\begin{pmatrix}
      z_1\\
      z_2\\
      \vdots\\
      z_{n-1}\\
      z_n
  \end{pmatrix}   \in \mathbb{R}^{nk\times 1}, \  \alpha_v = 
  \begin{pmatrix}
      \alpha_1\\
      \alpha_2\\
      \vdots\\
      \alpha_{T-1}\\
      \alpha_T
  \end{pmatrix}\in \mathbb{R}^{nT\times 1},\  (\Theta_t)_v=
  \begin{pmatrix}
     \Theta_{t,11}/\sqrt{2}\\[2pt]
     \vdots\\
     \Theta_{t,nn}/\sqrt{2}\\[2pt]
     \Theta_{t,12}\\[2pt]
     \Theta_{t,13}\\[2pt]
     \vdots\\
     \Theta_{t,n-1,n}
  \end{pmatrix}\in \mathbb R^{\frac{n(n+1)}{2} \times 1},
\end{align}
respectively. 
With a little abuse of notation, we let $\dot{L}_Z(Z,\alpha),\dot{L}_{\alpha}(Z,\alpha)$ and $\dot{L}_{\Theta_t}(Z,\alpha)$ denote the partial derivatives of $L(Z,\alpha)$   
with respect to vectors $Z_v$, $ {\alpha}_v$, and $(\Theta_t)_v$, respectively (see precise formulae in Section \ref{sec:derieffscore}). 
Also,  when there is no ambiguity, we also write  $\dot{L}_Z$, $\dot{L}_{\alpha}$, and $\dot{L}_{\Theta_t}$ for the simplicity of notation.
Similarly,
when considering partial derivatives (for both first and higher orders) of or with respect to $Z,\alpha,
\Theta_t$, such  vectorization is also implicitly applied unless otherwise specified. 
 
 \begin{remark}
 \label{rmk::vectorization_sqrt2}
 In the definition of $(\Theta_t)_v$, the diagonal elements $\Theta_{t,jj}$ in $\Theta_t$ is scaled by $1/\sqrt{2}$. 
The scaling is considered as it preserves the Frobenius norm, namely, $\|\Theta_t\|_{\mathrm{F}}^2 = 2 \|(\Theta_t)_v\|_2^2 $.
 \end{remark}

We next derive the formulae of the efficient score function $\Seff(Z,\alpha)$, the efficient information matrix $\Ieff(Z,\alpha)$, and the observed efficient information matrix $-\Heff(Z,\alpha)$ in Sections \ref{sec:derieffscore}--\ref{sec:obseffinf}, respectively. 

\subsubsection{Formula of the efficient score function $\Seff$} \label{sec:derieffscore}

Recall the efficient score function of the target parameters $Z$ is defined as 
\begin{align}
	\label{eq:def_eff_score}
    \Seff(Z,\alpha) =&~\dot{L}_Z - \mathbb{E}\big(\dot{L}_{Z}\, \dot{L}_{\alpha}^{\top}\big)\times \big\{\mathbb{E}\big(\dot{L}_{\alpha}\, \dot{L}_{\alpha}^{\top}\big)\big\}^{-1}\times \dot{L}_{\alpha} \ \in \  \mathbb{R}^{nk\times 1},
\end{align}
 where the expectations are taken over $\mathbf A$ under the proposed model with parameter values $(Z, \alpha)$. 
We next derive \eqref{eq:def_eff_score} by first deriving the first-order derivatives $\dot{L}_Z$ and $\dot{L}_{\alpha}$, and second, deriving the expectations  $\mathbb{E}\big(\dot{L}_{Z}\, \dot{L}_{\alpha}^{\top}\big)$ and $\mathbb{E}\big(\dot{L}_{\alpha}\, \dot{L}_{\alpha}^{\top}\big)$.

\paragraph*{Step 1: First-order derivatives.}
We derive partial derivatives  $\dot{L}_Z$ and $\dot{L}_{\alpha}$ with $(\Theta_t)_v$ as an intermediate variable. 
Specifically, by the chain rule of the first-order derivatives, 
\begin{align}
\dot{L}_Z = & ~\frac{\partial L}{\partial Z_v} =\sum_{t=1}^T \frac{\partial  (\Theta_t)_v^{\top}}{\partial Z_v}\frac{\partial L}{\partial (\Theta_t)_v} \in \mathbb{R}^{nk\times 1},\label{eq:dotL_Z1}\\[2pt]
\dot{L}_{\alpha}=&~
\begin{pmatrix}
	\dot{L}_{\alpha_1}\\
	\vdots\\
	\dot{L}_{\alpha_T}
\end{pmatrix}\in \mathbb{R}^{nT\times 1}, \quad \text{ with } \ \dot{L}_{\alpha_t} =\frac{\partial L}{\partial \alpha_t} = \frac{\partial  (\Theta_t)_v^{\top}}{\partial \alpha_t}\frac{\partial L}{\partial (\Theta_t)_v}\in \mathbb{R}^{n\times 1}. \label{eq:dotL_alpha1}
\end{align}
In particular, 
\begin{align}\label{eq:ldot_theta_t}
	\dot{L}_{\Theta_t}:=\frac{\partial L}{\partial (\Theta_{t})_v} =
\begin{pmatrix}
\frac{\partial L_t}{\partial (\Theta_{t,11}/\sqrt{2})}	\\[2pt]
\vdots\\[2pt]
\frac{\partial L_t}{\partial (\Theta_{t,nn}/\sqrt{2})} \\[8pt] 
\frac{\partial L_t}{\partial \Theta_{t,12}}	\\[2pt]
\vdots\\[2pt]
\frac{\partial L_t}{\partial \Theta_{n-1,n} } 
\end{pmatrix}=\begin{pmatrix}
\sqrt 2(A_{t,11}-\mu_{t,11})\\[4pt] \vdots \\[4pt] \sqrt2(A_{t,nn}-\mu_{t,nn})\\[4pt] A_{t,12} - \mu_{t,12}\\[4pt] \vdots \\[4pt] A_{t,n-1,n} - \mu_{t,n-1,n} 
\end{pmatrix} \in \mathbb{R}^{\frac{n(n+1)}{2}\times 1},
\end{align}
where we define
$$
\mu_{t,ij}=\exp(\Theta_{t,ij})=\exp\big(\alpha_{it}+\alpha_{jt}+\langle z_i, z_j\rangle \big)
$$ 
for $1\leqslant i,j\leqslant n$.
Moreover, the partial derivatives of $(\Theta_t)_v$ with respect to $Z_v$ and $\alpha_t$ are
\begin{align}
\frac{\partial (\Theta_t)_v^{\top}}{\partial Z_v}=&
		\begin{pmatrix}
	\frac{\partial (\Theta_{t,11}/\sqrt{2})}{\partial z_1} &  \cdots  & \frac{\partial (\Theta_{t,nn}/\sqrt{2})}{\partial z_1}  & \frac{\partial \Theta_{t,12}}{\partial z_1} & \cdots & \frac{\partial \Theta_{t,n-1,n}}{\partial z_1} \\
		\vdots &  \cdots  & \vdots  & \vdots & \cdots & \vdots \\
	\frac{\partial (\Theta_{t,11}/\sqrt{2})}{\partial z_n} & \cdots  & \frac{\partial (\Theta_{t,nn}/\sqrt{2})}{\partial z_n}  & \frac{\partial \Theta_{t,12}}{\partial z_n} & \cdots & \frac{\partial \Theta_{t,n-1,n}}{\partial z_n} 
	\end{pmatrix} \notag\\[5pt]
	=&
\begin{array}{c@{}c}
\left(
  \begin{BMAT}[2pt]{c:c:c}{c}
    \begin{BMAT}[3pt]{ccc}{cccccc}
   \sqrt{2}z_1   &\cdots  & 0\\
     0 &\cdots & 0\\
    0   &\cdots & 0 \\
    \vdots   & & \vdots  \\
    0  & & 0\\
     0 &\cdots &\sqrt{2}z_n
    \end{BMAT}
    & 
       \begin{BMAT}[3.1pt]{cccc}{cccccc}
  z_2   &z_3 &\cdots  & z_n\\
   z_1 &0 & \cdots & 0\\
    0   &z_1 & \cdots & 0 \\
    \vdots   & & &  \vdots  \\
    0  & 0 & \cdots &  0\\
     0 &0 & \cdots &z_1
    \end{BMAT}    
    &
  \begin{BMAT}[3pt]{ccc}{cccccc}
\ 0   & \cdots  & 0\\
\ z_3 & \cdots & 0\\
\ z_2 &  \cdots & 0 \\
    \vdots   & & \vdots  \\
  \  0  & & z_n\\
  \   0 &\cdots &z_{n-1}
    \end{BMAT}
  \end{BMAT} 
\right)  \\[-1ex]
\ \hexbrace{2.7cm}{n}\hexbrace{2.5cm}{n-1}\hexbrace{2.5cm}{\frac{n(n+1)}{2}-(2n-1)}
\end{array} \in \mathbb{R}^{nk\times \frac{n(n+1)}{2}},\label{eq:partialthetapartialzmatrix}
\end{align}
\begin{align}\label{eq:dthetaalphadefexpre}	
 \frac{\partial  (\Theta_t)_v^{\top}}{\partial \alpha_t} 
	=
\begin{array}{c@{}c}
\left(
  \begin{BMAT}[2pt]{c:c:c}{c}
    \begin{BMAT}[3.7pt]{ccc}{cccccc}
   \sqrt{2}   &\cdots  & 0\\
     0 &\cdots & 0\\
    0   &\cdots & 0 \\
    \vdots   & & \vdots  \\
    0  & & 0\\
     0 &\cdots &\sqrt{2}
    \end{BMAT}
    & 
       \begin{BMAT}[4.1pt]{cccc}{cccccc}
\  1   &1 &\cdots  & 1\\
\   1 &0 & \cdots & 0\\
\    0   &1 & \cdots & 0 \\
\    \vdots   & & &  \vdots  \\
\    0  & 0 & \cdots &  0\\
\     0 &0 & \cdots &1
    \end{BMAT}    
    &
  \begin{BMAT}[4.2pt]{ccc}{cccccc}
\ 0   & \cdots  & 0\\
\ 1 & \cdots & 0\\
\ 1 &  \cdots & 0 \\
    \vdots   & & \vdots  \\
  \  0  & & 1\\
  \   0 &\cdots & 1
    \end{BMAT}
  \end{BMAT} 
\right)  \\[-1ex]
\ \ \hexbrace{2.4cm}{n}\hexbrace{2.3cm}{n-1}\hexbrace{2.1cm}{\frac{n(n+1)}{2}-(2n-1)}
\end{array} \in \mathbb{R}^{n\times \frac{n(n+1)}{2}}. 
\end{align}

\begin{remark}
   In \eqref{eq:partialthetapartialzmatrix},  the definition of  ${\partial (\Theta_t)_v^{\top}}/{\partial Z_v}$ implicitly restricts   that no orthogonal transformation will be applied to $Z\in \mathbb{R}^{n\times k}$. Despite that,  the derived efficient information matrix of $Z$ still characterizes the quotient manifold induced by $Z$, as  suggested by Section  \ref{subsec:manifold}. Also, see more discussions on the structures of the efficient information matrix  of $Z$ in Section \ref{sec:ieffcombine}.  
\end{remark}

We point out the following facts that will be used in the following derivations. 
\begin{fact}\label{fact:partialderivfunction}\ 
\begin{enumerate}
    \item[(i)] $\displaystyle \frac{\partial (\Theta_t)_v^{\top}}{\partial Z_v}$ only depends on $Z$ and   does not depend on $\alpha$ and the observations $\mathbf{A}$. 
    \item[(ii)]$\displaystyle 
 \frac{\partial (\Theta_t)_v^{\top}}{\partial \alpha_t}$ only has constant entries and does not depend on $(Z,\alpha)$ and  $\mathbf{A}$. 
 \item[(iii)] Formulae of these two partial derivatives do not change with respect to  $t=1,\ldots, T$. 
\end{enumerate}
\end{fact}
\begin{notation}\label{notation:partialderiv}
As the above partial derivatives appear repeatedly in the theoretical derivations, to facilitate the presentation, we define the following notations: 
\begin{align*}
\DZt := & ~\frac{\partial \Theta_t}{\partial Z}^{\top}:= 	\frac{\partial (\Theta_t)_v^{\top}}{\partial Z_v}   \ \in \ \mathbb{R}^{nk\times \frac{n(n+1)}{2}},\\
 \myDonet:=&~\frac{\partial \Theta_t}{\partial \alpha_t}^{\top}:=\frac{\partial  (\Theta_t)_v^{\top}}{\partial \alpha_t}\ \in \ \mathbb{R}^{n\times \frac{n(n+1)}{2}}. 
\end{align*}
We point out that $\DZt$ and $ \myDonet$ are the same across $t=1,\ldots, T$, 
which implicitly uses the Fact \ref{fact:partialderivfunction} (iii).
Similarly, we denote the transposes of the above matrices as
\begin{align}
   \DZ :=&~ \frac{\partial \Theta_t}{\partial Z}:= 	\frac{\partial (\Theta_t)_v}{\partial Z_v^{\top}} \ \in \ \mathbb{R}^{ \frac{n(n+1)}{2} \times nk } , \label{eq:dthetazdef}\\
\myDone :=&~ \frac{\partial \Theta_t}{\partial \alpha_t}:= \frac{\partial (\Theta_t)_v}{\partial \alpha_t^{\top}}\ \in\ \mathbb{R}^{ \frac{n(n+1)}{2} \times n }. \label{eq:dthetaalphadef}
\end{align}
We give two remarks on the use of the notation. 
\begin{itemize}
    \item[(i)] Because  $\DZt$ and $\DZ$ are functions of $Z$ by the Fact \ref{fact:partialderivfunction} (i), we also write them as $\DZt(Z)$ and $\DZ(Z)$ when it is necessary to emphasize their dependence on $Z$. 
    \item[(ii)] On the other hand, $\myDonet$ and $\myDone $ would not depend on $(Z,\alpha)$. 
\end{itemize}
\end{notation}
In summary, by \eqref{eq:dotL_Z1}--\eqref{eq:dotL_alpha1}, the Fact \ref{fact:partialderivfunction}, and the Notation  \ref{notation:partialderiv}, we can write
\begin{align}
\dot{L}_Z(Z,\alpha) =& ~
\DZt(Z) \sum_t\dot{L}_{\Theta_t}(Z,\alpha)\  \in \  \mathbb{R}^{nk\times 1},\label{eq:dotL_Z} \\
\dot{L}_{\alpha_t}(Z,\alpha) =&~ 
 \myDonet \dot{L}_{\Theta_t}(Z,\alpha)\ \in \ \mathbb{R}^{n\times 1}.\label{eq:dotL_alpha_t}
\end{align}
\paragraph*{Step 2: Expectation terms.}
We next derive the expectations  $\mathbb{E}\big(\dot{L}_{Z}\, \dot{L}_{\alpha}^{\top}\big)$ and $\mathbb{E}\big(\dot{L}_{\alpha}\, \dot{L}_{\alpha}^{\top}\big)$ based on  \eqref{eq:dotL_Z}--\eqref{eq:dotL_alpha_t}.  In particular, 
\begin{align*}
   \mathbb{E}\big(\dot{L}_{Z}\, \dot{L}_{\alpha}^{\top}\big) = 
   \begin{pmatrix}
        \mathbb{E}\big(\dot{L}_{Z}  \dot{L}_{\alpha_1}^{\top}\big) &   \mathbb{E}\big(\dot{L}_{Z}\, \dot{L}_{\alpha_2}^{\top}\big) &  
 \cdots &  \mathbb{E}\big( \dot{L}_{Z}\, \dot{L}_{\alpha_T}^{\top}\big)
   \end{pmatrix}  \  \in \ \mathbb{R}^{nk \times nT}
\end{align*}
where for $t=1,\ldots, T$, 
\begin{align*}
   \mathbb{E}\big(\dot{L}_{Z}\, \dot{L}_{\alpha_t}^{\top}\big) =  \sum_{t_1, t_2=1}^T\mathbb E\left( \frac{\partial \Theta_{t_1}}{\partial Z}^\top \dot{L}_{\Theta_{t_1}}\dot{L}_{\Theta_{t_2}}^\top \frac{\partial \Theta_{t_2}}{\partial \alpha_t} \right)  =  \sum_{t_1, t_2=1}^T\DZt \mathbb E\left( \dot{L}_{\Theta_{t_1}}\dot{L}_{\Theta_{t_2}}^\top \right)\myDone   \ \in \ \mathbb R^{nk \times n}. 
\end{align*}
where the second equation follows by the Fact \ref{fact:partialderivfunction} and the Notation \ref{notation:partialderiv}. 
Moreover, by plugging \eqref{eq:ldot_theta_t}, we have for $ t_1\neq t_2 \in \{1,\ldots, T\}$, 
\begin{align*}
   \mathbb{E}\left(\dot{L}_{\Theta_{t_1}} \dot{L}_{\Theta_{t_2}}^\top \right) = 0 \ \in \ \mathbb{R}^{\frac{n(n+1)}2\times \frac{n(n+1)}2},  
\end{align*}
and for $t \in \{1,\ldots, T\}$, 
\begin{align*}
   \mathbb{E}\left(\dot{L}_{\Theta_t} \dot{L}_{\Theta_t}^\top\right) = \Dmut \ \in \ \mathbb{R}^{\frac{n(n+1)}2\times \frac{n(n+1)}2}, 
\end{align*}
where we define 
\begin{align}\Dmut = \text{diag}(2\mu_{t,11},\ldots,2\mu_{t,nn},\mu_{t,12},\ldots,\mu_{t,n-1,n} )   \ \in \ \mathbb{R}^{\frac{n(n+1)}2\times \frac{n(n+1)}2},\label{Dmutnotation}
\end{align} which is a diagonal matrix.
In summary, we can write
\begin{align}
\mathbb{E}\left(\dot{L}_{Z}\dot{L}_{\alpha}^\top\right) =&~ 
\begin{pmatrix}
 \DZt\mathcal{D}_{\mu_1}\myDone & 
\cdots & 
  \DZt\mathcal{D}_{\mu_T}\myDone
\end{pmatrix} 
\ \in \ \mathbb R^{nk \times nT},
\label{eq:EdotL_alphadotL_Ztop} 
\end{align}

Following a similar analysis, we have  
\begin{align} \label{eq:EdotL_alphadotL_alphatop} 
\mathbb{E}\left(\dot{L}_{\alpha}\dot{L}_{\alpha}^\top\right) =&~\text{diag}\left( \mathbb{E}\left(\dot{L}_{\alpha_t}\dot{L}_{\alpha_t}^\top\right)_{t=1, \ldots,T} \right) =\text{diag}\Big\{(\myDonet \Dmut \myDone)_{t=1, \ldots,T}\Big\} 
 \ \in \ \mathbb R^{nT \times nT},
 \end{align}
 which is a block-diagonal matrix  with $T$ blocks of $n\times n$ submatrices. 
\begin{remark}[Invertibility of $\mathbb{E}(\dot{L}_{\alpha}\dot{L}_{\alpha}^\top)$]
\label{rmk::invertable_ddotLalphaalpha}
Note that \eqref{eq:def_eff_score} involves the inverse of  $\mathbb{E}(\dot{L}_{\alpha}\dot{L}_{\alpha}^\top)$,
which is indeed invertible by \eqref{eq:EdotL_alphadotL_alphatop},  
Lemma \ref{lem:thm1.1_Hillar} and the explicit formula
\begin{align*}
\myDonet \Dmut \myDone = 
\begin{pmatrix}
4\mu_{t,11} + \sum_{j\neq 1} \mu_{t,1j} &\mu_{t,12} &\cdots &\mu_{t,1n}\\
\mu_{t,21} &4\mu_{t,22} + \sum_{j\neq 2} \mu_{t,2j} &\cdots &\mu_{t,2n}\\
\vdots &\vdots &\ddots & \vdots\\
\mu_{t,n1} &\mu_{t,n2} &\cdots & 4\mu_{t,nn} + \sum_{j\neq n} \mu_{t,nj} 
\end{pmatrix} \in \mathbb R^{n\times n}.
\end{align*}
\end{remark}

 Combining \eqref{eq:def_eff_score}, \eqref{eq:dotL_Z}, \eqref{eq:dotL_alpha_t}, \eqref{eq:EdotL_alphadotL_Ztop}, and \eqref{eq:EdotL_alphadotL_alphatop}, we can write
 \begin{align}\label{eq:seffformulaesupp}
     \Seff = \sum_{t=1}^T \left\{ \DZt \dot{L}_{\Theta_t} -  \DZt\Dmut\myDone (\myDonet \Dmut \myDone)^{-1} \myDonet \dot{L}_{\Theta_t}\right\}.
 \end{align}

\subsubsection{Formula of the efficient information matrix $\Ieff$}\label{sec:ieffformulae}
Now we calculate the efficient information matrix $ \Ieff=\mathbb{E}(\Seff \Seff^{\top})$. By  $\eqref{eq:def_eff_score}$, 
\begin{align*}
   \Ieff =&~\mathbb{E} \left[\dot{L}_Z - \mathbb{E}\big(\dot{L}_{Z}\, \dot{L}_{\alpha}^{\top}\big) \big\{\mathbb{E}\big(\dot{L}_{\alpha}\, \dot{L}_{\alpha}^{\top}\big)\big\}^{-1} \dot{L}_{\alpha}\right]\left[\dot{L}_Z - \mathbb{E}\big(\dot{L}_{Z}\, \dot{L}_{\alpha}^{\top}\big) \big\{\mathbb{E}\big(\dot{L}_{\alpha}\, \dot{L}_{\alpha}^{\top}\big)\big\}^{-1} \dot{L}_{\alpha}\right]^\top \\
    =&~ \mathbb{E}\left(\dot{L}_Z\dot{L}_{Z}^\top\right)
    - \mathbb{E}\left(\dot{L}_Z\dot{L}_{\alpha}^\top\right)\big\{\mathbb{E}\big(\dot{L}_{\alpha}\, \dot{L}_{\alpha}^{\top}\big)\big\}^{-1}\mathbb{E}\left(\dot{L}_\alpha\dot{L}_{Z}^\top\right).
\end{align*}
Following similar analysis to that of \eqref{eq:EdotL_alphadotL_Ztop}, we have 
\begin{align}
\label{eq:EdotL_ZdotL_Ztop}
\mathbb{E}\left(\dot{L}_Z\dot{L}_{Z}^\top\right) =&~ \sum_{t_1,t_2}\mathbb E\left( \frac{\partial \Theta_{t_1}}{\partial Z}^\top \dot{L}_{\Theta_{t_1}}\dot{L}_{\Theta_{t_2}}^\top \frac{\partial \Theta_{t_2}}{\partial Z} \right) = \sum_{t} \DZt  \Dmut \DZ \in\mathbb R^{nk \times nk}
\end{align}
and we note that $\mathbb{E}\big(\dot{L}_\alpha\dot{L}_{Z}^\top\big)=\mathbb{E}\big(\dot{L}_{Z}\dot{L}_\alpha^\top\big)^{\top}$. 
Moreover, by \eqref{eq:EdotL_alphadotL_Ztop} and \eqref{eq:EdotL_ZdotL_Ztop}, we can write
\begin{align}
    \Ieff = \sum_{t=1}^T \Big\{ \DZt\Dmut\DZ - \DZt\Dmut\myDone (\myDonet \Dmut \myDone)^{-1} \myDonet \Dmut \DZ   \Big\}.
    \label{eq:Ieff_Dmut_formula}
\end{align}

\subsubsection{Formula of the observed efficient information matrix $-\Heff$}  \label{sec:obseffinf}

Remark \ref{RMK:Heff_replace_Ieff} 
mentions that the efficient information matrix $\Ieff$ can also be replaced by the observed efficient information matrix $-\Heff $, whose formula is given by
\begin{align}
\Heff(Z,\alpha) = \ddot{L}_{Z,Z}(Z,\alpha) - \ddot{L}_{Z,\alpha} (Z,\alpha)\ddot{L}^{-1}_{\alpha,\alpha}(Z,\alpha)\ddot{L}_{\alpha,Z}(Z,\alpha) \ \in \ \mathbb{R}^{nk \times nk}.  
\label{eq:Heff_def_supple}
\end{align}
We next derive the second-order derivatives. 
In particular, 
\begin{align*}
 \ddot{L}_{Z,Z}(Z,\alpha) =   &  ~\frac{\partial^2 L(Z,\alpha)}{\partial Z_v\partial Z_v^{\top}}  
 =  \sum_t\left\{  \frac{\partial\DZt(Z)}{\partial Z_v^{\top}} \, \dot{L}_{\Theta_t}(Z,\alpha) +\DZt(Z) \, \frac{\partial\dot{L}_{\Theta_t}(Z,\alpha)}{\partial Z_v^{\top}} \right\}\ \in\  \mathbb{R}^{nk \times nk},
\end{align*}
where we use \eqref{eq:dotL_Z}  in the second equation. 
By \eqref{eq:partialthetapartialzmatrix}
 and  the chain rule, 
\begin{align*}
\frac{\partial\dot{L}_{\Theta_t}(Z,\alpha)}{\partial Z_v^{\top}} = &~\frac{\partial\dot{L}_{\Theta_t}(Z,\alpha)}{\partial (\Theta_t)_v^{\top}}\ \frac{\partial (\Theta_t)_v}{\partial Z_v^{\top}} 
= -\Dmut \DZ\ \in \ \mathbb{R}^{\frac{n(n+1)}{2}\times nk}, 
\end{align*}
where we use Notation \ref{notation:partialderiv} and 
\begin{equation}\label{eq:secondlthetatheta}
    \ddot{L}_{\Theta_t,\Theta_t} = \frac{\partial\dot{L}_{\Theta_t}}{\partial (\Theta_t)_v^{\top}} = -\Dmut.
\end{equation}
Moreover, 
by \eqref{eq:ldot_theta_t}  and \eqref{eq:partialthetapartialzmatrix}, 
\begin{align*}
   \frac{\partial\DZt(Z)}{\partial Z_v^{\top}} \, \dot{L}_{\Theta_t}(Z,\alpha) =  \Big\{\mathbf A_t - \bdmu_t+\text{diag}(A_{t,11}-\mu_{t,11}, \ldots,A_{t,nn}-\mu_{t,nn})\Big\}\otimes \mathrm I_k \ \in\ \mathbb{R}^{nk\times nk} 
\end{align*}
where $\bdmu_t=(\mu_{t,ij})_{n\times n} \in \mathbb{R}^{n\times n}$, and $\otimes$ denotes the matrix Kronecker product. 

\begin{notation}
\label{notation::NtN2t}
For the simplicity of notation, we define
    \begin{align}
    N_t =&~ \mathbf A_t - \bdmu_t \ \in\ \mathbb{R}^{n\times n} \label{eq::def_NtN2t}\\
    N_{2t} =&~ N_t + \operatorname{diag}(N_{t,11},\ldots,N_{t,nn}) \ \in\ \mathbb{R}^{n\times n} \label{eq::def_N2tdef}
\end{align}
Remark: $N_t$ is a function of $(Z,\alpha)$ as  $\bdmu_t=\exp( ZZ^{\top}+\alpha_t 1_n^{\top}+1_n\alpha_t^{\top} )$ depends on $(Z,\alpha)$. To emphasize that, we also write $N_t(Z,\alpha)$ and  $N_{2t}(Z,\alpha)$ in the following. 
\end{notation} 
In summary,
\begin{equation}
    \ddot{L}_{Z,Z}(Z,\alpha) =\sum_{t=1}^T \big\{-\DZt(Z) \Dmut(Z,\alpha) \DZ(Z) + N_{2t}(Z,\alpha) \otimes \mathrm I_k \big\} \in\mathbb R^{nk\times nk}.
\end{equation}

In addition, by $\alpha_v=(\alpha_1^{\top},\ldots, \alpha_T^{\top})^{\top}$, 
\begin{align*}
   \ddot{L}_{\alpha,\alpha}(Z,\alpha)=&~ \frac{\partial^2 L(Z,\alpha)}{\partial \alpha_v \partial \alpha_v^{\top} }  = \left(\ddot{L}_{\alpha_{t_1},\alpha_{t_2}} \right)_{1\leqslant t_1,t_2\leqslant T} \ \in \ \mathbb{R}^{nT\times nT} 
\end{align*}
where $\ddot{L}_{\alpha_{t_1},\alpha_{t_2}}=\frac{\partial^2 L}{\partial \alpha_{t_1} \partial \alpha_{t_2}^{\top} }\in \mathbb{R}^{n\times n}$.
When $t_1\neq t_2$, $ \ddot{L}_{\alpha_{t_1},\alpha_{t_2}}=0$, and by plugging \eqref{eq:dotL_alpha_t} and the Fact \ref{fact:partialderivfunction}, 
\begin{align}
\ddot{L}_{\alpha_t,\alpha_t} =&~ 
\frac{\partial(\myDonet \dot{L}_{\Theta_t})}{\partial \alpha_t^{\top}} =  \myDonet \frac{\partial \dot{L}_{\Theta_t}}{\partial \alpha_t^{\top}} \notag\\
=&~ \myDonet \ddot{L}_{\Theta_t,\Theta_t} \myDone = - \myDonet\Dmut\myDone
\in \mathbb{R}^{n\times n}, 
\label{eq:ddotL_alpha_t,alpha_t}
\end{align}
where the chain rule obtains the third equality, and the fourth equation  uses \eqref{eq:secondlthetatheta}. 
Therefore, 
\begin{align*}
  \ddot{L}_{\alpha,\alpha} = - \text{diag}\big\{(\myDonet \Dmut \myDone)_{t=1, \ldots,T}\big\} 
 \ \in \ \mathbb R^{nT \times nT},  
\end{align*}
which equals the negative of \eqref{eq:EdotL_alphadotL_alphatop}.


Moreover, following similar analysis, 
\begin{align}
        \ddot{L}_{Z,\alpha} =&~\frac{\partial^2 L}{\partial Z_v \partial \alpha_v^{\top}} =(\ddot{L}_{Z,\alpha_1},\ldots,\ddot{L}_{Z,\alpha_T}) \in \mathbb{R}^{nk\times nT}.\label{eq:ddotL_Z,alpha}
\end{align}
where 
by the chain rule and \eqref{eq:dotL_alpha_t}, 
\begin{align} 
\ddot{L}_{Z,\alpha_t}(Z,\alpha) =&~\frac{\partial L(Z,\alpha)}{\partial Z_v\partial \alpha_t^{\top}}= 
\frac{\partial \dot{L}_{\Theta_t}^\top}{\partial Z_v}\myDone \notag\\
=&~\DZ^\top\ddot{L}_{\Theta_t,\Theta_t}\myDone = - \DZt\Dmut\myDone\in \mathbb{R}^{nk\times n}.\label{eq:ddotL_Z,alpha_t} 
\end{align}
In summary, we can write
\begin{align*}
\ddot{L}_{Z,\alpha}(Z,\alpha)  =&~ 
-\begin{pmatrix}
 \DZt\mathcal{D}_{\mu_1}\myDone & 
\cdots & 
  \DZt\mathcal{D}_{\mu_T}\myDone
\end{pmatrix} 
\ \in \ \mathbb R^{nk \times nT},
\end{align*}

\begin{remark}[Relationship between $\Ieff$ and $\Heff$]
Comparing to Sections \ref{sec:ieffformulae} and \ref{sec:obseffinf}, 
\begin{equation}
\label{eq:HeffIeff}
   \Heff(Z,\alpha) = -\Ieff(Z,\alpha) + \sum_t N_{2t}(Z,\alpha)\otimes \mathrm{I}_k. 
\end{equation}
\end{remark}


\subsection{Eigen-Structure of the Efficient Information Matrix}
\label{subsec:formula_UH}
To establish the theorems, one major challenge is  the singularity of the  efficient  information matrix $\Ieff(\check{Z},\check{\alpha})$ in \eqref{eq:newtonsolpseudo}, as discussed in Remark \ref{lb:singularity}. 
We address this issue by studying the   eigenspace of $\Ieff$ in this section.

In particular, 
we first construct a restricted eigenspace in Section 
\ref{sec:eigenspace}  providing the basis for studying $\Ieff$. 
Moreover, to study $\Ieff$, we note that by \eqref{eq:Ieff_Dmut_formula}, 
\begin{align}
    \Ieff =&~ \DZt \left(\sum_{t=1}^T \myXt\right) \DZ.
     \label{eq:IeffHeffinXt}
\end{align}
where we define
\begin{align}\label{eq:myxtdef}
    \myXt=  \Dmut- \Dmut\myDone (\myDonet \Dmut \myDone)^{-1} \myDonet \Dmut \in\mathbb R^{{n(n+1)}/2\times {n(n+1)}/2}.
\end{align}
Thus, to study $\Ieff$, 
we  introduce the properties of $\mathcal{X}_t$ and $\DZ $ in Sections \ref{sec:ieffxt} and \ref{sec:ieffdthetaz}, respectively, and then examine $\Ieff$ in Section \ref{sec:ieffcombine}. 
We clarify that  $\myXt$, $\DZ$, and $\Ieff$ are functions of $(Z,\alpha)$. The following derivations  assume that a particular $(Z,\alpha)$ is given, so we do not emphasize the dependence on $(Z,\alpha)$ in the notation when there is no ambiguity. 


\subsubsection{Restricted Eigenspace} \label{sec:eigenspace}

Given any $Z\in \mathbb{R}_0^{n\times k}$, 
Definition \ref{def::UZA} below  gives an explicit characterization of the corresponding intrinsic restricted eigenspace based on the singular value decomposition of $Z$.  
In particular, 
let  $Z=\sum_{i=1}^k\lambda_i U_i  V_i^{\top} $,  
where $U_i\in \mathbb{R}^n$ are the left singular vectors, $V_i\in \mathbb{R}^k$ are the right singular vectors,  $\lambda_1\geqslant \cdots \geqslant \lambda_k > 0$.
We let $\{U_{k+1},\ldots, U_n\}$ be a set of orthonormal vectors that are orthogonal to $\mathrm{span}\{U_1, \ldots, U_{k}\}$, 
where $\mathrm{span}\{\cdot\}$ represents the linear space formed by the included vectors. 
As $1_n^{\top}Z=0$ for any $Z\in \mathbb{R}_0^{n\times k}$, 
we set, without loss of generality, that
\begin{align}\label{eq:uk+1def1n}
    U_{k+1}=1_n/\sqrt{n}. 
\end{align} 
Note that each $U_i$ may not be uniquely identified.
But the following conclusions and  proofs hold as long as $\mathrm{span}\{U_1,\ldots, U_k\}$ and $\mathrm{span}\{U_{k+2},\ldots, U_n\}$ are uniquely identified, respectively; also see Remark \ref{rm:uzunique}. 
Define 
\begin{align}\label{eq:betaijdefA}
    \beta_{ij} = U_i \otimes V_j \ \in \  \mathbb{R}^{nk}\quad \text{ for }\ 1\leqslant i \leqslant n,\ 1\leqslant j \leqslant k
\end{align}
where $\otimes$ denotes the Kronecker product.
By \eqref{eq:betaijdefA} and properties of $U_i$ and $V_j$, 
$\{\beta_{ij} : i=1, \ldots,n;j=1, \ldots,k \} $ form an orthonormal basis of $\mathbb R^{nk}$.
Note that $(U_i, V_i, \lambda_i)$'s and  $\beta_{ij}$'s are functions of $Z\in \mathbb{R}_0^{n\times k}$, which are not emphasized in the notation for the simplicity of notation here. 

\begin{definition}
\label{def::UZA}
Let  $\rankU = nk - k(k+1)/2$. 
Given $Z\in \mathbb{R}_0^{n\times k}$, define $\UZ\in \mathbb{R}^{nk \times \rankU}$ as a  matrix with column vectors in $  \mathcal{B}_{C,1}\cup   \mathcal{B}_{C,2}\cup   \mathcal{B}_{C,3}$, where
\begin{align*}
    \mathcal{B}_{C,1}=&~\{\beta_{ii} : i=1, \ldots,k\}, \hspace{2em} \mathcal{B}_{C,2}=\left\{\frac{\lambda_j \beta_{ij} + \lambda_i \beta_{ji}}{\sqrt{\lambda_i^2 + \lambda_j^2}} : 1\leqslant i < j \leqslant k \right\},\\
   \mathcal{B}_{C,3}=&~\{\beta_{ji} : i=1, \ldots, k \ \text{ and }\ j = k+ 2, \ldots, n\},
\end{align*}
and $\beta_{ij}$'s are defined as in \eqref{eq:betaijdefA}.  
Let $\operatorname{col}(\UZ)$ denote the column space of $\UZ$. 
\end{definition}

In the following proofs, 
we refer to 
 $\mathrm{col}(\UZ)$  as the intrinsic restricted eigenspace at a given $Z\in \mathbb{R}_0^{n\times k}$.  
In general, definitions of $\UZ$ and  $\beta_{ij}$'s are not  unique since  
$\{U_{i}:i=k+2,\ldots, n\}$ are unique up to an orthogonal transformation,
and $\{(U_{i},V_i):i=1,\ldots, k\}$ are unique up to sign change only if the non-zero singular values are distinct. 
But our proofs 
only rely on the space $\mathrm{col}(\UZ)$, which remains the same even when $\UZ$ is not unique;  see, e.g., Remark \ref{rm:uzunique}. 
We will establish  important properties of $\UZ$ and show its connections with $\Ieff(Z,\alpha)$ in Lemmas \ref{lem_DZtDz_eigenstructure} and  \ref{lem_presentingUZmain} below.

\bigskip
\begin{remark}\label{rm:takeuzcheck}
In the one-step estimator \eqref{newtonsolUZ}, 
$\check{\mathcal U}$ can be taken as any ${nk\times (nk-k(k+1)/2)}$ matrix such that 
$\mathrm{col}(\check{\mathcal U})=\mathrm{col}(\Ieff(\check{Z},\check{\alpha}))$,
and the proof remains the same. 
For the convenience of theoretical analysis, 
we will fix our choice to be $\mathcal{U}_{\check{Z}}$, i.e., the construction in Definition  \ref{def::UZA} at $\check{Z}$. Please also see Remark \ref{rm:uzunique}. 
\end{remark}

\subsubsection{Properties of $\mathcal{X}_t$}\label{sec:ieffxt}

By the definition of $\myXt$  in  \eqref{eq:myxtdef}, we know   $\myXt$  is a function of 
 $\bdmu_t = (\mu_{t,ij})_{1\leqslant i,j \leqslant n}\in \mathbb{R}^{n\times n}$.
 Lemma \ref{lem:property_myXt} establish properties of  $\myXt$ that can hold for any $\bdmu_t$ with positive entries, including the special cases with 
$\mu_{t,ij} = \exp(z_i^\top z_j +\alpha_{it} + \alpha_{jt}) $ or $\mu_{t,ij} = \exp(G_{ij} +\alpha_{it} + \alpha_{jt}) $  under our considered model.  
For any matrix $X\in\mathbb{R}^{a\times b} $, we denote its null space by $\text{null}(X) = \{ v\in\mathbb{R}^{b} : Xv = 0 \} $.
\begin{lemma}
\label{lem:property_myXt}
Given any $\bdmu_t\in \mathbb{R}^{n\times n}$ with positive entries, 
 $\myXt$ defined in \eqref{eq:myxtdef} satisfies the following properties. 
\begin{enumerate}
    \item $\myXt \Dmut^{-1} \myXt = \myXt $ and $\myXt\myDone = 0 $.
    \item  The null space of $\myXt$ equals the column space of $\myDone $, and $$\operatorname{rank}(\myXt) = {n(n+1)}/2-n.$$
    \item The smallest non-zero eigenvalue of $\myXt$ is lower bounded by $ \min_{1 \leqslant i \leqslant j \leqslant n} \{\mu_{t,ij}\} $, where $\mu_{t,ij}$ is the $(i,j)$-element of $\bdmu_t$.
\end{enumerate}
\end{lemma}

\begin{proof}
See Section \ref{sec:pfA1} on Page \pageref{sec:pfA1}.
\end{proof}




\subsubsection{Properties of $\DZ$}\label{sec:ieffdthetaz}

Let  $\DZ \in \mathbb{R}^{\frac{n(n+1)}{2}\times nk} $  be 
defined as in \eqref{eq:dthetazdef}. 
We next present the right singular vectors of $\DZ$, which are the eigenvectors of 
\begin{align*}
    \Xi_Z:=  \DZ^\top \DZ = \DZt \DZ \ \in \  \mathbb{R}^{nk\times nk}. 
\end{align*}
As $ \DZ$ is a function of $Z$ by  
Fact \ref{fact:partialderivfunction}, 
we know that $\Xi_Z$ is also a function of  $Z$. 
 The following Lemma \ref{lem_DZtDz_eigenstructure} explicitly characterizes the eigen-decomposition of $\Xi_Z$. 
 

\begin{lemma}
\label{lem_DZtDz_eigenstructure}
Given any $Z\in\mathbb R_0^{n\times k} =\{Z\in\mathbb R^{n\times k} : \operatorname{rank}(Z) = k, 1_n^\top Z = 0\} $, consider its singular value decomposition  $Z=\sum_{i=1}^k \lambda_i U_i V_i^{\top}$ specified as in Section
\ref{sec:eigenspace}.
Then define $\{\beta_{ij}: 1\leqslant i\leqslant n, 1\leqslant j\leqslant k \}$ as in Eq. \eqref{eq:betaijdefA},
 define $\mathcal{B}_{C,j}$ for $j=1,2,3$
as in Section \ref{sec:eigenspace}, and let  
\begin{align} \label{eq:bn1bn2setsupp} 
    \mathcal{B}_{N,1} = \{\beta_{k+1,i}: 1\leqslant i \leqslant k\}, \quad \text{and}\quad \mathcal{B}_{N,2} = \left\{ \frac{\lambda_i\beta_{ij}-\lambda_j\beta_{ji}}{\sqrt{\lambda_i^2+\lambda_j^2}}  : 1\leqslant i < j \leqslant k \right\}.  
\end{align} 
We have
\begin{enumerate}
	\item $1_n^\top U_i=0 $ for $i=1, \ldots,k$.
	\item Let $Z_v\in \mathbb{R}^{nk\times 1}$ be the  vectorization of $Z$ defined as in \eqref{eq:vecdef}. Then $Z_v$ is in the $k$-dimensional subspace spanned by vectors in $\mathcal{B}_{C,1}=\{\beta_{ii}: i=1,\ldots, k \}$. 
	\item  
  Vectors in $\cup_{j=1}^3\mathcal{B}_{C,j} \cup \mathcal{B}_{N,1}\cup \mathcal{B}_{N,2}$
   form an orthonormal basis of $\mathbb R^{nk} $, and they  
 form a set of (normalized) eigenvectors  of $\Xi_Z$.  
Specifically, the vectors satisfy 

\quad $\Xi_Z \beta_{ii} = 2\lambda_i^2 \beta_{ii}$ for $i=1, \ldots,k$; 

\quad $\Xi_Z \beta_{ji} = \lambda_i^2 \beta_{ji}$ for $j=k+1, \ldots,n$ and $i=1, \ldots,k$; 

\quad $\Xi_Z(\lambda_j\beta_{ij} + \lambda_i\beta_{ji})/\sqrt{\lambda_i^2+\lambda_j^2} = (\lambda_i^2+\lambda_j^2)(\lambda_j\beta_{ij} + \lambda_i\beta_{ji})/\sqrt{\lambda_i^2+\lambda_j^2} $ for $1\leqslant i < j \leqslant k$;

\quad $\Xi_Z(\lambda_i\beta_{ij} - \lambda_j\beta_{ji})/\sqrt{\lambda_i^2+\lambda_j^2} = 0 $ for $1\leqslant i < j \leqslant k$.

\end{enumerate} 
\end{lemma}

\begin{proof}
See Section \ref{sec:pfA2} on Page \pageref{sec:pfA2}.
\end{proof}

\subsubsection{Properties of $\Ieff$}\label{sec:ieffcombine}
In the following Lemma  \ref{lm:presentuzmain}, we present useful properties of $\Ieff$, which are derived based on Lemmas \ref{lem:property_myXt} and \ref{lem_DZtDz_eigenstructure}.

\begin{lemma}\label{lm:presentuzmain} \label{lem_presentingUZmain} 
For $Z\in \mathbb{R}_0^{n\times k}$ and  $\alpha \in \mathbb{R}^{n\times T}$,
define $\UZ$ as in Definition \ref{def::UZA},
and  define $\Ieff(Z,\alpha)$ as in Eq. \eqref{eq:ieffdef}. Then  we have

\vspace{2pt}




\noindent  (i) $   \mathrm{col}(\Ieff(Z,\alpha)) = \mathrm{col}(\UZ) = \mathrm{span}(\cup_{j=1}^3\mathcal{B}_{C,j})$,  $\operatorname{rank} (\Ieff(Z,\alpha)) = nk-k(k+1)/2$, and  
 $$\sigma_{\min}(\UZ^{\top} \Ieff(Z,\alpha) \UZ) \geqslant nT\, \sigma_{\operatorname{min}}( Z^\top  Z /n)\min_{{1\leqslant i,j\leqslant n, 1\leqslant t\leqslant T}}\exp(z_i^\top  z_j +\alpha_{it} +  \alpha_{jt}).$$

\vspace{2pt}

\noindent (ii) $\mathrm{null}(\Ieff(Z,\alpha)) = \mathrm{col}(\UZ)^{\perp}= \mathrm{span}(\mathcal{B}_{N,1}\cup \mathcal{B}_{N,2})$, where $\mathrm{col}(\UZ)^{\perp}$ denotes the orthogonal complement of $\mathrm{col}(\UZ)$. 

\vspace{2pt}

\noindent (iii)  $\UZ$ satisfies the equation $\UZ\UZ^\top \DZt\myXt = \DZt\myXt$, 
     and thus  
    $$\UZ\UZ^\top \Seff(Z,\alpha) =  \Seff(Z,\alpha), \quad \UZ\UZ^\top \Ieff(Z,\alpha) =  \Ieff(Z,\alpha),$$i.e., $\Seff$ lies in the eigenspace corresponding to the non-zero eigenvalues of $\Ieff$. 
\end{lemma}

\begin{proof}
See Section \ref{sec:pfA3} on Page \pageref{sec:pfA3}.
\end{proof}

\bigskip


Interestingly,
Lemma \ref{lem_presentingUZmain}-(ii) suggests that 
$\mathrm{null}(\Ieff(Z,\alpha))$ has  connections with two sources of singularity in   Remark \ref{lb:singularity}. 
For $ \mathcal{B}_{N,1}$, note 
 $\beta_{k+1,i} =  1_n \otimes  V_i /\sqrt{n}$ 
 and thus $\beta_{k+1,i}^{\top}Z_v=1_n^{\top}ZV_i/\sqrt{n}=0$, 
corresponding to the $k$-dimensional linear constraint $1_n^{\mytrans} Z=0$. 
Moreover, 
it will be shown in Lemma \ref{lem_HspaceDiff} that, 
$\mathcal{B}_{N,2}$ forms a set of basis of tangent space at $Z$ to the equivalence class $[Z]$  \citep{massart2020quotient}.
Intuitively, $\mathcal{B}_{N,2}$ accounts for the unidentifiability from the invariance with respect to orthonormal transformations. 
Also see Remark \ref{rm:internulluzieff}. 
In addition, 
although the efficient information matrix  $\Ieff(Z,\alpha)$ is a function of both $Z$ and $\alpha$, 
Lemma \ref{lm:presentuzmain} shows that its corresponding  restricted eigenspace $\mathrm{col}(\UZ)$ only 
depends on $Z$. This is consistent with our intuition that  
the null space of $\Ieff(Z,\alpha)$ is caused by the identifiability constraints of $Z$ only in  Remark \ref{lb:singularity}.

\bigskip

\begin{remark}\label{rm:internulluzieff}
This remark provides interpretation on $\mathrm{null}(\UZ \UZ^\mytrans)$ and $\mathrm{null}(\Ieff(Z,\alpha))$, which are the same by Lemma \ref{lem_presentingUZmain}. 
Lemma \ref{lem_presentingUZmain}-(ii) shows that $\mathrm{null}(\UZ \UZ^\mytrans)$ consists of  
 $\mathcal{B}_{N,1}$ and $\mathcal{B}_{N,2}$.
From the proof of Lemma \ref{lem_presentingUZmain}, we see that 
  $\mathcal{B}_{N,1}$ corresponds $\mathrm{null}(\myXt)$. 
Intuitively, Lemma \ref{lem:property_myXt} shows that $\mathrm{null}(\myXt)$ comes form the centering constraint  $1_n^{\top}Z = 0 $, 
and 
Lemmas \ref{lem_DZtDz_eigenstructure} and \ref{lem_HspaceDiff} indicate that $\mathcal{B}_{N,2}$  is due to the orthogonal group transformation invariance of $Z$. 
Therefore,  the singularity  of both $\UZ$ and  $\Ieff(Z,\alpha)$ can be  intuitively viewed as consequences of the orthogonal transformation invariance and centering constraint of $Z$. This is consistent with Remark \ref{lb:singularity}. 
\end{remark} 

\subsection{Log Profile Likelihood of \texorpdfstring{$G$}{} and Its Derivatives}
\label{subsec::prelim_G}
\paragraph*{Log Profile Likelihood of $G$}
As mentioned in Remark \ref{rm:techchanlplm},
we consider the profile likelihood of $G$ to establish the theoretical properties of the penalized MLE estimator. 
Specifically, following \cite{murphy2000profile}, 
we define the log profile likelihood of $G$ as $pl(G) = l(G,\hat\alpha(G)) $, where $l(G,\alpha)$ denotes the   log-likelihood function of 
  $(G,\alpha)$, and $\hat\alpha(G) = \argmax_{\alpha\in \mathbb{R}^{n\times T}} l(G,\alpha)  $.
Since 
$l(G,\alpha)=\sum_{t=1}^T l_t(G,\alpha_t)$ 
with  
\begin{equation*}
  l_t(G,\alpha_t)=\sum_{1\leqslant i\leqslant j\leqslant n} \left[A_{t,ij} (\alpha_{it}+\alpha_{jt}+G_{ij}) - \exp(\alpha_{it}+\alpha_{jt}+G_{ij}) \right],
\end{equation*} 
we have 
\begin{align*}
\hat\alpha(G) = \argmax_{\alpha_t\in \mathbb{R}^{n};\ t=1,\ldots, T}\sum_{t=1}^T l_t(G,\alpha_t) = \begin{pmatrix}
    \hat{\alpha}_1(G)&   \cdots &   \hat{\alpha}_T(G)
\end{pmatrix} \ \in \mathbb{R}^{n\times T}
\end{align*}
where for $t=1,\ldots, T$, we define
\begin{align}\label{eq:defalphahattg}
     \hat\alpha_t(G) = \argmax_{\alpha_t\in \mathbb{R}^n} l_t(G,\alpha_t).
\end{align}
It follows that $ pl(G) = \sum_{t=1}^T pl_t(G)$ with $pl_t(G)=l_t(G,\hat\alpha_t(G))$.


As we consider the Taylor expansion of $pl(G)$ in Section \ref{sec:pfthm3A}, 
we 
next derive the first-order and second-order derivatives of $pl(G)$ with respect to 
\begin{align}\label{eq:defgv}
  G_v = &~(G_{11}/\sqrt 2, G_{22}/\sqrt 2, \ldots, G_{nn}/\sqrt 2, G_{12}, G_{13}, \ldots, G_{n-1,n} )^\top \in \mathbb R^{\frac{n(n+1)}{2} \times 1}, 
\end{align}
which is a vectorization of $G\in \mathbb{R}^{n\times n}$, 
where $1/\sqrt{2}$ scaling  is used due to the reason similar to Remark \ref{rmk::vectorization_sqrt2}. 
By $ pl(G) = \sum_{t=1}^T pl_t(G)$, 
\begin{align}\label{eq:linearformpartialprofileg}
\frac{\partial pl(G)}{\partial G_v} =&~ \sum_{t=1}^T  \frac{\partial pl_t(G)}{\partial G_v},  \quad \text{and}\quad \quad
   \frac{\partial^2 \ pl(G)}{\partial G_v \partial G_v^\top}=\sum_{t=1}^T \frac{\partial^2 \ pl_t(G)}{\partial G_v \partial G_v^\top}.
\end{align}
Therefore, to derive the first-order and second-order derivatives of $pl(G)$, 
it suffices to derive the first-order and second-order derivatives of $pl_t(G)$, which are given Sections \ref{subsubsec::1st_deriv_profile} and \ref{subsubsec::2nd_deriv_profile}, respectively.



\subsubsection{First-Order Derivative of $pl_t(G)$} 
\label{subsubsec::1st_deriv_profile}

By $pl_t(G)=l_t(G,\hat{\alpha}_t(G))$, we have
\begin{align}
    \frac{\partial pl_t(G)}{\partial G_v}=&~ \frac{\partial l_t(G,\hat\alpha_t(G))}{\partial G_v}\notag \\
    =&~ \left.\left\{ \frac{\partial l_t(G,\alpha_t)}{\partial G_v} + \frac{\partial \hat{\alpha}_t(G)^{\mytrans}}{\partial G_v} \times \frac{\partial l_t(G,\alpha_t)}{\partial \alpha_t} \right\} \right|_{(G,\hat{\alpha}_t(G))},\label{eq::dpltG_dGv=dltG_dGv0} 
\end{align}
where the second equation follows by the chain rule, and
\begin{align}\label{eq:ltdot_Gv}
	\frac{\partial l_t(G,\alpha_t)}{\partial G_v} =&~
\begin{pmatrix}
\frac{\partial l_t}{\partial (G_{11}/\sqrt{2})}	\\[2pt]
\vdots\\[2pt]
\frac{\partial l_t}{\partial (G_{nn}/\sqrt{2})} \\[8pt] 
\frac{\partial l_t}{\partial G_{12}}	\\[2pt]
\vdots\\[2pt]
\frac{\partial l_t}{\partial G_{n-1,n} } 
\end{pmatrix}=\begin{pmatrix}
\sqrt 2(A_{t,11}-\mu_{t,11})\\[4pt] \vdots \\[4pt] \sqrt2(A_{t,nn}-\mu_{t,nn})\\[4pt] A_{t,12} - \mu_{t,12}\\[4pt] \vdots \\[4pt] A_{t,n-1,n} - \mu_{t,n-1,n} 
\end{pmatrix} \in \mathbb{R}^{\frac{n(n+1)}{2}\times 1},\\  
\frac{\partial l_t(G,\alpha_t)}{\partial \alpha_t} =&~ \begin{pmatrix}
          \frac{\partial l_t}{\partial \alpha_{1t}}&\ldots &   \frac{\partial l_t}{\partial \alpha_{nt}}
    \end{pmatrix}^{\top} \in \mathbb{R}^{n\times 1}.\notag
\end{align}
represent the derivatives of $l_t(G,\alpha_t)$ with respect to $G_v$ and $\alpha_t$, respectively. 
By \eqref{eq:defalphahattg}, we have the score equation 
\begin{align}\label{eq:scoreeqalphahat}
    \left.\left\{\frac{\partial l_t(G,\alpha_t)}{\partial \alpha_t}\right\}\right|_{(G,\hat\alpha_t(G))} =0.
\end{align}
Plugging \eqref{eq:scoreeqalphahat} in \eqref{eq::dpltG_dGv=dltG_dGv0}, we obtain
\begin{align}
    \frac{\partial pl_t(G)}{\partial G_v} =     &~ \left.\frac{\partial l_t(G,\alpha_t)}{\partial G_v}  \right|_{(G,\hat{\alpha}_t(G))}. \label{eq::dpltG_dGv=dltG_dGv}
\end{align}

\subsubsection{Second-Order Derivative of $pl_t(G)$} 
\label{subsubsec::2nd_deriv_profile}
By the first-order derivative in \eqref{eq::dpltG_dGv=dltG_dGv}, 
  \begin{align}
    \frac{\partial^2 \ pl_t(G)}{\partial G_v \partial G_v^\top} =&~ \frac{\partial\left\{ \left.\frac{\partial l_t(G,\alpha_t)}{\partial G_v}  \right|_{(G,\hat{\alpha}_t(G))}\right\}}{\partial G_v^{\top}}  \notag \\
    =&~\left.\left[ \frac{\partial\left\{\frac{\partial l_t(G,\alpha_t)}{\partial G_v}\right\}}{\partial G_v^{\top}}   +  \frac{\partial\left\{\frac{\partial l_t(G,\alpha_t)}{\partial G_v}\right\}}{\partial \alpha_t^{\top}} \times \frac{\partial \hat{\alpha}_t(G)}{\partial G_v^{\mytrans}}  \right]\right|_{(G,\hat{\alpha}_t(G))}  \notag\\
    =&~\left.\left\{ \frac{\partial^2 l_t(G,\alpha_t)}{\partial G_v\partial G_v^\top} +  \frac{\partial^2 l_t(G,\alpha_t)}{\partial G_v\partial \alpha_t^\top} \times  \frac{\partial \hat{\alpha}_t(G)}{\partial G_v^{\mytrans}}\right\}\right|_{(G,\hat{\alpha}_t(G))} \label{eq:2nd_deriv_pltG}
\end{align}  
where the second equation follows by the chain rule,
and in the third equation, 
we let 
${\partial^2 l_t(G,\alpha_t)}/{\partial G_v\partial G_v^\top}$ and $ {\partial^2 l_t(G,\alpha_t)}/{\partial G_v\partial \alpha_t^\top}$ denote the second-order partial derivatives of $l_t(G,\alpha_t)$ similarly to the definitions in \eqref{eq:ltdot_Gv}. 

By \eqref{eq:scoreeqalphahat}, we have
\begin{align*}
    0 = &~\frac{\partial\left\{ \left.\frac{\partial l_t(G,\alpha_t)}{\partial \alpha_t}  \right|_{(G,\hat{\alpha}_t(G))}\right\}}{\partial G_v^{\top}}  \notag\\
    =&~ \left.\left\{ \frac{\partial^2 l_t(G,\alpha_t)}{\partial \alpha_t \partial G_v^{\top}} +   \frac{\partial^2 l_t(G,\alpha_t)}{\partial \alpha_t\partial \alpha_t^\top} \times \frac{\partial \hat{\alpha}_t(G)}{\partial G_v^{\mytrans}}\right\}\right|_{(G,\hat{\alpha}_t(G))}
\end{align*}
where the second equation follows by the  chain rule, and ${\partial^2 l_t(G,\alpha_t)}/{\partial \alpha_t\partial \alpha_t^\top}$ and $ {\partial^2 l_t(G,\alpha_t)}\\/{\partial \alpha_t\partial G_v^\top}$ denote the second-order partial derivatives of $l_t(G,\alpha_t)$  similarly to the definition in \eqref{eq:ltdot_Gv}.  
Thus, 
\begin{align}\label{eq:partialalphahtformula}
   \frac{\partial \hat{\alpha}_t(G)}{\partial G_v^{\mytrans}} =-\left.\left[ \left\{\frac{\partial^2 l_t(G,\alpha_t)}{\partial \alpha_t\partial \alpha_t^\top}\right\}^{-1} \frac{\partial^2 l_t(G,\alpha_t)}{\partial \alpha_t\partial G_v^\top}\right]\right|_{(G,\hat{\alpha}_t(G))} 
\end{align}

Plugging \eqref{eq:partialalphahtformula} in \eqref{eq:2nd_deriv_pltG}, we obtain
\begin{align}
    \frac{\partial^2 \ pl_t(G)}{\partial G_v \partial G_v^\top} =&  ~  \left.\left[ \frac{\partial^2 l_t(G,\alpha_t)}{\partial G_v\partial G_v^\top} - \frac{\partial^2 l_t(G,\alpha_t)}{\partial G_v\partial \alpha_t^\top} \left\{\frac{\partial^2 l_t(G,\alpha_t)}{\partial \alpha_t\partial \alpha_t^\top}\right\}^{-1} \frac{\partial^2 l_t(G,\alpha_t)}{\partial \alpha_t\partial G_v^\top}\right]\right|_{(G,\hat{\alpha}_t(G))}.\label{eq::dpltG_dGv=dltG_dGvdGV}
\end{align}

\bigskip
\begin{remark}[Equivalent Expressions]\label{rm:partialgthetaequi}
By calculating the derivatives, and expressions in   \eqref{eq:ldot_theta_t}, \eqref{eq:dotL_alpha_t}, 
 \eqref{eq:secondlthetatheta},  and \eqref{eq:ddotL_alpha_t,alpha_t}, 
 we have
\begin{align*}
    \frac{\partial l_t(G,\alpha_t)}{\partial G_v} = &~\dot{L}_{\Theta_t}, \hspace{4.7em} \frac{\partial l_t(G,\alpha_t)}{\partial \alpha_t} = \dot{L}_{\alpha_t}, \hspace{4.2em} 
    \frac{\partial^2 l_t(G,\alpha_t)}{\partial \alpha_t \partial G_v^\top} =\ddot{L}_{\alpha_t,\Theta_t},\\
 \frac{\partial^2 l_t(G,\alpha_t)}{\partial G_v\partial G_v^\top} =&~\ddot{L}_{\Theta_t,\Theta_t},\hspace{3.2em} 
    \frac{\partial^2 l_t(G,\alpha_t)}{\partial \alpha_t \partial \alpha_t^\top} = \ddot{L}_{\alpha_t,\alpha_t}, \hspace{3.2em}\frac{\partial^2 l_t(G,\alpha_t)}{\partial G_v \partial \alpha_t^\top} = \ddot{L}_{\Theta_t,\alpha_t},
\end{align*}
where similarly to \eqref{eq:ddotL_alpha_t,alpha_t}, we define 
\begin{align}
    \ddot{L}_{\alpha_t,\Theta_t} = 
\frac{\partial \dot{L}_{\alpha_t}}{\partial (\Theta_t)_v^{\top}} =
\myDonet \frac{\partial \dot{L}_{\Theta_t}}{\partial (\Theta_t)_v^{\top}} =\myDonet\ddot{L}_{\Theta_t,\Theta_t} = - \myDonet\Dmut\in \mathbb{R}^{n\times \frac{n(n+1)}2}. 
\label{eq:ddotL_alpha_t,Theta_t}
\end{align}
The above identities suggest that taking derivatives with respect to $\Theta_t$ and $G_v$ gives the same form, which follows naturally by the chain rule and the fact $\partial (\Theta_t)_v / \partial G_v = \mathrm{I}_{{n(n+1)}/{2}}  $. 
Then by \eqref{eq::dpltG_dGv=dltG_dGv} and \eqref{eq::dpltG_dGv=dltG_dGvdGV}, we can also write
\begin{align*}
 \frac{\partial pl_t(G)}{\partial G_v}
    = \dot{L}_{\Theta_t}(G,\hat\alpha(G)) \quad  \text{and}\quad   \frac{\partial^2 \ pl_t(G)}{\partial G_v  \partial G_v^\top} = 
    \left.\left\{
    \ddot{L}_{\Theta_t,\Theta_t} -  \ddot{L}_{\Theta_t,\alpha_t} \ddot{L}_{\alpha_t,\alpha_t}^{-1} \ddot{L}_{\alpha_t,\Theta_t}
    \right\}\right|_{(G,\hat{\alpha}_t(G))}. 
\end{align*}
Moreover, by \eqref{eq:linearformpartialprofileg}, 
\begin{align}
 \frac{\partial pl(G)}{\partial G_v}
    = &~\sum_{t=1}^T  \dot{L}_{\Theta_t}(G,\hat\alpha(G)),   \label{eq::conclude_1s_deriv_plGall}\\
    \frac{\partial^2 \ pl(G)}{\partial G_v \partial G_v^\top} =&~\sum_{t=1}^T   
    \left.\left\{
    \ddot{L}_{\Theta_t,\Theta_t} -  \ddot{L}_{\Theta_t,\alpha_t} \ddot{L}_{\alpha_t,\alpha_t}^{-1} \ddot{L}_{\alpha_t,\Theta_t}
    \right\}\right|_{(G,\hat{\alpha}_t(G))}. \label{eq::conclude_2nd_deriv_plGall}    
\end{align}
\end{remark}

\begin{remark}
\label{rmk::Interpret_myXt}
  By \eqref{eq:secondlthetatheta},  
\eqref{eq:ddotL_alpha_t,alpha_t},  and \eqref{eq:ddotL_alpha_t,Theta_t}, $\myXt$ in \eqref{eq:myxtdef} can be equivalently expressed as 
\begin{equation}
\label{eq:myXt_formulae}
    \begin{aligned}
    \myXt =&~ - \left(\ddot{L}_{\Theta_t,\Theta_t} -\ddot{L}_{\Theta_t,\alpha_t} \ddot{L}^{-1}_{\alpha_t,\alpha_t}\ddot{L}_{\alpha_t,\Theta_t}\right).
    \end{aligned}
\end{equation}
By \eqref{eq:myXt_formulae} and  \eqref{eq::conclude_2nd_deriv_plGall}  , we have 
\begin{align}\label{eq:partiallggchitform}
     \frac{\partial^2 \ pl(G)}{\partial G_v \partial  G_v^{\mytrans}}=-\sum_{t=1}^T\myXt(G,\hat\alpha_t(G)).
\end{align}
We clarify that 
$(G,\hat\alpha_t(G))$ in the right hand side of \eqref{eq:partiallggchitform} represents that $\myXt$ is evaluated at $\mu_t:=G + \hat\alpha_t(G)1_n^{\top} + 1_n\hat\alpha_t(G)^{\top}$, and $\myXt$  depends on $(G,\hat\alpha_t(G))$  only through $\mu_t$.
As Lemma \ref{lem:property_myXt} obtains properties of $\myXt$ that do not depend on values of $\mu_t$, 
\eqref{eq:partiallggchitform} suggests we can apply Lemma \ref{lem:property_myXt} to ${\partial^2 \ pl(G)}/{\partial G_v \partial  G_v^{\mytrans}}$; see, e.g., Section \ref{sec:pfD2}.
\end{remark}

\subsection{Proof of Lemmas \ref{lem:property_myXt}--\ref{lm:presentuzmain} in
Section \ref{subsec:formula_UH}
}
\subsubsection{Proof of Lemma \ref{lem:property_myXt}} \label{sec:pfA1}

\paragraph*{Proof of Part 1.}
First, by \eqref{eq:myxtdef}, 
\begin{align*}
    \myXt\Dmut^{-1} \myXt= &~ \big\{ \Dmut- \Dmut\myDone (\myDonet \Dmut \myDone)^{-1} \myDonet \Dmut  \big\}\Dmut^{-1} \notag\\
    &~\times \big\{ \Dmut- \Dmut\myDone (\myDonet \Dmut \myDone)^{-1} \myDonet \Dmut  \big\} \notag\\
    =&~ \Dmut-2 \Dmut\myDone (\myDonet \Dmut \myDone)^{-1} \myDonet \Dmut\notag\\
    &~ +\Dmut\myDone (\myDonet \Dmut \myDone)^{-1} \myDonet  \Dmut\myDone (\myDonet \Dmut \myDone)^{-1} \myDonet \Dmut \notag\\
    =&~\Dmut- \Dmut\myDone (\myDonet \Dmut \myDone)^{-1} \myDonet = \myXt.\notag
\end{align*}
Second, by \eqref{eq:myxtdef}, 
\begin{align}
    \myXt\myDone =  &~ \Dmut\myDone- \Dmut\myDone (\myDonet \Dmut \myDone)^{-1} \myDonet \Dmut \myDone = 0.  \label{eq:myxtmydoneprod} 
\end{align}

\paragraph*{Proof of Part 2.}
On the one hand, 
\eqref{eq:myxtmydoneprod} suggests 
col$(\myDone) \subset$ null$(\myXt)$.
Then by rank$(\myDone) = n$, 
rank$(\myXt) \leqslant {n(n+1)}/{2} - n$. On the other hand, $$\operatorname{rank}(\myXt)\geqslant \text{rank}(\Dmut) - \text{rank}\{\Dmut\myDone (\myDonet \Dmut \myDone)^{-1} \myDonet \Dmut \}\geqslant \frac{n(n+1)}{2} - n,$$
where the  second inequality is obtained by rank$(\Dmut) = {n(n+1)}/{2}$ and rank$(\myDonet \Dmut \myDone) = n$ from Remark \ref{rmk::invertable_ddotLalphaalpha}. 
Combing the lower and upper bounds of $\text{rank}(\Dmut)$, we obtain rank$(\myXt) = {n(n+1)}/{2}-n$ and  $\text{null}(\myXt) = \text{col}(\myDone)$.

\paragraph*{Proof of Part 3.}
Suppose the eigendecomposition of $\myXt$ is $\myXt = Q_t \Lambda_t Q_t^\top $, where $Q_t$ is an orthonormal matrix and $\Lambda_t = \text{diag}(\tilde\Lambda_t, 0_n)$ and $\tilde\Lambda_t \in\mathbb R^{{n(n+1)}/{2}-n} $ is the vector of $\myXt$'s nonzero eigenvalues. Besides, let $\tilde Q_t$ denote the first $n(n+1)/2 -n$ columns of $Q_t$. Then from the relationship $\myXt \Dmut^{-1} \myXt = \myXt$ we get $\tilde Q_t^\top \Dmut^{-1} \tilde Q_t = \tilde\Lambda_t^{-1} $. Further denote $(Q_t)_{(ij,m)} $ to be the entry of $Q_t$ at row $ij$ and column $m$, where $i,j=1, \ldots,n$, $m=1, \ldots,{n(n+1)}/{2}$. Thus, for any nonzero eigenvalue $\lambda_m$, we have
$$\lambda_m^{-1} = \sum_{i<j} (Q_t)_{(ij,m)}^2 \mu_{t,ij}^{-1} + \sum_{i} (Q_t)_{(ii,m)}^2 (2\mu_{t,ij})^{-1} \leqslant \max(\mu_{t,ij}^{-1}), $$ i.e., $\lambda_m \geqslant \min(\mu_{t,ij}) $.

\subsubsection{Proof of Lemma \ref{lem_DZtDz_eigenstructure}}
\label{sec:pfA2}
While the proof is mainly direct verification of the statements in the Lemma, some basic facts that will be used repeatedly are first presented here:
\begin{equation}
\sum_{l=1}^n z_l U_{il} = Z^\top U_i = \left(\begin{smallmatrix}
V_1 & \cdots & V_k
\end{smallmatrix}\right)
\left(\begin{smallmatrix}
\text{diag}(\lambda_1, \ldots,\lambda_k) & 0_{k\times (n-k)} 
\end{smallmatrix}\right)
\left(\begin{smallmatrix}
U_1^\top \\ \vdots \\ U_n^\top
\end{smallmatrix}\right) U_i
= \left\{\begin{aligned}
&\lambda_i V_i, \quad & i\leqslant k;\\
&0, \quad & i>k.
\end{aligned}
\right.
\label{lem_dzeig_1}
\end{equation}
\begin{equation}
z_l^\top V_i = l\text{th row of }(ZV_i) = l\text{th row of }(\lambda_i U_i) = \lambda_i U_{il}. 
\label{lem_dzeig_2}
\end{equation}
\begin{equation}\begin{aligned}
&\DZt\DZ = \begin{pmatrix}
z_1z_1^\top + \sum_{l=1}^n z_lz_l^\top & z_2z_1^\top & \cdots & z_nz_1^\top \\
z_1z_2^\top & z_2z_2^\top + \sum_{l=1}^n z_lz_l^\top & \cdots & z_nz_2^\top \\
 & \vdots & \vdots & \\
 z_1z_n^\top  & z_2z_n^\top & \cdots & z_nz_n^\top + \sum_{l=1}^n z_lz_l^\top 
\end{pmatrix}\in \mathbb R^{nk\times nk}.
\end{aligned}
\label{lem_dzeig_3}
\end{equation}
Let $\Xi_Z = \DZt\DZ = (\Gamma_1 + \Gamma_2 ) $, where we define
\begin{align*}
\Gamma_1 = \begin{pmatrix}
\sum_{l=1}^n z_lz_l^\top &0 & \cdots & 0 \\
0 &  \sum_{l=1}^n z_lz_l^\top & \cdots &0 \\
 & \vdots & \vdots & \\
 0 & 0 & \cdots &  \sum_{l=1}^n z_lz_l^\top 
\end{pmatrix}, \quad
\Gamma_2 = \begin{pmatrix}
z_1z_1^\top  & z_2z_1^\top & \cdots & z_nz_1^\top \\
z_1z_2^\top & z_2z_2^\top  & \cdots & z_nz_2^\top \\
 & \vdots & \vdots & \\
 z_1z_n^\top  & z_2z_n^\top & \cdots & z_nz_n^\top 
\end{pmatrix}.
\end{align*}

We first prove statement 1 in the Lemma. The constraint on $Z$ could be written as $1_n^\top Z =1_n^\top U D V = 0$, which gives us $1_n^\top U_i = 0 $ for $i=1, \ldots,k$.

Next, we prove statement 2 in the Lemma. Since $\{\beta_{ij} : i=1, \ldots,n;j=1, \ldots,k \} $ is an orthonormal basis of $\mathbb R^{nk}$, it suffices to show that $Z_v^\top \beta_{ij} = 0 $ for $i\neq j$. Note that by \eqref{lem_dzeig_1},
\begin{align*}
Z_v^\top \beta_{ij} = \sum_{l=1}^n U_{il} z_l^\top V_j = \left\{\begin{aligned}
&\lambda_i V_i^\top V_j, \quad & i\leqslant k,\\
&0, \quad & i>k.
\end{aligned}
\right.
\end{align*}
As $V$ is an orthonormal matrix, the above equation shows that $Z_v^\top\beta_{ij} = \lambda_i  1_{(i=j)} $ where $ 1_{(i=j)}$ is the indicator function.

Finally, we prove statement 3 in the Lemma. From the fact that $\{\beta_{ij} : i=1, \ldots,n;j=1, \ldots,k \} $ form an orthonormal basis of $\mathbb R^{nk}$, we have that the $nk$ eigenvectors in statement 3 of the Lemma also form an orthonormal basis of $\mathbb R^{nk} $. It remains to check all the eigen equations given in statement 3. For $i=1, \ldots,k$,
\begin{align*}
\DZt\DZ \beta_{ii} &= \Gamma_1 \beta_{ii} + \Gamma_2 \beta_{ii} = \begin{pmatrix}
\sum_{l=1}^n z_lz_l^\top U_{i1} V_i \\ \vdots \\ \sum_{l=1}^n z_lz_l^\top U_{in} V_i
\end{pmatrix} + \begin{pmatrix}
\sum_{l=1}^n z_lz_1^\top U_{il}  V_i \\ \vdots \\ \sum_{l=1}^n  z_lz_n^\top U_{il} V_i
\end{pmatrix}\\
&= \begin{pmatrix}
\sum_{l=1}^n z_l \lambda_i U_{il} U_{i1} \\ \vdots \\ \sum_{l=1}^n z_l \lambda_i U_{il} U_{in} 
\end{pmatrix} + \begin{pmatrix}
\sum_{l=1}^n z_l \lambda_i U_{i1} U_{il} \\ \vdots \\ \sum_{l=1}^n z_l \lambda_i U_{in} U_{il} 
\end{pmatrix} =
2\begin{pmatrix}
\lambda_i^2 U_{i1}V_i \\ \vdots \\ \lambda_i^2 U_{in}V_i
\end{pmatrix} =2 \lambda_i^2 \beta_{ii},
\end{align*}
where in the third ``$=$'' we used \eqref{lem_dzeig_2} and in the fourth ``$=$'' we used \eqref{lem_dzeig_1}.
For $j=k+1, \ldots,n $ and $i=1, \ldots,k $,
\begin{align*}
\DZt\DZ \beta_{ji} &= \Gamma_1 \beta_{ji} + \Gamma_2 \beta_{ji} = \begin{pmatrix}
\sum_{l=1}^n z_lz_l^\top U_{j1} V_i \\ \vdots \\ \sum_{l=1}^n z_lz_l^\top U_{jn} V_i
\end{pmatrix} + \begin{pmatrix}
\sum_{l=1}^n z_lz_1^\top U_{jl}  V_i \\ \vdots \\ \sum_{l=1}^n  z_lz_n^\top U_{jl} V_i
\end{pmatrix}\\
&=\begin{pmatrix}
\sum_{l=1}^n z_l \lambda_i U_{il} U_{j1} \\ \vdots \\ \sum_{l=1}^n z_l \lambda_i U_{il} U_{jn} 
\end{pmatrix} + \begin{pmatrix}
\sum_{l=1}^n z_l \lambda_i U_{i1} U_{jl} \\ \vdots \\ \sum_{l=1}^n z_l \lambda_i U_{in} U_{jl} 
\end{pmatrix} =
\begin{pmatrix}
\lambda_i^2 U_{j1}V_i \\ \vdots \\ \lambda_i^2 U_{jn}V_i
\end{pmatrix} + 0 = \lambda_i^2 \beta_{ji}.
\end{align*}
For $1\leqslant i \neq j \leqslant k $,
\begin{align*}
\DZt\DZ \beta_{ij} &= \Gamma_1 \beta_{ij} + \Gamma_2 \beta_{ij} = \begin{pmatrix}
\sum_{l=1}^n z_lz_l^\top U_{i1} V_j \\ \vdots \\ \sum_{l=1}^n z_lz_l^\top U_{in} V_j
\end{pmatrix} + \begin{pmatrix}
\sum_{l=1}^n z_lz_1^\top U_{il}  V_j \\ \vdots \\ \sum_{l=1}^n  z_lz_n^\top U_{il} V_j
\end{pmatrix}\\
&= \begin{pmatrix}
\sum_{l=1}^n z_l \lambda_j U_{jl} U_{i1} \\ \vdots \\ \sum_{l=1}^n z_l \lambda_j U_{jl} U_{in} 
\end{pmatrix} + \begin{pmatrix}
\sum_{l=1}^n z_l \lambda_j U_{j1} U_{il} \\ \vdots \\ \sum_{l=1}^n z_l \lambda_j U_{jn} U_{il} 
\end{pmatrix} \\
&=
\begin{pmatrix}
\lambda_j^2 U_{i1}V_j \\ \vdots \\ \lambda_j^2 U_{in}V_j
\end{pmatrix} + \begin{pmatrix}
\lambda_i\lambda_j U_{j1} V_i \\ \vdots  \\ \lambda_i\lambda_j U_{jn} V_i
\end{pmatrix}
= \lambda_j^2 \beta_{ij} + \lambda_i\lambda_j\beta_{ji}.
\end{align*}
Thus we have
\begin{align*}
\DZt\DZ (\lambda_j \beta_{ij} + \lambda_i \beta_{ji}) = \lambda_j^3 \beta_{ij} + \lambda_i\lambda_j^2 \beta_{ji} + 
\lambda_i^3 \beta_{ji} + \lambda_j\lambda_i^2 \beta_{ij} 
= (\lambda_i^2+\lambda_j^2) (\lambda_j \beta_{ij} + \lambda_i \beta_{ji}) 
\end{align*}
and
\begin{align*}
\DZt\DZ (\lambda_i \beta_{ij} - \lambda_j \beta_{ji}) = \lambda_i\lambda_j^2 \beta_{ij} + \lambda_i^2\lambda_j \beta_{ji} - 
\lambda_i^2\lambda_j \beta_{ji} - \lambda_j^2\lambda_i \beta_{ij} 
= 0. 
\end{align*}

\subsubsection{Proof of Lemma \ref{lm:presentuzmain} }
\label{sec:pfA3}

\paragraph*{Proof of (i) and (ii).} 
To prove (i) and (ii), 
it suffices to show that 
(a) for any $v\in \text{col}(\UZ)^\perp$, $\Ieff v =0$, 
and (b) for any $v \in$ col$(\UZ)$, $v^{\top}\Ieff v$ has a positive lower bound given as in (i).

Part (a): By Definition \ref{def::UZA} in the main text and Lemma \ref{lem_DZtDz_eigenstructure}, we have
\begin{align}\label{eq:columnuzdecom}
    \text{col}(\UZ)^\perp = \text{span}(\mathcal{B}_{N,1}) \oplus \text{span}(\mathcal{B}_{N,2}),
\end{align} where $\oplus$ denotes orthogonal direct sum. 
To prove (a), we next show that $\mathrm{span}(\mathcal{B}_{N,1})\subseteq \text{null}(\Ieff)$
and $\mathrm{span}(\mathcal{B}_{N,2})\subseteq \text{null}(\Ieff) $, respectively. 
First, 
by the formula of $\Ieff$ in \eqref{eq:IeffHeffinXt}, 
it is obvious that 
$\text{null}(\DZ) \subseteq \text{null}(\Ieff)$. 
By Lemma \ref{lem_DZtDz_eigenstructure}, we know that 
$\mathrm{span}(\mathcal{B}_{N,2})=\text{null}(\DZ)$.
Therefore, $\mathrm{span}(\mathcal{B}_{N,2}) \subseteq \text{null}(\Ieff)$. 
Second, recall that we set  $U_{k+1} = 1_n/\sqrt{n} $ in \eqref{eq:uk+1def1n}.  
Then we have 
\begin{align}\label{eq:spanbn1col}
    \mathrm{span}(\mathcal{B}_{N,1})=\text{col}(\myDtri)
\end{align}
where we define
\begin{align*}
\myDtri = \begin{pmatrix}
e_1^{(k)} & e_2^{(k)} & \cdots & e_k^{(k)}\\
\vdots & \vdots & \cdots & \vdots\\
e_1^{(k)} & e_2^{(k)} & \cdots & e_k^{(k)}
\end{pmatrix} = \begin{pmatrix}
\mathrm{I}_{k\times k}\\\vdots\\\mathrm{I}_{k\times k}
\end{pmatrix} \in \mathbb R^{nk\times k}
\end{align*}
where $e_i^{(k)} \in \mathbb R^k $ is the vector whose $i$th element is 1 and other elements are 0. 
To prove $\mathcal{B}_{N,1}\subseteq \text{null}(\Ieff)$, 
it suffices to prove 
$\Ieff \myDtri = 0$. 
In particular, by \eqref{eq:IeffHeffinXt}, 
\begin{align}
    \Ieff \myDtri = \sum_{t=1}^T \DZt\myXt\DZ  \myDtri = \sum_{t=1}^T \DZt\myXt \myDone Z = 0
    \label{eq::IeffDI=0}
\end{align}
where the second equation follows by
the fact $$ \DZ  \myDtri  = \myDone Z,$$ which is obtained  
 by the formulae in Section \ref{sec:derieffscore},
 and the third equation follows by 
 $ \myXt  \myDone=0$ in Lemma \ref{lem:property_myXt}.

Part (b): 
Consider any $v \in$ col$(\UZ)$. 
We first show that 
\begin{align}\label{eq:ieffquadineq}
     v^{\top} \Ieff v =\sum_{t=1}^T v^\top \DZt \myXt \DZ v
    \geqslant T\min_{t,i,j}(\mu_{t,ij}) \|\DZ v\|_2^2
\end{align}
where the first equation follows by the definition of $\Ieff$. 
To prove the inequality in \eqref{eq:ieffquadineq}, 
by Lemma \ref{lem:property_myXt}, 
it suffices to show that $\DZ v \in \operatorname{null}(\myXt)^{\perp}$.
As $\operatorname{null}(\myXt) =  \text{col}(\myDone)$ by Lemma  \ref{lem:property_myXt}, 
it suffices to show that $\myDone^{\top}\DZ v =0.$ 
In particular, we have
\begin{align*}
    \myDone^{\top}\DZ v = \myDonet \DZ v =Z\myDtri^\top v =0 
\end{align*}
where in the second equation, we use the fact that
\begin{align*}
    \myDonet \DZ =  Z\myDtri^\top
\end{align*}
which follows from  formulae in Section \ref{sec:derieffscore}, and 
in the third equation, we use $ \myDtri^\top v = 0$ for $v\in \operatorname{col}(\UZ)$,
which follows from \eqref{eq:columnuzdecom} and \eqref{eq:spanbn1col}. 



We next derive a lower bound for $\|\DZ v\|_2^2$ in \eqref{eq:ieffquadineq}. 
For any $v\in \operatorname{col}(\UZ)$, 
\begin{align}\label{eq:dzvquadlowbd}
    \|\DZ v\|_2^2 = v^\top \DZt\DZ v \geqslant  \sigma_{\min}(Z^\top Z) \|v\|_2^2, 
\end{align}
where the inequality follows by 
 Lemma \ref{lem_DZtDz_eigenstructure}. 
 In summary, combining   \eqref{eq:ieffquadineq} and  \eqref{eq:dzvquadlowbd}, 
  $$\sigma_{\min}(\UZ^{\top} \Ieff(Z,\alpha) \UZ) \geqslant nT\, \sigma_{\operatorname{min}}( Z^\top  Z /n)\min_{{1\leqslant i,j\leqslant n, 1\leqslant t\leqslant T}}\exp(z_i^\top  z_j +\alpha_{it} +  \alpha_{jt}).$$

\paragraph*{Proof of (iii).}

Let $\UZ^\perp$ be an $nk\times k(k+1)/2$ matrix with columns vectors in 
$\mathcal{B}_{N,1}$ and $\mathcal{B}_{N,2}$ defined as in \eqref{eq:bn1bn2setsupp}. 
Then column vectors of  $\UZ^\perp$  are orthonormal and orthogonal to the columns of $\UZ$, i.e., 
$(\UZ^\perp)^\top \UZ^\perp = \mathrm{I}_{k(k+1)/2} $ and $\UZ^\top \UZ^\perp = 0 $. 
We first show that $\UZ^{\perp\top} \DZt \myXt=0$, so that
\begin{align}\label{eq:spaceprojres}
   \DZt\myXt  = \UZ\UZ^\top \DZt\myXt + \UZ^{\perp}\UZ^{\perp\top} \DZt\myXt =  \UZ\UZ^\top \DZt\myXt. 
\end{align}
By \eqref{eq:columnuzdecom},
to show $\UZ^{\perp\top} \DZt\myXt=0$, or equivalently $ \myXt \DZ \UZ^\perp = 0  $, it suffices to show $\text{span}(\mathcal{B}_{N,1})\oplus\text{span}(\mathcal{B}_{N,2}) \subseteq \text{null}(\myXt \DZ)$. By Lemma \ref{lem_DZtDz_eigenstructure}, we know that 
$\mathrm{span}(\mathcal{B}_{N,2})=\text{null}(\DZ) \subseteq \text{null}(\myXt \DZ)  $. By the proof same as in \eqref{eq::IeffDI=0}, we have $\myXt\DZ\myDtri = 0 $, so that $\mathrm{span}(\mathcal{B}_{N,1})\subseteq \text{null}(\myXt \DZ)  $
by \eqref{eq:spanbn1col}. Thus we have proved \eqref{eq:spaceprojres}.

Now recall that by  definitions in \eqref{eq:seffformulaesupp} and  \eqref{eq:IeffHeffinXt}, we have
\begin{align*}
     \Seff =
     \sum_{t=1}^T \DZt \myXt \Dmut^{-1} \dot{L}_{\Theta_t},\quad\quad \text{and}\quad \quad 
\Ieff = \sum_{t=1}^T  \DZt  \myXt \DZ. 
 \end{align*} 
By \eqref{eq:spaceprojres}, and the formulae of $\Seff$ and $\Ieff$ above, we obtain that
$$\UZ\UZ^\top \Seff(Z,\alpha) =  \Seff(Z,\alpha),\quad\quad \text{and}\quad \quad  \UZ\UZ^\top \Ieff(Z,\alpha) =  \Ieff(Z,\alpha).$$

\newpage 
\section{Technical Results for the One-Step  Estimator} \label{sec:theoryonestep}

In this section, we prove theoretical results for the one-step estimator in Section \ref{sec:newton} including Theorem \ref{thm_NRerr} and Remark \ref{sec:pfRMK}. 



\subsection{Proof of Theorem \ref{thm_NRerr}}\label{sec:pfthmnerr}
By the definitions,   
$\mathrm{dist}(\hat{Z}, Z^{\star}) \leqslant \|\hat Z_v -  (\ZQs)_v\|_2$ with $\check Q$ in \eqref{eq:defcheckq}.
To prove Theorem \ref{thm_NRerr},  
it suffices to establish an upper bound of $\|\hat Z_v -  (\ZQs)_v\|_2$.
However, directly studying the asymptotic expansion of $\hat Z_v -  (\ZQs)_v$ is challenged by the singularity  of  
 $\Ieff(\check{Z},\check{\alpha})$ in Remark \ref{lb:singularity}.  
To overcome the issue, 
we establish the proof on a properly constructed  subspace of $ \mathbb{R}^{nk}$. 
Specifically, 
given an initial estimate $\check{Z}\in \mathbb{R}_0^{n\times k}$, we  construct $\UZb\in \mathbb{R}^{nk\times \rankU}$ following Definition \ref{def::UZA} in Section \ref{sec:eigenspace}. We  use $\mathrm{col}(\UZb)\subset \mathbb{R}^{nk}$ as the restricted eigenspace,  
which turns out to contain all the intrinsic information for characterizing the estimation error $ \|\hat Z_v -  (\ZQs)_v\|_2$ (see, e.g., \eqref{pfthm1-3A}). 
We next present an outline of the proof and give more details afterwards. 

   \textit{Step 1:} We  show $\hat Z_v - (\ZQs)_v \, \in \,  \text{col}(\UZb)  $. 
This suggests that the error vector always lies in the restricted space $\text{col}(\UZb)$, and thus
\begin{align} \label{pfthm1-3A}
\big\| \hat Z_v -  (\ZQs)_v \big\|_2= \big\|\UZb^\top ( \hat Z_v -  (\ZQs)_v) \big\|_2.
\end{align}


\textit{Step 2:}
Define $\check{I}_{e,U}=\UZb^\top \Ieff(\check{Z}, \check{\alpha})\UZb /(nT)\in \mathbb{R}^{\rankU \times \rankU}$ as the  $\UZb$-restricted efficient information matrix,
and let $\sigma_{\min}(\check{I}_{e,U})$ represent its minimum eigenvalue.
We prove that when $n$ is sufficiently large, 
$\check{I}_{e,U}$ is positive-definite,
and $1/\sigma_{\min}(\check{I}_{e,U})$ has a constant upper bound.   


\textit{Step 3:} 
By Steps 1 and 2, we have
\begin{align*}
   \eqref{pfthm1-3A} = \big\| \check{I}_{e,U}^{-1}\ \check{I}_{e,U}\   \UZb^\top ( \hat Z_v -  (\ZQs)_v) \big\|_2
    \leqslant  \frac{1}{\sigma_{\min}( \check{I}_{e,U})}\big\|\check{I}_{e,U}\   \UZb^\top ( \hat Z_v -  (\ZQs)_v)\big\|_2.
\end{align*}
We prove that with high probability, the numerator $\big\|\check{I}_{e,U}\,   \UZb^\top ( \hat Z_v -  (\ZQs)_v)\big\|_2$ is of the order of 
$\max\{T^{-1/2}, n^{-1/2}\}$ 
up to logarithmic factors.


\textit{Step 4:}
We combine Steps 1--3 to obtain Theorem \ref{thm_NRerr}. 



\vspace{5pt}
We next present more details of Steps 1--4. 

\vspace{5pt}
 \textbf{{Step 1:}} 
For the proposed estimator \eqref{newtonsolUZ}, 
$\check{\mathcal U}$ can be taken as any ${nk\times \rankU}$ matrix such that 
$\mathrm{col}(\check{\mathcal U})=\mathrm{col}(\Ieff(\check{Z},\check{\alpha}))$. 
By Lemma \ref{lem_presentingUZmain}-(i), 
we take, without loss of generality,  
$\check{\mathcal{U}}$ in 
\eqref{newtonsolUZ} as $\UZb$ and then obtain
\begin{align}\label{eq:newtonsolUZ_diffA} 
     \hat Z_v - (\ZQs)_v = \check Z_v-  (\ZQs)_v + \UZb\left\{ \UZb^\top \Ieff(\check Z, \check\alpha) \UZb \right\}^{-1} \UZb^\top \Seff(\check Z,\check\alpha).
\end{align}
In \eqref{eq:newtonsolUZ_diffA}, we know ${\check Z}_v\in$ col$(\UZb)$ 
 by Lemma \ref{lem_DZtDz_eigenstructure}, 
 and $\UZb\{ \UZb^\top \Ieff(\check Z, \check\alpha) \UZb \}^{-1} \UZb^\top \Seff(\check Z,\check\alpha) \in $ col$(\UZb)$  by its formula. 
To prove Step 1, 
it remains to show  $(\ZQs)_v\in \mathrm{col}(\UZb)$. 

We next prove that for any vector 
$a\in \mathrm{col}(\UZb)^{\perp}$, $a^{\top} (\ZQs)_v =0$ by the bases representation in \eqref{eq:bn1bn2setsupp}. 
Define $\check{\beta}_{ij}$'s similarly as in \eqref{eq:betaijdefA} with  the singular value decomposition  $\check{Z}=\sum_{i=1}^n \check{\lambda}_i\check{U}_i\check{V}_i^{\top}$. 
By $\check{\beta}_{k+1,i}=1_n\otimes \check{V}_i/\sqrt{n}$ and the property of Kronecker product, we have
$
     \check{\beta}_{k+1,i}^{\top} (\ZQs)_v  =  1_n^{\top} \ZQs \check{V}_i /\sqrt{n} = 0,  
$ 
where the second equation follows by $1_n^{\top}Z^{\star}=0$. 
Moreover, by  $\check{\beta}_{ij}=\check{U}_i\otimes \check{V}_j$ and the property of Kronecker product, 
\begin{align}
    (\check{\lambda}_i\check{\beta}_{ij}-\check{\lambda}_j\check{\beta}_{ji})^{\top} (\ZQs)_v =&~  \check{\lambda}_i\check{U}_i^{\top}(\ZQs )\check{V}_j- \check{V}_i^{\top}(\ZQs)^{\top}\check{U}_j\check{\lambda}_j\notag\\
=&~\check{V}_i^{\top} (\check{Z}^{\top}Z^{\star}\check{Q}) \check{V}_j-\check{V}_i^{\top} (\check{Z}^{\top}Z^{\star}\check{Q})^{\top}\check{V}_j=0,\notag
\end{align}
where in the second equation, we use $ \check{\lambda}_i\check{U}_i=\check{Z}\check{V}_i$,
and in the third equation, we use the symmetricity of $\check{Z}^{\top}Z^{\star}\check{Q}$, which follows from   the  property of
the orthogonal Procrustes problem and the definition \eqref{eq:defcheckq}.


\vspace{5pt}
 {\textbf{Step 2:}} 
By the definition of $ \check{I}_{e,U}$ and Lemma \ref{lem_presentingUZmain}-(i),
we have 
\begin{align} \label{pfthm1-2A}
    \sigma_{\operatorname{min}}(\check{I}_{e,U})= \sigma_{\operatorname{min}}(\UZb^\top \Ieff(\check Z,\check\alpha) \UZb/(nT)) \geqslant \check M_{ Z,2} \exp(-\check{M}_{\Theta,2}),
\end{align}
where we define  $\check M_{ Z,2} = \sigma_{\operatorname{min}}(\check Z^\top \check Z /n)$
and $\check M_{\Theta,2}=\max_{1\leqslant i,j\leqslant n, 1\leqslant t\leqslant T}\{| \check z_i^\top \check z_j + \check \alpha_{it} + \check \alpha_{jt} |\}$. 
By Conditions \ref{cond:truvalueregularity} and \ref{cond_elem_init}, we have 
$\check M_{ Z,2} \geqslant M_{Z,2} - 3{M_{Z,1}^{1/2}} M_bn^{-1/2}\log^{\epsilon}(nT), $
and $ \check M_{\Theta,2} \leqslant M_{\Theta,2}+ 3(1+{M_{Z,1}^{1/2}})M_bn^{-1/2}\log^{\epsilon}(nT).$
When $n / \log^{2\epsilon}(nT)$  is sufficiently large,  $1/  \sigma_{\operatorname{min}}(\check{I}_{e,U})$ has a constant upper bound as $M$'s are constants.


\vspace{5pt}
 {\textbf{Step 3:}} By  the definition of $\check{I}_{e,U}$ and \eqref{eq:newtonsolUZ_diffA},  
\begin{align}
     \big\|\check{I}_{e,U}\   \UZb^\top ( \hat Z_v -  (\ZQs)_v)\big\|_2= &~\big\| \UZb^{\top}\big\{\Ieff(\check Z, \check\alpha) (\check Z_v-  (\ZQs)_v ) +   \Seff(\check{Z},\check{\alpha})\big\}\big\|_2 /(nT)\notag\\   
    \leqslant &~\big\| \Ieff(\check Z, \check\alpha) (\check Z_v-  (\ZQs)_v ) +   \Seff(\check{Z},\check{\alpha})\big\|_2/(nT)  \label{eq:multipliediff2A}
\end{align}
 where the second equation follows by $\|\UZb\|_{\operatorname{op}} = 1$.
Define
\begin{align}
    S_1 =&~   \Ieff(\check Z, \check\alpha) (\check Z_v-  (\ZQs)_v ), \hspace{2.6em} S_2=  \Seff(\check Z,\check\alpha)- \Seff(\ZQs,\check\alpha),\label{eq:thm1s1termA}\\
S_3=&~ \Seff( \ZQs,\check\alpha)- \Seff(\ZQs,\alpha^\star). \label{eq:thm1s3termA} 
\end{align}
Then we have
\begin{align*}
     \eqref{eq:multipliediff2A} =&~\|\Seff(Z^\star \check{Q},\alpha^\star)+S_1+S_2+S_3\|_2/(nT) \notag\\
\leqslant &~\left\|  \Seff(Z^\star \check{Q},\alpha^\star)\right\|_2/(nT)+\sqrt{nk}   \left\|  S_1+S_2\right\|_{\infty} /(nT) +\sqrt{nk}\left\|S_3 \right\|_\infty/(nT)  . 
\end{align*}
We can show that for any $s>0$, 
there exists a universal constant $C_s$ such that
\begin{align}
\|\Seff(\ZQs, \alpha^\star )\|_2 /(nT)
\leqslant &~
C_{s}  \, \big( T^{-1/2}\big)\, \log^{1/2}(nk),
\label{eq:Seff_starstar_main}
\\
\sqrt{nk}\|S_1+S_2\|_\infty/(nT) 
\leqslant  &~C_{s}  \,  \big(n^{-1/2}+ T^{-1/2}\big)\, \log^{2\varsigma }(nT), 
\label{eq:S1+S2_infty_main}
\\
\sqrt{nk}\|S_3\|_\infty /(nT)
\leqslant &~ C_{s} \,  \big(n^{-1/2}\big)\, \log^{2\varsigma }(nT),
\label{eq:S3_infty_main}
\end{align}
with probability $1-O(n^{-s})$  (see proofs in Section \ref{sec:pfthm1} in the Supplementary Material).
In summary, we obtain that for any $s>0$, there exists a universal constant $C_s$ such that 
\begin{align}\label{pfthm1-12A}
    \big\|\check{I}_{e,U}\   \UZb^\top ( \hat Z_v -  (\ZQs)_v)\big\|_2 \leqslant C_{s} \max\{   T^{-1/2} ,n^{-1/2}\}\log^{2\varsigma}(nT).  
\end{align}
 with probability $1-O(n^{-s}) $. 

\vspace{5pt}

{\textbf{Step 4:}} 
By \eqref{pfthm1-3A}, \eqref{pfthm1-2A},  and \eqref{pfthm1-12A}, we have that when  $n / \log^{2\varsigma}(T)$  is sufficiently large, for any $s>0$, there exists a universal constant $C_{s}>0$ such that
\begin{align*}
    \text{dist}^2(\hat Z, Z^\star) &\leqslant   \big\| \hat Z_v -  (\ZQs)_v \big\|_2^2\leqslant C_s\max\{T^{-1},n^{-1}\}\log^{4\varsigma}(nT) = \frac{1}{T} C_s  r_{n,T} 
\end{align*}
with probability $1 - O(n^{-s})$. 

\begin{remark}\label{rm:uzunique}
The proof in Section \ref{sec:pfthmnerr}  only relies on the restricted eigenspace $\mathrm{col}(\UZb)$ but not the original matrix $\UZb$. 
Therefore, even when  $\UZb$  changes, the above analysis remains the same as long as $\mathrm{col}(\UZb)$ does not change. 
For instance, 
permuting the columns of the matrix $\UZb$ would change the matrix but not its column space.   
\end{remark}

\subsection{Proof of Eq. \texorpdfstring{\eqref{eq:Seff_starstar_main}--\eqref{eq:S3_infty_main} in Section \ref{sec:pfthmnerr}}{}}\label{sec:pfthm1}  
To finish the proof of Theorem \ref{thm_NRerr},
it remains to prove 
Eq. \eqref{eq:Seff_starstar_main}--\eqref{eq:S3_infty_main} 
in Section \ref{sec:pfthmnerr}.
We next present proofs of Eq. \eqref{eq:S3_infty_main}, Eq. \eqref{eq:S1+S2_infty_main}, and Eq. \eqref{eq:Seff_starstar_main}  in Sections \ref{sec:pfB2}, \ref{sec:pfB3},  and \ref{sec:pfB4}, respectively.

\subsubsection{Proof of Eq. \eqref{eq:S3_infty_main}} \label{sec:pfB2}
For the ease of presentation, 
let $\ZQssub=\ZQs$ and $(\ZQssub)_v=(\ZQs)_v$. Define $W(Z,\alpha) =  \ddot{L}_{Z,\alpha}(Z,\alpha) \ddot{L}^{-1}_{\alpha,\alpha}(Z,\alpha)$.
We have 
\begin{align*}
S_3 &= \dot{L}_Z( \ZQssub,\check\alpha) - W(\ZQssub,\check\alpha) \dot{L}_\alpha( \ZQssub,\check\alpha) - \dot{L}_Z(\ZQssub,\alpha^\star)+ W(\ZQssub,\alpha^\star) \dot{L}_\alpha(\ZQssub,\alpha^\star)\\
&= S_{31} - S_{32} - S_{33}
\end{align*}
where we define
\begin{align*}
    S_{31} =  &~ \dot{L}_Z(\ZQssub,\check\alpha)  - \dot{L}_Z(\ZQssub,\alpha^\star),\\
   S_{32} = &~  \big\{W(\ZQssub,\check\alpha) - W(\ZQssub,\alpha^\star) \big\}\dot{L}_\alpha(\ZQssub,\check\alpha), \\
   S_{33} = &~  W(\ZQssub,\alpha^\star)  \big\{\dot{L}_\alpha(\ZQssub,\check\alpha)-\dot{L}_\alpha(\ZQssub,\alpha^\star) \big\}.
\end{align*}
\paragraph*{Proof Outline}
Deriving the upper bound of $S_3$ consists of two steps. 

\noindent \textit{Step 1} derives an upper bound of $S_{31}-S_{33}$.  
Intuitively, by Taylor expansion, 
\begin{align}
S_{31} =&~ \ddot{L}_{Z,\alpha}(\ZQssub,\alpha^\star)(\check{\alpha}_v- \alpha_v^{\star})+ \text{higher-order terms},
\label{eq:S_31_heuristic} \\
S_{33} =&~ \ddot{L}_{Z,\alpha}(\ZQssub,\alpha^\star) \ddot{L}_{\alpha,\alpha}^{-1}(\ZQssub,\alpha^\star) \big[\ddot{L}_{\alpha,\alpha}(\ZQssub,\alpha^\star)(\check{\alpha}_v- \alpha_v^{\star}) + \text{higher-order terms} \big]\notag\\
=&~ \ddot{L}_{Z,\alpha}(\ZQssub,\alpha^\star)(\check{\alpha}_v- \alpha_v^{\star}) + \text{higher-order terms},
\label{eq:S_33_heuristic} 
\end{align}
 where we plug in $W(\ZQssub,\alpha^\star)=\ddot{L}_{Z,\alpha}(\ZQssub,\alpha^\star) \ddot{L}_{\alpha,\alpha}^{-1}(\ZQssub,\alpha^\star)$. 
Then \eqref{eq:S_31_heuristic}--\eqref{eq:S_33_heuristic}
imply $S_{31}-S_{33}=0 + \text{higher-order terms}$. 
The proof below will show that 
the first-order terms in $S_{31}-S_{33}$ would be 0 and then examine the higher-order terms. 


\medskip
\noindent  \textit{Step 2} derives an upper bound of $S_{32}$. 
It consists of two substeps: 
\textit{Step 2-1} derives an upper bound of $W(\ZQssub,\check\alpha) - W(\ZQssub,\alpha^\star)$, and 
\textit{Step 2-2} derives an upper bound of $\dot{L}_\alpha(\ZQssub,\check\alpha)$.

 

\paragraph*{Step 1} In this step we consider $S_{31} - S_{33} $.
Recall we denote $\Theta_t = ZZ^\top + \alpha_t 1_n^\top + 1_n \alpha_t^\top$. Besides, 
By the Fact \ref{fact:partialderivfunction}, 
\begin{align*}
S_{31}=&~\dot{L}_Z(\ZQssub,\check\alpha)  - \dot{L}_Z(\ZQssub,\alpha^\star) = \sum_t\left.\left( \frac{\partial \Theta_t}{\partial Z}^\top \dot{L}_{\Theta_t}\right)\right|_{(\ZQssub{},\check\alpha ) } - \sum_t\left.\left( \frac{\partial \Theta_t}{\partial Z}^\top \dot{L}_{\Theta_t} \right)\right|_{(\ZQssub{},\alpha^\star ) } 
\\ =&~ \DZt(\ZQssub) \sum_t\big[ \myUt(\ZQssub,\check\alpha) - \myUt(\ZQssub,\alpha^\star) \big].
\end{align*}
Let $d_{t,ij} = \alpha_{it}+\alpha_{jt},\Delta{\alpha_{it}}= \check\alpha_{it} - \alpha_{it}^\star $, and $\Delta d_{t,ij} = \Delta{\alpha_{it}}+\Delta{\alpha_{jt}} $. Let 
\begin{equation} \label{eq::def_varrhoij}
    \varrho_{ij} = 1+(\sqrt2-1)1_{(i=j)}
\end{equation}
so that the $(i,j)$th element of $\myUt\in \mathbb R^{{n(n+1)}/{2}}$ is $\varrho_{ij} (A_{t,ij} - \mu_{t,ij})$. By Taylor expansion, each element of $\myUt(\ZQssub,\check\alpha) - \myUt(\ZQssub,\alpha^\star)$ can be decomposed as
\begin{align}\label{eq:mu_Taylorexp}
    \dot{L}_{\Theta_{t,ij}} {(\ZQssub{},\check\alpha ) } -\dot{L}_{\Theta_{t,ij}} {(\ZQssub{},\alpha^\star ) } 
    =&- \varrho_{ij} \Big\{\exp(G_{ij}^{\star}+ \check{d}_{t,ij} ) - \exp(G_{ij}^{\star} + d_{t,ij}^{\star} )\Big\}\notag\\
    =&-\varrho_{ij} \Big\{\exp(G_{ij}^{\star} + d_{t,ij}^{\star})\Delta d_{t,ij} +\exp( G_{ij}^{\star} + \tilde{d}_{t,ij}) (\Delta d_{t,ij})^2/2\Big\},
\end{align}
where $\tilde{d}_{t,ij}$ denotes a value between $ d_{t,ij}^{\star} $ and $\check d_{t,ij}$. 
Let $$\myUt(\ZQssub,\check\alpha) - \myUt(\ZQssub,\alpha^\star) = -(\Delta L_{1t} + \Delta L_{2t} ) ,$$ where $\Delta L_{1t}$ is an $n(n+1)/2$ dimensional vector with $(i,j)$th element being $\varrho_{ij} \exp(G_{ij}^{\star} + d_{t,ij}^{\star})\Delta d_{t,ij} $, and $\Delta L_{2t}$ is an $n(n+1)/2$ dimensional vector with $(i,j)$th element being $\varrho_{ij} \exp(G_{ij}^{\star} + \tilde d_{t,ij})(\Delta d_{t,ij})^2/2$.
We define $S_{31} = -(S_{311} + S_{312} ) $ with
\begin{align*}
S_{311}:= \DZt(\ZQssub) \sum_t \Delta L_{1t}
\quad\text{and}\quad S_{312}:= \DZt(\ZQssub) \sum_t \Delta L_{2t}.
\end{align*}
Here $\Delta L_{1t}$ corresponds to $\ddot{L}_{\Theta_t,\alpha_t} \times \Delta \alpha_t $
so that $S_{311}$ corresponds to the first order term $\ddot{L}_{Z,\alpha} \times \Delta\alpha $ and $S_{312}$ corresponds to the smaller order term in the heuristic argument \eqref{eq:S_31_heuristic}.
To see an elementwise bound for $S_{312} $, we consider the matrix $l_\infty$ norm defined by $\|A\|_\infty = \sup_{x\neq 0} \|Ax\|_\infty / \|x\|_\infty = \max_{1\leqslant i\leqslant m} \sum_{j=1}^n |a_{ij}| $ for $A\in \mathbb R^{m\times n} $.
We have $$\|\Delta L_{2t}\|_\infty \leqslant Ce^{-M_{\Theta,1}} (M_b\log^{\epsilon}(nT) /\sqrt{n})^2, $$ 
where $C$ is a universal constant whose value may change from line to line, and we use 
\begin{equation}
\label{eq:checkM_bd}
    \exp(G^\star_{ij} + \tilde d_{t,ij}) \leqslant \max_{t,i,j}\Big[ \max\big\{\exp(G^\star_{ij} +  d^\star_{t,ij}), \exp(G^\star_{ij} + \check d_{t,ij}) \big\} \Big]\leqslant Ce^{-M_{\Theta,1}}
\end{equation} in which the second inequality is due to  Condition \ref{cond_elem_init}.
Hence
\begin{align*}
\|S_{312}\|_\infty \leqslant&~ \sum_t \|\DZt(\ZQssub)\|_\infty \|\Delta L_{2t}\|_\infty
\\\leqslant&~
C T n \sqrt{M_{Z,1}} e^{-M_{\Theta,1}} (M_b\log^{\epsilon}(nT)/\sqrt n)^2 = CM_{Z,1}^{\frac12} M_b^2 e^{-M_{\Theta,1}}T\log^{2\epsilon}(nT) 
\end{align*}

Next we rigorously show the intuition in \eqref{eq:S_31_heuristic} and \eqref{eq:S_33_heuristic} that the first order term of $S_{31} $, i.e., $S_{311}$, is canceled with the first order term in $S_{33}$. 
We have
\begin{align*}
S_{33} =&~ W(\ZQssub,\alpha^\star) \times \big\{\dot{L}_\alpha(\ZQssub,\check\alpha)-\dot{L}_\alpha(\ZQssub,\alpha^\star) \big\} \\
=&~  \sum_t \ddot{L}_{Z,\alpha_t}(\ZQssub,\alpha^\star) \ddot{L}^{-1}_{\alpha_t,\alpha_t}(\ZQssub,\alpha^\star)  \bigg\{\bigg.\bigg( \frac{\partial \Theta_t}{\partial \alpha_t}^\top \dot{L}_{v_{\Theta_t}} \bigg)\bigg|_{(\ZQssub{},\check\alpha ) } - \sum_t\bigg.\bigg( \frac{\partial \Theta_t}{\partial \alpha_t}^\top \dot{L}_{v_{\Theta_t}} \bigg)\bigg|_{(\ZQssub{},\alpha^\star ) } \bigg\} \\
=&~ \sum_t \ddot{L}_{Z,\alpha_t}(\ZQssub,\alpha^\star) \ddot{L}^{-1}_{\alpha_t,\alpha_t}(\ZQssub,\alpha^\star) \myDonet \big[ \myUt(\ZQssub,\check\alpha) - \myUt(\ZQssub,\alpha^\star) \big] \\
=&~ - \sum_t \ddot{L}_{Z,\alpha_t}(\ZQssub,\alpha^\star) \ddot{L}^{-1}_{\alpha_t,\alpha_t}(\ZQssub,\alpha^\star) \myDonet \big\{ \Delta L_{1t} + \Delta L_{2t} \big\}.
\end{align*}
Define $S_{33} = -(S_{331} + S_{332} ) $ with
\begin{align*}
S_{331} :=&~ \sum_t \ddot{L}_{Z,\alpha_t}(\ZQssub,\alpha^\star) \ddot{L}^{-1}_{\alpha_t,\alpha_t}(\ZQssub,\alpha^\star) \myDonet  \Delta L_{1t}, \\
S_{332} :=&~ \sum_t \ddot{L}_{Z,\alpha_t}(\ZQssub,\alpha^\star) \ddot{L}^{-1}_{\alpha_t,\alpha_t}(\ZQssub,\alpha^\star) \myDonet \Delta L_{2t},
\end{align*}
so that analogous to the decompostion of $S_{31} $, $S_{331}$ corresponds to the first order term and $S_{332}$ corresponds to the smaller order term in the heuristic argument \eqref{eq:S_33_heuristic}.

It is straightforward to verify that $\Delta L_{1t}$ has the matrix multiplication expression
\begin{align*}
\Delta L_{1t}
= -\myCt(Z^\star,\alpha^\star) \myDone \times \Delta{\alpha_t}
\end{align*}
where $\Delta{\alpha_t}$ is an $n$-dimensional vector whose $i$th element is $\Delta{\alpha_{it}}$.
Therefore, we have 
\begin{align*}
S_{331} =&~ \sum_t \ddot{L}_{Z,\alpha_t}(\ZQssub,\alpha^\star) \ddot{L}^{-1}_{\alpha_t,\alpha_t}(\ZQssub,\alpha^\star) \myDonet  \Delta L_{1t}, \\
=&~ -\sum_t \ddot{L}_{Z,\alpha_t}(\ZQssub,\alpha^\star) \ddot{L}^{-1}_{\alpha_t,\alpha_t}(\ZQssub,\alpha^\star)\myDonet  \myCt(Z^\star,\alpha^\star) \myDone \times \Delta{\alpha_t}\\
=&~ -\sum_t\ddot{L}_{Z,\alpha_t}(\ZQssub,\alpha^\star)  \times \Delta{\alpha_t}\\
=&~ -\sum_t  \DZt(\ZQssub)\myCt(Z^\star,\alpha^\star) \myDone \times \Delta{\alpha_t}\\
=&~ \DZt(\ZQssub) \sum_t \Delta L_{1t} = S_{311},
\end{align*}
where in the third and fourth lines, we use the 
relationships
\eqref{eq:ddotL_alpha_t,alpha_t} and \eqref{eq:ddotL_Z,alpha_t} 
from Section \ref{subsec:formula_efficient}, respectively.

Thus we have $S_{331} = S_{311} $, and hence $S_{31} - S_{33} = -S_{312} + S_{332} $. 
Similar to the case of $S_{312}$ we elementwisely bound $S_{332} $.
By applying Lemma \ref{lem:thm1.1_Hillar} we get
\begin{align}
\label{eq:myAt_star_inv_inftynorm}
    \|\myAt^{-1}(Z^\star,\alpha^\star)\|_\infty \leqslant 2 e^{M_{\Theta,2}}/n. 
\end{align} Thus,
\begin{align*}
\|S_{332}\|_\infty &\leqslant
\sum_t \left\|\ddot{L}_{Z,\alpha_t}(\ZQssub,\alpha^\star) \myAt^{-1}(Z^\star,\alpha^\star)  \myDonet  \Delta L_{2t} \right\|_\infty\\ 
&\leqslant \sum_t \|\DZt(\ZQssub)\|_\infty \|\myCt(Z^\star,\alpha^\star)\|_\infty \|\myDone\|_\infty \|\myAt^{-1}(Z^\star,\alpha^\star) \|_\infty \|\myDonet \|_\infty \|\Delta L_{2t} \|_\infty  \\
&\leqslant T\times 2 n\sqrt{M_{Z,1}} \times 2e^{-M_{\Theta,1}} \times 2 \times 2 e^{M_{\Theta,2}}/n \times 2 n \times2e^{-M_{\Theta,1}} (M_b\log^{\epsilon}(nT)/\sqrt{n})^2\\
&\leqslant C M_{Z,1}^{\frac12} e^{M_{\Theta,2}-2M_{\Theta,1}} M_b^2 T\log^{2\epsilon}(nT).
\end{align*}

\paragraph*{Step 2}
Next we consider $S_{32} $. Since $W(Z,\alpha) = \ddot{L}_{Z,\alpha} \ddot{L}_{\alpha,\alpha}^{-1} $ and $\ddot{L}_{Z,\alpha_t} = \DZt \ddot{L}_{\Theta_t,\alpha_t} $, we write
\begin{align*}
S_{32} = \big[W(\ZQssub,\check\alpha) - W(\ZQssub,\alpha^\star) \big]\dot{L}_\alpha(\ZQssub,\check\alpha)
= \DZt(\ZQssub) \sum_t 
\Delta W_{t1}
\dot{L}_{\alpha_t}(Z^\star,\check\alpha),
\end{align*}
where
$\Delta W_{t1} = \myBt^\top(Z^\star,\check\alpha) \myAt^{-1}(Z^\star,\check\alpha) - \myBt^\top(Z^\star,\alpha^\star) \myAt^{-1}(Z^\star,\alpha^\star) $. 

\medskip
\noindent \textit{Step 2-1:} We decompose $\Delta W_{t1} = \Delta W_{t11} + \Delta W_{t12} $ where
\begin{align*}
\Delta W_{t11} =&~ \{\myBt^\top(Z^\star,\check\alpha)-\myBt^\top(Z^\star,\alpha^\star)\} \myAt^{-1}(Z^\star,\check\alpha),\\
\Delta W_{t12} =&~\myBt^\top(Z^\star,\alpha^\star)\{ \myAt^{-1}(Z^\star,\check\alpha)- \myAt^{-1}(Z^\star,\alpha^\star)\}.
\end{align*}
Again by applying $\ddot{L}_{\alpha_t,\Theta_t} = \myDonet \ddot{L}_{\Theta_t,\Theta_t} $
and Condition \ref{cond_elem_init}, we get 
\begin{align*}
\|\Delta W_{t11} \|_\infty \leqslant&~ \|\myCt(Z^\star,\check\alpha) - \myCt(Z^\star,\alpha^\star)\|_\infty \|\myDone\|_\infty \| \myAt^{-1}(Z^\star,\check\alpha) \|_\infty \\
\leqslant&~ 2\max_{i,j} | \mu_{t,ij}(Z^\star,\check\alpha) - \mu_{t,ij}(Z^\star,\alpha^\star) |\times2 \times 2e^{M_{\Theta,2}}/n\\
\leqslant&~ Ce^{-M_{\Theta,1}}M_{Z,1}^{\frac12} M_b\log^{\epsilon}(nT)/\sqrt n \times 2 \times 2e^{M_{\Theta,2}}/n \\
=&~ C e^{M_{\Theta,2}-M_{\Theta,1}}M_{Z,1}^{\frac12} M_b n^{-\frac32}\log^{\epsilon}(nT) 
\end{align*}
where in the second inequality we use
\begin{align}
    \|\myAt^{-1}(Z^\star,\check\alpha) \|_\infty \leqslant&~ 2\exp\{\max_{i,j,t}( \langle z_i^\star, z_j^\star \rangle + \check\alpha_{it}+ \check\alpha_{jt}) \}/n \notag \\
    \leqslant&~ 2\exp(M_{\Theta,2}) \exp(2M_b n^{-1/2} \log^{\epsilon}(nT) ) / n \leqslant 4 \exp(M_{\Theta,2})/n \label{eq:myAt_check_inv_inftynorm}
\end{align}
by Lemma \ref{lem:thm1.1_Hillar}, and in the third inequality we use \eqref{eq:checkM_bd} to bound $\mu_{t,ij} $ at mean values.

For $\Delta W_{t12}$, we apply a special case of Woodbury matrix identity 
$$(E-F)^{-1} = \sum_{k=0}^\infty (E^{-1}F)^k E^{-1}$$
with $E = \myAt(Z^\star,\alpha^\star) , F = \myAt(Z^\star,\alpha^\star) - \myAt(Z^\star,\check\alpha)  $.
 Note that 
 \begin{align*}
\|F\|_\infty &= \|\myDonet \{\myCt(Z^\star,\alpha^\star) -\myCt(Z^\star,\check\alpha)\} \myDone\|_\infty \\
&\leqslant 2\|\myDonet\|_\infty \max_{i,j,t}| \mu_{t,ij}(Z^\star,\alpha^\star) - \mu_{t,ij}(Z^\star,\check\alpha)|\|\myDone\|_\infty \\ &\leqslant C e^{-M_{\Theta,1}}M_b\sqrt{n}\log^{\epsilon}(nT) 
\end{align*}
 by applying mean value theorem similarly to \eqref{eq:checkM_bd}.
Besides, $E^{-1} = \myAt^{-1}(Z^\star,\alpha^\star) $ has matrix $l_\infty$ norm bounded by $2e^{M_{\Theta,2}}/n $ as shown previously.
Thus,
\begin{align*}
&\quad\ \|\myAt(Z^\star,\alpha^\star)^{-1} - \myAt(Z^\star,\check\alpha)^{-1} \|_\infty 
= \| E^{-1} - (E-F)^{-1}\|_\infty
\\
&=\bigg\| \sum_{k=1}^\infty [E^{-1}F]^k E^{-1}\bigg\|_\infty  
\leqslant \sum_{k=1}^\infty (\|E^{-1}\|_\infty \|F\|_\infty)^k \|E^{-1}\|_\infty\\ 
&= \sum_{k=1}^\infty (Ce^{M_{\Theta,2}-M_{\Theta,1}}M_{Z,1}^{\frac12}M_b\log^{\epsilon}(nT)/\sqrt{n} )^k (2e^{M_{\Theta,2}}/n)\\
&\leqslant \frac{ C e^{2M_{\Theta,2}-M_{\Theta,1}}M_{Z,1}^{\frac12}M_bn^{-\frac32}\log^{\epsilon}(nT) }{1 - Ce^{M_{\Theta,2}-M_{\Theta,1}}M_{Z,1}^{\frac12}M_b\log^{\epsilon}(nT)/\sqrt{n}  } \\
&\leqslant C e^{2M_{\Theta,2}-M_{\Theta,1}}M_{Z,1}^{\frac12}M_bn^{-\frac32}\log^{\epsilon}(nT)
\end{align*}
as long as $n/\log^{2\varsigma}(nT)$  (and hence $n/\log^{2\epsilon}(nT)$) is sufficiently large. 
Consequently, we obtain the matrix $l_\infty $ norm bound for $\Delta W_{t12} $: 
\begin{align*}
\|\Delta W_{t12}\|_\infty  \leqslant&~ \|\myCt(Z^{\star},\alpha^{\star})\|_\infty \|\myDone\|_\infty \|\myAt(Z^\star,\alpha^\star)^{-1} - \myAt(Z^\star,\check\alpha)^{-1} \|_\infty \\ \leqslant&~ C e^{-M_{\Theta,1}} \times 2 \times e^{2M_{\Theta,2}-M_{\Theta,1}}M_{Z,1}^{\frac12}M_bn^{-\frac32}\log^{\epsilon}(nT)
\\ \leqslant&~ C  e^{2M_{\Theta,2}-2M_{\Theta,1}}M_{Z,1}^{\frac12}M_bn^{-\frac32}\log^{\epsilon}(nT).
\end{align*}

\medskip
\noindent \textit{Step 2-2:} 
Next we consider the bound on the vector $l_\infty$ norm for $\dot{L}_{\alpha_t}(Z^\star,\check\alpha) \in \mathbb R^{n} $. Decomposing it into $\dot{L}_{\alpha_t}(Z^\star,\check\alpha) = \dot{L}_{\alpha_t}(Z^\star,\alpha^\star) + (\dot{L}_{\alpha_t}(Z^\star,\check\alpha) - \dot{L}_{\alpha_t}(Z^\star,\alpha^\star) )  $, we can respectively bound the vector $l_\infty$ norm of the two parts: for the first part,
\begin{align*}
\| \dot{L}_{\alpha_t}(Z^\star,\alpha^\star) \|_\infty & \leqslant \max_{i,t} \Big| \sum_{j=1}^n (A_{t,ij} -\mu_{t,ij}(Z^\star,\alpha^\star))(1+1_{(j=i)})\Big|. 
\end{align*}
Noting that $A_{t,ij} - \mu_{t,ij}(Z^\star,\alpha^\star) $ are independent centered Poisson random variables,
by Lemma \ref{lem:bern_pois}, we have for any $x\geqslant 0$,
\begin{align} \label{eq:bern_dotL}
    \Pr\bigg( \Big| \sum_{j=1}^n (A_{t,ij} -\mu_{t,ij}(Z^\star,\alpha^\star))(1+1_{(j=i)})\Big| \geqslant x \bigg) \leqslant \exp\left( -\frac{x^2}{(n+3)e^{-M_{\Theta,1}} + 10 x } \right).
\end{align}
Then, by a union bound, we have 
\begin{align} \label{eq:unionbern_dotL}
  \| \dot{L}_{\alpha_t}(Z^\star,\alpha^\star) \|_\infty \leqslant 2(s+1) \sqrt{(n+3)e^{-M_{\Theta,1}} \log(nT)}   
\end{align}
with probability $1 - O((nT)^{-s} ) $.
For the second part,
\begin{align*}
    \|\dot{L}_{\alpha_t}(Z^\star,\check\alpha) - \dot{L}_{\alpha_t}(Z^\star,\alpha^\star) \|_\infty &\leqslant 2\max_{i,t} \sum_j |\mu_{t,ij}(Z^\star,\check\alpha) - \mu_{t,ij}(Z^\star,\alpha^\star) | \\&\leqslant Ce^{-M_{\Theta,1}}M_b \sqrt n \log^{\epsilon}(nT) 
\end{align*} by a mean value thoerem as in \eqref{eq:checkM_bd}.
Combining the above two bounds, we have
 \begin{equation}
 \label{eq:L_alphat_Zstar_alphacheck_inftynorm}
     \|\dot{L}_{\alpha_t}(Z^\star,\check\alpha)  \|_\infty \leqslant C(s+1) e^{-0.5M_{\Theta,1}} \sqrt{n}\max(\sqrt{\log(nT) },M_b \log^{\epsilon}(nT) ).
 \end{equation}

Concluding the results in Step 2, we have 
\begin{align*}
\|S_{32}\|_\infty \leqslant&~ \|\DZt(\ZQssub)\|_\infty \sum_t \|\Delta W_{t1}\|_\infty \|\dot{L}_{\alpha_t}(Z^\star,\check\alpha)\|_\infty \\
\leqslant&~ C n  {M_{Z,1}^{\frac12}} T  e^{2M_{\Theta,2}-2.5M_{\Theta,1}}M_b   n^{-\frac32}\log^{\epsilon}(nT)   (s+1) \sqrt{n}\max(\sqrt{\log(nT) },M_b \log^{\epsilon}(nT) ) \\
\leqslant&~ C (s+1)e^{2M_{\Theta,2}-2.5 M_{\Theta,1}}M_b^2 M_{Z,1} T   \log^{\epsilon+\max(\frac12,\epsilon)}(nT) 
\end{align*}
with probability $1-O((nT)^{-s}) $,
where $C$ is a  constant.

\medskip

Summing up the bounds in Steps 1 and 2, we get with probability $1-O((nT)^{-s}) $,
\begin{align*}
\|S_3\|_\infty \leqslant&~ \|S_{312}\|_\infty + \|S_{332}\|_\infty + \|S_{32}\|_\infty \\
\leqslant&~ C M_{Z,1}^{\frac12} M_b^2 e^{-M_{\Theta,1}} T \log^{2\epsilon}(nT) + C M_{Z,1}^{\frac12} M_b^2 e^{M_{\Theta,2}-2M_{\Theta,1}} T\log^{2\epsilon}(nT) \\&~+ C (s+1)e^{2M_{\Theta,2}-2.5 M_{\Theta,1}}M_b^2 M_{Z,1} T   \log^{\epsilon+\max(\frac12,\epsilon)}(nT)\\
\leqslant&~ C_{M,s} T \log^{\epsilon+\max(\frac12,\epsilon)}(nT).
\end{align*}
where $C_{M,s}$ is a constant that only depends  on $M_{\Theta,1}, M_{\Theta,2}, M_{Z,1}$, $M_b$ and $s$. 
Thus we have proved \eqref{eq:S3_infty_main}: For any $s>0$, there exists a universal constant $C_s>0$ such that
$$
\sqrt{nk}\|S_3\|_\infty /(nT)
\leqslant  C_{s} \,  \big(n^{-1/2}\big)\, \log^{2\varsigma }(nT)
$$
with probability $1-O((nT)^{-s}) $.

\subsubsection{Proof of Eq. \eqref{eq:S1+S2_infty_main}} \label{sec:pfB3}
Similarly to the analysis of $S_{3} $, we decompose $S_{2} $ into
\begin{align*}
S_2 =&~ \dot{L}_Z(\check Z,\check\alpha) - W(\check Z,\check\alpha) \dot{L}_\alpha(\check Z,\check\alpha) - \dot{L}_Z(\ZQssub,\check\alpha)+ W(\ZQssub,\check\alpha) \dot{L}_\alpha(\ZQssub,\check\alpha)\\
=&~ S_{21} - S_{22} - S_{23},
\end{align*}
where we define
\begin{align*}
S_{21} =&~ \dot{L}_Z(\check Z,\check\alpha)  - \dot{L}_Z(\ZQssub,\check\alpha),\\
S_{22} =&~ \big[ W(\check Z,\check\alpha) -W(\ZQssub,\check\alpha) \big]\dot{L}_\alpha(\ZQssub,\check\alpha),\\
S_{23} =&~ W(\check Z,\check \alpha)\big[\dot{L}_\alpha(\check Z,\check\alpha) - \dot{L}_\alpha(\ZQssub,\check\alpha) \big].
\end{align*}

\paragraph*{Proof Outline}
By the above decomposition, we have  $S_1+S_2 =(S_1+S_{21} -S_{23}) -S_{22} $.
We derive an upper bound for $S_1+S_2$ by  two steps. 

\medskip
\noindent \textit{Step 1} derives an upper bound of $S_{22}$,  similarly to the $S_{32}$ term in the proof of \eqref{eq:S3_infty_main}. 
In particular, we first derive an upper bound for $W(\check Z,\check\alpha) -W(\ZQssub,\check\alpha)$ and then use the upper bound of $\dot{L}_\alpha(\ZQssub,\check\alpha)$ derived in the proof of \eqref{eq:S3_infty_main}. 

\medskip
\noindent \textit{Step 2} 
derives an upper bound of $S_1+S_{21} -S_{23}$ similarly to $S_{31}-S_{33}$ in the proof of \eqref{eq:S3_infty_main}. Intuitively,  through Taylor expansion with respect to $Z_v$, we have
\begin{align}
S_{21}=&~ \ddot{L}_{Z,Z}(\check Z,\check\alpha)(\check Z_v-(\ZQssub)_v )  + \text{higher-order term}
\label{eq:S_21_heuristic} \\
S_{23} =&~  \ddot{L}_{Z,\alpha}(\check Z,\check \alpha) \ddot{L}_{\alpha,\alpha}^{-1}(\check Z,\check \alpha) \ddot{L}_{\alpha,Z}(\check Z,\check \alpha) (\check Z_v-(\ZQssub)_v ) + \text{higher-order term},
\label{eq:S_23_heuristic}
\end{align}
where we plug in $W(\check{Z}, \check{\alpha})=\ddot{L}_{Z,\alpha}(\check Z,\check \alpha) \ddot{L}_{\alpha,\alpha}^{-1}(\check Z,\check \alpha)$. 
On the other hand, by \eqref{eq:HeffIeff}, we have
\begin{align*}
S_1  =&~ \Ieff(\check Z,\check \alpha) (\check Z_v-(\ZQssub)_v )  \\ 
 =&~ \left[- \big( \ddot{L}_{Z,Z}(\check Z,\check\alpha) - \ddot{L}_{Z,\alpha}(\check Z,\check \alpha) \ddot{L}_{\alpha,\alpha}^{-1}(\check Z,\check \alpha) \ddot{L}_{\alpha,Z}(\check Z,\check \alpha) \big) + \text{perturbation term} \right](\check Z_v-(\ZQssub)_v ) \\
 =&~ -\left[\ddot{L}_{Z,Z}(\check Z,\check\alpha) - \ddot{L}_{Z,\alpha}(\check Z,\check \alpha) \ddot{L}_{\alpha,\alpha}^{-1}(\check Z,\check \alpha) \ddot{L}_{\alpha,Z}(\check Z,\check \alpha)  \right](\check Z_v-(\ZQssub)_v ) + \text{perturbation term}.
\end{align*}
Hence $ S_1 + S_{21} - S_{23} = 0 + \text{higher-order terms}+ \text{perturbation term}$.
In the following, we give detailed calculations for the above heuristic argument.

\paragraph*{Step 1}
We first consider $S_{22}$. 
By the Fact \ref{fact:partialderivfunction},
\begin{align*}
S_{22} &= \big[W(\check Z,\check\alpha) - W(\ZQssub,\check\alpha) \big]\dot{L}_\alpha(\ZQssub,\check\alpha)\\
&=\sum_t \big[  \ddot{L}_{Z,\alpha_t}(\check Z,\check\alpha) \ddot{L}^{-1}_{\alpha_t,\alpha_t}(\check Z,\check\alpha) - \ddot{L}_{Z,\alpha_t}(\ZQssub,\check\alpha) \ddot{L}^{-1}_{\alpha_t,\alpha_t}(\ZQssub,\check\alpha)  \big]\dot{L}_{\alpha_t}(\ZQssub,\check\alpha)
\\
&=\sum_t \big[ \DZt(\check Z) \myCt(\check Z,\check\alpha) \myDone \myAt^{-1}(\check Z,\check\alpha) \\
&\quad\quad\quad- \DZt(\ZQssub) \myCt(Z^\star,\check\alpha)\myDone \myAt^{-1}(Z^\star,\check\alpha) \big]\dot{L}_{\alpha_t}(Z^\star,\check\alpha)\\
&=: \sum_t [\Delta W_{t2} + \DZt(\check Z)\Delta W_{t3}] \dot{L}_{\alpha_t}(Z^\star,\check\alpha)
\end{align*}
where
\begin{align*}
\Delta W_{t2} =&~ \{\DZt(\check Z) - \DZt(\ZQssub)\}\myCt(Z^\star,\check\alpha)\myDone \myAt^{-1}(Z^\star,\check\alpha) \\
\text{and}\quad
\Delta W_{t3} =&~  \myCt(\check Z,\check\alpha)\myDone \myAt^{-1}(\check Z,\check\alpha)  -\myCt(Z^\star,\check\alpha) \myDone\myAt^{-1}(Z^\star,\check\alpha) .
\end{align*}
Applying the same technique as in bounding $\Delta W_{t1}$, we have $$\|\Delta W_{t3}\|_\infty \leqslant C e^{2M_{\Theta,2}-2M_{\Theta,1}}M_b n^{-\frac32} \log^{\epsilon}(nT) .$$ 
Similarly we could bound $\Delta W_{t2} $ by 
\begin{align*}
\|\Delta W_{t2}\|_\infty \leqslant&~ \| \DZt(\check Z) - \DZt(\ZQssub) \|_\infty \| \myCt(Z^\star,\check\alpha)\|_\infty \|\myDone\|_\infty \|\myAt^{-1}(Z^\star,\check\alpha)\|_\infty \\ 
\leqslant&~ C M_b e^{M_{\Theta,2}-M_{\Theta,1}} n^{-\frac12}\log^{\epsilon}(nT). 
\end{align*} 
Besides, we have {already derived} in \eqref{eq:L_alphat_Zstar_alphacheck_inftynorm} that $$\|\dot{L}_{\alpha_t}(\ZQssub,\check\alpha)\|_\infty \leqslant C(s+1) \sqrt{ne^{-M_{\Theta,1}}}M_b \log^{\max(\epsilon,\frac12)}(nT)$$ with probability $1-O((nT)^{-s}) $ in the proof for $S_{32}$. Combining those results we have 
\begin{align*}
\|S_{22}\|_\infty \leqslant&~ \sum_t (\|\Delta W_{t2}\|_\infty + \|\DZt(\check Z)\|_\infty \|\Delta W_{t3}\|_\infty)\|\dot{L}_{\alpha_t}(Z^\star,\check\alpha)\|_\infty 
\\
\leqslant&~ C (s+1)e^{2M_{\Theta,2}-2.5 M_{\Theta,1}}M_b (M_{Z,1}^{\frac12}+1) T \log^{\epsilon+\max(\epsilon,\frac12)}(nT)
\end{align*}
with probability $1-O((nT)^{-s}) $ for some constant $C$.

\paragraph*{Step 2} In step 2 we bound $S_1+S_{21} -S_{23}$.
First we consider $S_{21}$.
\begin{align*}
S_{21} &= \dot{L}_Z(\check Z,\check\alpha)  - \dot{L}_Z(\ZQssub,\check\alpha)\\
&=\sum_t \big[ \DZt(\check Z)\myUt(\check Z,\check\alpha) -  \DZt(\ZQssub)\myUt(\ZQssub,\check\alpha)\big]\\
&= \DZt(\check Z)\sum_t\big[\myUt(\check Z,\check \alpha)- \myUt(\ZQssub,\check \alpha) \big] + (\DZt(\check Z) - \DZt(\ZQssub) )\sum_t\myUt(\ZQssub,\check\alpha).
\end{align*}
The $(i,j) $th entry of $\big[\myUt(\check Z,\check \alpha)- \myUt(\ZQssub,\check \alpha)\big]$ is 
$$\varrho_{ij} \{\exp( (z_i^\star)^\top z_j^\star + \check\alpha_{it} + \check\alpha_{jt})-\exp(\check z_i^\top \check z_j + \check\alpha_{it} + \check\alpha_{jt})\},$$ on which we apply Taylor expansion similarly to \eqref{eq:mu_Taylorexp}:
\begin{equation}
\begin{aligned}
 &\exp( (z_i^\star)^\top z_j^\star + \check\alpha_{it} + \check\alpha_{jt})-\exp(\check z_i^\top \check z_j + \check\alpha_{it} + \check\alpha_{jt})\\ 
 =&~ \frac{\partial \exp(\check z_i^\top \check z_j + \check\alpha_{it} + \check\alpha_{jt}) }{\partial z_i^\top } (\check{Q}^{\mytrans} z_i^\star- \check z_i) 
 + \frac{\partial \exp(\check z_i^\top \check z_j + \check\alpha_{it} + \check\alpha_{jt}) }{\partial z_j^\top } (\check{Q}^{\mytrans} z_j^\star- \check z_j) \\
 &+ \frac{1}{2}(\check{Q}^{\mytrans} z_i^\star- \check z_i)^\top \frac{\partial^2 \exp(\tilde z_i^{(ij)\top} \tilde z_j^{(ij)} + \check\alpha_{it} + \check\alpha_{jt}) }{\partial z_i \partial z_i^\mytrans } (\check{Q}^{\mytrans} z_i^\star- \check z_i)
 \\&+ \frac{1}{2}(\check{Q}^{\mytrans} z_j^\star- \check z_j)^\top \frac{\partial^2 \exp(\tilde z_i^{(ij)\top} \tilde z_j^{(ij)} + \check\alpha_{it} + \check\alpha_{jt}) }{\partial z_j \partial z_j^\mytrans } (\check{Q}^{\mytrans} z_j^\star- \check z_j)\\
 &+ (\check{Q}^{\mytrans} z_i^\star- \check z_i)^\top \frac{\partial^2 \exp(\tilde z_i^{(ij)\top} \tilde z_j^{(ij)} + \check\alpha_{it} + \check\alpha_{jt}) }{\partial z_i \partial z_j^\mytrans} (\check{Q}^{\mytrans} z_j^\star- \check z_j)
\end{aligned}
\label{eq:mu_Taylorexp_wrt_z}
\end{equation}
where $\tilde z_i^{(ij)}$ is a midpoint between $\check{Q}^{\mytrans} z_i^\star$ and $\check z_i$. (We use the subscript $^{(ij)} $ because for a pair of $(i,j)$ the midpoint between $\check{Q}^{\mytrans} z_i^\star$ and $\check z_i$ in the expansion of $\mu_{t,ij}(\check Z,\check\alpha) - \mu_{t,ij}(\ZQssub,\check\alpha) $ depends on both $z_i$ and $z_j$.) 
Similar to the definition of $\Delta L_{1t} $ and $\Delta L_{2t} $ we let $\myUt(\check Z,\check \alpha)- \myUt(\ZQssub,\check \alpha) =\Delta L_{3t} + \Delta L_{4t}$, where $\Delta L_{3t}$ is an $n(n+1)/2 $ dimensional vector whose $(i,j)$th element is
\begin{align*}
&\bigg\{\frac{\partial \exp(\check z_i^\top \check z_j + \check\alpha_{it} + \check\alpha_{jt}) }{\partial z_i^\top } (\check{Q}^{\mytrans} z_i^\star- \check z_i) 
 + \frac{\partial \exp(\check z_i^\top \check z_j + \check\alpha_{it} + \check\alpha_{jt}) }{\partial z_j^\top } (\check{Q}^{\mytrans} z_j^\star- \check z_j)\bigg\}\varrho_{ij} \\
 =&~\left\{
\exp(\check z_i^\top \check z_j + \check\alpha_{it} + \check\alpha_{jt})\check z_j^\top (\check{Q}^{\mytrans} z_i^\star- \check z_i) +  \exp(\check z_i^\top \check z_j + \check\alpha_{it} + \check\alpha_{jt})\check z_i^\top (\check{Q}^{\mytrans} z_j^\star- \check z_j)\right\}\varrho_{ij},
\end{align*} 
and $\Delta L_{4t}$ is an $n(n+1)/2$ dimensional vector whose $(i,j)$th element is given by  
\begin{align*}
 &\bigg\{\frac{1}{2}(\check{Q}^{\mytrans} z_i^\star- \check z_i)^\top \frac{\partial^2 \exp(\tilde z_i^{(ij)\top} \tilde z_j^{(ij)} + \check\alpha_{it} + \check\alpha_{jt}) }{\partial z_i \partial z_i^\mytrans } (\check{Q}^{\mytrans} z_i^\star- \check z_i) \bigg.
 \\&\bigg. +~\frac{1}{2} (\check{Q}^{\mytrans} z_j^\star- \check z_j)^\top \frac{\partial^2 \exp(\tilde z_i^{(ij)\top} \tilde z_j^{(ij)} + \check\alpha_{it} + \check\alpha_{jt}) }{\partial z_j \partial z_j^\mytrans } (\check{Q}^{\mytrans} z_j^\star- \check z_j)\bigg.\\
 &\bigg.+~ (\check{Q}^{\mytrans} z_i^\star- \check z_i)^\top \frac{\partial^2 \exp(\tilde z_i^{(ij)\top} \tilde z_j^{(ij)} + \check\alpha_{it} + \check\alpha_{jt}) }{\partial z_i \partial z_j^\mytrans } (\check{Q}^{\mytrans} z_j^\star- \check z_j)\bigg\}\varrho_{ij} \\
 =&~ \left\{\frac{1}{2}(\check{Q}^{\mytrans} z_i^\star- \check z_i)^\top\tilde z_j^{(ij)}\exp(\tilde z_i^{(ij)\top} \tilde z_j^{(ij)} + \check\alpha_{it} + \check\alpha_{jt})\tilde z_j^{(ij)\top} (\check{Q}^{\mytrans} z_i^\star- \check z_i)\right. 
 \\&\left.+ ~\frac{1}{2}(\check{Q}^{\mytrans} z_j^\star- \check z_j)^\top\tilde z_i^{(ij)}\exp(\tilde z_i^{(ij)\top} \tilde z_j^{(ij)} + \check\alpha_{it} + \check\alpha_{jt})\tilde z_i^{(ij)\top} (\check{Q}^{\mytrans} z_j^\star- \check z_j) \right.
 \\
 &\left.+ ~(\check{Q}^{\mytrans} z_j^\star- \check z_j)^\top\tilde z_i^{(ij)}\exp(\tilde z_i^{(ij)\top} \tilde z_j^{(ij)} + \check\alpha_{it} + \check\alpha_{jt})\tilde z_j^{(ij)\top} (\check{Q}^{\mytrans} z_i^\star- \check z_i) \right.
 \\
 &\left.+ ~(\check{Q}^{\mytrans} z_j^\star- \check z_j)^\top\exp(\tilde z_i^{(ij)\top} \tilde z_j^{(ij)} + \check\alpha_{it} + \check\alpha_{jt}) (\check{Q}^{\mytrans} z_i^\star- \check z_i)\right\}\varrho_{ij}.
\end{align*}

Let $S_{21} =S_{211} +S_{212} + S_{213} $ with
\begin{equation}
\label{eq::def_S21_123}
    \begin{aligned}
S_{211} :=&~ \DZt(\check Z)\sum_t\Delta L_{3t},\\
S_{212} :=&~ \DZt(\check Z)\sum_t \Delta L_{4t},\\
S_{213} :=&~ (\DZt(\check Z) - \DZt(\ZQssub) ) \sum_t\myUt(\ZQssub,\check\alpha).
\end{aligned}
\end{equation}
Similar to the case of $\Delta L_{1t} $, 
 $\Delta L_{3t}$ can be equivalently written in matrix form as
\begin{align*}
\Delta L_{3t}
 = \myCt(\check Z,\check\alpha) \DZ(\check Z) ({\check Z}_v - {(\ZQssub)}_v).
\end{align*}
Then 
$$S_{211} = \sum_t\big[ \DZt(\check Z) \ddot{L}_{\Theta_t,\Theta_t}(\check Z,\check\alpha) \DZ(\check Z)  ({\check Z}_v - {(\ZQssub)}_v) \big]. 
$$
In the following we first bound $S_{212} $ and $S_{213} $, and then give rigorous arguments on how $S_{211}$ cancel with terms in $S_{23} $ and $S_1$.

Using the same technique as bounding $S_{312} $, and from Condition \ref{cond_elem_init} on the error of $\check Z$, $S_{212} \in \mathbb R^{nk}$ is elementwisely bounded by $
C e^{-M_{\Theta,1}}M_b^2(M_{Z,1}+1)^{\frac32}T \log^{2\epsilon}(nT) $ from the bound $$\|\Delta L_{4t}\|_\infty \leqslant Ce^{-M_{\Theta,1}} (M_{Z,1}+1 )(M_b\log^{\epsilon}(nT) /\sqrt{n})^2 .$$ For $S_{213} $, using mean value theorem similarly to \eqref{eq:checkM_bd},
we first have 
\begin{equation}
\label{eq::DZcheck-DZstar}
    \|\DZt(\check Z) - \DZt(\ZQssub) \|_\infty \leqslant C M_b\log^{\epsilon}(nT)  \sqrt{n} .
\end{equation}
Then for  $\sum_t\myUt(\ZQssub,\check\alpha)=\sum_t\myUt(Z^\star,\check\alpha) $, we decompose it into $$\sum_t\myUt(Z^\star,\alpha^\star) + \sum_t\big[\myUt(Z^\star,\check\alpha)- \myUt(Z^\star,\alpha^\star)\big].$$ 
Similarly applying a Bernstein bound as in \eqref{eq:bern_dotL} and then a union bound as in \eqref{eq:unionbern_dotL}, we have that for the first part $$\|\sum_t\myUt(Z^\star,\alpha^\star)\|_\infty \leqslant 4(s+1)\sqrt{Te^{-M_{\Theta,1}}\log(nT)} $$ with probability $1-O((nT)^{-s} ) $. In the second part, we have 
\begin{align*}
    \|\sum_t\myUt(Z^\star,\check\alpha)-\myUt(Z^\star,\alpha^\star)\|_\infty &\leqslant C \sum_t\max_{t,i,j} | \mu_{t,ij}(Z^\star,\check\alpha) -\mu_{t,ij}(Z^\star,\alpha^\star) | \\&\leqslant CTe^{-M_{\Theta,1}}M_b\log^{\epsilon}(nT) /\sqrt n
\end{align*} 
where we again apply the mean value theorem as in \eqref{eq:checkM_bd}.
Summing up, we have 
$$\|S_{213}\|_\infty \leqslant 
C\big[(s+1)e^{-0.5M_{\Theta,1}} M_b \sqrt{nT}\log^{\epsilon+\frac12}(nT) + e^{-M_{\Theta,1}} M_b^2 T\log^{2\epsilon}(nT) \big] $$ with probability $1-O((nT)^{-s} ) $.

We next decompose $S_{23}$ following the heuristic argument of \eqref{eq:S_23_heuristic}, and
show that $S_{211}$ 
together with the first order part in $S_{23} $
could be canceled with $S_{1}$. 
\begin{align*}
S_{23} =&~ W(\check Z,\check \alpha)\big[\dot{L}_\alpha(\check Z,\check\alpha) - \dot{L}_\alpha(\ZQssub,\check\alpha) \big]\\
=&~ \sum_t  \ddot{L}_{Z,\alpha_t}(\check Z,\check\alpha) \ddot{L}^{-1}_{\alpha_t,\alpha_t}(\check Z,\check\alpha)   \big[\dot{L}_{\alpha_t}(\check Z,\check\alpha)-\dot{L}_{\alpha_t}(\ZQssub,\check\alpha) \big] \\
=&~ \sum_t  \ddot{L}_{Z,\alpha_t}(\check Z,\check\alpha) \myAt^{-1}(\check Z,\check\alpha) \myDonet \big[ \myUt(\check Z,\check\alpha) - \myUt(\ZQssub,\check\alpha) \big]. 
\end{align*}
Using $\myUt(\check Z,\check \alpha)- \myUt(\ZQssub,\check \alpha) =\Delta L_{3t} + \Delta L_{4t}$ again, we
let $S_{23} = S_{231} + S_{232} $, where
\begin{align*}
S_{231} :=&~ \sum_t  \ddot{L}_{Z,\alpha_t}(\check Z,\check\alpha) \myAt^{-1}(\check Z,\check\alpha)\myDonet \Delta L_{3t}, \\
S_{232} :=&~ \sum_t  \ddot{L}_{Z,\alpha_t}(\check Z,\check\alpha) \myAt^{-1}(\check Z,\check\alpha)\myDonet \Delta L_{4t}.
\end{align*}
Since $$\myDonet\Delta L_{3t} = \myDonet\myCt(\check Z,\check\alpha) \DZ(\check Z) (\check Z_v-(\ZQssub)_v ) = \ddot{L}_{\alpha_t,Z}(\check Z,\check \alpha)(\check Z_v-(\ZQssub)_v ),$$
we could see that
$S_{231}$ corresponds to the first-order term $\ddot{L}_{Z,\alpha}(\check Z,\check \alpha) \ddot{L}_{\alpha,\alpha}^{-1}(\check Z,\check \alpha) \ddot{L}_{\alpha,Z}(\check Z,\check \alpha) (\check Z_v-(\ZQssub)_v )$ in \eqref{eq:S_23_heuristic}.
By our heuristic argument in \eqref{eq:S_21_heuristic} and \eqref{eq:S_23_heuristic} we should have $S_{211} - S_{231}$ canceling with $S_1$. Indeed,
\begin{align*}
S_{211} - S_{231} &=  \DZt(\check Z) \sum_t \big\{ \myCt(\check Z,\check\alpha) \DZ(\check Z) ({\check Z}_v - {(\ZQs)}_v)  \big\}
\\ & \quad- \sum_t \ddot{L}_{Z,\alpha_t}(\check Z,\check\alpha) \myAt^{-1}(\check Z,\check\alpha)\myDonet \ddot{L}_{\alpha_t,Z}(\check Z,\check\alpha) ({\check Z}_v - {(\ZQs)}_v) \\
&=\Big\{-\sum_t \DZt(\check Z) \mathcal{D}_{\mu_t}(\check Z,\check\alpha) \DZ(\check Z)  \\
&\quad\quad- \ddot{L}_{Z,\alpha}(\check Z,\check \alpha) \ddot{L}^{-1}_{\alpha,\alpha}(\check Z,\check \alpha)\ddot{L}_{\alpha,Z}(\check Z,\check \alpha)\Big\} ({\check Z}_v - {(\ZQs)}_v)
\\&=-
\Ieff(\check Z,\check\alpha)({\check Z}_v - {(\ZQs)}_v) = -S_1.
\end{align*}
Hence, $S_1 + S_{211} - S_{231} 
= 0$. 

Similar to the case of $S_{332}$ and $S_{212}$, we could  elementwisely bound $S_{232}$ by $$\|S_{232}\|_{\infty} \leqslant C e^{M_{\Theta,2}-2M_{\Theta,1}} M_b^2(M_{Z,1}+1)^{\frac32} T\log^{2\epsilon}(nT). $$

\medskip

To sum up the bounds in Steps 1 and 2, we can bound $S_1+S_2$ elementwisely with probability $1-O((nT)^{-s}) $ by
\begin{equation}
    \label{eq::finalbd_S1+S2_infty}
    \begin{aligned}
\|S_1+S_2\|_\infty \leqslant &~ \|S_{22}\|_\infty + \|S_{212}\|_\infty + \|S_{213}\|_\infty + \|S_{232}\|_\infty + \|S_1 + S_{211} - S_{231}\|_\infty \\
\leqslant&~ C (s+1)e^{2M_{\Theta,2}-2.5 M_{\Theta,1}}M_b (M_{Z,1}^{\frac12}+1) T \log^{\epsilon+\max(\epsilon,\frac12)}(nT) 
\\&+ C e^{-M_{\Theta,1}}M_b^2(M_{Z,1}+1)^{\frac32}T \log^{2\epsilon}(nT) 
\\&+ C[(s+1)e^{-0.5M_{\Theta,1}} M_b \sqrt{nT}\log^{\epsilon+\frac12}(nT) + e^{-M_{\Theta,1}} M_b^2 T\log^{2\epsilon}(nT) ] 
\\&+ C e^{M_{\Theta,2}-2M_{\Theta,1}} M_b^2(M_{Z,1}+1)^{\frac32} T\log^{2\epsilon}(nT)  
\\
\leqslant&~ C_{M,s} \left(T \log^{\epsilon+\max(\epsilon,\frac12)}(nT) +  \sqrt{nT}\log^{\epsilon+\frac12}(nT) \right)
\end{aligned}
\end{equation}
where $C_{M,s}$ is a constant that only depends  on $M_{\Theta,1},M_{\Theta,2},M_{Z,1}$, $M_b$, and $s$. 
Note that the $\sqrt{nT}$ order term (ignoring the $\log$ factors) comes only from the first part of term $S_{213} $.
Thus we have proved \eqref{eq:S1+S2_infty_main}: For any $s>0$, there exists a universal constant $C_s>0$ such that
$$
\sqrt{nk}\|S_1+S_2\|_\infty/(nT) 
\leqslant  C_{s}  \,  \big(n^{-1/2}+ T^{-1/2}\big)\, \log^{2\varsigma }(nT)
$$
with probability $1-O((nT)^{-s}) $.

\subsubsection{Proof of Eq. \eqref{eq:Seff_starstar_main}} \label{sec:pfB4}

Note that $\|\Seff(\ZQs,\alpha^\star)\|_2 = \|\Seff(Z^\star, \alpha^\star )\|_2$ as
\begin{align*}
\Seff(\ZQs,\alpha^\star) &= \sum_t \left.(\DZt\myXt\myCt^{-1}\myUt)\right|_{(\ZQs,\alpha^\star)} \\
&= [\mathrm{I}_{n}\otimes \check Q^\top]\sum_t \left.(\DZt\myXt\myCt^{-1}\myUt)\right|_{(Z^\star,\alpha^\star)} =[\mathrm{I}_{n}\otimes \check Q^\top] \Seff(Z^\star, \alpha^\star ),
\end{align*}
where 
\begin{align*}
\Seff(Z^\star,\alpha^\star) =&~ \left.(\dot{L}_Z - \ddot{L}_{Z,\alpha} \ddot{L}^{-1}_{\alpha,\alpha}\dot{L}_{\alpha})\right|_{(Z^\star,\alpha^\star)} \\
=&~ \sum_t \left.\left\{\DZt(\mathrm{I}_{{n(n+1)}/{2}} -\ddot{L}_{\Theta_t,\alpha_t} \myAt^{-1} \myDonet)\myUt\right\}\right|_{(Z^\star,\alpha^\star)}.
\end{align*}
Therefore,  $\|\Seff(\ZQs,\alpha^\star)\|_2 /\sqrt{nk}=\|\Seff(Z^\star, \alpha^\star )\|_2/\sqrt{nk}\leqslant \|\Seff(Z^\star, \alpha^\star )\|_{\infty}$. To derive an upper bound for $\|\Seff(\ZQs,\alpha^\star)\|_2 /\sqrt{nk}$, we  derive an upper bound of $\|\Seff(Z^\star, \alpha^\star )\|_{\infty}$ below. 

In the proof of this lemma $(Z,\alpha)$ only take value at $(Z^\star,\alpha^\star)$, so we may omit the argument $(Z^\star,\alpha^\star) $ of the functions.
Recall in Section \ref{subsec:formula_UH}
we denote
$$\myXt = -(\myCt - \ddot{L}_{\Theta_t,\alpha_t} \myAt^{-1} \ddot{L}_{\alpha_t,\Theta_t}) \in \mathbb R^{\frac{n(n+1)}{2} \times \frac{n(n+1)}{2}}$$ so that 
\begin{equation}
\label{eq:Seff_matrixform_D13}
    \Seff =- \DZt \sum_{t=1}^T \myXt\myCt^{-1}\dot{L}_{\Theta_t}.
\end{equation}
We index the rows/columns of $\myXt$ by $(i,j),1\leqslant i\leqslant j \leqslant n$; for notational simplicity,  the $(j,i)$th row also means the $(i,j)$th row for $i\leqslant j$.
Recall $ \dot{L}_{\Theta_{t,ij}}(Z^\star,\alpha^\star)=\varrho_{ij} (A_{t,ij} - \mu^\star_{t,ij}) $. Since it is used repeatedly in the proofs throughout this manuscript, we denote
\begin{equation}
\label{eq::def_Mtij}
    M_{t,ij} := N_{t,ij}(Z^\star,\alpha^\star)= A_{t,ij} - \mu^\star_{t,ij}
\end{equation}
where $N_{t} $ is defined in  \eqref{eq::def_NtN2t}.
Now for $\Seff
\in \mathbb R^{nk}$, we could write out its $(j,k')$th element ($1\leqslant j\leqslant n,1\leqslant k' \leqslant k $) as

\begin{align}(\Seff)_{(j,k')}=&~-\sum_{t=1}^T
\sum_{l=1}^n \varrho_{jl} z_{l k'} \left\{\myXt\myCt^{-1}\myUt \right\}_{(j,l)} \notag \\
=&~\sum_{t=1}^T \sum_{l=1}^n \varrho_{jl} z_{lk'}\Big\{\sum_{1 \leqslant i_1\leqslant i_2 \leqslant n} (\myXt)_{(j,l)(i_{1},i_{2})} 
(\varrho_{i_1i_2}^2 \mu_{t,i_{1}i_{2}})^{-1} 
 (\varrho_{i_1i_2}M_{t,i_1 i_2} ) 
\Big\} \notag \\
=&~
\sum_{t=1}^T \sum_{1 \leqslant i_1\leqslant i_2 \leqslant n}\Big\{\sum_{l=1}^n\varrho_{jl} z_{lk'} (\myXt)_{(j,l)(i_{1},i_{2})} 
(\varrho_{i_1i_2}^2 \mu_{t,i_{1}i_{2}})^{-1} 
 (\varrho_{i_1i_2}M_{t,i_1 i_2} ) 
\Big\} \label{eq:Seff_elem} 
\end{align}
where $\{\myXt\myCt^{-1}\myUt\}_{(j,l)} $ denotes the $(j,l) $th element of the $n(n+1)/2$ dimensional vector $\{\myXt\myCt^{-1}\myUt \}$, and similarly $(\myXt)_{(j,l)(i_{1},i_{2})}$ denotes the entry at the $(j,l)$th row and the $(i_1,i_2)$th column of the ${n(n+1)}/{2} \times {n(n+1)}/{2}$ matrix $\myXt$.
Due to the non-randomness of  $\varrho, z, \mu_{t,ij} $ and $\myXt $, and the independence among $M_{t,ij}$, we change the order of summation to write $(S_{eff})_{(j,k')} $ as   
a sum of independent centered Poisson variables in the last line of  \eqref{eq:Seff_elem}.
Thus, we can apply Bernstein's inequality  in Lemma \ref{lem:bernin} to bound each $(\Seff)_{(j,k')} $ as follows:
\begin{align} \label{eq::c7_Bern}
    \Pr\Big( \left|(\Seff)_{(j,k')}\right|  >  (s+1) \left(\sqrt{V_{(j,k')}} + L\right)\log(nk)\Big)  \leqslant O((nk)^{-(s+1)}) 
\end{align}
for any $s>0$, in which 
$$
V_{(j,k')} = \operatorname{var}\Big\{(\Seff)_{(j,k')}\Big\}
$$
and 
$$
L = \max_{1\leqslant j \leqslant n, 1\leqslant k' \leqslant K, 1\leqslant t\leqslant T, 1\leqslant i_1\leqslant i_2 \leqslant n} \Big\{\sum_{l=1}^n\varrho_{jl} z_{lk'} (\myXt)_{(j,l)(i_{1},i_{2})} 
(\varrho_{i_1i_2}^2 \mu_{t,i_{1}i_{2}})^{-1} 
 (\varrho_{i_1i_2} ) 
\Big\}.
$$
First we bound $L$, for which we need a bound on $ \| \myXt \|_\infty$:
\begin{align*}
     \| \myXt \|_\infty \leqslant&~  \| \Dmut \|_\infty + \| \Dmut \|_\infty \| \myDone \|_\infty \| (\myDonet\Dmut\myDone)^{-1} \|_\infty \| \myDonet \|_\infty \| \Dmut \|_\infty \\
     \leqslant&~ e^{-M_{\Theta,1}} + e^{-M_{\Theta,1}} \times 2 \times  2 e^{M_{\Theta,2}} / n \times (n+\sqrt{2}-1) \times e^{-M_{\Theta,1}} \\
     \leqslant&~ e^{-M_{\Theta,1}} + 6e^{M_{\Theta,2}-2M_{\Theta,1}} 
\end{align*}
where we used $\| (\myDonet\Dmut\myDone)^{-1} \|_\infty \leqslant 2 e^{M_{\Theta,2}} / n $ by Lemma \ref{lem:thm1.1_Hillar} and Remark \ref{rmk::invertable_ddotLalphaalpha}.
Then we can bound $L$ by
\begin{align}
    L \leqslant&~ 2 \cdot \max_{1\leqslant l \leqslant n, 1\leqslant k' \leqslant K} |z_{lk'}| \cdot \max_{1\leqslant i_1\leqslant i_2 \leqslant n} |(\mu_{t,i_{1}i_{2}})^{-1} | \cdot \| \myXt \|_\infty \notag \\
    \leqslant&~ 2 M_{Z,1}^{\frac{1}{2}} e^{M_{\Theta,2}} (e^{-M_{\Theta,1}} + 6e^{M_{\Theta,2}-2M_{\Theta,1}} ) \label{eq::c7_L}
\end{align}

Next we bound the $V_{(j,k')}= \operatorname{var}\Big\{(\Seff)_{(j,k')}\Big\}$.

\begin{align} \label{eq:var_plstar}
V_{(j,k')}
&=\operatorname{var}\left[\sum_{t=1}^T \sum_{1 \leqslant i_1\leqslant i_2 \leqslant n}\Big\{\sum_{l=1}^n\varrho_{jl} z_{lk'} (\myXt)_{(j,l)(i_{1},i_{2})} 
(\varrho_{i_1i_2}^2 \mu_{t,i_{1}i_{2}})^{-1} 
 (\varrho_{i_1i_2}M_{t,i_1 i_2} ) 
\Big\}\right]\notag\\
&\leqslant 2 \sum_{t=1}^T\sum_{1 \leqslant i_1 \leqslant i_2 \leqslant n}\Big(\sum_{l=1}^n z_{lk'} (\myXt)_{(j,l)(i_{1},i_{2})} 
\Big)^2(\mu_{t,i_{1}i_{2}})^{-1}\notag\\
&\leqslant 4\sum_{t=1}^T \sum_{1 \leqslant i_1, i_2 \leqslant n}\Big(\sum_{l=1}^n z_{lk'}(\myCt)_{(j,l)(i_{1},i_{2})}
\Big)^2(\mu_{t,i_{1}i_{2}})^{-1}\notag
\\&\quad+ 4 \sum_{t=1}^T \sum_{1 \leqslant i_1, i_2 \leqslant n}\Big(\sum_{l=1}^n z_{lk'}(\ddot{L}_{\Theta_t,\alpha_t} \myAt^{-1}\myBt)_{(j,l)(i_{1},i_{2})}
\Big)^2(\mu_{t,i_{1}i_{2}})^{-1}. \notag
\\&=: V_1 + V_2
\end{align}
Since $\myCt$ is diagonal, in the first summation 
$$ \sum_{1 \leqslant i_1, i_2 \leqslant n}\Big(\sum_{l=1}^n z_{lk'}(\myCt)_{(j,l)(i_{1},i_{2})}
\Big)^2(\mu_{t,i_{1}i_{2}})^{-1}$$
the summand is non-zero only when $i_1=j$ and $l=i_2$. Thus the first term in \eqref{eq:var_plstar} 
$$V_1 = 4\sum_{t=1}^T \sum_{i_2=1}^n (z_{i_2k'}\varrho_{ji_2} \mu_{t,ji_2} )^2 (\varrho_{ji_2}\mu_{t,ji_2})^{-1} \leqslant 8 M_{Z,1} e^{-M_{\Theta,1}} nT .$$

Next, we consider the second term in \eqref{eq:var_plstar}, for which we derive a bound on 
$$(\myBt^\top\myAt^{-1}\myBt)_{(j,l)(i_{1},i_{2})}.$$ For the ease of notation we let 
\begin{align}
\label{eq:Wtstar_def}
    \mathcal W_t^\star = (\myAt(Z^\star,\alpha^\star))^{-1}
\end{align}
and its $(i,j)$th element is denoted by $\omega^{\star}_{t,ij}$. From \eqref{eq:myAt_star_inv_inftynorm} we have $\max_{i,t} \sum_{j=1}^n |\omega^{\star}_{t,ij}| \leqslant \|\myAt^{-1}(Z^\star,\alpha^\star)\|_\infty \leqslant 2 e^{M_{\Theta,2}}/n $.  Then
\begin{align*}
(\myBt^\top\myAt^{-1}\myBt)_{(j,l)(i_{1},i_{2})} &= [\myCt \myDone \mathcal W_t^\star \myDonet \myCt]_{(j,l)(i_{1},i_{2})}\\
&= \varrho_{jl} \mu_{t,jl} (e_j^{(n)}+e_l^{(n)})^\top \mathcal W^\star_t (e_{i_1}^{(n)}+e_{i_2}^{(n)})\varrho_{i_1i_2} \mu_{t,i_1i_2} \\
&= (\omega^\star_{t,ji_1}+ \omega^\star_{t,ji_2}+ \omega^\star_{t,li_1}+ \omega^\star_{t,li_2})\mu_{t,jl} \mu_{t,i_1i_2}\varrho_{jl}\varrho_{i_1i_2}
\end{align*}
where $e_i^{(n)}\in \mathbb R^n $ deontes the vector whose $i$th element is $1$ and others are $0$. Thus,
\begin{align*}
V_2 = &~4\sum_{t=1}^T \sum_{1 \leqslant i_1 , i_2 \leqslant n}\Big(\sum_{l=1}^n z_{lk'} (\myBt^\top\myAt^{-1}\myBt)_{(j,l)(i_{1},i_{2})}
\Big)^2(\mu_{t,i_{1}i_{2}})^{-1} \\
\leqslant&~ C M_{Z,1} e^{-3M_{\Theta,1}} \Big[ \sum_t \sum_{1 \leqslant i_1 , i_2 \leqslant n}(\sum_{l=1}^n\omega^\star_{t,li_1} )^2 +\sum_t \sum_{1 \leqslant i_1 , i_2 \leqslant n}(\sum_{l=1}^n\omega^\star_{t,li_2} )^2\Big.
\\&~\Big.+\sum_t \sum_{1 \leqslant i_1 , i_2 \leqslant n}(\sum_{l=1}^n\omega^\star_{t,ji_1} )^2 
+\sum_t \sum_{1 \leqslant i_1 ,i_2 \leqslant n}(\sum_{l=1}^n\omega^\star_{t,ji_2} )^2 \Big] \\
\leqslant&~ CM_{Z,1}e^{-3M_{\Theta,1}} \Big[ \sum_t n \sum_{i_1=1}^n (\sum_{l=1}^n \omega^\star_{t,li_1})^2 + \sum_t n^3 \sum_{i_1=1}^n (\omega^\star_{t,ji_1})^2 \Big] 
\\ \leqslant&~ CM_{Z,1}e^{-3M_{\Theta,1}} \Big[ \sum_t n \sum_{i_1=1}^n e^{2M_{\Theta,2}}/n^2 + \sum_t n^3 e^{2M_{\Theta,2}}/n^2 \Big] \\
\leqslant&~ CM_{Z,1}e^{2M_{\Theta,2}-3M_{\Theta,1}} nT.
\end{align*}

Substituting the above bound into \eqref{eq:var_plstar}, we get
\begin{align}
V_{(j,k')} = 
\operatorname{var}\Big\{(\Seff)_{(j,k')} \Big\}
&=\operatorname{var}\left[\sum_{t=1}^T \sum_{1 \leqslant i_1, i_2 \leqslant n}\Big\{\sum_{l=1}^n z_{lk'} (\myXt)_{(j,l)(i_{1},i_{2})} 
(\mu_{t,i_{1}i_{2}})^{-1} 
 (M_{t,i_1 i_2} ) 
\Big\}\right]\notag
\\&\leqslant CM_{Z,1}e^{2M_{\Theta,2}-3M_{\Theta,1}} nT. \label{eq::c7_V}
\end{align}
Plugging \eqref{eq::c7_L} and \eqref{eq::c7_V} into the Bernstein bound \eqref{eq::c7_Bern}, we have $$|(\Seff)_{(j,k')}| < C_{M,s}\sqrt{nT\log(nk)} $$ with probability at least $O((nk)^{-(s+1)}) $, where $C_{M,s} $ is a constant depending only on $M_{\Theta,1}$, $M_{\Theta,2}$, $M_{Z,1}$, and $s$. Then by a union bound, we have $\|(\Seff)\|_\infty < C_{M,s}\sqrt{nTlog(nk)} $ with probability at least $O((nk)^{-s}) $.


Thus, we obtain that for any $s>0$, there exists a universal constant $C_s>0$ such that 
\begin{align*}
\|\Seff(\ZQs, \alpha^\star )\|_2 /(nT)
\leqslant 
C_{s}  \, \big( T^{-1/2}\big)\, \log^{1/2}(nk) 
\end{align*}
with probability $1-O((nk)^{-s}) $.

\subsection{Proof of Remark \ref{RMK:Heff_replace_Ieff}} \label{sec:pfRMK}

The proof consists of three steps similar to that in Section \ref{sec:pfthmnerr} with $\Ieff(\check Z, \check\alpha)$ replaced by $\Heff(\check Z, \check\alpha)$. 
\begin{itemize}
    \item[] \textit{Step 1}:  Note Eq. \eqref{pfthm1-3A} in the main text still holds by the same proof. 
Replacing $\Ieff$ in the proof of Theorem \ref{thm_NRerr} with $\Heff$,
we now consider    
\begin{align}\label{eq:pfthm1-bdratiohess}
\big\|\UZb^\top ( \hat Z_v -  (\ZQs)_v) \big\|_2  \leqslant \frac{ \|-\UZb^\top \Heff(\check Z, \check\alpha) \UZb \UZb^{\top}(   \hat Z_v - (\ZQs)_v) \|_2 }{\sigma_{\min}\big(-\UZb^\top \Heff(\check Z, \check\alpha) \UZb \big)};
\end{align}
    \item[] \textit{Step 2}:   derive an upper bound of  $\|-\UZb^\top \Heff(\check Z, \check\alpha) \UZb \UZb^{\top}(   \hat Z_v - (\ZQs)_v) \|_2$;
    \item[]\textit{Step 3:} obtain a  lower bound of $\sigma_{\min}\big(-\UZb^\top \Heff(\check Z, \check\alpha) \UZb \big)$. 
\end{itemize}
We next present the details of \textit{Steps 2} and \textit{3}.

\bigskip
\noindent \textit{Step 2:}  
Similarly to \eqref{eq:multipliediff2A},  we have
\begin{align}
    -\UZb^\top \Heff(\check Z, \check\alpha) \UZb \UZb^{\top}\big\{  \hat Z_v - (\ZQs)_v\big\}
=&~\UZb^\top \big\{\tilde{S}_1  + S_2+S_3+\Seff\big(\ZQs,\alpha^\star\big)\big\},
\label{eq:multipliediff2hss}
\end{align}
where $S_2$ and $S_3$ are defined in \eqref{eq:thm1s1termA} and \eqref{eq:thm1s3termA}, and we define
\begin{align}
\label{eq:rmk2s1tildeterm}
    \tilde{S}_1=-\Heff(\check Z, \check\alpha) (\check Z_v-  (\ZQs)_v )
\end{align}
which is similar to $S_1$ in \eqref{eq:thm1s1termA} except that $\Ieff(\check Z, \check\alpha) $ is replaced by $\Heff(\check Z, \check\alpha) $. 
By \eqref{eq:multipliediff2hss}, we have
\begin{align} 
&\quad\ \Big\| -\UZb^\top \Heff(\check Z, \check\alpha) \UZb\UZb^\top  ( \hat Z_v -  (\ZQs)_v) \Big\|_2  \notag\\ 
&\leqslant \sqrt{nk} \left(  \|  \tilde{S}_1+S_2\|_{\infty} +\left\|S_3 \right\|_\infty+\left\|  \Seff(Z^\star \check{Q},\alpha^\star)\right\|_2/\sqrt{nk}\right).  \label{pfthm1-1hess}
\end{align}
 We can see that 
 the bounds in \eqref{eq:Seff_starstar_main} and \eqref{eq:S3_infty_main}
 can still be applied to $\| \Seff(Z^\star \check{Q},\alpha^\star)\|_2/\sqrt{nk}$ and $\|S_3\|_{\infty}$ in \eqref{pfthm1-1hess}. 
Furthermore, we can
show the following high probability bound on 
$\|\tilde{S}_1+S_2\|$, similarly 
to \eqref{eq:S1+S2_infty_main}:
for any $s>0$, there exists a universal constant $C_s>0$, such that
\begin{align}
\Pr\left[
\|\tilde S_1+S_2\|_\infty 
\geqslant 
C_{s}T\log^{2\varsigma}(nT)
\right] = O((nT)^{-s}).
\label{eq:tildeS1+S2_infty_D}
\end{align}
(The proof of \eqref{eq:tildeS1+S2_infty_D} in given in Section \ref{subsubsec::pf_lemmaB5} on Page \pageref{subsubsec::pf_lemmaB5}.)

Combining the bounds on the terms in \eqref{pfthm1-1hess}, we obtain that with probability $1-O((nk)^{-s}) $, 
\begin{align}
\eqref{pfthm1-1hess}&\leqslant 
C_{s}\sqrt{nk} \left\{T \log^{2\varsigma}(nT)+ \sqrt{nT\log(nT)} \right\}.
\label{pfthm1-12hess}   
\end{align}

\bigskip
\noindent \textit{Step 3:}
Based on the analyses in \eqref{eq:HeffIeff}  and Lemma \ref{lem_presentingUZmain},  
we give a high probability bound on the smallest eigenvalue of $-\UZb^\top \Heff(\check Z, \check\alpha) \UZb$ in the following lemma.
\begin{lemma}
\label{lem::sigmamin_UZbHeffUZb}
Under the conditions of Theorem \ref{thm_NRerr}, we have for any $s>0$, there exists a universal constant $C_s>0$, such that
\begin{align}
\label{eq::lem_sigmamin_-UZbHeffUZb}
\Pr\left[
\sigma_{\min}\big(-\UZb^\top \Heff(\check Z, \check\alpha) \UZb \big) \leqslant C_{s} nT
\right] = O(n^{-s}).
\end{align}
\end{lemma}
\begin{proof}
See Section \ref{pf::lem::sigmamin_UZbHeffUZb} on Page \pageref{pf::lem::sigmamin_UZbHeffUZb}.
\end{proof}

\bigskip
\noindent \textit{Combining Steps 1-3:} By \eqref{eq:pfthm1-bdratiohess}, \eqref{pfthm1-12hess} and \eqref{eq::lem_sigmamin_-UZbHeffUZb}, when $n/ \log^{2\varsigma}(nT)$ is sufficiently large, we have that for any $s>0$,  there exists a universal constant $C_s>0$, such that
\begin{align*}
\big\| \hat Z_v -  (\ZQs)_v \big\|_2 \leqslant 
C_{s} \left(\frac{\log^{2\varsigma}(nT)\sqrt k}{\sqrt n}  + \frac{\log^{\frac12}(nT) \sqrt{k}}{\sqrt{T}}\right)
\end{align*}
with probability 
$1-O(n^{-s}) $.
The rest of the proof follows from the last step in the proof of Theorem \ref{thm_NRerr}.
\\


\subsubsection{Proof of Eq. \eqref{eq:tildeS1+S2_infty_D}}
\label{subsubsec::pf_lemmaB5}
Recall from \eqref{eq:IeffHeffinXt} and \eqref{eq:HeffIeff} that 
\begin{align*}
    \Ieff = \DZt \big(\sum_t \myXt\big) \DZ, \quad \Heff = -\DZt \big(\sum_t \myXt\big) \DZ + \sum_t N_{2t}(Z,\alpha)\otimes \mathrm{I}_k,
\end{align*}
where $N_{2t}$ is defined in \eqref{eq::def_NtN2t}.
We could see that $-\Heff$ differs from $\Ieff$ only by a mean-zero smaller order term: 
\begin{equation}
\label{eq::def_S4=Ieff+Heff}
    -\Heff(\check Z,\check \alpha) = \Ieff(\check Z,\check \alpha) - S_4 \quad \text{where} \quad S_4 :=\sum_t N_{2t}(\check Z,\check \alpha) \otimes \mathrm{I}_{k} .
\end{equation}
Thus, from the definition of $S_1$ and $\tilde S_1$ in \eqref{eq:thm1s1termA} and \eqref{eq:rmk2s1tildeterm}, we have $\tilde S_1 = S_1 -S_4(\check Z_v - (\ZQs)_v) $. Recall from the proof of \eqref{eq:S1+S2_infty_main} that to bound $\|S_1+S_2\|_\infty$ we decompose it into (with probability $1- O((nT)^{-s}) $)
\begin{equation}
    \label{eq::finalbd_S1+S2_infty_rmk2}
    \begin{aligned}
\|S_1+S_2\|_\infty \leqslant &~ \|S_{22}\|_\infty + \|S_{212}\|_\infty + \|S_{213}\|_\infty + \|S_{232}\|_\infty + \|S_1 + S_{211} - S_{231}\|_\infty \\
\leqslant&~ C_{M,s} T \log^{\epsilon+\max(\epsilon,\frac12)}(nT) +\|S_{213}\|_\infty 
\end{aligned}
\end{equation}
where in the second inequality, we plug in all the bounds in \eqref{eq::finalbd_S1+S2_infty} except for the bound on $\|S_{213}\|_\infty $. (Recall that in \eqref{eq::finalbd_S1+S2_infty} the $\sqrt{nT}$ order term (ignoring the $\log$ factors) comes only from the first part of term $S_{213} $.)
Following \eqref{eq::finalbd_S1+S2_infty_rmk2}, we could bound $\|\tilde S_1 + S_2\|_\infty $ by
\begin{align}
    \|\tilde S_1 + S_2\|_\infty =&~ \|S_1 -S_4(\check Z_v - (\ZQs)_v)+ S_2\|_\infty\notag
    \\ 
    \leqslant&~ C_{M,s} T \log^{\epsilon+\max(\epsilon,\frac12)}(nT) +\|S_{213}-S_4(\check Z_v - (\ZQs)_v)\|_\infty 
    \label{eq::finalbd_tildeS1+S2_infty_rmk2}
\end{align}
 with probability $1- O((nT)^{-s}) $.

Next we verify via direct calculation that $S_4(\check Z_v - (\ZQs)_v) $ cancels with part of the term $S_{213}$ (defined in \eqref{eq::def_S21_123}) in the proof of \eqref{eq:S1+S2_infty_main}. 
We decompose $S_{213} $ into $S_{2131} + S_{2132} $, where
\begin{align*}
S_{2131} :=&~ (\DZt(\check Z) - \DZt(\ZQssub) ) \sum_t\myUt(\check Z,\check\alpha)\\
\text{and} \quad S_{2132} :=&~ (\DZt(\check Z) - \DZt(\ZQssub) ) \big[\sum_t\myUt(\ZQssub,\check\alpha)-\myUt(\check Z,\check\alpha)\big].
\end{align*}
First
consider the term $S_{2132} $.
Same as in \eqref{eq::DZcheck-DZstar}, from Condition \ref{cond_elem_init} on the error of $\check Z$ and by a mean value theorem we have 
$$\|\DZt(\check Z) - \DZt(\ZQssub) \|_\infty \leqslant  2M_b\sqrt{n}\log^{\epsilon}(nT) .$$ Besides, by plugging in the form of $\dot{L}_{\Theta_t} $ in \eqref{eq:ldot_theta_t},
\begin{align*}
    \Big\|\sum_t\big\{\myUt(\ZQssub,\check\alpha)-\myUt(\check Z,\check\alpha)\big\}\Big\|_\infty &\leqslant 2\max_{i,j} \sum_t | \mu_{t,ij}(\ZQssub,\check\alpha) -\mu_{t,ij}(\check Z,\check\alpha) | \\
    &\leqslant 4Te^{-M_{\Theta,1}}M_{Z,1}^{\frac12} (M_b\log^{\epsilon}(nT)/\sqrt n).
\end{align*}
Thus, we have 
\begin{align*}
    \|S_{2132}\|_\infty \leqslant&~ \|\DZt(\check Z) - \DZt(\ZQssub) \|_\infty\Big\|\sum_t\big\{\myUt(\ZQssub,\check\alpha)-\myUt(\check Z,\check\alpha)\big\}\Big\|_\infty \\ \leqslant&~ 8^{-M_{\Theta,1}}M_{Z,1}^{\frac12} M_b^2T\log^{2\epsilon}(nT) . 
\end{align*}
Next we show $S_{2131}$ cancels with $S_{4}(\check Z_v - (\ZQs)_v) $:
\begin{align*}
S_{2131} &= (\DZt(\check Z) - \DZt(\ZQs) ) \sum_t\myUt(\check Z,\check\alpha) \\ 
&=\sum_t \begin{pmatrix}
\sum_{i=1}^n \varrho_{1i}^2 (A_{t,1i}  - \mu_{t,1i}(\check Z, \check\alpha) ) (\check z_i - z_i^\star)\\
\vdots \\
\sum_{i=1}^n \varrho_{ni}^2 (A_{t,ni}  - \mu_{t,ni}(\check Z, \check\alpha) ) (\check z_i - z_i^\star)
\end{pmatrix}\\
&=\sum_t \left[N_{2t}(\check Z,\check \alpha) \otimes \mathrm{I}_{k}\right] (\check Z_v - (\ZQs)_v ) = S_{4}(\check Z_v - (\ZQs)_v),
\end{align*}
where recall $\varrho_{ij} $ is defined in \eqref{eq::def_varrhoij}  so that the $(i,j)$th element of $\myUt\in \mathbb R^{{n(n+1)}/{2}}$ is $\varrho_{ij} (A_{t,ij} - \mu_{t,ij})$.

To conclude, from \eqref{eq::finalbd_tildeS1+S2_infty_rmk2} we have with probability $1- O((nT)^{-s}) $
\begin{align*}
    \|\tilde S_1 + S_2\|_\infty \leqslant&~ C_{M,s} T \log^{\epsilon + \max(\epsilon,\frac{1}{2})}(nT) + 
    \|S_{2131} + S_{2132} - S_4(\check Z_v - (Z^\star\check Q)_v\|_\infty \\
    =&~ C_{M,s} T \log^{\epsilon + \max(\epsilon,\frac{1}{2})}(nT) + 
    \| S_{2132} \|_\infty \\ \leqslant&~ C_{M,s} T \log^{\epsilon + \max(\epsilon,\frac{1}{2})}(nT) \leqslant~ C_{M,s} T \log^{2\varsigma}(nT),
\end{align*}
where $C_{M,s}$ is a constant that only depends on $M_{\Theta,1},M_{\Theta,2},M_{Z,1},M_b $ and $s$.
Thus, we have proved the statement of \eqref{eq:tildeS1+S2_infty_D}:
for any $s>0$, there exists a universal constant $C_s>0$, such that
\begin{align*}
\Pr\left[
\|\tilde S_1+S_2\|_\infty 
\geqslant 
C_{s}T\log^{2\varsigma}(nT)
\right] = O((nT)^{-s}).
\end{align*}

\subsubsection{Proof of Lemma \ref{lem::sigmamin_UZbHeffUZb}}
\label{pf::lem::sigmamin_UZbHeffUZb}
From \eqref{eq::def_S4=Ieff+Heff}, we have
\begin{align}
\label{eq::sigmamin_Heff_decomp}
    \sigma_{\min}\big(-\UZb^\top \Heff(\check Z, \check\alpha) \UZb \big) \geqslant \sigma_{\min}\big(\UZb^\top \Ieff(\check Z, \check\alpha) \UZb \big) - \|S_4\|_{\text{op}}.
\end{align}
The lower bound on $\sigma_{\min}\big(\UZb^\top \Ieff(\check Z, \check\alpha)\UZb\big)$ has been obtained in \eqref{pfthm1-2A}; next we derive the upper bound on the operator norm of $S_4$.
We decompose  
\begin{equation}
\label{eq::decomp_N2tcheck}
S_4=\sum_t N_{2t}(\check Z,\check \alpha) \otimes \mathrm{I}_{k} 
=\sum_t M_{2t}\otimes \mathrm{I}_{k}  -\sum_t \{\bdmu_{2t}( \check Z,\check \alpha) - \bdmu_{2t}( Z^\star,\alpha^\star) \}\otimes \mathrm{I}_{k} ,
\end{equation}
where $N_{2t} $ is defined in  \eqref{eq::def_NtN2t}, and 
\begin{align}
    M_{2t} :=&~ N_{2t}(Z^\star,\alpha^\star) \ \in\ \mathbb{R}^{n\times n}, \label{eq::def_M2t}
    \\
    \bdmu_{2t}(Z,\alpha) :=&~ \bdmu_t(Z,\alpha) + \operatorname{diag}(\mu_{t,11}(Z,\alpha),\ldots,\mu_{t,nn}(Z,\alpha))\ \in\ \mathbb{R}^{n\times n}.  \label{eq::def_bdmu2t}
\end{align}
For the first score at true $(Z^\star,\alpha^\star)$ part in \eqref{eq::decomp_N2tcheck}, Lemma \ref{Lemma-concen}
gives that  
$$\Big\|\sum_t  M_{2t} \otimes \mathrm{I}_{k}  \Big\|_{\text{op}} \leqslant C_s \sqrt{nTe^{-M_{\Theta,1}} + \log(nk) }$$ with probability $1-n^{-s}$.
For the second part, from
 Condition \ref{cond_elem_init} we have 
$$\Big\|\sum_t  \{\bdmu_{2t}( \check Z,\check \alpha) - \bdmu_{2t}( Z^\star,\alpha^\star) \}\otimes \mathrm{I}_{k}  \Big\|_{\mathrm{F}} \leqslant C(1+\sqrt{M_{Z,1}})e^{-M_{\Theta,1}}M_b T\sqrt{nk}\log^{\epsilon}(nT). $$ 
Plugging the above two bounds and \eqref{pfthm1-2A} into \eqref{eq::sigmamin_Heff_decomp}, the smallest eigenvalue of $-\UZb^\top \Heff(\check Z,\check\alpha) \UZb$ is bounded from below by
$$nT\check M_{Z,2} e^{-\check M_{\Theta,2}} - C_s\sqrt{nTe^{-M_{\Theta,1}} + \log(nk) }- C(1+\sqrt{M_{Z,1}})e^{-M_{\Theta,1}}M_b T\sqrt{nk}\log^{\epsilon}(nT) $$ with probability $1-n^{-s} $,
in which $\check M_{Z,2} $ and $\check{M}_{\Theta,2} $ are defined after \eqref{pfthm1-2A}.
Under the condition that $n/\log^{2\varsigma}(nT) $ is sufficiently large, we have 
$$\sigma_{\operatorname{min}}(-\UZb^\top \Heff(\check Z,\check\alpha) \UZb) \geqslant C_{M,s} nT$$
with probability $1-O(n^{-s}) $, where $C_{M,s}$ is a constant that only depends on $M_{\Theta,1}$, $M_{\Theta,2}$, $M_{Z,1}$, $M_{Z,2}$, $M_b$, and $s$.


\newpage 
\section{Technical Results for the Initial Estimator}\label{sec:techinitial}

\subsection{Details for the Initialization Algorithm} \label{sec:details}

\begin{remark}
We point out that \eqref{eq:grandzformulaalo}--\eqref{eq:grandalphaformulaalo} in Stage 2 are slight modifications from the gradients of $L(Z,\alpha)$ with respect to $(Z,\alpha)$ and can be viewed as gradients of a pseudo likelihood. 
 We use the ``pseudo gradients''  in Algorithm \ref{algor:init} only for    the convenience of proving Theorem \ref{thm:initialerror}. 
In particular, define the pseudo likelihood 
\begin{align*}
    \tilde{L}(Z,\alpha)=\sum_{t=1}^T\sum_{1\leqslant i\leqslant j\leqslant n}\left[A_{t,ij} (\alpha_{it}+\alpha_{jt}+z_i^\top z_j) - \exp(\alpha_{it}+\alpha_{jt}+z_i^\top z_j) \right] w_{ij}
\end{align*}
with $w_{ij}=1$ if $i\neq j$ and $w_{ii}=1/2$.  
We have 
$g_{\alpha_t}(Z,\alpha)=\frac{\partial \tilde{L}}{\partial \alpha_t}(Z,\alpha)$ and 
$g_{Z}(Z,\alpha) = \frac{\partial \tilde{L}}{\partial Z}(Z ,\alpha ) $. 
In addition, \eqref{eq:grandzformulaalo2}--\eqref{eq:grandalphaformulaalo2} in Stage 2 are gradients of the following pseudo likelihood
\begin{align}\label{eq:halffixedlzalpha}
    \sum_{t=1}^T\sum_{1\leqslant i, j\leqslant n} \left[A_{t,ij} (\alpha_{it}+\alpha_{jt}^{R_1}+z_i^\top z_j^{R_1}) - \exp(\alpha_{it}+\alpha_{jt}^{R_1}+z_i^\top z_j^{R_1}) \right] 
\end{align}
evaluated at $(Z^r,\alpha^r)$.
\end{remark}

\begin{remark}[\textbf{Implementation Details}]\label{rm:e2}
We next present the choice of tuning parameters and implementation details of Algorithm \ref{algor:init} in simulations. 
\begin{enumerate}
\setlength{\itemsep}{0pt} 
    \item[(i)]
For thresholds $\tau_t$, 
we set   
$\tau_{t} = (n\widehat{p}_{t})^{1/2}$ with  $\widehat{p}_{t} = \sum_{ij} A_{t,ij} /n^2$. 
\item[(ii)] For $\eta_Z$ and $\eta_{\alpha}$, we use Barzilai-Borwein step sizes \citep{barzilai1988two}, instead of fixed step sizes,  to further improve computational efficiency. Specifically, in each iteration $r\geqslant 1$, we use step sizes defined as
$$
\eta_{Z}^{r}=\frac{\left\langle Z^{r}-Z^{r-1}, g_{Z}^{r}-g_{Z}^{r-1}\right\rangle}{\left\langle g_{Z}^{r}-g_{Z}^{r-1}, g_{Z}^{r}-g_{Z}^{r-1}\right\rangle} \quad \quad\text{and}\quad \quad \eta_{\alpha}^{r}=\frac{\left\langle\alpha^{r}-\alpha^{r-1}, g_{\alpha}^{r}-g_{\alpha}^{r-1}\right\rangle}{\left\langle g_{\alpha}^{r}-g_{\alpha}^{r-1}, g_{\alpha}^{r}-g_{\alpha}^{r-1}\right\rangle} ,
$$
where we let $g_{Z}^{r}=- g_Z \left(Z^{r}, \alpha^r\right)$ and $g_{\alpha}^{r} = - (g_{\alpha_1}\left(Z^{r},\alpha^r\right), \ldots, g_{\alpha_T}\left(Z^{r},\alpha^r\right))$.
\item[(iii)] For the iteration numbers, 
in principle, $R_1$ and $R_2$ should be sufficiently large such that the likelihood does not change too much. 
In simulations, we find that it suffices to set $R_1 = 1000$ in Stage 2-1 and skip Stage 2-2.  
Theoretically, we propose Stage 2-2 for the ease of establishing elementwise error bounds in Theorem \ref{thm:initialerror}. 
\item[(iv)] For the projection in lines 17--18 and  23--24,
we simplify them as 
$Z^{r+1}= J\widetilde{Z}^{r+1}$
and $\alpha^{r+1}=\widetilde{\alpha}^{r+1}$, i.e., ignoring the constraints \eqref{eq:czsetproj}--\eqref{eq:czprimesetproj} that depend on unknown parameters $(M_{Z,1},M_{\alpha})$.  This leads to equally good performance in simulations and is similar to the simplification in Eq. (13) of  \cite{ma2020universal}. 
\end{enumerate} 
\end{remark}

\subsection{Proof of Theorem \ref{thm:initialerror}} \label{sec:pfthm2}

In this section, we establish upper bounds of estimation errors of   $\check{\alpha}_{it}$ and $\check{z}_i$ from 
Algorithm \ref{algor:init} for $1\leqslant i\leqslant n$ and $1\leqslant t\leqslant T$.

%
%
%



\vspace{1em}

The proof of Theorem \ref{thm:initialerror} consists of two parts. 
In Section \ref{sec:pfstage1}, 
we first establish deterministic error upper bounds of $\check{\alpha}_{it}$ and $\check{z}_i$ (i.e., \eqref{eq:alphadeterminbound} and \eqref{eq:zicheckdetermbd}) through the following random quantities 
  \begin{align}
 \zeta_{\infty} = \max_{1\leqslant t\leqslant T} \zeta_t, \quad \quad 
\zeta_{i,\infty}^2 = \max_{1\leqslant t\leqslant T} \zeta_{it}^2,  \quad \text{ and }\quad 
\zeta_{\infty,\infty}^2 =  \max_{1\leqslant i\leqslant n} \zeta_{i,\infty}^2, \label{eq:zetadef}
  \end{align}   
where we define that for $1\leqslant i\leqslant n$ and $1\leqslant t\leqslant T$,  
\begin{align*}
    \zeta_{t} =&~  \|\mathbf A_{t} - \exp(\Theta_{t}^{\star})\|_{\mathrm{op}} + \sqrt{n}, \\
       \zeta_{it}^2 = &~ \Big(\sum_{j=1}^n M_{t,ij}\Big)^2 + \Big\| \sum_{j=1}^n  M_{t,ij}z_{j}^{\star}\Big\|_2^2+ \sum_{j=1}^n M_{t,ij}^2 + n, \  \text{ with }\ M_{t,ij}=A_{t,ij} -  \exp(\Theta_{t,ij}^{\star}).
\end{align*}
Then, in Section \ref{sec:pfstage12prob}, we establish probabilistic upper bounds of the random quantities in \eqref{eq:zetadef}, which leads to the conclusion in  Theorem \ref{thm:initialerror}. 

\subsubsection{Deterministic Upper Bounds}\label{sec:pfstage1}
The proof of \eqref{eq:alphadeterminbound} and \eqref{eq:zicheckdetermbd} consists of 6 steps:
{Step 1} establishes  a deterministic error upper bound of $\Theta_t^0=\mathring{\Theta}_t=\mathring{\alpha}_t1_n^{\top}+1_n\mathring{\alpha}_t^{\top} +\mathring{Z}\mathring{Z}^{\top}$,  where $(\mathring{\alpha},\mathring{Z})$ is the ``initial of initial'' estimator from Stage 1.  
{Steps 2--6} examine the $(Z^r,\alpha^r)$ in Stage 2 of Algorithm \ref{algor:init}. In particular, Steps 2--6 establish upper bounds of the following terms
\begin{align}
   &~ e^{r} =2n\left\|{\alpha^{r}}-\alpha^{\star}\right\|_{\mathrm{F}}^{2} + T\left\|Z^{\star}\right\|_{\mathrm{op}}^{2}\left\|Z^{r} - Z^{\star}Q_r\right\|_{\mathrm{F}}^{2},  \notag\\ 
 &~e_t^r =2n \|\alpha_t^r - \alpha_t^{\star}\|_2^2, \notag\\
 &~e_{i\cdot}^r =2n \|\alpha_{i\cdot}^r - \alpha_{i\cdot}^{\star}\|_2^2 + T\|Z^{\star}\|_{\mathrm{op}}^2 \|z_{i}^r - Q_{R_1}^{\top} z_{i}^{\star}\|_2^2, \notag\\
 &~e_{it}^r =  2n (\alpha_{it}^r - \alpha_{it}^{\star})^2, \notag\\
&~\mathrm{dist}_i^2(\check{z}_i, z_i^{\star})= \big\|\check{z}_{i} - \check{Q}^\top z^\star_i \big\|_2^2, \notag
\end{align}
respectively, 
where
$Q_{r}=\argmin_{Q \in \mathcal{O}(k)}\left\|Z^{r}-Z^{\star} Q\right\|_{\mathrm{F}}$, $\alpha_{i\cdot}^{\star}=(\alpha_{i1}^{\star},\ldots, \alpha_{iT}^{\star})^{\top}\in \mathbb{R}^T$,  and $\alpha_{i\cdot}^{r}=(\alpha_{i1}^{r},\ldots, \alpha_{iT}^{r})^{\top} \in \mathbb{R}^T$. 


\bigskip
 \textit{Step 1:} Lemma \ref{lem:in of in} establishes an upper bound of $\sum_{t = 1}^T \left\|\Theta_{t}^0 - \Theta_{t}^{\star}\right\|_{\mathrm{F}}^2$. 
\begin{lemma} \label{lem:in of in}
    Assume Condition \ref{cond:truvalueregularity} and Condition \ref{cond:tuning}. 
For $t=1,\ldots, T$,  if we set $\tau_t = c \zeta_t$ for a  constant $c > 1$,  there exists a positive constant $C$ such that 
	$$
	\sum_{t = 1}^T \left\|\Theta_{t}^0 - \Theta_{t}^{\star}\right\|_{\mathrm{F}}^2 \leqslant C e^{2M_{\Theta,2}}M_{\Theta,2}Tn^{\frac{3}{2}-\frac{1}{k+3}} \zeta_{\infty},
	$$
where $\zeta_{\infty}$ is defined in \eqref{eq:zetadef} and $M_{\Theta, 2}$ is defined in Condition \ref{cond:truvalueregularity}. 
\end{lemma}
\begin{proof}
See Section \ref{sec:pf_inofin} on Page \pageref{sec:pf_inofin}.
\end{proof}


\bigskip

Steps 2--6 assume
Condition \ref{cond:zstareigenbd} below, which assumes a lower bound of $\|Z^{\star}\|_{\mathrm{op}}$ through  $\zeta_{\infty}$ in \eqref{eq:zetadef}. We will show Condition \ref{cond:zstareigenbd} holds with a high probability  in Section \ref{sec:pfstage12prob}. 
\begin{condition}\label{cond:zstareigenbd}
There exists a large constant $C_1$ such that 
$$\left\|Z^{\star}\right\|_{\mathrm{op}}^{4} \geqslant C_{1}  e^{4M_{\Theta,2}} \kappa^{4} k \max\Big\{ \zeta_{\infty}^2,\, M_{\Theta,2} n^{\frac{3}{2}- \frac{1}{k+3}} \zeta_{\infty} \Big\},$$
where $\zeta_{\infty}$ is defined in \eqref{eq:zetadef} and $\kappa$ represents the condition number of $Z^{\star}$. 
\end{condition}


\bigskip
 \textit{Step 2:} Lemma \ref{thm:linear} establishes linear convergence of $e^r$ for $1\leqslant r\leqslant R_1$.  

\begin{lemma} \label{thm:linear}
Assume conditions in Lemma \ref{lem:in of in} 
 and Condition \ref{cond:zstareigenbd}. 
  There exist positive constants $\rho$ and $C$ such that for any $r \leqslant R_1$, 
	$$e^{r} \leqslant 2\left(1-\frac{\eta  \rho}{e^{M_{\Theta,2}} \kappa^{2}}\right)^{r} e^{0}+\frac{C e^{2M_{\Theta,2}} \kappa^{2} k T \zeta_{\infty}^2}{\rho},$$
 where $\zeta_{\infty}$ are defined in \eqref{eq:zetadef}. 
    Since we choose $R_1$ to be sufficiently large, there exists an integer $R_{0.5} < R_1$ such that for any $r = R_{0.5}, \ldots, R_1$,  
    $$e^{r} \leqslant {C e^{2M_{\Theta,2}}\kappa^{2}k  T\zeta_{\infty}^2},$$
    where $C$ denotes a universal constant. 
\end{lemma}

\begin{proof}
See Section \ref{sec:pf_linear1} on Page \pageref{sec:pf_linear1}.
\end{proof}


\bigskip
 \textit{Step 3:} Lemma \ref{thm:linear2} establishes linear convergence  of $e^r_t$ for $R_{0.5}\leqslant r\leqslant R_1$, where $R_{0.5}$ is specified in Lemma \ref{thm:linear}.

\begin{lemma} \label{thm:linear2}
	Under the conditions of Lemma \ref{thm:linear}, 
	there exist positive constants $\rho$ and $C$ such that for all $t = 1, \ldots, T$ and $r = R_{0.5}, \ldots, R_1$, 
	$$e_{t}^{r} \leqslant \left(1-\frac{\eta \rho}{e^{M_{\Theta,2}} } \right)^{r - R_{0.5}} e_{t}^{R_{0.5}}+\frac{C e^{2M_{\Theta,2}}}{\rho}\left(\zeta_t^2 + e^{4M_{\Theta,2}}\kappa^{2}k \zeta_{\infty}^2 \right),$$
  where $\zeta_{t}$ and $\zeta_{\infty}$ are defined in \eqref{eq:zetadef}. 
	Then, with a large enough $R_1$, we have
    $$\max_{1\leqslant t\leqslant T} e_{t}^{R_1} \leqslant Ce^{6M_{\Theta,2}}\kappa^2k\zeta_{\infty}^2.$$
\end{lemma}
\begin{proof}
See Section \ref{sec:pf_linear2} on Page \pageref{sec:pf_linear2}.
\end{proof}

\bigskip
 \textit{Step 4:} Lemma \ref{thm:linear3} establishes linear convergence of $e^r_{i\cdot}$ for $r\geqslant R_1$. 


\begin{lemma} \label{thm:linear3}
	Under the conditions of Lemma \ref{thm:linear}, 
 there exist positive constants $\rho$ and $C$ such that for all $r \geqslant R_1$ and $i = 1, \ldots, n$,
	$$e_{i\cdot}^{r} \leqslant 2\left(1-\frac{\eta \rho}{e^{M_{\Theta,2}} \kappa^{2}} \right)^{r - R_1} e_{i\cdot}^{R_1}+\frac{C e^{2M_{\Theta,2}}\kappa^{4}T}{\rho} \left( \frac{\zeta_{i,\infty}^2 + e^{4M_{\Theta,2}}  M_{\Theta,2} \kappa^2 k \zeta_{\infty}^2}{\min\{n,\left\|Z^{\star}\right\|_{\mathrm{op}}^{2}\}} + \frac{e^{2M_{\Theta,2}} \kappa^2 k \zeta_{\infty}^2\zeta_{i,\infty}^2 }{\left\|Z^{\star}\right\|_{\mathrm{op}}^{4}}\right),$$
  where $\zeta_{\infty}$ and $\zeta_{i,\infty}$ are defined in \eqref{eq:zetadef}. 
    It follows that there exists an integer $R_1< R_{1.5} <R_2$ such that for any $R_{1.5}\leqslant r \leqslant R_2$, we have 
    $$\max_{1\leqslant i\leqslant n} e_{i\cdot}^{r} \leqslant {C e^{6M_{\Theta,2}} M_{\Theta,2} \kappa^{6} k T} \  \frac{\zeta_{\infty}^2\zeta_{\infty,\infty}^2 }{\min\{n^2,\left\|Z^{\star}\right\|_{\mathrm{op}}^{4}\}},$$
      where $\zeta_{\infty}$ and $\zeta_{\infty,\infty}$ are defined in \eqref{eq:zetadef}. 
\end{lemma}

\begin{proof}
See Section \ref{sec:pf_linear3} on Page \pageref{sec:pf_linear3}.
\end{proof}

\bigskip
 \textit{Step 5:} Lemma \ref{thm:linear4} establishes  linear convergence of $e^r_{it}$ for $R_{1.5}\leqslant r \leqslant R_2$, where $R_{1.5}$ is specified in Lemma \ref{thm:linear3}. 


\begin{lemma} \label{thm:linear4}
	Under the conditions of Lemma \ref{thm:linear}, there exist positive constants $\rho$ and $C$ such that for all $i = 1, \ldots, n$, $t = 1, \ldots, T$ and $r \geqslant R_{1.5}$
	$$e_{it}^{r} \leqslant \left(1-\frac{\eta \rho}{e^{M_{\Theta,2}} } \right)^{r - R_{1.5}} e_{it}^{R_{1.5}}+\frac{C e^{2M_{\Theta,2}}}{\rho} \left(\frac{\zeta_{it}^2}{n} + \frac{e^{8M_{\Theta,2}}\kappa^2k\zeta_{\infty}^2}{\min\{n,\left\|Z^{\star}\right\|_{\mathrm{op}}^{2}\}} + \frac{e^{8M_{\Theta,2}} M_{\Theta,2} \kappa^6 k \zeta_{\infty}^2 \zeta_{\infty,\infty}^2 }{\min\{n^2,\left\|Z^{\star}\right\|_{\mathrm{op}}^{4}\}}\right),$$
           where $\zeta_{\infty}$ and $\zeta_{\infty,\infty}$ are defined in \eqref{eq:zetadef}. 
	With a large enough $R_2$, we have
    $$\max_{1\leqslant i\leqslant n, 1\leqslant t\leqslant T} e_{it}^{R_2} \leqslant Ce^{10M_{\Theta,2}}M_{\Theta,2}\kappa^6k   \frac{\zeta_{\infty}^2\zeta_{\infty,\infty}^2 }{\min\{n^2,\left\|Z^{\star}\right\|_{\mathrm{op}}^{4}\}}.$$
\end{lemma}

\begin{proof}
See Section \ref{sec:pf_linear4} on Page \pageref{sec:pf_linear4}.
\end{proof}


As $\check{\alpha}=\alpha^{R_2}$, Lemma \ref{thm:linear4} implies that there exists a positive constant $C$ such that 
\begin{align}\label{eq:alphadeterminbound}
    \max_{1\leqslant i\leqslant n, 1\leqslant t\leqslant T} |\check{\alpha}_{it}-\alpha_{it}^{\star}|^2 \leqslant Ce^{10M_{\Theta,2}}M_{\Theta,2}\kappa^6k   \frac{\zeta_{\infty}^2\zeta_{\infty,\infty}^2 }{n \cdot \min\{n^2,\left\|Z^{\star}\right\|_{\mathrm{op}}^{4}\}}.
\end{align}

\medskip
 \textit{Step 6:} Lemma \ref{thm:linear5} establishes an upper bound of $  \max_{1\leqslant i\leqslant n} \|\check{z}_{i} - \check{Q}^\top z^\star_i \|_2^2  $. 
\begin{lemma} \label{thm:linear5}
	Under the conditions of Lemma \ref{thm:linear}, there exists a positive constant $C$ such that 
 \begin{align}\label{eq:zicheckdetermbd}
    \max_{1\leqslant i\leqslant n} \big\|\check{z}_{i} - \check{Q}^\top z^\star_i \big\|_2^2  \leqslant {C e^{6M_{\Theta,2}} M_{\Theta,2}^2 \kappa^{8} k}  \, \frac{n\zeta_{\infty}^2\zeta_{\infty,\infty}^2 }{\min\{n^4,\left\|Z^{\star}\right\|_{\mathrm{op}}^{8}\}}
\end{align}
           where $\zeta_{\infty}$ and $\zeta_{\infty,\infty}$ are defined in \eqref{eq:zetadef}. 
 \end{lemma}
\begin{proof}
See Section \ref{sec:pf_linear5} on Page \pageref{sec:pf_linear5}.
\end{proof}


\subsubsection{Probabilistic Upper Bounds}\label{sec:pfstage12prob}
The deterministic upper bounds in \eqref{eq:alphadeterminbound} and \eqref{eq:zicheckdetermbd} depend on random quantities  $\zeta_{\infty}$ and $\zeta_{\infty, \infty}$.
Lemma \ref{lem:concen} below establishes probabilistic upper bounds of $\zeta_{\infty}$ and $\zeta_{\infty, \infty}$.

\begin{lemma}\label{lem:concen}
Under Condition \ref{cond:truvalueregularity},  there exists a constant $C>0 $ such that  for any   $s >0$
\begin{align*}
   &\Pr\left( \zeta_{\infty}^2 \geqslant C(s+1)^2 n \log^2(nT)  \right) \leqslant (nT)^{-s},\\     &\Pr\left( \zeta_{\infty, \infty}^2  \geqslant C(s+2)^2 e^{M_{\Theta,2}} nk \log^2(nT) \right) \leqslant (k+2)(nT)^{-s}.
\end{align*}
\end{lemma}

\begin{proof}
See Section \ref{sec:pf_concen} on Page \pageref{sec:pf_concen}.
\end{proof}

We next illustrate the additional conditions required by Lemmas \ref{lem:in of in}--\ref{thm:linear5} are satisfied with high probability according to Lemma \ref{lem:concen} and Lemma \ref{Lemma-concen}. Firstly, we have $\zeta_t \asymp \sqrt{n}$ with high probability by Lemma \ref{Lemma-concen} and the definition of $\zeta_t$ in \eqref{eq:zetadef}. Combining with Condition \ref{cond:tuning}(i), the thresholds are set as $\tau_t = c \zeta_t$ for $c > 1$, which is the additional condition of Lemma \ref{lem:in of in}. Secondly, by Condition \ref{cond:truvalueregularity} we have
$$
nM_{Z,2} \leqslant \sigma_k^2(Z^{\star}) \leqslant \left\|Z^{\star}\right\|_{\mathrm{op}}^{2} \leqslant \left\|Z^{\star}\right\|_{\mathrm{F}}^{2} \leqslant nM_{Z,1} \quad\text{and}\quad 1 \leqslant \kappa \leqslant M_{Z,1}/M_{Z,2},
$$
which indicates $\left\|Z^{\star}\right\|_{\mathrm{op}}^2 \asymp n$ and $\kappa \asymp 1$. Then when $n/\log^{k+3}(T)$ is sufficiently large, we know Condition \ref{cond:zstareigenbd} is satisfied for any constant $s >0 $ with probability $1 - (nT)^{-s}$ since we consider all $M$ and $k$ as constants.

In summary,
combining \eqref{eq:alphadeterminbound}, \eqref{eq:zicheckdetermbd}, Lemma \ref{lem:concen},
under the conditions of Theorem \ref{thm:initialerror}, we have
\begin{align*}
    {\max_{1\leqslant i\leqslant n,1 \leqslant t \leqslant T}} \left( \left|\check\alpha_{it} - \alpha_{it}^\star\right|^2 +
    \|\check z_{i} - \check{Q}^\top z^\star_i \|_2^2 \right) \leqslant 
    C  (s+2)^4 e^{12M_{\Theta,2}} \kappa^8 k^2 \,  \frac{n^3 \log^4(nT)}{\min\{n^4,\left\|Z^{\star}\right\|_{\mathrm{op}}^{8}\}}
\end{align*}
happens with probability $1 - (k+3) \, (nT)^{-s}$, where we used $M_{\Theta,2} \leqslant e^{M_{\Theta,2}}$. Treated $M, k, s$ as constants and substitute $\left\|Z^{\star}\right\|_{\mathrm{op}}^2 \asymp n$ and $\kappa \asymp 1$ into the above inequality, we have
$$
\Pr\left[{\max_{1\leqslant i\leqslant n,1 \leqslant t \leqslant T}} \left( \left|\check\alpha_{it} - \alpha_{it}^\star\right|^2 +
    \|\check z_{i} - \check{Q}^\top z^\star_i \|_2^2 \right) > 
     \frac{  C_{M,s,k}\log^4(nT) }{{n}} \right] \leqslant (k+3) (nT)^{-s}
$$
for any $s > 0$, where $C_{M,s,k}$ is a constant that only depends on $M_{\Theta,2}, M_{Z,1}, M_{Z,2}, s$, and $k$.

\subsection{Proof of Lemmas \ref{lem:in of in}--\ref{lem:concen} in Section \ref{sec:pfthm2}} 

\subsubsection{Proof of Lemma \ref{lem:in of in}} \label{sec:pf_inofin}

By Lemma \ref{Lemma-3.5}, when we select $\tau_t = c \zeta_t$ for any $c > 1$, we have
\begin{align} \label{lem1-1}
\big\|\widetilde{E}_t - \exp(\Theta^{\star}_t)\big\|_{\mathrm{F}}^2
\leqslant C  \zeta_t \left\|\exp(\Theta^{\star}_t)\right\|_{\mathrm{*}}.
\end{align}
Then we want to find a low-rank approximation of $\exp(\Theta_t^\star)$.
Let $ \mathcal{N}_\alpha$ be the smallest $\delta\text{-net}$ of $\{\alpha \in  \mathbb{R}:|\alpha|\leqslant{M_{\alpha}}\}$, 
 $\mathcal{N}_Z$ be the smallest $3 M_{Z,1}^{-1/2} \cdot \delta \text{-net}$ of $\{z\in \mathbb{R}^k:\|z\|_{\mathrm{2}}^2\leqslant{M_{Z,1}}\}$.
When $\delta$ is small enough, we have $$|\mathcal{N}_\alpha|\leqslant \frac{M_{\alpha}}{\delta},\quad |\mathcal{N}_Z|\leqslant\left( \frac{M_{Z,1}}{\delta}\right)^k.$$
Select $\alpha_{it}^\dagger$ as the closest element to $\alpha_{it}^\star$ in $\mathcal{N}_{\alpha}$, and $z_i^\dagger$ as the closest element to $z_i^\star$ in $\mathcal{N}_Z$.
Let $\Theta_{t,ij}^\dagger=\alpha_{it}^\dagger+\alpha_{jt}^\dagger+(z_i^\dagger)^\top(z_j^\dagger),$ then
$$\operatorname{rank}(\exp(\Theta_t^\dagger))\leqslant|\mathcal{N}_\alpha|\cdot|\mathcal{N}_Z|\leqslant\left(\frac{M_{\Theta,2}}{\delta}\right)^{k+1}.$$ 
Since $\Theta_{t,ij}^{\star} \leqslant -M_{\Theta,1}$ and $\delta$ is small enough, it is reasonable to assume $\Theta_{t,ij}^{\dagger} \leqslant 0$. Then 
\begin{align*}
&\quad \ \big|\exp(\Theta_{t,ij}^\star)-\exp(\Theta_{t,ij}^\dagger)\big|
\leqslant\big|\Theta_{t,ij}^\star-\Theta_{t,ij}^\dagger\big|\\
& \leqslant |\alpha_{it}^\star-\alpha_{it}^\dagger|
+|\alpha_{jt}^\star-\alpha_{jt}^\dagger|
+\|z_i^\star-z_i^\dagger\|_{\mathrm{2}}\cdot\|z_j^\star\|_{\mathrm{2}}
+\|z_j^\star-z_j^\dagger\|_{\mathrm{2}}\cdot\|z_i^\star\|_{\mathrm{2}}\\
&\leqslant 8\delta. 
\end{align*}
Therefore, we can obtain
$$
\left\| \exp(\Theta_t^\star)\right\|_{\mathrm{*}}
\leqslant \big\|\exp(\Theta_t^\star)-\exp(\Theta_t^\dagger)\big\|_{\mathrm{*}}
+\big\|\exp(\Theta_t^\dagger)\big\|_{\mathrm{*}}
\leqslant 8n^{1.5}\delta+\sqrt{\frac{M_{\Theta,2}^{k+1}}{\delta^{k+1}}}\cdot n.
$$
Set $\delta=c_1M_{\Theta,2}n^{-{1}/{(k+3)}}$ with a small constant $c_1$, then 
\begin{align} \label{lem1-2}
\left\|\exp(\Theta_t^\star)\right\|_{\mathrm{*}}\leqslant CM_{\Theta,2}n^{\frac{3}{2}-\frac{1}{k+3}}.
\end{align}
Combining \eqref{lem1-1} and \eqref{lem1-2}, we have
\begin{align*}
\big\|\widetilde{E}_t - \exp(\Theta_t^\star)\big\|_{\mathrm{F}}^2
\leqslant C {M_{\Theta,2}n^{\frac{3}{2}-\frac{1}{k+3}}} \zeta_t.
\end{align*}
Because $\log(\cdot)$ is $e^{M_{\Theta,2}}\text{-Lip}$ in $[e^{-M_{\Theta,2}},1]$, then we can obtain that
\begin{align}
   \big\|\mathring{\Theta}_t - \Theta_t^\star\big\|_{\mathrm{F}}^2
    &\leqslant {e^{2M_{\Theta,2}}}\big\|\mathring{E}_t-\exp(\Theta_t^\star)\big\|_{\mathrm{F}}^2\notag\\
    &\leqslant {e^{2M_{\Theta,2}}}\big\|\widetilde{E}_t-\exp(\Theta_t^\star)\big\|_{\mathrm{F}}^2\notag\\
    &\leqslant C {e^{2M_{\Theta,2}}}  {M_{\Theta,2}n^{\frac{3}{2}-\frac{1}{k+3}}} \zeta_t.\label{th2}
\end{align}

Note that we have the following equality.
$$
\mathring{\alpha}_t=
\big(n \mathrm{I}_n+1_n1_n^\top\big)^{-1}\mathring{\Theta}_t1_n
=\alpha_t^\star+\frac{1}{n}\big(\mathrm{I}_n-\frac{1}{2n}1_n1_n^\top\big)\big(\mathring{\Theta}_t-\Theta_t^\star\big)1_n,
$$
and therefore
\begin{align*}
\left\|\mathring{\alpha}_t-\alpha^\star_t\right\|_{\mathrm{2}}
&=\frac{1}{n}\left\|\big(\mathrm{I}_n-\frac{1}{2n}1_n1_n^\top\big)\big(\mathring{\Theta}_t-\Theta_t^\star\big)1_n\right\|_{\mathrm{2}}\\
&\leqslant\frac{1}{n}\left\|\mathrm{I}_n-\frac{1}{2n}1_n1_n^\top\right\|_{\mathrm{op}}\big\|\big(\mathring{\Theta}_t-\Theta_t^\star\big)1_n\big\|_{\mathrm{2}}
\leqslant\frac{1}{\sqrt{n}}\big\|\mathring{\Theta}_t-\Theta_t^\star\big\|_{\mathrm{F}},\\
\|\mathring{\alpha}-\alpha^\star\|_{\mathrm{F}}^2
&=\sum_{t=1}^T\|\mathring{\alpha}_t-\alpha^\star_t\|_{\mathrm{2}}^2
\leqslant\frac{1}{n}\sum_{t=1}^T\big\|\mathring{\Theta}_t-\Theta_t^\star\big\|_{\mathrm{F}}^2.
\end{align*}
Because $G^\star\in \mathbb{S}_+^n,$ then we obtain
\begin{align}
    \big\|\mathring{G}-G^\star\big\|_{\mathrm{F}}
    &\leqslant\big\|\mathring{G}-J\widetilde{G}J\big\|_{\mathrm{F}}+\big\|{G^\star}-J\widetilde{G}J\big\|_{\mathrm{F}} \notag\\
    &\leqslant2\big\|{G^\star}-J\widetilde{G}J\big\|_{\mathrm{F}}\notag\\
    &\leqslant\frac{2}{T}\sum_{t=1}^T\big\|{G^\star}-J\mathring{\Theta}_t J\big\|_{\mathrm{F}}\notag\\
    &\leqslant\frac{2}{T}\sum_{t=1}^T\big\|{\Theta_t^\star}-\mathring{\Theta}_t \big\|_{\mathrm{F}}. \label{eq:gringbd1}
\end{align}
We know that $\operatorname{rank}(G^\star)=k$, therefore
\begin{align*}
    \big\|\mathring Z \mathring Z^\top-G^\star\big\|_{\mathrm{F}}
    \leqslant\big\|\mathring Z \mathring Z^\top-\mathring{G}\big\|_{\mathrm{F}}+\big\|\mathring{G}-G^\star\big\|_{\mathrm{F}}
    \leqslant2\big\|\mathring{G}-G^\star\big\|_{\mathrm{F}}
    \leqslant\frac{4}{T}\sum_{t=1}^T\big\|\Theta_t^\star-\mathring{\Theta}_t \big\|_{\mathrm{F}},
\end{align*}
then $\|\mathring Z \mathring Z^\top-G^\star\|_{\mathrm{F}}^2\leqslant{16}\sum_{t=1}^T\|\Theta_t^\star-\mathring{\Theta}_t \|_{\mathrm{F}}^2/{T}$.

Combining the above bounds of $\mathring \alpha$ and $\mathring Z$, we have
\begin{align*}
    \sum_{t = 1}^T \big\|\Theta_{t}^0 - \Theta_{t}^{\star}\big\|_{\mathrm{F}}^2
    \leqslant2T\big\|\mathring Z \mathring Z^\top-G^\star\big\|_{\mathrm{F}}^2
    +4n\left\|{\mathring \alpha}-\alpha^{\star}\right\|_{\mathrm{F}}^{2}
    \leqslant 36 \sum_{t=1}^T\big\|\Theta_t^\star-\mathring{\Theta}_t \big\|_{\mathrm{F}}^2.
\end{align*}
Recall equation \eqref{th2}, we obtain that 
\begin{align*}
\sum_{t = 1}^T \big\|\Theta_{t}^0 - \Theta_{t}^{\star}\big\|_{\mathrm{F}}^2 \leqslant C {e^{2M_{\Theta,2}}}  {M_{\Theta,2}Tn^{\frac{3}{2}-\frac{1}{k+3}}} \zeta_{\infty}.
\end{align*}


\subsubsection{Proof of Lemma \ref{thm:linear}} \label{sec:pf_linear1}

For the convenience of analysis, we will instead analyze the following quantity,
$$
{\epsilon}^{r} =  2n\left\|\Delta_{\alpha^{r}}\right\|_{\mathrm{F}}^{2} + T\left\|Z^{0}\right\|_{\mathrm{op}}^{2}\left\|\Delta_{Z^{r}}\right\|_{\mathrm{F}}^{2},
$$
where we denote $ \Delta_{\alpha^r} = \alpha^r - \alpha^{\star}$ and $\Delta_{Z^{r}}=Z^{r}-Z^{\star} Q_{r}$ for simplicity.
Note that the only difference between $\epsilon^r$ and $e^r$ is that the $Z^{\star}$ in $e^r$ has been replaced by $Z^{0}$. Intuitively, $\epsilon^r$ and $e^r$ are close since we obtain appropriate $Z^0$ through Stage 1 in Algorithm \ref{algor:init}.
Combining the conclusion of Lemma \ref{lem:in of in} and the condition $$\left\|Z^{\star}\right\|_{\mathrm{op}}^{4} \geqslant C_{1}  e^{4M_{\Theta,2}} M_{\Theta,2}\kappa^{4}n^{1.5- 1/(k+3) }\zeta_{\infty}, $$ we obtain another deterministic error bound for the ``initial of initial'' estimator
\begin{align*}
    \sum_{t = 1}^T \big\|\Theta_{t}^0 - \Theta_{t}^{\star}\big\|_{\mathrm{F}}^2 \leqslant c e^{-M_{\Theta,2}} T \left\|Z^{\star}\right\|_{\mathrm{op}}^{4}/\kappa^4,
\end{align*}
where $c$ is a sufficiently small constant. Denote $ \Delta_{G^r} = G^r - G^{\star}$ and $\Delta_{\Theta_t^r} = \Theta_t^r - \Theta_t^{\star}$.
Note that $e^0$ has a similar bound,
\begin{align} \label{lem2-1}
e^0 \leqslant 2n \left\|\Delta_{\alpha^0}\right\|_{\mathrm{F}}^2 + \kappa^2 T\left\|\Delta_{G^0} \right\|_{\mathrm{F}}^2 \leqslant \kappa^2 \sum_{t = 1}^T \big\|\Delta_{\Theta_{t}^0} \big\|_{\mathrm{F}}^2 \leqslant c e^{- 2M_{\Theta,2}}T\left\|Z^{\star}\right\|_{\mathrm{op}}^{4} / \kappa^{2},
\end{align}
where the first inequality follows from Lemma \ref{Lemma-28} and the second inequality follows from \eqref{l1-a1}.
The above inequality \eqref{lem2-1} shows that
$$
\left|\|Z^0\|_{\mathrm{op}} - \|Z^{\star}\|_{\mathrm{op}}\right| \leqslant
\left\|\Delta_{Z^{0}}\right\|_{\mathrm{op}} \leqslant \left\|\Delta_{Z^{0}}\right\|_{\mathrm{F}} \leqslant \delta^{\prime}\left\|Z^{\star}\right\|_{\mathrm{op}}
$$
for some sufficiently small positive constant $\delta^{\prime}$. As a result,
$$
\left|{\epsilon}^{r} - e^r \right|\leqslant T \left|\|Z^0\|_{\mathrm{op}}^2 - \|Z^{\star}\|_{\mathrm{op}}^2\right| \left\|\Delta_{Z^{r}}\right\|_{\mathrm{F}}^2 \leqslant \delta T \|Z^{\star}\|_{\mathrm{op}}^2 \left\|\Delta_{Z^{r}}\right\|_{\mathrm{F}}^2 \leqslant
\delta e^r,
$$
where $\delta = \delta^{\prime} (2+\delta^\prime)$ is a sufficiently small positive constant. We would consider $\delta = 0.1$ in the following analysis. 

To establish Lemma \ref{thm:linear}, it suffices to prove 
\begin{align} \label{eq:epilinear}
   \epsilon^{r} \leqslant \left(1-\frac{\eta }{e^{M_{\Theta,2}} \kappa^{2} } \rho\right)^{r} \epsilon^{0}+\frac{C e^{2M_{\Theta,2}}\kappa^{2}k T\zeta_{\infty}^2}{\rho }.
\end{align}
The proof of \eqref{eq:epilinear} relies on the following two lemmas.

\begin{lemma} \label{lemma-1}
	Assume the conditions in Lemma \ref{thm:linear}. If $\left\|\Delta_{Z^{r}}\right\|_{\mathrm{F}}^2 \leqslant c_{0} e^{-2M_{\Theta,2}}\left\|Z^{\star}\right\|_{\mathrm{op}}^2 / \kappa^{2}$ for any fixed $0 \leqslant r <R_1$ and a sufficiently small constant $c_0$, then  for any sufficiently small $\eta > 0$, there
	exist constants $\rho$ and $C$,
	$$
	\epsilon^{r+1} \leqslant\left(1-\frac{\eta }{e^{M_{\Theta,2}} \kappa^{2}  } \rho\right) \epsilon^r+  {\eta Ce^{M_{\Theta,2}}k T \zeta_{\infty}^2}.
	$$
\end{lemma}
\begin{proof}
See Section \ref{sec:pf_linear11} on Page \pageref{sec:pf_linear11}.
\end{proof}

\begin{lemma} \label{lemma-2}
Assume the conditions in Lemma \ref{thm:linear}.
    If $\epsilon^0 \leqslant c_{0}e^{-2M_{\Theta,2}}T\left\|Z^{\star}\right\|_{\mathrm{op}}^{4} / (4 \kappa^{2})$ for a sufficiently small constant $c_0$ , then for all $r = 0, 1, \ldots, R_1$, 
	$$
	\left\|\Delta_{Z^{r}}\right\|_{\mathrm{F}}^2 \leqslant \frac{c_{0}}{e^{2M_{\Theta,2}} \kappa^{2}}\left\|Z^{\star}\right\|_{\mathrm{op}}^2.
	$$
\end{lemma}
\begin{proof}
See Section \ref{sec:pf_linear12} on Page \pageref{sec:pf_linear12}.
\end{proof}

The proof of \eqref{eq:epilinear} is straightforward with these two lemmas. Firstly, \eqref{lem2-1} indicates the additional condition of Lemma \ref{lemma-2} is satisfied. Lemma \ref{lemma-2} shows that for all $0 \leqslant r  < R_1$, we have
$$
\left\|\Delta_{Z^{r}}\right\|_{\mathrm{F}}^2 \leqslant \frac{c_{0}}{e^{2M_{\Theta,2}} \kappa^{2}}\left\|Z^{\star}\right\|_{\mathrm{op}}^2, 
$$
where $c_0$ is a sufficiently small constant depend on $c$ in \eqref{lem2-1}. Then by Lemma \ref{lemma-1}, there exists positive constants $\rho$ and $C$ such that for all $r < R_1$,
$$
\epsilon^{r+1} \leqslant\left(1-\frac{\eta}{e^{M_{\Theta,2}} \kappa^{2}} \rho\right) \epsilon^{r}+\eta Ce^{M_{\Theta,2}}k T \zeta_{\infty}^2.
$$
As a result, for any $r \leqslant R_1$
$$
\begin{aligned}
\epsilon^{r} &\leqslant \left(1-\frac{\eta }{e^{M_{\Theta,2}} \kappa^{2}  } \rho\right)^{r} \epsilon^0+\sum_{i=1}^{r}{\eta Ce^{M_{\Theta,2}}k T\zeta_{\infty}^2}\left(1-\frac{\eta }{e^{M_{\Theta,2}} \kappa^{2}} \rho\right)^{i-1}\\
&\leqslant \left(1-\frac{\eta }{e^{M_{\Theta,2}} \kappa^{2} } \rho\right)^{r} \epsilon^{0}+\frac{C e^{2M_{\Theta,2}}\kappa^{2}k T\zeta_{\infty}^2}{\rho }.
\end{aligned}
$$

Finally, as $0.9 e^r \leqslant \epsilon^r \leqslant 1.1 e^r$, by \eqref{eq:epilinear} we have 
$$
e^r \leqslant 2\left(1-\frac{\eta }{e^{M_{\Theta,2}} \kappa^{2}  } \rho\right)^{r} e^{0}+\frac{2C e^{2M_{\Theta,2}}\kappa^{2}k  T \zeta_{\infty}^2}{\rho }.
$$
\subsubsection{Proof of Lemma \ref{lemma-1}}
\label{sec:pf_linear11}

For $t = 1, \ldots, T$, define 
$
h_{t}(\Theta)= \sum_{i, j=1}^n  (\exp (\Theta_{i j})- A_{t, i j} \Theta_{i j}),
$ 
and its gradient is $\nabla h_{t}(\Theta)= \exp (\Theta) -\mathbf A_{t} \in \mathbb{R}^{n\times n}$. By \eqref{eq:grandzformulaalo} and \eqref{eq:grandalphaformulaalo}, 
we can rewrite  
$$ \nabla_Z L(Z^r,\alpha^r)= - \sum_{t=1}^{T} \nabla h_{t}(\Theta^r_t) Z^r \quad
\text{and} \quad \nabla_{\alpha_t} {L}(Z^r,\alpha^r)=- \nabla h_{t}\left(\Theta_{t}^r\right) 1_{n}, \ t=1, \ldots, T.$$ 

Given $Z^r$ and $\widetilde{Z}^r$ in the $r$-th iteration of Algorithm \ref{algor:init}, 
let $Q_{r}=\arg \min_{Q \in \mathcal{O}(k)}\left\|Z^{r}-Z^{\star} Q\right\|_{\mathrm{F}}$ and 
 $\widetilde{Q}_{r}=\arg \min _{Q \in \mathcal{O}(k)}\|\widetilde{Z}^{r}-Z^{\star} Q\|_{\mathrm{F}}$. 
By the definition of $Q_{r+1}$, 
\begin{align*}
   \left\|Z^{r+1}-Z^{\star} Q_{r+1}\right\|_{\mathrm{F}}^2
\leqslant &~ \left\|Z^{r+1}-Z^{\star} \widetilde{Q}_{r+1}\right\|_{\mathrm{F}}^2\notag\\
\leqslant  &~\left\|\widetilde{Z}^{r+1}-Z^{\star} \widetilde{Q}_{r+1}\right\|_{\mathrm{F}}^2
\leqslant \left\|\widetilde{Z}^{r+1}-Z^{\star}{Q}_{r}\right\|_{\mathrm{F}}^2,
\end{align*}
where the second inequality follows by the projection step in line 23 of Algorithm \ref{algor:init},
and the third inequality follows by the definition of $\widetilde{Q}_{r+1}$.
Recall that $\Delta_{Z^{r+1}} =Z^{r+1}-Z^{\star}Q_{r+1}$. 
By	$\widetilde{Z}^{r+1}=Z^{r} - \eta_{Z} \sum_{t=1}^T \nabla h_{t}\left(\Theta_{t}^r\right) Z^r$ in line 21 of Algorthm \ref{algor:init}, 
\begin{align}\label{l1-1}
    \left\| \Delta_{Z^{r+1}} \right\|_{\mathrm{F}}^2
    \leqslant \left\|{Z}^{r}-Z^{\star}{Q}_{r}\right\|_{\mathrm{F}}^2
    +\eta_Z^2 \left\|\sum_{t=1}^T \nabla h_{t}\left(\Theta_{t}^r\right)Z^r \right\|_{\mathrm{F}}^2
    -2\eta_Z\Big\langle \sum_{t=1}^T \nabla h_{t}\left(\Theta_{t}^r\right)Z^r,{Z}^{r}-Z^{\star}{Q}_{r}\Big\rangle.
\end{align}
For $t=1, \ldots, T$, 
\begin{align}\label{l1-2}
 \left\langle \nabla h_{t}\left(\Theta_{t}^r\right)Z^r,{Z}^{r}-Z^{\star}Q_r\right\rangle\notag
&=\Big\langle \nabla h_{t}\left(\Theta_{t}^r\right),({Z}^{r}-Z^{\star}Q_r)(Z^r)^\top\Big\rangle \notag\\
&=\Big\langle \nabla h_{t}\left(\Theta_{t}^r\right),\frac{1}{2}({Z}^{r}({Z}^{r})^\top-Z^{\star}(Z^{\star})^\top)\Big\rangle \notag\\
&\quad+\Big\langle \nabla h_{t}\left(\Theta_{t}^r\right),\frac{1}{2}({Z}^{r}({Z}^{r})^\top+Z^{\star}(Z^{\star})^\top)-Z^{\star}Q_r(Z^r)^\top\Big\rangle \notag\\
&=\frac{1}{2}\left\langle \nabla h_{t}\left(\Theta_{t}^r\right),\Delta_{G^r}\right\rangle + \frac{1}{2}\big\langle \nabla h_{t}\left(\Theta_{t}^r\right),\Delta_{Z^r}\Delta_{Z^r}^\top\big\rangle.
\end{align}
Therefore, by \eqref{l1-1} and \eqref{l1-2},
\begin{align}\label{l1-3}
 \left\| \Delta_{Z^{r+1}} \right\|_{\mathrm{F}}^2
    &\leqslant \left\|\Delta_{Z^r}\right\|_{\mathrm{F}}^2
    +\eta_Z^2 \left\|\sum_{t=1}^{T} \nabla h_{t}\left(\Theta_{t}^r\right)Z^r \right\|_{\mathrm{F}}^2 \notag\\
    &\quad-\eta_Z\sum_{t=1}^{T}\left\langle \nabla h_{t}\left(\Theta_{t}^r\right),\Delta_{G^r}\right\rangle 
     -\eta_Z\sum_{t=1}^{T}\big\langle \nabla h_{t}\left(\Theta_{t}^r\right),\Delta_{Z^r}\Delta_{Z^r}^\top\big\rangle.
\end{align}
Similarly to \eqref{l1-3}, 
we also have
\begin{align} \label{l1-3.5}
 &\left\|\alpha_t^{r+1}-\alpha_t^{\star} \right\|_{\mathrm{2}}^2
    \leqslant \left\|\widetilde{\alpha}_t^{r+1}-\alpha_t^{\star} \right\|_{\mathrm{2}}^2 \notag\\
    =&\left\|\alpha_t^{r}-\alpha_t^{\star} \right\|_{\mathrm{2}}^2
    +\eta_\alpha^2 \left\|\nabla h_{t}\left(\Theta_{t}^r\right)1_n \right\|_{\mathrm{2}}^2 
   -2\eta_\alpha\left\langle \nabla h_{t}\left(\Theta_{t}^r\right),(\alpha_t^r-\alpha_t^{\star})1_n^\top \right\rangle.
\end{align}
Taking the summation over $t$ in \eqref{l1-3.5}, we get
\begin{align}\label{l1-4}
 \left\| \Delta_{\alpha^{r+1}} \right\|_{\mathrm{F}}^2
    \leqslant \left\|\Delta_{\alpha^r}\right\|_{\mathrm{F}}^2
    +\eta_\alpha^2 \sum_{t=1}^{T} \left\|\nabla h_{t}\left(\Theta_{t}^r\right)1_n \right\|_{\mathrm{F}}^2
    -2\eta_\alpha\sum_{t=1}^{T}\left\langle \nabla h_{t}\left(\Theta_{t}^r\right),\Delta_{\alpha_t^r}1_n^\top \right\rangle,
\end{align}
where $ \Delta_{\alpha_t^r} = \alpha_t^r - \alpha_t^{\star}$.

Recall that $\epsilon^{r+1}=T\|Z^0\|_{\operatorname{op}}^2\|\Delta_{Z^{r+1}}\|_{\mathrm{F}}^2 +2n\|\Delta_{\alpha^{r+1}}\|_{\mathrm{F}}^2$ and combine with \eqref{l1-3} and \eqref{l1-4}
\begin{align}
\epsilon^{r+1}
&\leqslant\epsilon^{r}-\eta\sum_{t=1}^T\left\langle \nabla h_{t}\left(\Theta_{t}^r\right),\Delta_{G^r}
+\Delta_{Z^r}\Delta_{Z^r}^\top+2\Delta_{\alpha_t^r}1_n^\top\right\rangle \notag\\
&\quad\ +\frac{\eta^2}{T\left\| Z^0 \right\|_{\mathrm{op}}^2}\left\|\sum_{t=1}^T\nabla h_{t}\left(\Theta_{t}^r\right)Z^r \right\|_{\mathrm{F}}^2
+\frac{\eta^2}{2n}\sum_{t=1}^T \left\|\nabla h_{t}\left(\Theta_{t}^r\right)1_n \right\|_{\mathrm{2}}^2.\label{eq:l1-9-1}
\end{align}
For $t=1,\ldots, T$, 
define $H_t(\Theta) = \EXPT[h_t(\Theta)]$, where the expectation is taken given the true parameters $(Z^{\star},\alpha^{\star})$.
Then $H_t(\Theta)=\sum_{i,j=1}^n\{\exp(\Theta_{ij})-\exp(\Theta_{t,ij}^{\star})\Theta_{ij}\}$,
and its gradient is $\nabla H_t(\Theta) = \exp(\Theta) - \exp(\Theta^{\star}_t) \in \mathbb{R}^{n\times n}$.
In  \eqref{eq:l1-9-1},
by $$\nabla h_{t}(\Theta_{t}^r)=\nabla H_{t}(\Theta_{t}^r) + \{\nabla h_{t}(\Theta_{t}^r)-\nabla H_{t}(\Theta_{t}^r)\} \quad\text{and}\quad \langle A,\Delta_{\Theta_t^r}\rangle =\langle A,\Delta_{G^r}+2\Delta_{\alpha_t^r}1_n^{\top}\rangle$$ for any matrix $A\in \mathbb{R}^{n\times n}$, we obtain
\begin{align} \label{l1-9}
    \epsilon^{r+1}
\leqslant\epsilon^{r} - \eta D_1 + \eta D_2+\eta D_3+\eta^2D_4 
\end{align}
where
\begin{align*}
    D_1 = &~ \sum_{t=1}^T\left\langle \nabla H_{t}\left(\Theta_{t}^r\right),\Delta_{\Theta_t^r}\right\rangle ,\\
    D_2 = &~ \sum_{t=1}^T\left|\left\langle \nabla h_{t}\left(\Theta_{t}^r\right)-\nabla H_{t}\left(\Theta_{t}^r\right),\Delta_{\Theta_t^r}\right\rangle\right| ,\\
    D_3 = &~ \sum_{t=1}^T\left|\left\langle \nabla h_{t}\left(\Theta_{t}^r\right),\Delta_{Z^r}\Delta_{Z^r}^\top\right\rangle\right| ,\\
    D_4 = &~ \frac{1}{T\left\| Z^0 \right\|_{\mathrm{op}}^2}\left\|\sum_{t=1}^T\nabla h_{t}\left(\Theta_{t}^r\right)Z^r \right\|_{\mathrm{F}}^2 +\frac{1}{2n}\sum_{t=1}^T \left\|\nabla h_{t}\left(\Theta_{t}^r\right)1_n \right\|_{\mathrm{2}}^2.
\end{align*}

By the constraint set $\mathcal{S}_C$ of $(Z,\alpha)$ in Section \ref{sec:initest}, 
it suffices to consider $\Theta$ in the following set 
$$
\mathcal{F}_{\Theta}=\{\Theta = \alpha1_n^\top+1_n\alpha^\top+ZZ^\top:JZ=Z,\max_{1 \leqslant i \leqslant n}\left\|z_i\right\|_2^2 \leqslant {M_{Z,1}},\max_{1\leqslant i\leqslant n, 1\leqslant t \leqslant T} |\alpha_{it}|\leqslant M_{\alpha}\}. $$ 
Note that for any ${\Theta\in\mathcal{F}_{\Theta}}$ we have
$$
e^{-M_{\Theta,2}} \mathrm{I}_{n^2\times n^2}\preceq \nabla^2 H_{t}(\Theta)=\text{diag}( \text{vec}(\exp{(\Theta)}))\preceq  e^{M_{\Theta,2}} \mathrm{I}_{n^2\times n^2},$$
where  for two matrices $A$ and $B$ of the same size, $A\preceq B$ represents that $B-A$ is  positive semidefinite. 
Hence $H_{t}(\cdot)$ is $e^{-M_{\Theta,2}}$-strongly convex and $e^{M_{\Theta,2}}$-smooth in the $\mathcal{F}_{\Theta}$. Further
notice that $\nabla H_{t}(\Theta_{t}^{\star})=0$ and $\Theta_t^r \in \mathcal{F}_{\Theta}$, then by Lemma \ref{Lemma-30}

\begin{align*}
\left\langle \nabla H_{t}\left(\Theta_{t}^r\right),\Delta_{\Theta_t^r}\right\rangle
&\geqslant \frac{1}{e^{-M_{\Theta,2}}+e^{M_{\Theta,2}}}\left\|\Delta_{\Theta_t^r}\right\|_{\mathrm{F}}^2+\frac{1}{e^{-M_{\Theta,2}}+e^{M_{\Theta,2}}}\left\|\nabla H_t(\Theta_t^r)\right\|_{\mathrm{F}}^2.
\end{align*}
Calculating the sum over $t$, we obtain
\begin{align}\label{l1-5}
    D_1 \geqslant \frac{1}{e^{-M_{\Theta,2}}+e^{M_{\Theta,2}}}\sum_{t=1}^T\left\|\Delta_{\Theta_t^r}\right\|_{\mathrm{F}}^2
+\frac{1}{e^{-M_{\Theta,2}}+e^{M_{\Theta,2}}}\sum_{t=1}^T\left\|\nabla H_t(\Theta_t^r)\right\|_{\mathrm{F}}^2.
\end{align}
Note that for $t = 1, \ldots, T$,
\begin{align} \label{l1-a1}
    \left\|\Delta_{\Theta_t^r}\right\|_{\mathrm{F}}^2 &=  \left\|\Delta_{G^r} + \Delta_{\alpha_t^r} 1_n^{\top} + 1_n \Delta_{\alpha_t^r}^{\top}\right\|_{\mathrm{F}}^2 \notag\\
    &= \left\|\Delta_{G^r}\right\|_{\mathrm{F}}^2 + 2 \left\|\Delta_{\alpha_t^r} 1_n^{\top}\right\|_{\mathrm{F}}^2 + \text{tr}(\Delta_{\alpha_t^r} 1_n^{\top}\Delta_{\alpha_t^r} 1_n^{\top}) \notag\\
    &= \left\|\Delta_{G^r}\right\|_{\mathrm{F}}^2 + 2 \left\|\Delta_{\alpha_t^r} 1_n^{\top}\right\|_{\mathrm{F}}^2 +  (1_n^{\top}\Delta_{\alpha_t^r})^2 \notag\\
    &\geqslant \left\|\Delta_{G^r}\right\|_{\mathrm{F}}^2 + 2 \left\|\Delta_{\alpha_t^r} 1_n^{\top}\right\|_{\mathrm{F}}^2.
\end{align}
Then combining \eqref{l1-5} and \eqref{l1-a1} shows that
\begin{align} \label{l1-a2}
    D_1\geqslant \frac{1}{e^{-M_{\Theta,2}}+e^{M_{\Theta,2}}}\left(T\left\|\Delta_{G^r}\right\|_{\mathrm{F}}^2+2n\left\|\Delta_{\alpha^r}\right\|_{\mathrm{F}}^2\right)+\frac{1}{e^{-M_{\Theta,2}}+e^{M_{\Theta,2}}}\sum_{t=1}^T\left\|\nabla H(\Theta_t^r)\right\|_{\mathrm{F}}^2.
\end{align}

For all $t=1, \ldots, T$ and any ${c_2}>0$,
\begin{align*}
& \quad \  \big|\left\langle \nabla h_{t}\left(\Theta_{t}^r\right)-\nabla H_{t}\left(\Theta_{t}^r\right),\Delta_{\Theta_t^r}\right\rangle\big|
=\left|\left\langle \exp{(\Theta_t^\star)}-\mathbf A_t,\Delta_{\Theta_t^r} \right\rangle\right|\\
&\leqslant  \left|\left\langle \exp{(\Theta_t^\star)}-\mathbf A_t,\Delta_{G^r} \right\rangle\right|
+2\left|\left\langle \exp{(\Theta_t^\star)}-\mathbf A_t,\Delta_{\alpha_t^r}1_n^\top \right\rangle\right|\\
& \leqslant  \zeta_t \left\|\Delta_{G^r}\right\|_{\mathrm{*}}+2\zeta_t\left\|\Delta_{\alpha^r_t}1_n^\top\right\|_{\mathrm{*}}\\
& \leqslant  \zeta_t \sqrt{2k}\left\|\Delta_{G^r}\right\|_{\mathrm{F}}+2\zeta_t\left\|\Delta_{\alpha^r_t}1_n^\top\right\|_{\mathrm{F}}\\
& =  \frac{\zeta_t \sqrt{2k}}{\sqrt{c_2}}\sqrt{c_2}\left\|\Delta_{G^r}\right\|_{\mathrm{F}}+\frac{\sqrt{2}\zeta_t }{\sqrt{c_2}}\sqrt{2c_2}\left\|\Delta_{\alpha^r_t}1_n^\top\right\|_{\mathrm{F}}\\
& \leqslant \sqrt{c_2 \left\|\Delta_{G^r}\right\|_{\mathrm{F}}^2 +2c_2\left\|\Delta_{\alpha^r_t}1_n^\top\right\|_{\mathrm{F}}^2}\cdot\sqrt{\frac{2k\zeta_t^2+2\zeta_t^2}{c_2}}\\
& \leqslant c_2 \left(\left\|\Delta_{G^r}\right\|_{\mathrm{F}}^2 +2\left\|\Delta_{\alpha^r_t}1_n^\top\right\|_{\mathrm{F}}^2\right)+\frac{4}{c_2}\zeta_t^2k.
\end{align*}
Calculate the sum over $t$, then we get
\begin{align}\label{l1-6}
D_2 \leqslant c_2\left(T\left\|\Delta_{G^r}\right\|_{\mathrm{F}}^2 +2n\left\|\Delta_{\alpha^r}\right\|_{\mathrm{F}}^2\right)+\frac{4kT}{c_2} \zeta_{\infty}^2.
\end{align}
By Lemma \ref{Lemma-28}, 
\begin{align*}
    T\left\|\Delta_{G^r}\right\|_{\mathrm{F}}^2 +2n\left\|\Delta_{\alpha^r}\right\|_{\mathrm{F}}^2
    \geqslant \frac{\sigma_k^2(Z^\star)}{2\left\|Z^\star\right\|_{\mathrm{op}}^2} T\left\|Z^\star\right\|_{\mathrm{op}}^2\left\|\Delta_{Z^r}\right\|_{\mathrm{F}}^2+2n\left\|\Delta_{\alpha^r}\right\|_{\mathrm{F}}^2
    \geqslant\frac{1}{2\kappa^2}e^r.
\end{align*}
Therefore, combining the above two equations \eqref{l1-a2} and \eqref{l1-6},
\begin{align}\label{l1-10}
    D_1-D_2&\geqslant\left(\frac{1}{e^{-M_{\Theta,2}}+e^{M_{\Theta,2}}}-c_2\right)\left(T\left\|\Delta_{G^r}\right\|_{\mathrm{F}}^2 +2n\left\|\Delta_{\alpha^r}\right\|_{\mathrm{F}}^2\right)
    \notag \\
    &\quad+\frac{1}{e^{-M_{\Theta,2}}+e^{M_{\Theta,2}}}\sum_{t=1}^T\left\|\nabla H_t{(\Theta_t^r)}\right\|_{\mathrm{F}}^2-\frac{4kT}{c_2}\zeta_{\infty}^2\notag\\
    &\geqslant \frac{1}{2\kappa^2}\left(\frac{1}{e^{-M_{\Theta,2}}+e^{M_{\Theta,2}}}-c_2\right)e^r+\frac{1}{e^{-M_{\Theta,2}}+e^{M_{\Theta,2}}}\sum_{t=1}^T\left\|\nabla H_t{(\Theta_t^r)}\right\|_{\mathrm{F}}^2-\frac{4kT}{c_2}\zeta_{\infty}^2.
\end{align}

For all $t=1, \ldots, T$ and any $c_3>0$, we have
\begin{align*}
    \left|\left\langle \nabla h_{t}\left(\Theta_{t}^r\right),\Delta_{Z^r}\Delta_{Z^r}^\top\right\rangle\right|
    &\leqslant \left\|\nabla h_t{(\Theta_t^r)}\right\|_{\mathrm{op}}
    \cdot \left\|\Delta_{Z^r}\Delta_{Z^r}^\top \right\|_{\mathrm{*}}\notag\\
    &=\left\|\nabla h_t{(\Theta_t^r)}\right\|_{\mathrm{op}}
    \cdot \left\|\Delta_{Z^r} \right\|_{\mathrm{F}}^2\notag\\
    &\leqslant\left\|\nabla h_t{(\Theta_t^r)}-\nabla H_t{(\Theta_t^r)}\right\|_{\mathrm{op}}\cdot \left\|\Delta_{Z^r} \right\|_{\mathrm{F}}^2
    +\left\|\nabla H_t{(\Theta_t^r)}\right\|_{\mathrm{op}}
    \cdot \left\|\Delta_{Z^r} \right\|_{\mathrm{F}}^2\notag\\
    &\leqslant\zeta_t\cdot\left\|\Delta_{Z^r} \right\|_{\mathrm{F}}^2+\left\|\nabla H_t{(\Theta_t^r)}\right\|_{\mathrm{F}}
    \cdot \left\|\Delta_{Z^r} \right\|_{\mathrm{F}}^2\\
    &\leqslant c_3 \zeta_{t}^2 + c_3 \left\|\nabla H_t{(\Theta_t^r)}\right\|_{\mathrm{F}}^2 + \frac{1}{2c_3} \left\|\Delta_{Z^r} \right\|_{\mathrm{F}}^4\\
    &\leqslant  c_3 \zeta_{t}^2 + c_3 \left\|\nabla H_t{(\Theta_t^r)}\right\|_{\mathrm{F}}^2 + \frac{c_0}{2c_3 e^{2M_{\Theta,2}}\kappa^{2}}\left\|\Delta_{Z^r}\right\|_{\mathrm{F}}^2 \left\|{Z^\star}\right\|_{\mathrm{op}}^2.
\end{align*}
The last inequality is due to the assumption $\left\|\Delta_{Z^{r}}\right\|_{\mathrm{F}}^2 \leqslant c_{0} e^{-2M_{\Theta,2}}\left\|Z^{\star}\right\|_{\mathrm{op}}^2 / \kappa^{2}.$
Calculating the sum over $t$ shows that
\begin{align}\label{l1-11}
D_3\leqslant \frac{c_0}{2c_3 e^{2M_{\Theta,2}}\kappa^{2}}e^r+c_3\sum_{t=1}^T \left\|\nabla H_t{(\Theta_t^r)}\right\|_{\mathrm{F}}^2 + c_3 T\zeta_{\infty}^2.
\end{align}

For all $t=1, \ldots, T$, we have
\begin{align*}
\left\|\nabla h_{t}\left(\Theta_{t}^r\right)Z^r \right\|_{\mathrm{F}}^2
&\leqslant 2\left\|(\nabla h_{t}\left(\Theta_{t}^r\right)-\nabla H_{t}\left(\Theta_{t}^r\right))Z^r \right\|_{\mathrm{F}}^2+2\left\|\nabla H_{t}\left(\Theta_{t}^r\right)Z^r \right\|_{\mathrm{F}}^2\\
&\leqslant 2\left\|\exp{(\Theta_{t}^\star)}-\mathbf A_t \right\|_{\mathrm{op}}^2\left\|Z^r \right\|_{\mathrm{F}}^2
+2\left\|\nabla H_{t}\left(\Theta_{t}^r\right) \right\|_{\mathrm{F}}^2\left\|Z^r \right\|_{\mathrm{op}}^2\\
&\leqslant2\zeta_t^2\cdot k\left\|Z^r \right\|_{\mathrm{op}}^2
+2\left\|\nabla H_{t}\left(\Theta_{t}^r\right) \right\|_{\mathrm{F}}^2\left\|Z^r \right\|_{\mathrm{op}}^2, \\
\left\|\nabla h_{t}\left(\Theta_{t}^r\right)1_n \right\|_{\mathrm{2}}^2
&\leqslant 2\left\|(\nabla h_{t}\left(\Theta_{t}^r\right)-\nabla H_{t}\left(\Theta_{t}^r\right))1_n \right\|_{\mathrm{2}}^2+2\left\|\nabla H_{t}\left(\Theta_{t}^r\right)1_n \right\|_{\mathrm{2}}^2\\
&\leqslant2n\zeta_t^2+2n\left\|\nabla H_{t}\left(\Theta_{t}^r\right) \right\|_{\mathrm{F}}^2.
\end{align*}
Combining the above two equations
\begin{align*}
D_4
&\leqslant \frac{1}{\left\| Z^0 \right\|_{\mathrm{op}}^2}\sum_{t=1}^T\left\|\nabla h_{t}\left(\Theta_{t}^r\right)Z^r \right\|_{\mathrm{F}}^2+\frac{1}{2n}\sum_{t=1}^T \left\|\nabla h_{t}\left(\Theta_{t}^r\right)1_n \right\|_{\mathrm{2}}^2\\
&\leqslant \left(\frac{\left\| Z^r \right\|_{\mathrm{op}}^2}{\left\| Z^0 \right\|_{\mathrm{op}}^2}\cdot2k T\zeta_{\infty}^2+T\zeta_{\infty}^2\right)
+\left(\frac{2\left\| Z^r \right\|_{\mathrm{op}}^2}{\left\| Z^0 \right\|_{\mathrm{op}}^2}+1\right)\cdot\sum_{t=1}^T \left\|\nabla H_{t}\left(\Theta_{t}^r\right) \right\|_{\mathrm{F}}^2.
\end{align*}
Note that \eqref{lem2-1} indicates $\left\|\Delta_{Z^{0}}\right\|_{\mathrm{F}}^2 \leqslant ce^{-2M_{\Theta,2}}\left\|Z^{\star}\right\|_{\mathrm{op}}^2 / \kappa^{2}$ for a sufficiently small constant $c$. Combining with the assumption $\left\|\Delta_{Z^{r}}\right\|_{\mathrm{F}}^2 \leqslant c_{0} e^{-2M_{\Theta,2}}\left\|Z^{\star}\right\|_{\mathrm{op}}^2 / \kappa^{2}$, we have $(1-c^{\prime})\left\|{Z^{0}}\right\|_{\mathrm{op}}^2\leqslant\left\|{Z^{r}}\right\|_{\mathrm{op}}^2\leqslant(1+c^{\prime})\left\|{Z^{0}}\right\|_{\mathrm{op}}^2$ for another sufficiently small constant $c^{\prime}$.
Then, there exists a constant $C_4$ such that  
\begin{align}\label{l1-12}
    D_4\leqslant C_4\left(kT\zeta_{\infty}^2+\sum_{t=1}^T \left\|\nabla H_{t}\left(\Theta_{t}^r\right)\right\|_{\mathrm{F}}^2\right).
\end{align}

Finally, we combine \eqref{l1-9}, \eqref{l1-10}, \eqref{l1-11} and \eqref{l1-12}, then obtain
\begin{align*}
    \epsilon^{r+1} & \leqslant   \epsilon^{r}-\eta\left(\frac{1}{2\kappa^2}\left(\frac{1}{e^{-M_{\Theta,2}}+e^{M_{\Theta,2}}}-c_2\right)-\frac{c_0}{2c_3 e^{2M_{\Theta,2}}\kappa^{2}}\right)e^r+\eta\left(\frac{4k}{c_2}+c_3\right)T\zeta_{\infty}^2\\
    & \quad\  -\eta\left(\frac{1}{e^{-M_{\Theta,2}}+e^{M_{\Theta,2}}}-c_3-C_4\eta\right)\cdot\sum_{t=1}^T \left\|\nabla H_{t}\left(\Theta_{t}^r\right) \right\|_{\mathrm{F}}^2
    +\eta^2 C_4kT\zeta_{\infty}^2,
\end{align*}
where $c_2,c_3$ are arbitrary constants, $c_0$ is a sufficiently small constant, and $C_4$ is a sufficiently large constant. Let $c_2=ce^{-M_{\Theta,2}}$, $c_3 = 2c_0 e^{-M_{\Theta,2}}$ and $c,\eta$ are small enough such that
\begin{align*}
    \frac{1}{e^{-M_{\Theta,2}}+e^{M_{\Theta,2}}}-c_2-\frac{c_0}{c_3 e^{2M_{\Theta,2}}}&>2 \widetilde{\rho}\ e^{-M_{\Theta,2}}, \quad\text{and} \\
    \frac{1}{e^{-M_{\Theta,2}}+e^{M_{\Theta,2}}}-c_3-C_4\eta&\geqslant0,
\end{align*}
for some positive constant $\widetilde{\rho}$. Then there exists a universal constant $C>0$ such that
$$
\epsilon^{r+1}\leqslant \epsilon^{r}-\frac{\eta \widetilde{\rho}}{e^{M_{\Theta,2}}\kappa^2}e^r+C\cdot{\eta e^{M_{\Theta,2}}} kT\zeta_{\infty}^2.
$$
Recall that $e^r\geqslant(1-\delta)\epsilon^{r}$, the proof is completed by setting $\rho=(1-\delta)\widetilde{\rho}$
$$
\epsilon^{r+1}\leqslant \left(1-\frac{\eta \rho}{e^{M_{\Theta,2}}\kappa^2}\right)\epsilon^{r}+{\eta C e^{M_{\Theta,2}}kT\zeta_{\infty}^2}.
$$

\subsubsection{Proof of Lemma \ref{lemma-2}}
\label{sec:pf_linear12}
At initialization, we have
\begin{align*}
\|\Delta_{Z^0}\|_{\mathrm{F}}^2
\leqslant\frac{\epsilon^0}{T\|Z^0\|_{\mathrm{op}}^2}
\leqslant\frac{c_0}{4e^{2M_{\Theta,2}}\kappa^2}\cdot\frac{\|Z^\star\|_{\mathrm{op}}^4}{\|Z^0\|_{\mathrm{op}}^2}
\leqslant\frac{c_0}{e^{2M_{\Theta,2}}\kappa^2}\|Z^\star\|_{\mathrm{op}}^2,
\end{align*}
where the last inequality is obtained from
\begin{align*}
2\|Z^0\|_{\mathrm{op}}^2\geqslant\|Z^\star\|_{\mathrm{op}}^2 - 2\|\Delta_{Z^0}\|_{\mathrm{op}}^2
\geqslant\left(1-\frac{c_0}{2e^{2M_{\Theta,2}}\kappa^2}\right)\|Z^\star\|_{\mathrm{op}}^2
\geqslant\frac{2}{3}\|Z^\star\|_{\mathrm{op}}^2,
\end{align*}
where the second and the last inequalities are due to the additional condition of Lemma \ref{lemma-2}.

Suppose the claim is true for all $r\leqslant r_0$, where $r_0$ can be any positive integer smaller than $R_1$. Then by Lemma \ref{lemma-1},
\begin{align*}
\epsilon^{r_0 + 1}
&\leqslant \left(1-\frac{\eta \rho}{e^{M_{\Theta,2}} \kappa^{2} } \right)^{r_0 + 1} \epsilon^{0}+\frac{Ce^{2M_{\Theta,2}} \kappa^{2}k  T \zeta_{\infty}^2}{\rho }\\
&\leqslant \epsilon^{0}+\frac{C e^{2M_{\Theta,2}}\kappa^{2}k  T \zeta_{\infty}^2}{\rho}\\
&\leqslant \frac{c_0T}{4e^{2M_{\Theta,2}}\kappa^2}\|Z^\star\|_{\mathrm{op}}^4+\frac{C e^{2M_{\Theta,2}} \kappa^{2} k  T \zeta_{\infty}^2}{\rho}\\
&\leqslant \frac{c_0T}{e^{2M_{\Theta,2}}\kappa^2}\|Z^\star\|_{\mathrm{op}}^4
\left(\frac{1}{4}+\frac{C e^{4M_{\Theta,2}}\kappa^{4}k   \zeta_{\infty}^2}{c_0 \rho \|Z^\star\|_{\mathrm{op}}^4  }\right)\\
&\leqslant \frac{c_0T}{e^{2M_{\Theta,2}}\kappa^2}\|Z^\star\|_{\mathrm{op}}^4
\left(\frac{1}{4}+\frac{C}{c_0 \rho C_1}\right),
\end{align*}
where the last inequality follows from Condition \ref{cond:zstareigenbd}.

Choosing $C_1$ large enough such that ${C}/({c_0 \rho C_1})\leqslant{1}/{12},$ then
\begin{align*}
    \epsilon^{r_0 + 1}\leqslant \frac{c_0T}{3e^{2M_{\Theta,2}}\kappa^2}\|Z^\star\|_{\mathrm{op}}^4
\end{align*}
and therefore
\begin{align*}
    \|\Delta_{Z^{r_0+1}}\|_{\mathrm{F}}^2\leqslant\frac{\epsilon^{r_0 + 1}}{T\|Z^0\|_{\mathrm{op}}^2}
    \leqslant\frac{c_0}{3e^{2M_{\Theta,2}}\kappa^2}\|Z^\star\|_{\mathrm{op}}^2\frac{\|Z^\star\|_{\mathrm{op}}^2}{\|Z^0\|_{\mathrm{op}}^2}\leqslant\frac{c_0}{e^{2M_{\Theta,2}}\kappa^2}\|Z^\star\|_{\mathrm{op}}^2.
\end{align*}
This completes the proof.

\subsubsection{Proof of Lemma \ref{thm:linear2}} \label{sec:pf_linear2}
For any $t = 1, \ldots, T$ and $r = R_{0.5}, \ldots, R_1$, $2n \times$ \eqref{l1-3.5}  gives that
\begin{align} \label{s2:00}
    e_t^{r+1} &\leqslant e_t^r 
    -2\eta\left\langle \nabla h_{t}\left(\Theta_{t}^r\right),(\alpha_t^r-\alpha_t^{\star})1_n^\top \right\rangle
    +\frac{\eta^2}{2n} \left\|\nabla h_{t}\left(\Theta_{t}^r\right)1_n \right\|_{\mathrm{2}}^2 .
\end{align}
Define $\Theta_t^{r;\star} = \alpha_t^r 1_n^{\top} + 1_n (\alpha_t^r)^{\top} + Z^{\star}(Z^{\star})^{\top}$.
In \eqref{s2:00}, we decompose $\nabla h_{t}(\Theta_{t}^r)=\nabla H_{t}(\Theta_{t}^{r;\star}) + \{\nabla h_{t}(\Theta_{t}^r)-\nabla H_{t}(\Theta_{t}^{r;\star})\} $  and obtain
\begin{align} \label{s2:0}
 e^{r+1}_t \leqslant e_t^r -\eta D_{21} +\eta D_{22}+\eta^2 D_{23}
\end{align}
where
\begin{align*}
    D_{21} = &~ \left\langle \nabla H_{t}\left(\Theta_{t}^{r;\star}\right),2(\alpha_t^r-\alpha_t^{\star})1_n^\top \right\rangle  ,\\
    D_{22} = &~ \left|  \left\langle \nabla h_{t}\left(\Theta_{t}^{r}\right) - \nabla H_{t}\left(\Theta_{t}^{r;\star}\right),2(\alpha_t^r-\alpha_t^{\star})1_n^\top \right\rangle\right| ,\\
    D_{23} = &~ \frac{1}{2n} \left\|\nabla h_{t}\left(\Theta_{t}^r\right)1_n \right\|_{\mathrm{2}}^2.
\end{align*}

We will use similar techniques as in Lemma \ref{lemma-1} to bound these terms. Firstly by Lemma \ref{Lemma-30},
\begin{align} \label{s2:3}
    D_{21} &
    = \left\langle \nabla H_{t}\left(\Theta_{t}^{r;\star}\right) - \nabla H_{t}\left(\Theta_{t}^{\star}\right), \Theta_{t}^{r;\star} - \Theta_{t}^{\star}\right\rangle \notag\\
    &\geqslant \frac{1}{e^{-M_{\Theta,2}}+e^{M_{\Theta,2}}}\left\|\Theta_{t}^{r;\star} - \Theta_{t}^{\star}\right\|_{\mathrm{F}}^{2}+\frac{1}{e^{-M_{\Theta,2}}+e^{M_{\Theta,2}}}\left\|\nabla H_{t}\left(\Theta_{t}^{r;\star}\right)\right\|_{\mathrm{F}}^{2} \notag\\
    &\geqslant \frac{1}{e^{-M_{\Theta,2}}+e^{M_{\Theta,2}}} \cdot  e_t^r+\frac{1}{e^{-M_{\Theta,2}}+e^{M_{\Theta,2}}}\left\|\nabla H_{t}\left(\Theta_{t}^{r;\star}\right)\right\|_{\mathrm{F}}^{2}.
\end{align}
Then for any $c_{22} > 0$,
\begin{align} \label{s2:1}
    D_{22} &\leqslant 2 \left\|\nabla h_{t}\left(\Theta_{t}^{r}\right) - \nabla H_{t}\left(\Theta_{t}^{r;\star}\right)\right\|_{\mathrm{op}} \left\|\Delta_{\alpha_{t}^{r}} 1_n^\top\right\|_{*} \notag\\
    &\leqslant 2\zeta_t \left\|\Delta_{\alpha_{t}^{r}} 1_{n}^{\top}\right\|_{*} + 2\left\|\nabla H_{t}\left(\Theta_{t}^{r}\right) - \nabla H_{t}\left(\Theta_{t}^{r;\star}\right)\right\|_{\mathrm{op}} \left\|\Delta_{\alpha_{t}^{r}} 1_{n}^{\top}\right\|_{*} \notag\\
    &\leqslant 2\zeta_t \left\|\Delta_{\alpha_{t}^{r}} 1_{n}^{\top}\right\|_{\mathrm{F}} + 2\left\|\nabla H_{t}\left(\Theta_{t}^{r}\right) - \nabla H_{t}\left(\Theta_{t}^{r;\star}\right)\right\|_{\mathrm{F}} \left\|\Delta_{\alpha_{t}^{r}} 1_{n}^{\top}\right\|_{\mathrm{F}} \notag\\
    &\leqslant c_{22} \left(\zeta_t^2 + \left\|\nabla H_{t}\left(\Theta_{t}^{r}\right) - \nabla H_{t}\left(\Theta_{t}^{r;\star}\right)\right\|_{\mathrm{F}}^2\right) + \frac{1}{c_{22}} \left\|\Delta_{\alpha_{t}^{r}} 1_{n}^{\top}\right\|_{\mathrm{F}}^2.
\end{align}
Because $\nabla H_t(\cdot)$ is $e^{M_{\Theta,2}}$-Lip in $[-M_{\Theta,2}, M_{\Theta,2}]$, then we can obtain that
\begin{align} \label{s2:2}
    \left\|\nabla H_{t}\left(\Theta_{t}^{r}\right) - \nabla H_{t}\left(\Theta_{t}^{r;\star}\right)\right\|_{\mathrm{F}}^2 \leqslant e^{2M_{\Theta,2}} \left\|\Delta_{G^r}\right\|_{\mathrm{F}}^2 \leqslant 9e^{2M_{\Theta,2}} \left\|Z^{\star}\right\|_{\mathrm{op}}^2 \left\| \Delta_{Z^r}\right\|_{\mathrm{F}}^2
    \leqslant Ce^{4M_{\Theta,2}} \kappa^2 k \zeta_{\infty}^2,
\end{align}
where the second inequality follows from Lemma \ref{Lemma-29} and the last inequality follows from Lemma \ref{thm:linear}. Combining \eqref{s2:1} and \eqref{s2:2}, we have
\begin{align} \label{s2:4}
    D_{22} \leqslant c_{22} \left(\zeta_t^2 + Ce^{4M_{\Theta,2}}\kappa^2k\zeta_{\infty}^2\right) + \frac{1}{2c_{22}} \cdot e_t^r.
\end{align}
Finally,
\begin{align} \label{s2:5}
D_{23} &
\leqslant \frac{1}{2} \left\|\nabla h_{t}\left(\Theta_{t}^{r}\right) \right\|_{\mathrm{op}}^{2} \notag\\
& \leqslant 2\left\|\nabla h_{t}\left(\Theta_{t}^{r}\right)-\nabla H_{t}\left(\Theta_{t}^{r}\right)\right\|_{\mathrm{op}}^{2}
+ 2\left\|\nabla H_{t}\left(\Theta_{t}^{r}\right)-\nabla H_{t}\left(\Theta_{t}^{r;\star}\right) \right\|_{\mathrm{F}}^{2} + 
\left\|\nabla H_{t}\left(\Theta_{t}^{r;\star}\right) \right\|_{\mathrm{F}}^{2} \notag\\
& \leqslant 2 \left(\zeta_{t}^{2}+ Ce^{4M_{\Theta,2}}\kappa^2k \zeta_{\infty}^2\right) +  
 \left\|\nabla H_{t}\left(\Theta_{t}^{r;\star}\right)\right\|_{\mathrm{F}}^{2}.
\end{align}

Then we combine \eqref{s2:0}, \eqref{s2:3}, \eqref{s2:4} and \eqref{s2:5}, and obtain
\begin{align*}
    e_{t}^{r+1} &\leqslant e_t^r - \eta \left( \frac{1}{e^{-M_{\Theta,2}} + e^{M_{\Theta,2}}} - \frac{1}{2c_{22}} \right)e_t^r  + \eta (c_{22} + 2\eta)  \left(\zeta_t^2 + Ce^{4M_{\Theta,2}}\kappa^2k\zeta_{\infty}^2\right) \\
    & \quad \ - \eta \left( \frac{1}{e^{-M_{\Theta,2}}+e^{M_{\Theta,2}}} - \eta \right)  \left\|\nabla H_{t}\left(\Theta_{t}^{r;\star}\right)\right\|_{\mathrm{F}}^{2},
\end{align*}
where $c_{22}$ is an arbitrary constant. Let $c_{22} = 2e^{M_{\Theta,2}}$ and $\eta$ is small enough, we obtain 
\begin{align*}
    e_t^{r+1} \leqslant \left( 1- \frac{\eta  \rho}{e^{M_{\Theta,2}}}\right) e_t^r + \eta C e^{M_{\Theta,2}} \left( \zeta_t^2 + e^{4M_{\Theta,2}} \kappa^2 k \zeta_{\infty}^2\right)
\end{align*}
for some universal positive constants $\rho, C$ and any $r = R_{0.5}, \ldots, R_1 - 1$. 

As a result, for any $r = R_{0.5}, \ldots, R_{1}$ and $t = 1, \ldots, T$
\begin{align*}
e_t^{r} & \leqslant\left(1-\frac{\eta \rho}{e^{M_{\Theta,2}} } \right)^{r - R_{0.5}} e_t^{R_{0.5}}+\sum_{i=1}^{r - R_{0.5}} \eta C e^{M_{\Theta,2}} \left( \zeta_t^2 + e^{4M_{\Theta,2}} \kappa^2 k \zeta_{\infty}^2\right) \left(1-\frac{\eta \rho}{e^{M_{\Theta,2}} } \right)^{i-1} \\
& \leqslant\left(1-\frac{\eta \rho}{e^{M_{\Theta,2}} } \right)^{r - R_{0.5}} e_t^{R_{0.5}}+\frac{C e^{2 M_{\Theta,2}}}{\rho} \left( \zeta_t^2 + e^{4M_{\Theta,2}} \kappa^2 k \zeta_{\infty}^2\right).
\end{align*}

\subsubsection{Proof of Lemma \ref{thm:linear3}} 
\label{sec:pf_linear3}
We will show the linear convergence of the vector-wise errors after $R_1$ iterations. Note that the rotation between $Z^r$ and $Z$ is fixed at $Q_{R_1}$ in the error terms we defined. 

Similar to in Lemma \ref{thm:linear}, we instead analyze the following errors,
\begin{align*}
    \epsilon_{i\cdot}^r =2 n\left\|\alpha_{i\cdot}^{r}-\alpha_{i \cdot}^{\star}\right\|_{2}^{2}+T\left\|Z^{0}\right\|_{\mathrm{op}}^{2}\left\|z_{i}^{r}-Q_{R_{1}}^{\top} z_{i}^{\star}\right\|_{2}^{2}.
\end{align*}
And we also have $(1-\delta)e_{i\cdot}^r \leqslant \epsilon_{i\cdot}^r \leqslant (1+\delta)e_{i\cdot}^r$ for a sufficiently small positive constant $\delta$. 

Due to $\widetilde{\alpha}_{it}^{r+1} = \alpha_{it}^r - \eta_{\alpha} e_{i}^{\top} \nabla h_t(\Theta_t^{r,R_1})1_n$ for any $i,t$ and $r \geqslant R_1$, we have 
\begin{align} \label{s3:1}
    &\left\|\alpha_{i\cdot}^{r+1}-\alpha_{i\cdot}^{\star}\right\|_{2}^{2} \leqslant \left\|\widetilde{\alpha}_{i .}^{r+1}-\alpha_{i \cdot}^{\star}\right\|_{2}^{2} \notag\\
    =& \left\|\alpha_{i .}^{r}-\alpha_{i \cdot}^{\star}\right\|_{2}^{2} + \eta_{\alpha}^2 \sum_{t=1}^T \Big(e_{i}^{\top} \nabla h_t(\Theta_t^{r,R_1})1_n\Big)^2 - 2\eta_{\alpha} \sum_{t=1}^T \left(e_{i}^{\top} \nabla h_t(\Theta_t^{r,R_1})1_n\right)(\alpha_{it}^r - \alpha_{it}^{\star}).
\end{align}
Similarly, we also obtain
\begin{align} \label{s3:2}
    &\left\|z_i^{r+1} - Q_{R_1}^{\top} z_i^{\star} \right\|_2^2 \leqslant \left\|\widetilde{z}_i^{r+1} - Q_{R_1}^{\top} z_i^{\star} \right\|_2^2 \notag\\
    =& \left\|{z}_i^{r} - Q_{R_1}^{\top} z_i^{\star} \right\|_2^2  + \eta_{Z}^2 \left\| \sum_{t=1}^T e_{i}^{\top} \nabla h_t(\Theta_t^{r,R_1})Z^{R_1} \right\|_2^2 
    - 2 \eta_Z \sum_{t=1}^T e_{i}^{\top} \nabla h_t(\Theta_t^{r,R_1})Z^{R_1} \big({z}_i^{r} - Q_{R_1}^{\top} z_i^{\star} \big).
\end{align}
For simplicity we denote $\Delta_{z_i^r} = z_i^r - Q_{R_1}^{\top} z_i^{\star}$, $\Delta_{\alpha_{i\cdot}^r} = \alpha_{i\cdot}^r - \alpha_{i\cdot}^{\star}$, and $\Delta_{\alpha_{it}^r} = \alpha_{it}^r - \alpha_{it}^{\star}$,
then $2n \times \eqref{s3:1} + T\left\|Z^{0}\right\|_{\mathrm{op}}^{2} \times \eqref{s3:2}$ shows that for any $i = 1, \ldots, n$ and $r \geqslant R_1$
\begin{align} \label{s3:88}
    \epsilon_{i\cdot}^{r+1} &\leqslant \epsilon_{i\cdot}^{r} -2\eta \sum_{t=1}^T\left(e_i^{\top} \nabla h_t(\Theta_t^{r,R_1})\right)\left( \Delta_{\alpha_{it}^r} 1_n + Z^{R_1} \Delta_{z_i^r} \right) \notag\\
    &\quad\ +  \frac{\eta^2 }{2n}  \sum_{t=1}^T \left(e_{i}^{\top} \nabla h_t(\Theta_t^{r,R_1})1_n\right)^2 + \frac{\eta^2}{T\left\|Z^{0}\right\|_{\mathrm{op}}^{2}} \left\| \sum_{t=1}^T e_{i}^{\top} \nabla h_t(\Theta_t^{r,R_1})Z^{R_1} \right\|_2^2 .
\end{align}
In \eqref{s3:88}, we use the decomposition 
$$\nabla h_{t}(\Theta_{t}^{r,R_1})=\nabla \widetilde{H}_{t}(\alpha_t^r,Z^r) + \left\{\nabla h_{t}(\Theta_{t}^{r,R_1})-\nabla H_{t}(\Theta_{t}^{r,R_1})\right\} +\left\{\nabla H_{t}(\Theta_{t}^{r,R_1})-\nabla \widetilde{H}_{t}(\alpha_t^r,Z^r)\right\},$$ 
where we denote
$$\nabla \widetilde{H}_t(\alpha,Z) = \exp \big( \alpha 1_n^{\top} + 1_n (\alpha_t^{\star})^{\top} +Z(Z^{R_1})^{\top}\big) -  \exp\big( \alpha_t^{\star} 1_n^{\top} + 1_n (\alpha_t^{\star})^{\top} +Z^{\star}Q_{R_1}(Z^{R_1})^{\top}\big).$$ 
Then we obtain
\begin{align} \label{s3:8}
    \epsilon_{i\cdot}^{r+1}
\leqslant \epsilon_{i\cdot}^r - \eta D_{31} + \eta D_{32} + \eta D_{33} + \eta^2 D_{34}
\end{align}
where
\begin{align*}
    D_{31} = &~ \sum_{t=1}^T 2\left(e_i^{\top} \nabla \widetilde{H}_t(\alpha_t^r,Z^r)\right)\left( \Delta_{\alpha_{it}^r} 1_n + Z^{R_1} \Delta_{z_i^r} \right),\\
    D_{32} = &~ \sum_{t=1}^T 2\left| e_{i}^{\top} \left(\nabla h_t(\Theta_t^{r,R_1}) - \nabla H_t(\Theta_t^{r,R_1}) \right) \left( \Delta_{\alpha_{it}^r} 1_n + Z^{R_1} \Delta_{z_i^r} \right)\right| ,\\
    D_{33} = &~ \sum_{t=1}^T 2\left| e_{i}^{\top} \Big(\nabla H_t(\Theta_t^{r,R_1}) - \nabla \widetilde{H}_t(\alpha_t^r,Z^r) \Big) \left( \Delta_{\alpha_{it}^r} 1_n + Z^{R_1} \Delta_{z_i^r} \right)\right| ,\\
    D_{34} = &~  \frac{1}{2n} \sum_{t=1}^T \left(e_{i}^{\top} \nabla h_t(\Theta_t^{r,R_1})1_n\right)^2 + \frac{1}{T\left\|Z^{0}\right\|_{\mathrm{op}}^{2}} \left\| \sum_{t=1}^T e_{i}^{\top} \nabla h_t(\Theta_t^{r,R_1})Z^{R_1} \right\|_2^2 .
\end{align*}

For any $t = 1, \ldots, T$,
\begin{align} \label{s3:3}
    &\quad\ \left(e_i^{\top} \nabla \widetilde{H}_t(\alpha_t^r,Z^r)\right)\left( \Delta_{\alpha_{it}^r} 1_n + Z^{R_1} \Delta_{z_i^r} \right)  \notag\\
    &= \sum_{j = 1}^n \left(\exp\big(\alpha_{it}^r + \alpha_{jt}^{\star} + (z_{i}^r)^{\top} z_{j}^{R_1}\big) - \exp\big(\alpha_{it}^{\star} + \alpha_{jt}^{\star} + (z_{i}^{\star})^{\top} Q_{R_1} z_{j}^{R_1}\big)\right) \left( \Delta_{\alpha_{it}^r} + \Delta_{z_i^r}^{\top} z_j^{R_1}\right) \notag \\
    & \geqslant \frac{1}{e^{-M_{\Theta,2}}+e^{M_{\Theta,2}}} \left\|\Delta_{\alpha_{it}^r} 1_n + Z^{R_1} \Delta_{z_i^r}\right\|_2^2 + \frac{1}{e^{-M_{\Theta,2}}+e^{M_{\Theta,2}}} \left\|e_i^{\top} \nabla \widetilde{H}_t(\alpha_t^r,Z^r) \right\|_2^2,
\end{align}
where the last inequality follows from Lemma \ref{Lemma-30}. With the constraint $1_n^{\top} Z^{R_1} = 0$ we have
\begin{align} \label{s3:4}
    \left\|\Delta_{\alpha_{it}^r} 1_n + Z^{R_1} \Delta_{z_i^r}\right\|_2^2 = \left\|\Delta_{\alpha_{it}^r} 1_n\right\|_2^2 + \left\|Z^{R_1} \Delta_{z_i^r} \right\|_2^2 \geqslant n(\Delta_{\alpha_{it}^r})^2 + \left\|Z^{\star}  \right\|_{\mathrm{op}}^2 \left\|\Delta_{z_i^r}\right\|_2^2/(2\kappa^2),
\end{align}
where the last inequality follows from $0.9 \|Z^{\star}\|_{\mathrm{op}} \leqslant \|Z^{R_1}\|_{\mathrm{op}} \leqslant 1.1 \|Z^{\star}\|_{\mathrm{op}}$ (see Lemma \ref{lemma-2}). Combining \eqref{s3:3} and \eqref{s3:4}, and summing over $t$ shows that
\begin{align} \label{s3:9}
    D_{31} &\geqslant \frac{1}{e^{-M_{\Theta,2}}+e^{M_{\Theta,2}}} \left( n\left\|\Delta_{\alpha_{i\cdot}^{r}}\right\|_{2}^{2}+\frac{T}{2\kappa^2}\left\|Z^{\star}\right\|_{\mathrm{op}}^{2}\left\|\Delta_{z_{i}^{r}}\right\|_{2}^{2} + \sum_{t=1}^T \left\|e_i^{\top} \nabla \widetilde{H}_t(\alpha_t^r,Z^r) \right\|_2^2 \right) \notag \\
    &\geqslant \frac{1}{e^{-M_{\Theta,2}}+e^{M_{\Theta,2}}} \left( \frac{1}{2\kappa^2} \cdot \epsilon_{i\cdot}^r + \sum_{t=1}^T \left\|e_i^{\top} \nabla \widetilde{H}_t(\alpha_t^r,Z^r) \right\|_2^2 \right).
\end{align}

For any $c_{32} > 0$,
\begin{align} \label{s3:5}
    &\quad\  \Big| e_{i}^{\top} \left(\nabla h_t(\Theta_t^{r,R_1}) - \nabla H_t(\Theta_t^{r,R_1}) \right) \left( \Delta_{\alpha_{it}^r} 1_n + Z^{R_1} \Delta_{z_i^r} \right)\Big| \notag\\
    &\leqslant \left| e_{i}^{\top} \left(\nabla h_t(\Theta_t^{r,R_1}) - \nabla H_t(\Theta_t^{r,R_1}) \right) 1_n  \right| \left| \Delta_{\alpha_{it}^r}\right| + \left\| e_{i}^{\top} \left(\nabla h_t(\Theta_t^{r,R_1}) - \nabla H_t(\Theta_t^{r,R_1}) \right) Z^{R_1}  \right\|_2 \left\|\Delta_{z_i^r}\right\|_2 \notag\\
    &\leqslant \frac{c_{32}}{2} \left(n \left(\Delta_{\alpha_{it}^r}\right)^2 + \left\|Z^{\star}\right\|_{\mathrm{op}}^{2}\left\|\Delta_{z_{i}^{r}}\right\|_{2}^{2}\right) 
     + \frac{1}{2c_{32}} \Bigg(\frac{\zeta_{it}^{2}}{n} + \frac{\left\| e_{i}^{\top} \left(\nabla h_t(\Theta_t^{r,R_1}) - \nabla H_t(\Theta_t^{r,R_1}) \right) Z^{R_1}  \right\|_2^2}{\left\|Z^{\star}\right\|_{\mathrm{op}}^{2}}\Bigg).
\end{align}
Note that 
\begin{align} \label{s3:6}
    &\quad\ \Big\| e_{i}^{\top} \left(\nabla h_t(\Theta_t^{r,R_1}) - \nabla H_t(\Theta_t^{r,R_1}) \right) Z^{R_1}  \Big\|_2^2 \notag\\
    &\leqslant 2 \left\| e_{i}^{\top} \left(\nabla h_t(\Theta_t^{r,R_1}) - \nabla H_t(\Theta_t^{r,R_1}) \right) Z^{\star}  \right\|_2^2 + 2 \left\| e_{i}^{\top} \left(\nabla h_t(\Theta_t^{r,R_1}) - \nabla H_t(\Theta_t^{r,R_1}) \right) \Delta_{Z^{R_1}}  \right\|_2^2 \notag\\
    &\leqslant 2 \left\| e_{i}^{\top} \left(\nabla h_t(\Theta_t^{r,R_1}) - \nabla H_t(\Theta_t^{r,R_1}) \right) Z^{\star}  \right\|_2^2 + 2 \left\| e_{i}^{\top} \left(\nabla h_t(\Theta_t^{r,R_1}) - \nabla H_t(\Theta_t^{r,R_1}) \right) \right\|_2^2 \left\| \Delta_{Z^{R_1}}  \right\|_{\mathrm{F}}^2 \notag\\
    &\leqslant 2\zeta_{it}^2 + C \zeta_{it}^2 \cdot \frac{e^{2M_{\Theta,2}} \kappa^2 k \zeta_{\infty}^2}{\left\|Z^{\star}\right\|_{\mathrm{op}}^{2}},
\end{align}
where the last inequality follows from Lemma \ref{thm:linear}. Then combining the above two equations \eqref{s3:5} and \eqref{s3:6}, and summing over $t$
\begin{align} \label{s3:10}
    D_{32} &\leqslant c_{32} \left(n \left\|\Delta_{\alpha_{i\cdot}^r}\right\|_2^2 +T \left\|Z^{\star}\right\|_{\mathrm{op}}^{2}\left\|\Delta_{z_{i}^{r}}\right\|_{2}^{2}\right) + \frac{CT}{c_{32}} \left( \frac{\zeta_{i,\infty}^2}{n} + \frac{\zeta_{i,\infty}^2}{\left\|Z^{\star}\right\|_{\mathrm{op}}^{2}} + \frac{e^{2M_{\Theta,2}}\kappa^2k \zeta_{\infty}^2 \zeta_{i,\infty}^2}{\left\|Z^{\star}\right\|_{\mathrm{op}}^{4}} \right) \notag\\
    &\leqslant  2 c_{32} \cdot \epsilon_{i\cdot}^r + \frac{CT}{c_{32}} \left( \frac{\zeta_{i,\infty}^2}{\min\{n,\left\|Z^{\star}\right\|_{\mathrm{op}}^{2}\}} + \frac{e^{2M_{\Theta,2}}\kappa^2k \zeta_{\infty}^2 \zeta_{i,\infty}^2}{\left\|Z^{\star}\right\|_{\mathrm{op}}^{4}} \right).
\end{align}

For any $c_{33} > 0$,
\begin{align} \label{s3:7}
    D_{33} &\leqslant \sum_{t=1}^T 2 \left\| e_{i}^{\top} \left(\nabla H_t(\Theta_t^{r,R_1}) - \nabla \widetilde{H}_t(\alpha_t^r,Z^r) \right)\right\|_2 \left\|  \Delta_{\alpha_{it}^r} 1_n + Z^{R_1} \Delta_{z_i^r}  \right\|_2 \notag\\
    &\leqslant c_{33}  \sum_{t=1}^T \left\|  \Delta_{\alpha_{it}^r} 1_n + Z^{R_1} \Delta_{z_i^r}  \right\|_2^2 + \frac{1}{c_{33}} \sum_{t=1}^T \left\| e_{i}^{\top} \left(\nabla H_t(\Theta_t^{r,R_1}) - \nabla \widetilde{H}_t(\alpha_t^r,Z^r) \right)\right\|_2^2 \notag \\
    &\leqslant c_{33} \left(n \left\|\Delta_{\alpha_{i\cdot}^r}\right\|_2^2 +T \left\|Z^{R_1}\right\|_{\mathrm{op}}^{2}\left\|\Delta_{z_{i}^{r}}\right\|_{2}^{2}\right) + \frac{2e^{2M_{\Theta,2}}}{c_{33}}  \left( \left\|\Delta_{\alpha^{R_1}}\right\|_{\mathrm{F}}^2 + T \left\|z_i^{\star}\right\|_2^2 \left\| \Delta_{Z^{R_1}} \right\|_{\mathrm{F}}^2  \right) \notag\\
    &\leqslant c_{33} \left(n \left\|\Delta_{\alpha_{i\cdot}^r}\right\|_2^2 +T \left\|Z^{R_1}\right\|_{\mathrm{op}}^{2}\left\|\Delta_{z_{i}^{r}}\right\|_{2}^{2}\right) + \frac{CTe^{4M_{\Theta,2} }\kappa^2 k}{c_{33}} \left( \frac{\zeta_{\infty}^2}{n} +  \frac{M_{\Theta,2} \zeta_{\infty}^2}{\left\|Z^{\star}\right\|_{\mathrm{op}}^{2}}\right) \notag\\
    &\leqslant c_{33} \left(n \left\|\Delta_{\alpha_{i\cdot}^r}\right\|_2^2 +T \left\|Z^{R_1}\right\|_{\mathrm{op}}^{2}\left\|\Delta_{z_{i}^{r}}\right\|_{2}^{2}\right) + \frac{CT}{c_{33}} \cdot 
    \frac{e^{4M_{\Theta,2} }M_{\Theta,2}\kappa^2 k \zeta_{\infty}^2}{ \min\{n,\left\|Z^{\star}\right\|_{\mathrm{op}}^{2}\}} ,
\end{align}
where the third inequality follows from $\exp(\cdot)$ is $e^{M_{\Theta,2}}$-Lip in $[-M_{\Theta,2},M_{\Theta,2}]$, and the fourth inequality follows from Lemma \ref{thm:linear}. For all $t=1, \ldots, T$, we have
\begin{align*}
\left(e_{i}^{\top} \nabla h_t(\Theta_t^{r,R_1})1_n\right)^2
&\leqslant 3 \left(e_{i}^{\top} \nabla \widetilde{H}_t(\alpha_t^r,Z^r)1_n\right)^2 
+3\left( e_{i}^{\top} \left(\nabla h_t(\Theta_t^{r,R_1}) - \nabla H_t(\Theta_t^{r,R_1}) \right) 1_n  \right)^2 
\\
&\quad\ + 3\left( e_{i}^{\top} \left(\nabla H_t(\Theta_t^{r,R_1}) - \nabla \widetilde{H}_t(\alpha_t^r,Z^r) \right) 1_n  \right)^2 \\
&\leqslant 3n\left(\left\|  e_{i}^{\top} \nabla \widetilde{H}_t(\alpha_t^r,Z^r) \right\|_2^2 +  \frac{\zeta_{it}^2}{n} +  \left\| e_{i}^{\top} \left(\nabla H_t(\Theta_t^{r,R_1}) - \nabla \widetilde{H}_t(\alpha_t^r,Z^r) \right)\right\|_2^2\right)\\
& \leqslant 3n\left(\left\|  e_{i}^{\top} \nabla \widetilde{H}_t(\alpha_t^r,Z^r) \right\|_2^2 +  \frac{\zeta_{it}^2}{n}
 + \frac{Ce^{4M_{\Theta,2} }M_{\Theta,2}\kappa^2 k \zeta_{\infty}^2}{\min\{n,\left\|Z^{\star}\right\|_{\mathrm{op}}^{2}\}} \right),
\end{align*}
where the last inequality follows from the latter half of \eqref{s3:7}.
Similarly, we have
\begin{align*}
&\quad\ \Big\|e_{i}^{\top} \nabla h_t(\Theta_t^{r,R_1})Z^{R_1}\Big\|_2^2 
\\
&\leqslant 3 \left\|e_{i}^{\top} \nabla \widetilde{H}_t(\alpha_t^r,Z^r)Z^{R_1}\right\|_2^2 
+3\left\| e_{i}^{\top} \left(\nabla h_t(\Theta_t^{r,R_1}) - \nabla H_t(\Theta_t^{r,R_1}) \right) Z^{R_1}  \right\|_2^2 
\\
&\quad\ + 3\left\| e_{i}^{\top} \left(\nabla H_t(\Theta_t^{r,R_1}) - \nabla \widetilde{H}_t(\alpha_t^r,Z^r) \right) Z^{R_1}  \right\|_2^2 \\
&\leqslant 3 \left\|Z^{R_1}\right\|_{\mathrm{op}}^2 \left\|  e_{i}^{\top} \nabla \widetilde{H}_t(\alpha_t^r,Z^r) \right\|_2^2 + 
 C \zeta_{it}^2 \left(1+ \frac{e^{2M_{\Theta,2}} \kappa^2 k \zeta_{\infty}^2}{\left\|Z^{\star}\right\|_{\mathrm{op}}^{2}}\right)\\
&\quad\  +  3 \left\|Z^{R_1}\right\|_{\mathrm{op}}^2 \left\| e_{i}^{\top} \left(\nabla H_t(\Theta_t^{r,R_1}) - \nabla \widetilde{H}_t(\alpha_t^r,Z^r) \right)\right\|_2^2\\
& \leqslant C\left\|Z^{\star}\right\|_{\mathrm{op}}^2 \left(\left\|  e_{i}^{\top} \nabla \widetilde{H}_t(\alpha_t^r,Z^r) \right\|_2^2 +   \frac{\zeta_{it}^2}{\left\|Z^{\star}\right\|_{\mathrm{op}}^{2}} + \frac{e^{2M_{\Theta,2}}\kappa^2k \zeta_{\infty}^2 \zeta_{it}^2}{\left\|Z^{\star}\right\|_{\mathrm{op}}^{4}} 
  +\frac{e^{4M_{\Theta,2} }M_{\Theta,2}\kappa^2 k \zeta_{\infty}^2}{\min\{n,\left\|Z^{\star}\right\|_{\mathrm{op}}^{2}\}} \right),
\end{align*}
where the second inequality follows from \eqref{s3:6} and the last inequality follows from $0.9 \|Z^{\star}\|_{\mathrm{op}} \leqslant \|Z^{R_1}\|_{\mathrm{op}} \leqslant 1.1 \|Z^{\star}\|_{\mathrm{op}}$.

Finally, with the above two inequalities,
\begin{align} \label{s3:11}
    D_{34}
    &\leqslant \frac{1}{2n} \sum_{t=1}^T \left(e_{i}^{\top} \nabla h_t(\Theta_t^{r,R_1})1_n\right)^2 + \frac{1}{\left\|Z^{0}\right\|_{\mathrm{op}}^{2}} \sum_{t=1}^T\left\|  e_{i}^{\top} \nabla h_t(\Theta_t^{r,R_1})Z^{R_1} \right\|_2^2 \notag\\
    &\leqslant C  \sum_{t=1}^T\left\|  e_{i}^{\top} \nabla \widetilde{H}_t(\alpha_t^r,Z^r) \right\|_2^2 + CT \left(  \frac{\zeta_{i,\infty}^2 + e^{4M_{\Theta,2} }M_{\Theta,2}\kappa^2 k \zeta_{\infty}^2}{ \min\{n,\left\|Z^{\star}\right\|_{\mathrm{op}}^{2}\}} + \frac{e^{2M_{\Theta,2}}\kappa^2k \zeta_{\infty}^2 \zeta_{i,\infty}^2 }{\left\|Z^{\star}\right\|_{\mathrm{op}}^{4}} \right),
\end{align}
where the last inequality partially follows from  \eqref{lem2-1}. Combining \eqref{s3:8}, \eqref{s3:9}, \eqref{s3:10}, \eqref{s3:7} and \eqref{s3:11}, we have
\begin{align*}
    \epsilon_{i\cdot}^{r + 1} &\leqslant \epsilon_{i\cdot}^{r} - \eta\left(\frac{1}{2(e^{-M_{\Theta,2}}+e^{M_{\Theta,2}}) \kappa^2} - 2c_{32} - c_{33}\right) \epsilon_{i\cdot}^{r}\\
    &\quad\ -  \eta \left(\frac{1}{e^{-M_{\Theta,2}}+e^{M_{\Theta,2}}} - \eta C\right) \sum_{t=1}^T\left\|  e_{i}^{\top} \nabla \widetilde{H}_t(\alpha_t^r,Z^r) \right\|_2^2\\
    &\quad\ + \eta CT \left(\frac{1}{c_{32}} + \frac{1}{c_{33}} + \eta\right)\left(  \frac{\zeta_{i,\infty}^2 + e^{4M_{\Theta,2} }M_{\Theta,2}\kappa^2 k \zeta_{\infty}^2}{ \min\{n,\left\|Z^{\star}\right\|_{\mathrm{op}}^{2}\}} + \frac{e^{2M_{\Theta,2}}\kappa^2k \zeta_{\infty}^2 \zeta_{i,\infty}^2 }{\left\|Z^{\star}\right\|_{\mathrm{op}}^{4}} \right),
\end{align*}
where $c_{32}$ and $c_{33}$ are arbitrary constants. Let $c_{32} = c_{33} = ce^{-M_{\Theta,2}}/\kappa^2$ and $c,\eta$ are small enough, and we have 
\begin{align*}
    \epsilon_{i\cdot}^{r+1} \leqslant \left(1-\frac{\eta \rho}{e^{M_{\Theta,2}}\kappa^2}  \right)\epsilon_{i\cdot}^r + \eta C e^{M_{\Theta,2}} \kappa^2 T \left(  \frac{\zeta_{i,\infty}^2 + e^{4M_{\Theta,2} }M_{\Theta,2}\kappa^2 k \zeta_{\infty}^2}{ \min\{n,\left\|Z^{\star}\right\|_{\mathrm{op}}^{2}\}} + \frac{e^{2M_{\Theta,2}}\kappa^2k \zeta_{\infty}^2 \zeta_{i,\infty}^2 }{\left\|Z^{\star}\right\|_{\mathrm{op}}^{4}} \right)
\end{align*}
for some universal positive constants $\rho,C$ and any $r \geqslant R_1$.

As a result, for any $r \geqslant R_1$ and $i = 1, \ldots, n$, 
\begin{align*}
    \epsilon_{i\cdot}^{r} \leqslant \left(1-\frac{\eta \rho}{e^{M_{\Theta,2}} \kappa^{2}} \right)^{r - R_1} \epsilon_{i\cdot}^{R_1}+\frac{C e^{2M_{\Theta,2}}\kappa^{4}T}{\rho} \left( \frac{\zeta_{i,\infty}^2 + e^{4M_{\Theta,2}}  M_{\Theta,2} \kappa^2 k \zeta_{\infty}^2}{\min\{n,\left\|Z^{\star}\right\|_{\mathrm{op}}^{2}\}} + \frac{e^{2M_{\Theta,2}} \kappa^2 k \zeta_{\infty}^2\zeta_{i,\infty}^2 }{\left\|Z^{\star}\right\|_{\mathrm{op}}^{4}}\right) 
\end{align*}
follows from iterating over the above equation.

\subsubsection{Proof of Lemma \ref{thm:linear4}} \label{sec:pf_linear4}
For any $i,t$ and $r \geqslant R_{1.5}$ we have
\begin{align*}
    \left(\alpha_{it}^{r+1} - \alpha_{it}^{\star}\right)^2 \leqslant \left(\alpha_{it}^{r} - \alpha_{it}^{\star}\right)^2 + \eta_{\alpha}^2  \left(e_{i}^{\top} \nabla h_t(\Theta_t^{r,R_1})1_n\right)^2 - 2\eta_{\alpha}  \left(e_{i}^{\top} \nabla h_t(\Theta_t^{r,R_1})1_n\right)(\alpha_{it}^r - \alpha_{it}^{\star}).
\end{align*}
Equivalently,
\begin{align} \label{s4:33}
    e_{it}^{r+1} &\leqslant e_{it}^r - 2\eta \left(e_{i}^{\top} \nabla h_t(\Theta_t^{r,R_1})1_n\right)(\alpha_{it}^r - \alpha_{it}^{\star}) + \frac{\eta^2}{2n}\left(e_{i}^{\top} \nabla h_t(\Theta_t^{r,R_1})1_n\right)^2 .
\end{align}
We decompose $\nabla h_{t}(\Theta_{t}^{r,R_1})= \nabla H_t(\Theta_t^{r,\star;\star}) + \{\nabla h_t(\Theta_t^{r,R_1}) - \nabla H_t(\Theta_t^{r,\star;\star})\} $ in \eqref{s4:33}, where $\Theta_t^{r,\star;\star} = \alpha_t^r 1_n^{\top} + 1_n (\alpha_t^{\star})^{\top} + Z^{\star}(Z^{\star})^{\top}$. Then we obtain
\begin{align} \label{s4:3}
 e^{r+1}_{it} \leqslant e_{it}^r - \eta D_{41} + \eta D_{42} + \eta^2 D_{43}
\end{align}
where
\begin{align*}
    D_{41} = &~ 2 \left(e_{i}^{\top} \nabla H_t(\Theta_t^{r,\star;\star})1_n\right)(\alpha_{it}^r - \alpha_{it}^{\star})  ,\\
    D_{42} = &~ 2 \left|e_i^{\top} \left(\nabla h_t(\Theta_t^{r,R_1}) - \nabla H_t(\Theta_t^{r,\star;\star})\right) 1_n\right| \cdot\left|\alpha_{it}^r - \alpha_{it}^{\star}\right| ,\\
    D_{43} = &~ \frac{1}{2n}\left(e_{i}^{\top} \nabla h_t(\Theta_t^{r,R_1})1_n\right)^2.
\end{align*}

By Lemma \ref{Lemma-30} we have
\begin{align} \label{s4:4}
    D_{41} &= 2 \sum_{j=1}^n \left(\exp\big(\alpha_{it}^r + \alpha_{jt}^{\star} + (z_i^{\star})^{\top}z_j^{\star}\big) - \exp\big(\alpha_{it}^{\star} + \alpha_{jt}^{\star} + (z_i^{\star})^{\top}z_j^{\star}\big) \right) (\alpha_{it}^r  - \alpha_{it}^{\star}) \notag\\
    &\geqslant \frac{2n}{e^{-M_{\Theta,2}}+e^{M_{\Theta,2}}} (\alpha_{it}^r  - \alpha_{it}^{\star})^2 + \frac{2}{e^{-M_{\Theta,2}}+e^{M_{\Theta,2}}}\left\|e_i^{\top} \nabla H_t(\Theta_t^{r,\star;\star}) \right\|_2^2.
\end{align}
For any $c_{42} > 0$,
\begin{align} \label{s4:1}
    D_{42} &\leqslant c_{42} (\alpha_{it}^r - \alpha_{it}^{\star})^2 + \frac{1}{c_{42}} \left(e_i^{\top} \left(\nabla h_t(\Theta_t^{r,R_1}) - \nabla H_t(\Theta_t^{r,\star;\star})\right) 1_n\right)^2\notag\\
    &\leqslant c_{42} (\alpha_{it}^r - \alpha_{it}^{\star})^2 + \frac{2}{c_{42}} \zeta_{it}^2 +  \frac{2}{c_{42}}\left(e_i^{\top} \left(\nabla H_t(\Theta_t^{r,R_1}) - \nabla H_t(\Theta_t^{r,\star;\star})\right) 1_n\right)^2.
\end{align}
Note that for any $r \geqslant R_{1.5}$,
\begin{align} \label{s4:2}
    &\quad\  \Big(e_i^{\top} \left(\nabla H_t(\Theta_t^{r,R_1}) - \nabla H_t(\Theta_t^{r,\star;\star})\right) 1_n\Big)^2 \notag\\
    &\leqslant n \left\|e_i^{\top} \left(\nabla H_t(\Theta_t^{r,R_1}) - \nabla H_t(\Theta_t^{r,\star;\star})\right) \right\|_2^2 \notag\\
    &\leqslant Cne^{2M_{\Theta,2}}\left(\left\|\alpha^{R_1}_t - \alpha_t^{\star}\right\|_{2}^2 + \left\|Z^{R_1} z_i^r  - Z^{\star} z_i^{\star} \right\|_2^2\right)\notag\\
    &\leqslant Cne^{2M_{\Theta,2}}\left(\left\|\alpha^{R_1}_t - \alpha_t^{\star}\right\|_{2}^2 + \left\|Z^{R_1} - Z^{\star}Q_{R_1}\right\|_{\mathrm{F}}^2 \left\| z_i^r \right\|_2^2 + \left\|Z^{\star}\right\|_{\mathrm{op}}^2 \left\|Q_{R_1}z_i^r -  z_i^{\star} \right\|_2^2\right)\notag\\
    &\leqslant Cne^{2M_{\Theta,2}}\left( \frac{e^{6M_{\Theta,2}}\kappa^2k\zeta_{\infty}^2}{n} + \frac{e^{2M_{\Theta,2}}M_{\Theta,2}\kappa^2k\zeta_{\infty}^2}{\left\|Z^{\star}\right\|_{\mathrm{op}}^{2}} + \frac{e^{6M_{\Theta,2}} M_{\Theta,2} \kappa^6 k \zeta_{\infty}^2 \zeta_{\infty,\infty}^2 }{\min\{n^2,\left\|Z^{\star}\right\|_{\mathrm{op}}^{4}\}}\right)\\
    &\leqslant Cne^{2M_{\Theta,2}}\left( \frac{e^{6M_{\Theta,2}}\kappa^2k\zeta_{\infty}^2}{\min\{n,\left\|Z^{\star}\right\|_{\mathrm{op}}^{2}\}} + \frac{e^{6M_{\Theta,2}} M_{\Theta,2} \kappa^6 k \zeta_{\infty}^2 \zeta_{\infty,\infty}^2 }{\min\{n^2,\left\|Z^{\star}\right\|_{\mathrm{op}}^{4}\}}\right),
\end{align}
where the fourth inequality follows from Lemma \ref{thm:linear}, Lemma \ref{thm:linear2} and Lemma \ref{thm:linear3}. Then combining two above inequalities \eqref{s4:1} and \eqref{s4:2}, we obtain
\begin{align} \label{s4:5}
    D_{42} \leqslant c_{42} (\alpha_{it}^r - \alpha_{it}^{\star})^2 + \frac{Cn}{c_{42}} \left(\frac{\zeta_{it}^2}{n} +  \frac{e^{8M_{\Theta,2}}\kappa^2k\zeta_{\infty}^2}{\min\{n,\left\|Z^{\star}\right\|_{\mathrm{op}}^{2}\}}  + \frac{e^{8M_{\Theta,2}} M_{\Theta,2} \kappa^6 k \zeta_{\infty}^2 \zeta_{\infty,\infty}^2 }{\min\{n^2,\left\|Z^{\star}\right\|_{\mathrm{op}}^{4}\}}\right).
\end{align}
For the third term, we have
\begin{align} \label{s4:6}
    D_{43} &\leqslant \frac{1}{n}\left(e_i^{\top}  \nabla H_t(\Theta_t^{r,\star;\star}) 1_n\right)^2 + \frac{1}{n}\left(e_i^{\top} \left(\nabla h_t(\Theta_t^{r,R_1}) - \nabla H_t(\Theta_t^{r,\star;\star})\right) 1_n\right)^2 \notag\\
    &\leqslant \left\|e_i^{\top}  \nabla H_t(\Theta_t^{r,\star;\star})\right\|_2^2 + C\left(\frac{\zeta_{it}^2}{n} +  \frac{e^{8M_{\Theta,2}}\kappa^2k\zeta_{\infty}^2}{\min\{n,\left\|Z^{\star}\right\|_{\mathrm{op}}^{2}\}} + \frac{e^{8M_{\Theta,2}} M_{\Theta,2} \kappa^6 k \zeta_{\infty}^2 \zeta_{\infty,\infty}^2 }{\min\{n^2,\left\|Z^{\star}\right\|_{\mathrm{op}}^{4}\}}\right),
\end{align}
where the last inequality is derived similarly with \eqref{s4:1} and \eqref{s4:2}. 

Combining \eqref{s4:3}, \eqref{s4:4}, \eqref{s4:5} and \eqref{s4:6}, we obtain 
\begin{align*}
    e_{it}^{r+1} &\leqslant e_{it}^r - \eta \left(\frac{1}{e^{-M_{\Theta,2}}+e^{M_{\Theta,2}}} - \frac{c_{42}}{2n} \right)  e_{it}^r - \eta \left(\frac{2}{e^{-M_{\Theta,2}}+e^{M_{\Theta,2}}} - \eta\right) \left\|e_i^{\top}  \nabla H_t(\Theta_t^{r,\star;\star})\right\|_2^2\\
    &\quad\  + \eta C \left(\frac{n}{c_{42}} + \eta\right) \left(\frac{\zeta_{it}^2}{n} +  \frac{e^{8M_{\Theta,2}}\kappa^2k\zeta_{\infty}^2}{\min\{n,\left\|Z^{\star}\right\|_{\mathrm{op}}^{2}\}} + \frac{e^{8M_{\Theta,2}} M_{\Theta,2} \kappa^6 k \zeta_{\infty}^2 \zeta_{\infty,\infty}^2 }{\min\{n^2,\left\|Z^{\star}\right\|_{\mathrm{op}}^{4}\}}\right),
\end{align*}
where $c_{42}$ is an arbitrary constant. Let $c_{42} = cne^{-M_{\Theta,2}}$ and $c, \eta$ are small enough, we have 
\begin{align*}
    e_{it}^{r+1} \leqslant \left(1- \frac{\eta \rho}{e^{M_{\Theta,2}}} \right)e_{it}^r + \eta Ce^{M_{\Theta,2}}   \left(\frac{\zeta_{it}^2}{n} +  \frac{e^{8M_{\Theta,2}}\kappa^2k\zeta_{\infty}^2}{\min\{n,\left\|Z^{\star}\right\|_{\mathrm{op}}^{2}\}} + \frac{e^{8M_{\Theta,2}} M_{\Theta,2} \kappa^6 k \zeta_{\infty}^2 \zeta_{\infty,\infty}^2 }{\min\{n^2,\left\|Z^{\star}\right\|_{\mathrm{op}}^{4}\}}\right)
\end{align*}
for some universal positive constants $\rho, C$ and any $r \geqslant R_{1.5}$.

Then for any $i = 1, \ldots, n$, $t = 1, \ldots, T$ and $r \geqslant R_{1.5}$,
\begin{align*}
    e_{it}^{r} \leqslant \left(1-\frac{\eta \rho}{e^{M_{\Theta,2}} } \right)^{r - R_{1.5}} e_{it}^{R_{1.5}}+\frac{C e^{2M_{\Theta,2}}}{\rho} \left(\frac{\zeta_{it}^2}{n} + \frac{e^{8M_{\Theta,2}}\kappa^2k\zeta_{\infty}^2}{\min\{n,\left\|Z^{\star}\right\|_{\mathrm{op}}^{2}\}} + \frac{e^{8M_{\Theta,2}} M_{\Theta,2} \kappa^6 k \zeta_{\infty}^2 \zeta_{\infty,\infty}^2 }{\min\{n^2,\left\|Z^{\star}\right\|_{\mathrm{op}}^{4}\}}\right)
\end{align*}
follows from iterating over the above equation.

\subsubsection{Proof of Lemma \ref{thm:linear5}}\label{sec:pf_linear5}
For any $1\leqslant i\leqslant n$, $$\check{z}_i = z_i^{R_2} - \frac{1}{n} \sum_{j = 1}^n z_j^{R_2}  = z_i^{R_2} - \frac{1}{n}\sum_{j = 1}^n (z_j^{R_2} - Q_{R_1}^\top z^\star_j).$$ 
Then
\begin{align*}
    \big\|\check{z}_{i} - Q_{R_1}^\top z^\star_i \big\|_2 \leqslant  \big\|{z}_{i}^{R_2} - Q_{R_1}^\top z^\star_i \big\|_2 + \big\|\frac{1}{n}\sum_{j = 1}^n (z_j^{R_2} - Q_{R_1}^\top z^\star_j) \big\|_2 \leqslant 2 \max_j \big\|{z}_{j}^{R_2} - Q_{R_1}^\top z^\star_j \big\|_2.
\end{align*}
Combining with Lemma \ref{thm:linear3}, we have 
\begin{align} \label{t1-1}
\max_i \big\|\check{z}_{i} - Q_{R_1}^\top z^\star_i \big\|_2^2  \leqslant {C e^{6M_{\Theta,2}} M_{\Theta,2} \kappa^{6} k} \cdot \frac{\zeta_{\infty}^2\zeta_{\infty,\infty}^2 }{\min\{n^3,\left\|Z^{\star}\right\|_{\mathrm{op}}^{6}\}}.
\end{align}
Finally, we will show $\|\check{z}_{i} - \check{Q}^\top z^\star_i \|_2^2$ has a similar bound. The definition of $\check{Q}$ shows that
$\|\check{Z} - Z^{\star}\check{Q}\|_{\mathrm{F}} \leqslant \|\check{Z} - Z^{\star}Q_{R_1}\|_{\mathrm{F}}$, then we have
\begin{align*}
    \frac{\left\|Z^{\star}\right\|_{\mathrm{op}}}{\kappa} \,  \left\|Q_{R_1} - \check{Q}\right\|_{\mathrm{op}} \leqslant
    \sigma_k(Z^{\star}) \left\|Q_{R_1} - \check{Q}\right\|_{\mathrm{F}} \leqslant \left\|Z^{\star} (Q_{R_1} - \check{Q})\right\|_{\mathrm{F}} \leqslant 2\left\|\check{Z} - Z^{\star}Q_{R_1}\right\|_{\mathrm{F}}.
\end{align*}
As a result, 
\begin{align*}
    \big\|\check{z}_{i} - \check{Q}^\top z^\star_i \big\|_2^2 
    &\leqslant 2\big\|\check{z}_{i} - Q_{R_1}^\top z^\star_i \big\|_2^2 + 2 \left\|Q_{R_1}- \check{Q}\right\|_{\mathrm{op}}^2 \left\|z_i^{\star}\right\|_2^2 \\
    &\leqslant \Bigg( \frac{8 M_{\Theta,2} \kappa^2  n}{\left\|Z^{\star}\right\|_{\mathrm{op}}^2} + 2 \Bigg) \, \max_j \big\|\check{z}_{j} - Q_{R_1}^\top z^\star_j \big\|_2^2.
\end{align*}
Combining with \eqref{t1-1} we have
\begin{align*}
    \max_{1\leqslant i\leqslant n} \big\|\check{z}_{i} - \check{Q}^\top z^\star_i \big\|_2^2  \leqslant {C e^{6M_{\Theta,2}} M_{\Theta,2}^2 \kappa^{8} k} \,  \frac{n\zeta_{\infty}^2\zeta_{\infty,\infty}^2 }{\min\{n^4,\left\|Z^{\star}\right\|_{\mathrm{op}}^{8}\}}.
\end{align*}

\subsubsection{Proof of Lemma \ref{lem:concen}} \label{sec:pf_concen}
We first establish the probabilistic upper bound of $\zeta_{\infty}^2 = \max_t \zeta_{t}^2$. For $t = 1,\ldots, T$, we can write the $\mathbf A_{t} - \exp(\Theta^{\star}_t)$ as
\begin{align*}
   \mathbf A_{t} - \exp(\Theta^{\star}_t)
    =\sum_{1\leqslant i < j \leqslant n}\left(A_{t,ij} - \exp(\Theta^{\star}_{t,ij})\right)(\delta_{ij}+\delta_{ji}) + \sum_{i=1}^n\left(A_{t,ii} - \exp(\Theta^{\star}_{t,ii})\right)\delta_{ii}, 
\end{align*}
where we let $\delta_{ij}:=e_ie_j^\top$, and $e_i$ is an indicator vector with $i$-th position being 1 and remaining entries being 0. Here $(A_{t,ij} - \exp(\Theta^{\star}_{t,ij}))(\delta_{ij}+\delta_{ji})$ and $(A_{t,ii} - \exp(\Theta^{\star}_{t,ii}))\delta_{ii}$ are independent, symmetric random matrices with zero means, and
\begin{align*}
    \EXPT\left[\left(A_{t,ij} - \exp(\Theta^{\star}_{t,ij})\right)(\delta_{ij}+\delta_{ji})\right]^{\ell}
    =\left\{
    \begin{aligned}
    \EXPT\left(A_{t,ij} - \exp(\Theta^{\star}_{t,ij})\right)^{\ell}\cdot(\delta_{ij}+\delta_{ji}),\quad  \text{$\ell$ is odd,}\ \\
    \EXPT\left(A_{t,ij} - \exp(\Theta^{\star}_{t,ij})\right)^{\ell}\cdot(\delta_{ii}+\delta_{jj}),\quad \text{$\ell$ is even.}
    \end{aligned}
    \right.
\end{align*}

Without loss of generality, we assume $e^{-M_{\Theta,1}} < 1/2$, or $A_{t,ij}$ can be decomposed into independent Poisson random variables with smaller means. Then using Lemma \ref{Lemma-A2}, we have 
$$
\EXPT\left(A_{t,ij} - \exp(\Theta^{\star}_{t,ij})\right)^{\ell}\leqslant \frac{\ell!}{2}\left(\frac{e^2+1}{1-e^{-M_{\Theta,1}}}\right)^{\ell-2}e^{-M_{\Theta,1}}.
$$
Set $B=2(e^2+1),A_{ij}^2=e^{-M_{\Theta,1}}(\delta_{ii}+\delta_{jj})$, for any $\ell\geqslant2$ we obtain
$$
\EXPT\left[\left(A_{t,ij} - \exp(\Theta^{\star}_{t,ij})\right)(\delta_{ij}+\delta_{ji})\right]^{\ell}
\preceq \frac{\ell!}{2} B^{\ell-2}A_{ij}^2.
$$
Similarly, we have 
$$
\EXPT\left[\left(A_{t,ii} - \exp(\Theta^{\star}_{t,ii})\right)\delta_{ii}\right]^{\ell}
\preceq \frac{\ell!}{2} B^{\ell-2}A_{ii}^2,
$$
for any $\ell \geqslant 2$, where $A_{ii}^2 = e^{-M_{\Theta,1}} \delta_{ii}$.
Because of $\sigma^2= \|\sum_{i \leqslant j } A_{ij}^2\|_{\mathrm{op}}\leqslant n e^{-M_{\Theta,1}}$, we use Lemma \ref{Lemma-MB} and know that for any $x \geqslant 0$, 
\begin{align*}
    \Pr\left(\left\|\mathbf A_{t} - \exp(\Theta^{\star}_t)\right\|_{\mathrm{op}}\geqslant x\right)\leqslant 2n \, \exp\left(\frac{-x^2/2}{\sigma^2 + Bx}\right), \quad t= 1, \ldots, T.
\end{align*}
As a result, 
\begin{align*}
    \Pr\left(\max_{1\leqslant t\leqslant T} \zeta_t \geqslant x\right) \leqslant \sum_{t=1}^T \Pr\left(\zeta_t \geqslant x\right) \leqslant 2nT\, \exp\left(\frac{-x^2/2}{\sigma^2 + Bx}\right).
\end{align*}
Set $x=8(s+1)\cdot\max\{\sigma\sqrt{\log (nT)},B\log (nT)\}$, then we have
\begin{align*}
    \Pr\left(\max_{1\leqslant t\leqslant T} \zeta_{t}^2 \geqslant C (s+1)^2 n \log(nT) \,  \max\left\{ e^{-M_{\Theta,1}}, \frac{\log(nT)}{n}\right\}\right) \leqslant (nT)^{-s}.
\end{align*}
Since $M_{\Theta,1} \geqslant 0$, we can simplify the above probabilistic inequality to 
$$\Pr\left(\zeta_{\infty}^2 
\geqslant C(s+1)^2n \log^2(nT)\right)\leqslant (nT)^{-s}. $$

For $\zeta_{\infty,\infty}^2 = \max_{i,t} \zeta_{it}^2$, we decomposed it into $\max_{i,t} (\zeta_{it,1}^2 +  \zeta_{it,2}^2 +  \zeta_{it,3}^2) + n$, where
\begin{align*}
    \zeta_{it,1}^2 =  \Big(\sum_{j=1}^n M_{t,ij}\Big)^2, \quad \zeta_{it,2}^2 = \Big\| \sum_{j=1}^n  M_{t,ij}z_{j}^{\star}\Big\|_2^2, \quad \text{and} \quad \zeta_{it,3}^2 = \sum_{j=1}^n M_{t,ij}^2     
\end{align*}
with $M_{t,ij}=A_{t,ij} -\exp(\Theta_{t,ij}^{\star}).$
By Lemma \ref{lem:bern_pois} and similar arguments we used for $\zeta_{\infty}$, we have the following probabilistic upper bound:
\begin{align*}
    &\Pr\left( \max_{1 \leqslant i\leqslant n,1\leqslant t\leqslant T} \zeta_{it,1}^2 \geqslant C(s+1)^2 n \log(nT)  \max\{e^{-M_{\Theta,1}}, \frac{\log(nT)}{n}\}\right) \leqslant (nT)^{-s},\\
    &\Pr\left( \max_{1 \leqslant i\leqslant n,1\leqslant t\leqslant T} \zeta_{it,2}^2 \geqslant C(s+1)^2 nkM_{\Theta,2} \log(nT)  \max\{e^{-M_{\Theta,1}}, \frac{\log(nT)}{n}\}\right) \leqslant k(nT)^{-s},\\
    &\Pr\left( \max_{1 \leqslant i\leqslant n,1\leqslant t\leqslant T} \zeta_{it,3}^2 \geqslant C(s+2)^2 n\log^2(nT) \right) \leqslant (nT)^{-s}.
\end{align*}
Combining those inequalities, we have for any $s > 0$,
$$
\Pr\left( \zeta_{\infty, \infty}^2  \geqslant C(s+2)^2 e^{M_{\Theta,2}} nk \log^2(nT) \right) \leqslant (k+2)(nT)^{-s}.
$$

\newpage

\subsection{Simulation Results for the Initial Estimator}\label{sec:initialsim}


In this section, we 
 discuss the relationship between the initial estimator and the one-step estimator and compare their numerical performances. 

The projected gradient descent with an appropriate initialization in Algorithm 
 \ref{algor:init} is a commonly adopted method in the related literature \citep{ma2020universal,zhang2020flexible}. 
Both the one-step estimator and gradient descent algorithm are similar in the sense that they are motivated from solving the maximum likelihood estimator for $Z$. 
They can be viewed as different methods for solving the same problem. 
However, the performance of the gradient descent algorithm depends on tuning parameters such as the step size and the stopping criteria in practice. 
On the other hand, the one-step estimator does not depend on any tuning parameters and provide a useful refinement tool that  further improves the output from the projected gradient descent regardless of the choice of tuning parameters.  

\begin{figure}[!htbp]
\centering
\begin{subfigure}{0.32\textwidth}
	\centering
    \includegraphics[width=1\linewidth]{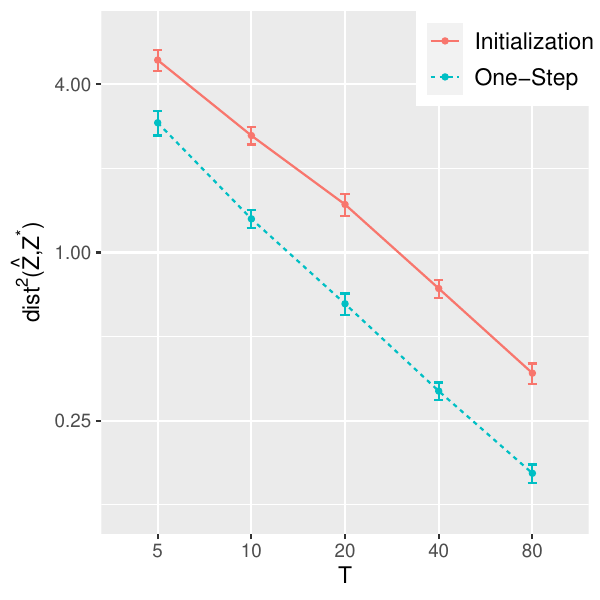}
    \caption{$t_Z= 10^{-1}$}
\end{subfigure}
\begin{subfigure}{0.32\textwidth}
	\centering
    \includegraphics[width=1\linewidth]{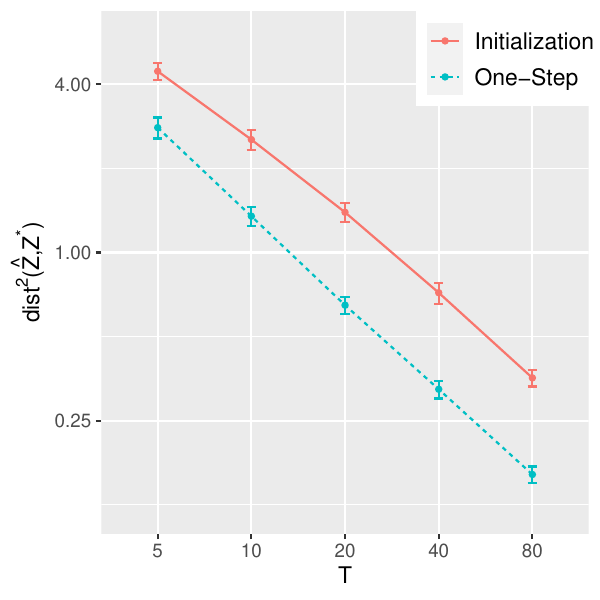}
    \caption{$t_Z= 10^{-2}$}
\end{subfigure}
\begin{subfigure}{0.32\textwidth}
	\centering
    \includegraphics[width=1\linewidth]{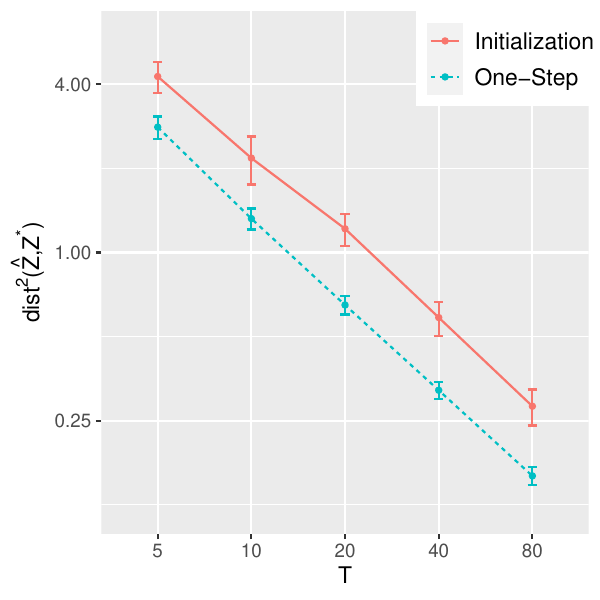}
    \caption{$t_Z= 10^{-3}$}
\end{subfigure}
\caption{Estimation errors of the initial estimator with different thresholds $t_Z$ and the one-step estimator: $\mathrm{dist}^2(\hat{Z},Z^{\star})$ versus $T$ under Case (I) with $n=200$.}	\label{fig:pgdvsonestep_supp}
\end{figure}

As an illustration, 
We empirically examine the estimation errors of both estimators versus $T$ under Case (I) with $n=200$ over 50 repetitions. 
For the initial estimator, we  set the projected gradient descent algorithm to stop if the 
distance between  estimates of two consecutive iterations 
is smaller than $10^{-t_Z}$ for $t_Z\in \{1,2,3\}$.   
The results are presented in Figure \ref{fig:pgdvsonestep_supp}, showing that the one-step update can effectively improve the output from the initial gradient descent algorithm regardless of the stopping criteria.

\begin{remark}
Theoretically, for the projected gradient descent algorithm, 
it is unclear how to establish, if they exist, the desired oracle properties for such a method, due largely to the lack of theoretical tools. 
On the other hand, 
 the one-step approach makes use of 
  the semiparametric efficient score equation for the target parameter that is more manageable since it bypasses the direct estimation of the  nuisance parameter. Taking advantage of the explicit score equation form, we are able to establish the oracle properties for this approach. 	
 Furthermore, we would like to emphasize that the proposed one-step update formula is a general construction  and does not rely on a specific initialization.  
We can provide different initialization estimators into the one-step update formula, and it would serve as a useful one-shot refinement method for improving input estimates. 
\end{remark}

\newpage

\section{Technical Results for the Penalized Maximum Likelihood Estimator}

\subsection{Proof of Theorem \ref{MLE_mainthm}} \label{sec:pfthm3A}

As discussed in Remark \ref{rm:techchanlplm}, 
to derive a sharp estimation error rate for $\hat{G}$,  
we establish an analysis based on profile  likelihood of the target parameters $G$.  
In particular, 
we define the log profile likelihood of $G$ as $pl(G) = l(G,\hat\alpha(G)) $, where $l(G,\alpha)$ denotes the   log-likelihood function of 
  $(G,\alpha)$, and $\hat\alpha(G) = \argmax_{\alpha\in \mathbb{R}^{n\times T}} l(G,\alpha)  $.

\vspace{5pt}
\textbf{{Step 1:}}
By definitions of $\hat{G}$ and $\hat{\alpha}(G)$, 
\begin{align}
0\leqslant  &~ l(\hat{G},\hat{\alpha}(\hat{G})) - \lambda_{n,T} \|\hat G\|_* - l(G^{\star}, \hat{\alpha}(G^{\star})) +  \lambda_{n,T} \|G^{\star}\|_*\notag\\
=&~ pl(\hat{G}) - pl(G^{\star}) + \lambda_{n,T} (\|G^{\star}\|_* - \|\hat{G}\|_*),  \label{basic_ineqA}
\end{align}
where in the second equation, we use $pl(G) = l(G,\hat{\alpha}(G)) $. 
Given $G\in \mathbb{R}^{n\times n}$, define its vectorization as $G_v = (G_{11}/\sqrt 2, \ldots, G_{nn}/\sqrt 2, G_{12}, \ldots, G_{n-1,n} )^\top \in \mathbb R^{\frac{n(n+1)}{2} \times 1}.$ Let 
\begin{align}\label{eq:deltagdefinitionA}
    \Delta_G = \hat G - G^\star \in \mathbb{R}^{n\times n} \hspace{2em} \text{and}\hspace{2em} \Delta G_v = \hat G_v - G_v^\star \in \mathbb R^{\frac{n(n+1)}{2} \times 1} 
\end{align}
satisfying $2\|\Delta G_v\|_2^2 = \|\Delta_G\|_{\mathrm{F}}^2 $.
We derive an upper bound of $\|\Delta_G\|_{\mathrm{F}}^2$ through the Taylor expansion of $pl(\hat{G})-pl(G^{\star})$ in \eqref{basic_ineqA} with respect to $G_v$. 
In particular,   the Right Hand Side of 
$\eqref{basic_ineqA} = I + II+ \lambda_{n,T} (\|G^{\star}\|_* - \|\hat{G}\|_*), 
$ 
where we define 
\begin{align*}
 I=&~(\Delta G_v)^\top \frac{\partial\ pl(G)}{\partial G_v}\biggr|_{G^{\star}}, \hspace{4em}
    II= \frac{1}{2}(\Delta G_v)^{\mytrans}\frac{\partial^2 \ pl(G)}{\partial G_v G_v^{\mytrans}}\biggr|_{\tilde{G}}(\Delta G_v), 
\end{align*}
in which formulae of 
${\partial pl}/{\partial G_v}$ and ${\partial^2 pl}/{\partial G_v\partial G_v^{\mytrans}}$ are derived in Section \ref{subsec::prelim_G} of the supplementary material,
and in $II$, $\tilde{G}$ is a convex combination of $\hat{G}$ and $G^{\star}$. 

\vspace{5pt}
\textbf{{Step 2:}}
We next derive an upper bound of $I$.
By  the formulae in \eqref{eq:ltdot_Gv} and \eqref{eq::dpltG_dGv=dltG_dGv}, 
 $I=\frac{1}{2}\langle I_a+I_b, \Delta_{G}\rangle$, where $I_a$ and $I_b$ are $n\times n$ matrices with their $(i,j)$-th element for $1\leqslant i,j\leqslant n$ defined as
 $ (I_a)_{ij}  = \sum_{t=1}^T \{A_{t,ij} - \exp(\alpha^{\star}_{it} + \alpha^{\star}_{jt} + G^{\star}_{ij})\}\{1+1_{(j=i)}\},$
 and $(I_b)_{ij} = \sum_{t=1}^T \{\exp(\alpha^{\star}_{it} + \alpha^{\star}_{jt} + G^{\star}_{ij}) -\exp( \hat{\alpha}_{it}(G^{\star}) + \hat{\alpha}_{jt}(G^{\star}) + G^{\star}_{ij} ) \}\{1+1_{(j=i)}\}$, respectively. 
Then by Von Neumann's trace inequality and Cauchy-Schwarz inequality, 
 \begin{align} \label{pfthm3-1A}
    I \leqslant \frac{1}{2}\|I_a \|_{\operatorname{op}} \|\Delta_G\|_* + \frac{1}{2}\| I_b\|_{\mathrm{F}} \| \Delta_G\|_{\mathrm{F}}. 
\end{align}

Note that $I_a$ is a $n\times n$ symmetric matrix with independent mean-zero  random variables in the entries.
We apply a random matrix concentration inequality (Lemma  \ref{Lemma-concen} in the Supplementary Material) 
to $I_a$ and 
 then obtain that  
 for any constant $s>0$, there exists $ C_s>0$ such that
\begin{align} 
\|I_a\|_{\operatorname{op}} \leqslant C_s \sqrt{nT} 
\label{lm:thm3boundiaA} 
\end{align}
with probability $1-O(n^{-s})$. 
In addition, we  show that for any constant $s>0$, 
there exists a constant $C_s>0$ such that 
\begin{align}
     \|I_b\|_{\mathrm{F}} \leqslant &~ C_{s} \sqrt{nT}\max\left\{1, \sqrt{\frac{T}{n}}\right\}\log(nT)=C_s\sqrt{nT}\sqrt{r'_{n,T}} 
    \label{lem_termI_inbasicineqA}  
\end{align}
with probability $1-O((nT)^{-s})$. 
We can interpret 
$  \|I_b\|_{\mathrm{F}}$
as an aggregated error from nuisance parameters
by noticing  
$(I_b)_{ij}=\sum_{t=1}^T \{\exp \left(\Delta \alpha_{i t}+\Delta \alpha_{j t}\right)-1\}\exp(\alpha^{\star}_{it} + \alpha^{\star}_{jt} + G^{\star}_{ij})$, 
where $\Delta \alpha_{i t}=\hat{\alpha}_{i t}(G^{\star})-\alpha_{i t}^{\star}$ is the difference between the profiled nuisance estimator $\hat{\alpha}_{it}(G^{\star})$ and true   $\alpha^{\star}_{it}$.
Intuitively,  
  \eqref{lem_termI_inbasicineqA} 
 shows that the aggregated error from nuisance parameters $  \|I_b\|_{\mathrm{F}}$ can be  well-controlled. 
The proof of \eqref{lem_termI_inbasicineqA} requires involved calculation, so we defer the details to Section \ref{sec:pfD1}  of the Supplementary Material.  


\vspace{5pt}
\textbf{{Step 3:}}
We next derive an upper bound of $II$.
In particular, we can show that for any $s>0$,  there exists a universal constant $C_s>0$ such that
\begin{align}\label{eq_IIinbasicineq}
    II \leqslant - C_{s} T \|\Delta_G\|_{\mathrm{F}}^2
\end{align}
with probability $1-O((nT)^{-s})$ 
(see the proof of \eqref{eq_IIinbasicineq} in Section \ref{sec:pfD2} of the Supplementary Material). 
Note that $II$ takes a quadratic form.
But similarly to Step 2 in Section \ref{sec:pfthmnerr}, 
proving \eqref{eq_IIinbasicineq} is challenged by the fact that 
the Hessian matrix $H(G):={\partial^2 pl(G)}/{\partial G_v \partial G_v^{\top}}$ is singular. 
To address the issue, 
we  characterize the eigen space of  $H(G)$ explicitly and show $\Delta G_v$ belongs to the column space of  $H(G)$. 
The details are provided in  Section \ref{sec:pfD2}  of the Supplementary Material.



\vspace{5pt}
\textbf{{Step 4:}}
By Steps 1--3, we next establish an upper bound of $\|\Delta_G\|_{\mathrm{F}}^2$. 
In particular, by \eqref{basic_ineqA}, \eqref{pfthm3-1A}, and \eqref{eq_IIinbasicineq},
we obtain that  
\begin{align}
C_{s} T \|\Delta_G \|_{\mathrm{F}}^2 \leqslant  &~ \|I_a\|_{\operatorname{op}} \|\Delta_G\|_* +  \|I_b\|_{\mathrm{F}} \| \Delta_G\|_{\mathrm{F}} + \lambda_{n,T} ( \|G^{\star}\|_* - \|\hat G\|_* )\label{basic_ineq3A_original}\\
\leqslant &~  \|I_a\|_{\operatorname{op}} \|\Delta_G\|_* +  \|I_b\|_{\mathrm{F}} \| \Delta_G\|_{\mathrm{F}} + \lambda_{n,T}\|\Delta_G\|_*\label{basic_ineq3A}
\end{align}
with probability $1-O((nT)^{-s}) $, 
where \eqref{basic_ineq3A} follows by the triangle inequality.

When $\lambda_{n,T} \geqslant 2 \|I_a\|_{\operatorname{op}} $, 
we can obtain 
\begin{align}\label{lem_basicineq_step3A}
\|\Delta_G\|_* \leqslant 4\sqrt{2k} \|\Delta_{G}\|_{\mathrm{F}} 
+ \frac{2}{\lambda_{n,T}} \|I_b\|_{\mathrm{F}}  \| \Delta_G\|_{\mathrm{F}}. 
\end{align}
The proof of \eqref{lem_basicineq_step3A} is given in Section \ref{sec:pfD3}  of the Supplementary Material. 
Plugging \eqref{lem_basicineq_step3A} 
into \eqref{basic_ineq3A}, we obtain $C_{s} T \|\Delta_G \|_{\mathrm{F}}^2 
\leqslant  
(6\sqrt{2k}\lambda_{n,T} + 4 \|I_b\|_{\mathrm{F}}   ) \|\Delta_G\|_{\mathrm{F}},$ which is a quadratic inequality of $\|\Delta_G \|_{\mathrm{F}}$. Solving the inequality gives   
\begin{align}
\|\Delta_G \|_{\mathrm{F}}^2 \leqslant  &~ C_{s}  \frac{\|I_b\|_{\mathrm{F}}^2 + k \lambda_{n,T}^2}{T^2}.\label{eq:pfthm3-bddeltag0A}
\end{align}


Take $\lambda_{n,T} \asymp \sqrt{nT}\log(nT) $. 
By  \eqref{lm:thm3boundiaA},  
we know that 
$\lambda_{n,T}\geqslant\|I_a\|_{\operatorname{op}}$ with probability  $1-O(n^{-s})$ 
when $nT$ is sufficiently large. 
Then   \eqref{lem_termI_inbasicineqA} and \eqref{eq:pfthm3-bddeltag0A} give
$\|\Delta_G\|_{\operatorname{F}}^2/n \leqslant C_sT^{-1} r'_{n,T}$
with probability $1-O(n^{-s})$.

\subsection{Proof of Eq. \texorpdfstring{\eqref{lem_termI_inbasicineqA}, \eqref{eq_IIinbasicineq}, and \eqref{lem_basicineq_step3A} in Section \ref{sec:pfthm3A}}{}}\label{sec:lem_of_thm}

To finish the proof of  Theorem \ref{MLE_mainthm},
it remains to prove Eq.  \eqref{lem_termI_inbasicineqA}, \eqref{eq_IIinbasicineq}, and \eqref{lem_basicineq_step3A} in Section \ref{sec:pfthm3A}. 
We next present a technical Lemma \ref{lem_bdd} to be used below,
and then we 
prove Eq. \eqref{lem_termI_inbasicineqA}, \eqref{eq_IIinbasicineq}, and \eqref{lem_basicineq_step3A}  in Sections \ref{sec:pfD1},  \ref{sec:pfD2}, and  \ref{sec:pfD3}, respectively. 



\bigskip
\begin{lemma} \label{lem_bdd} 
Consider $G \in \mathbb S^n$ satisfying $|G_{ij}| \leqslant M_{Z,1}.$
For any $s>0$,  we have that under the conditions of Theorem \ref{MLE_mainthm},
$$\Pr \left(\max_{1\leqslant i \leqslant n,1\leqslant t\leqslant  T}|\hat\alpha_{it}(G)| \geqslant M_{\alpha}'\right)
 = O((nT)^{-s}), $$ 
 where $M_{\alpha}'$  depends only on $M_{\Theta,1},M_{Z,1},M_{Z,2}$, and $s$. 
Moreover, let   $M_{\Theta,2}' = M_{Z,1} + 2M_{\alpha}' $, and then
$$\Pr \left(\max_{1\leqslant i,j \leqslant n , 1\leqslant  t\leqslant T }|G_{ij}+\hat\alpha_{it}(G)| \geqslant M_{\Theta,2}'\right)
 = O((nT)^{-s}).$$ 
\end{lemma}
\begin{proof}
See Section \ref{sec:pfD4} on Page \pageref{sec:pfD4}.
\end{proof}



\subsubsection{Proof of Eq. \eqref{lem_termI_inbasicineqA}}
\label{sec:pfD1}
Recall that 
$$
(I_b)_{ij} = \sum_t \{\exp(\alpha^{\star}_{it} + \alpha^{\star}_{jt} + G^{\star}_{ij}) -\exp( \hat{\alpha}_{it}(G^{\star}) + \hat{\alpha}_{jt}(G^{\star}) + G^{\star}_{ij} ) \}\{1+1_{(j=i)}\}.
$$
Also recall the definitions of $M_t $ and $M_{2t}$ from \eqref{eq::def_Mtij} and \eqref{eq::def_M2t}, as we will use them frequently in this section.
For simplicity, we define 
$$\Delta M_{t,ij} := \mu_{t,ij}^{\star} \times (\exp(\Delta {\alpha_{it}} + \Delta {\alpha_{jt}}  )-1), \quad\Delta \alpha_{it} := \hat\alpha_{it}(G^{\star}) - \alpha_{it}^{\star}  ,$$ and similarly to $M_{2t}$ define $\Delta M_{2t} \in  \mathbb S^n$ to have diagonal element $2\Delta M_{t,ii} $ and off-diagonal element $\Delta M_{t,ij} $.
Thus $I_b = \sum_t \Delta M_{2t} $.

In order to bound the Frobenius norm of $\Delta M_{2t}$,
we exploit the longitudinal structure of the model. 
In particular, for $1\leqslant i \leqslant j\leqslant n$,  we decompose $\sum_t \Delta M_{t,ij} $ to a form similar to sum of independent components:
\begin{align}
\sum_t \Delta M_{t,ij}
=&~\sum_t \mu_{t,ij}^{\star}\times \big(e^{\Delta \alpha_{it}+\Delta \alpha_{jt}}-1\big) \notag \\
=&~ \sum_t \mu_{t,ij}^{\star}\times \big((e^{\Delta \alpha_{it}}-1)+(e^{\Delta \alpha_{jt}}-1)+(e^{\Delta \alpha_{it}}-1)(e^{\Delta \alpha_{jt}}-1)\big). \label{eq:DeltaM_tij}
\end{align}

Now we aim to bound the Frobenius norm of \eqref{eq:DeltaM_tij}, for which we need to study $\Delta \alpha_{it} = \hat\alpha_{it}(G^{\star}) - \alpha_{it}^{\star}$.
Given $G^{\star}$, we have the estimating equations
\begin{align}
\label{eq::EstimEq_alpha}
0&= \left.\frac{\partial l_t}{\partial \alpha_{it}} \right|_{(G^\star,\hat\alpha_t(G^\star))} =\left(2A_{t,ii} - 2 e^{G_{ii}^{\star} + 2\hat\alpha^{\star}_{it} }\right) + \sum_{j\neq i}\left(A_{t,ij}- e^{G_{ij}^{\star}+\hat{\alpha}^{\star}_{it}+\hat{\alpha}^{\star}_{jt}} \right) \\
&= 2M_{t,ii}- 2e^{G_{ii}^{\star}}\left\{2e^{\alpha^{\star}_{it}}\left(e^{\hat{\alpha}^{\star}_{it}}-e^{\alpha^{\star}_{it}}\right) +\left(e^{\hat{\alpha}^{\star}_{it}}-e^{\alpha^{\star}_{it}}\right)^2\right\} \notag
\\&+
\sum_{j\neq i}\left[M_{t,ij}- e^{G_{ij}^{\star}}\left\{e^{\alpha^{\star}_{it}}\left(e^{\hat{\alpha}^{\star}_{jt}}-e^{\alpha^{\star}_{jt}}\right) + e^{\alpha^{\star}_{jt}}\left(e^{\hat{\alpha}^{\star}_{it}}-e^{\alpha^{\star}_{it}}\right)+\left(e^{\hat{\alpha}^{\star}_{it}}-e^{\alpha^{\star}_{it}}\right) \left(e^{\hat{\alpha}^{\star}_{jt}}-e^{\alpha^{\star}_{jt}}\right)\right\}\right] \notag
\end{align}
where $\hat\alpha^{\star}_{it}$ is an abbreviation for $\hat\alpha_{it}(G^{\star})$, and $M_{t,ij} = A_{t,ij} - e^{ G_{ij}^\star + \alpha_{it}^\star + \alpha_{jt}^\star } $ as is defined in \eqref{eq::def_Mtij}.
Denote 
\begin{equation} \label{eq::def_M2ticdot}
    M_{2t,i\cdot}= 2M_{t,ii}+ \sum_{j\neq i} M_{t,ij}.
\end{equation}
Recall the form of $\ddot{L}_{\alpha_t,\alpha_t} = - \myDonet\Dmut\myDone$ from Remark \ref{rmk::invertable_ddotLalphaalpha} 
(and  \eqref{eq:ddotL_alpha_t,alpha_t}).
By writing out the above equation for each $i=1,\ldots,n$ into a system and rearranging terms, we have
\begin{align*}
\begin{pmatrix}
M_{2t,1\cdot}\\
\vdots \\
M_{2t,n\cdot} 
\end{pmatrix} &= \ddot{L}_{\alpha_t, \alpha_t}
\begin{pmatrix}
e^{\hat{\alpha}^{\star}_{1t}-\alpha_{1t}^{\star}}-1\\
\vdots\\
e^{\hat{\alpha}^{\star}_{nt}-\alpha_{nt}^{\star}}-1 
\end{pmatrix}\\&+
\begin{pmatrix}
2e^{G_{11}^{\star}+2\alpha_{1t}^{\star}}(e^{\hat{\alpha}^{\star}_{1t}-\alpha_{1t}^{\star}}-1)^2+
\sum_{j\neq i} e^{G_{1j}^{\star}+\alpha_{1t}^{\star}+\alpha_{jt}^{\star}}(e^{\hat{\alpha}^{\star}_{1t}-\alpha_{1t}^{\star}}-1)(e^{\hat{\alpha}^{\star}_{jt}-\alpha_{jt}^{\star}}-1)\\
\vdots \\
2e^{G_{nn}^{\star}+2\alpha_{nt}^{\star}}(e^{\hat{\alpha}^{\star}_{nt}-\alpha_{nt}^{\star}}-1)^2+
\sum_{j\neq n} e^{G_{nj}^{\star}+\alpha_{nt}^{\star}+\alpha_{jt}^{\star}}(e^{\hat{\alpha}^{\star}_{nt}-\alpha_{nt}^{\star}}-1)(e^{\hat{\alpha}^{\star}_{jt}-\alpha_{jt}^{\star}}-1)   
\end{pmatrix},
\end{align*}
i.e.,
\begin{equation}
\begin{aligned}
\begin{pmatrix}
e^{\Delta \alpha_{1t}}-1\\
\vdots\\
e^{\Delta \alpha_{nt}}-1  
\end{pmatrix}=
\ddot{L}_{\alpha_t, \alpha_t}^{-1}  \begin{pmatrix}
M_{2t,1\cdot}\\
\vdots \\
M_{2t,n\cdot} 
\end{pmatrix}
-
\ddot{L}_{\alpha_t, \alpha_t}^{-1} \begin{pmatrix}
\mu_{t,11}^\star (e^{\Delta \alpha_{1t}}-1)^2 +
\sum_{j=1}^n\mu_{t,1j}^{\star} (e^{\Delta \alpha_{1t}}-1)(e^{\Delta \alpha_{jt}}-1)\\
\vdots \\
\mu_{t,nn}^\star (e^{\Delta \alpha_{nt}}-1)^2 +
\sum_{j=1}^n\mu_{t,nj}^{\star} (e^{\Delta \alpha_{nt}}-1)(e^{\Delta \alpha_{jt}}-1)  
\end{pmatrix}.
\end{aligned}
\label{eq:Deltaalpha_it}
\end{equation}
Here $\ddot{L}_{\alpha_t, \alpha_t}$ takes value at $(G^{\star}, \alpha^{\star})$.

We characterize the independence across $t$ of $\Delta M_{t,ij}$ by plugging \eqref{eq:Deltaalpha_it} into \eqref{eq:DeltaM_tij}.
Roughly speaking, this will give us 
$$\sum_t \Delta M_{t,ij}= \sum_t \big[ M_t \text{ term} + (e^{\Delta \alpha_{l_1t}}-1)(e^{\Delta \alpha_{l_2t}}-1)  \text{ term}\big],$$ i.e., true score which could concentrate, plus second-order terms. Again, plugging \eqref{eq:Deltaalpha_it} into the second-order terms, we could push the concentration to an even higher order. In the following, we put the above intuition into detail.

We first consider the $\sum_t \mu_{t,ij}^{\star}(e^{\Delta \alpha_{it}}-1)$ term in \eqref{eq:DeltaM_tij}. Recall from \eqref{eq:Wtstar_def} that we let $\mathcal W_t^{\star} := \ddot{L}_{\alpha_t, \alpha_t}^{-1}(G^{\star},\alpha^{\star})$, and denote its $(i,j)$th element by $\omega^{\star}_{t,ij}$. Then by \eqref{eq:Deltaalpha_it} we have
$e^{\Delta \alpha_{it}}-1 =K_{i1t}  - K_{i2t}$ where
\begin{align*}
K_{i1t} = &~\sum_{j=1}^n\omega^{\star}_{t,ij}M_{2t,j\cdot}, \\ 
K_{i2t} = &~\sum_{l=1}^n\omega^{\star}_{t,il}\sum_{j=1}^n(1 + 1_{(j=l)})\mu_{t,lj}^{\star}(e^{\Delta \alpha_{jt}}-1)(e^{\Delta \alpha_{lt}}-1). \notag
\end{align*}
Here $1_{(j=l)} $ denotes the indicator function of $j=l$.
As analyzed above, $K_{i1t}$ is the score term at true parameter value which could concentrate; $K_{i2t}$ is the second order term. Further plugging \eqref{eq:Deltaalpha_it} into $K_{i2t}$, we get
\begin{align*}
K_{i2t}\notag
=&~\sum_{l=1}^n\omega^{\star}_{t,il}\sum_{j=1}^n(1+1_{(j=l)})\mu_{t,lj}^{\star} \\
&~\times\left\{ \sum_{k_1=1}^n\omega^{\star}_{t,jk_1} M_{2t,k_1\cdot} - \sum_{k_1=1}^n \omega^{\star}_{t,jk_1} \sum_{k_2 =1}^n(1+1_{(k_2=k_1)}) \mu_{t,k_1k_2}^{\star}(e^{\Delta \alpha_{k_2t}}-1)(e^{\Delta \alpha_{k_1t}}-1)\right\}\notag\\
&~\times \left\{\sum_{k_3=1}^n\omega^{\star}_{t,lk_3} M_{2t,k_3\cdot} - \sum_{k_3=1}^n \omega^{\star}_{t,lk_3} \sum_{k_4=1}^n(1+1_{(k_4=k_3)}) \mu_{t,k_3k_4}^{\star}(e^{\Delta \alpha_{k_4t}}-1)(e^{\Delta \alpha_{k_3t}}-1)\right\}\notag\\
=:&~  K_{i2t1}- K_{i2t2}- K_{i2t3}+K_{i2t4}.
\end{align*}
In the above decomposition 
\begin{align*}
K_{i2t1} =&~\sum_{l=1}^n\omega^{\star}_{t,il}\sum_{j=1}^n(1+1_{(j=l)})\mu_{t,lj}^{\star}\sum_{k_1=1}^n\omega^{\star}_{t,jk_1} M_{2t,k_1\cdot} \sum_{k_3=1}^n\omega^{\star}_{t,lk_3} M_{2t,k_3\cdot}, \\
K_{i2t2} =&~\sum_{l=1}^n\omega^{\star}_{t,il}\sum_{j=1}^n(1+1_{(j=l)})\mu_{t,lj}^{\star}\sum_{k_1=1}^n\omega^{\star}_{t,jk_1} M_{2t,k_1\cdot}\\&~\times \sum_{k_3=1}^n \omega^{\star}_{t,lk_3} \sum_{k_4=1}^n(1+1_{(k_4=k_3)}) \mu_{t,k_3k_4}^{\star}(e^{\Delta \alpha_{k_4t}}-1)(e^{\Delta \alpha_{k_3t}}-1),\\
K_{i2t3} =&~\sum_{l=1}^n\omega^{\star}_{t,il}\sum_{j=1}^n(1+1_{(j=l)})\mu_{t,lj}^{\star}\sum_{k_3=1}^n\omega^{\star}_{t,lk_3} M_{2t,k_3\cdot} \\&~\times \sum_{k_1=1}^n \omega^{\star}_{t,lk_1} \sum_{k_2=1}^n(1+1_{(k_2=k_1)}) \mu_{t,k_1k_2}^{\star}(e^{\Delta \alpha_{k_2t}}-1)(e^{\Delta \alpha_{k_1t}}-1),\\
K_{i2t4} =&~\sum_{l=1}^n\omega^{\star}_{t,il}\sum_{j=1}^n(1+1_{(j=l)})\mu_{t,lj}^{\star} \sum_{k_1=1}^n \omega^{\star}_{t,lk_1} \sum_{k_2=1}^n(1+1_{(k_2=k_1)}) \mu_{t,k_1k_2}^{\star}(e^{\Delta \alpha_{k_2t}}-1)(e^{\Delta \alpha_{k_1t}}-1)
\\&~\times\sum_{k_3=1}^n \omega^{\star}_{t,lk_3} \sum_{k_4=1}^n(1+1_{(k_4=k_3)}) \mu_{t,k_3k_4}^{\star}(e^{\Delta \alpha_{k_4t}}-1)(e^{\Delta \alpha_{k_3t}}-1).
\end{align*}
We present their respective bounds in the following lemmas.
First, some preliminary bounds are given to facilitate the proof.
\begin{lemma}
\label{lem_hillar_union}
Under the conditions of Theorem \ref{MLE_mainthm}, we have 

 (a) 
 $\max_{i,t} \sum_j |\omega^{\star}_{t,ij}| \leqslant 2\exp(M_{\Theta,2}) / n $;

(b) for any $s > 0$, 
$$\Pr \left(\max_{i,t} |M_{2t,i\cdot}| \geqslant 
C_s\sqrt{ne^{-M_{\Theta,1}} \log(nT) }\right) = O((nT)^{-s})$$
 and  
 $$\Pr \left(\max_{i,t} |e^{\Delta{\alpha_{it}}}-1| \geqslant C_{s}e^{M_{\Theta,2}'-\frac12 M_{\Theta,1}}\sqrt{\log(nT)/ n}\right) = O((nT)^{-s}),$$ where $M_{\Theta,2}'$ is the constant defined in Lemma \ref{lem_bdd}.
\end{lemma}

\begin{proof}
See Section \ref{sec:pfD5} on Page \pageref{sec:pfD5}.
\end{proof}

Then we present bounds for the terms $K_{i1t}, K_{i2t2}, K_{i2t3}$ and $ K_{i2t4} $, respectively.
\begin{lemma}
\label{lem_Ki1t}
Under the conditions of Theorem \ref{MLE_mainthm}, we have for any $s > 0$,
\begin{align*}
&\Pr \left(\max_{i,j'} \Big|\sum_t \mu_{t,ij'}^{\star} K_{i1t}\Big|\geqslant C_{M,s} \sqrt{T \log(nT)/n}\right) = O((nT)^{-s}),\\
&\Pr \left(\max_{i,j'} \Big| \sum_t \mu_{t,ij'}^{\star} K_{i2t1}  \Big|
\geqslant C_{M,s} T\log(nT)/n\right) = O((nT)^{-s}), \\
&\Pr \left(\max_{i,j'} \Big|\sum_t \mu_{t,ij'}^{\star} (K_{i2t2}+K_{i2t3})  \Big|
\geqslant C_{M,s}  T n^{-3/2}\log^{3/2}(nT)\right) = O((nT)^{-s}),\\
&\Pr \left(\max_{i,j'} \Big|\sum_t \mu_{t,ij'}^{\star} K_{i2t4} \Big|
\geqslant C_{M,s} T n^{-2} \log^2(nT)\right) = O((nT)^{-s}).
\end{align*}
 Here $C_{M,s}$ denotes a constant that depends only on $M_{\Theta,1},M_{\Theta,2},M_{Z,1}, s$, and its value could change from line to line.
\end{lemma}
\begin{proof}
See Section \ref{sec:pfD6} on Page \pageref{sec:pfD6}.
\end{proof}

Recall that 
$$
\max_{i,j}
\Big| \sum_{t } \mu_{t,ij}^{\star} (e^{\Delta \alpha_{it}} - 1 )\Big| = \max_{i,j}\Big| \sum_{t } \mu_{t,ij}^{\star} (K_{i1t} - K_{i2t1} + K_{i2t2} + K_{i2t3} - K_{i2t4}) \Big|,
$$
then combining the results in Lemma \ref{lem_Ki1t}, we conclude
\begin{align} \label{I2-1}
\max_{i,j}
\Big| \sum_{t } \mu_{t,ij}^{\star} (e^{\Delta \alpha_{it}} - 1 )\Big| \leqslant C_{M,s} \Big[\sqrt{\frac{T}n\log(nT)}  + \frac{T}{n} \log(nT)\Big]
\end{align}
with probability $1- O((nT)^{-s}) $.

Next we consider the second order term $\sum_{t } \mu_{t,ij}^{\star} (e^{\Delta \alpha_{it}} - 1 ) (e^{\Delta \alpha_{jt}} - 1 ) $ in \eqref{eq:DeltaM_tij}. We plug in the result of Lemma \ref{lem_hillar_union} to get with probability $1-O((nT)^{-s}) $,
\begin{align} \label{I2-2}
\max_{i,j}
\Big|\sum_t \mu_{t,ij}^\star(e^{\Delta \alpha_{it}} - 1 ) (e^{\Delta \alpha_{jt}} - 1 ) \Big|\leqslant&~
T e^{-M_{\Theta,1}} \max_{i,j,t} \Big|(e^{\Delta \alpha_{it}} - 1 ) (e^{\Delta \alpha_{jt}} - 1 )\Big| \notag\\ \leqslant&~ C_{M,s} T\log(nT)/n,
\end{align}
where $C_{M,s}$ denotes a constant that depends only on $M_{\Theta,1},M_{\Theta,2}$, $M_{Z,1}$, and $s$.

Combining \eqref{I2-1} and \eqref{I2-2}
we get 
\begin{align}\label{I2}
\max_{i,j}
\Big| \sum_t \Delta M_{t,ij} \Big| &\leqslant 2\max_{i,j}
\Big| \sum_{t } \mu_{t,ij}^{\star} (e^{\Delta \alpha_{it}} - 1 )\Big| + \max_{i,j}
\Big|\sum_t \mu_{t,ij}^\star(e^{\Delta \alpha_{it}} - 1 ) (e^{\Delta \alpha_{jt}} - 1 ) \Big|\notag \\
&\leqslant C_{M,s} \Big[\sqrt{\frac{T}n\log(nT)}  + \frac{T}{n} \log(nT)\Big]
\end{align}
with probability $1- O((nT)^{-s}) $. As a result, 
$$
\|I_b\|_{\mathrm{F}}^2 \leqslant 4n^2 \cdot \max_{i,j}
\Big( \sum_t \Delta M_{t,ij} \Big)^2\leqslant  C_{M,s} \left[{nT\log(nT)} +{T^2} \log^2(nT)\right]
$$
with probability $1- O((nT)^{-s}) $.
Thus, we have proved \eqref{lem_termI_inbasicineqA}: for any $s>0$, there exists a universal constant $C_s>0$ such that
\begin{align*}
     \|I_b\|_{\mathrm{F}} \leqslant &~ C_{s} \sqrt{nT}\max\Big\{1, \sqrt{\frac{T}{n}}\Big\}\log(nT)=C_s\sqrt{nT}\sqrt{r'_{n,T}} 
\end{align*}with probability $1- O((nT)^{-s}) $.

\subsubsection{Proof for Eq. \eqref{eq_IIinbasicineq}} \label{sec:pfD2}
 By the definition of $II$ 
 and  \eqref{eq:partiallggchitform}, 
\begin{align}\label{eq:iiinequichiform}
    -II= &~ \frac{1}{2}\sum_{t=1}^T(\Delta G_v)^\top \myXt\Big|_{(\tilde G,\hat{\alpha}_t(\tilde G))} (\Delta G_v). 
\end{align} 
Similarly to Remark \ref{rmk::Interpret_myXt}, $\myXt$ depends on $(\tilde G,\hat{\alpha}_t(\tilde G))$ only through $\tilde{\mu}_t:=\exp(\tilde{G}+\hat{\alpha}_t(\tilde G) 1_n^{\top}+1_n \hat{\alpha}_t(\tilde G)^{\top})\in \mathbb{R}^{n\times n}$.
For the simplicity of notation, 
we use $\myXt|_{\tilde{\mu}_t}$ to denote 
$\myXt|_{(\tilde G,\hat{\alpha}_t(\tilde G))}$ below.
Applying Lemma \ref{lem:property_myXt} to  $\myXt|_{\tilde{\mu}_t}$, 
 we obtain
\begin{itemize}
    \item[(X-i)] the null space of $\myXt|_{\tilde{\mu}_t}$  equals the column space of $\myDone$; 
    \item[(X-ii)] the smallest non-zero eigenvalue of $\myXt|_{\tilde{\mu}_t}$ $\geqslant \min_{i,j=1,\ldots,n} \tilde{\mu}_{t,ij} $, where $\tilde{\mu}_{t,ij}$ is the $(i,j)$-th element in $\tilde{\mu}_t$.
\end{itemize}

\bigskip
We next show that $\Delta G_v  \in \text{col}(\myXt|_{\tilde{\mu}_t})$ so that by the fact (X-ii) in the above paragraph, we have
\begin{align}
(\Delta G_v)^\top \myXt|_{\tilde{\mu}_{t}} (\Delta G_v) \geqslant \min_{i,j=1,\ldots,n}\tilde{\mu}_{t,ij} \times \|\Delta G_v\|_2^2 . \label{eq:nonzeroeigenspacesq}
\end{align}
To prove  $\Delta G_v  \in \text{col}(\myXt|_{\tilde{\mu}_t})$, 
by the fact (X-i) in the above paragraph, 
it suffices to show  $(\Delta G)^{\top}\mathcal{D}_{\Theta \alpha}=0$. 
By \eqref{eq:deltagdefinitionA}, \eqref{eq:defgv}, \eqref{eq:dthetaalphadef}, and \eqref{eq:dthetaalphadefexpre}, 
\begin{align*}
(\Delta G_v)^\top \myDone=\begin{pmatrix}
    \sum_{j=1}^n(G_{1j}^{\star}-\hat{G}_{1j}) & \cdots & \sum_{j=1}^n(G_{nj}^{\star}-\hat{G}_{nj})
\end{pmatrix}    = \begin{pmatrix}
0&\cdots &0
\end{pmatrix}
\end{align*}
where the second equation is obtained as  $\sum_{j=1}^n G_{ij}^{\star}=0$ by Condition   \ref{cond:truvalueregularity} and $\sum_{j=1}^n\hat{G}_{ij}=0$ as $\hat{G}$ is solved from the optimization problem \eqref{eq:penalizedMLEproblem}.

Moreover, for $i,j\in \{1,\ldots, n\}$,  we have  $|\tilde G_{ij}| \leqslant M_{Z,1}$, which follows by the definition of $\tilde G$ being a convex combination of $\hat{G}$ and $G^{\star}$,  $|G_{ij}^\star| \leqslant M_{Z,1}$ in Condition  \ref{cond:truvalueregularity}, and $|\hat G_{ij}| \leqslant M_{Z,1}$ in \eqref{eq:penalizedMLEproblem}.  
Therefore, by Lemma \ref{lem_bdd}, we have that with probability $\initerrprob$, $\tilde{\mu}_{t,ij}
\geqslant \tau $, where we define
$$\tau = \exp(-M_{\Theta,2}') = \exp(-M_{Z,1} - 2M_{\alpha}' ).$$
Then by \eqref{eq:iiinequichiform} and \eqref{eq:nonzeroeigenspacesq}, 
$$-II \geqslant \frac{T\tau}{2} \|\Delta G_v \|_2^2 = \frac{T\tau}{4} \|\hat G - G^{\star}\|_{\mathrm{F}}^2$$
with probability at least $\initerrprob$. Thus we have proved \eqref{eq_IIinbasicineq}: for any $s>0$,  there exists a universal constant $C_s>0$ such that
\begin{align*}
    II \leqslant - C_{s} T \|\Delta_G\|_{\mathrm{F}}^2
\end{align*}
with probability $1-O((nT)^{-s})$. 
\subsubsection{Proof of Eq. \eqref{lem_basicineq_step3A}} \label{sec:pfD3}
Recall $G^\star = Z^{\star} (Z^{\star})^\top $.
Define
\begin{align*}
 \mM :=& \left\{ H \in \mathbb R^{n\times n} \left| \text{row}(H) \subset \text{col}(Z^{\star}), \text{col}(H) \subset \text{col}(Z^{\star}) \right.\right\},\\
\bmMp :=& \left\{ H \in \mathbb R^{n\times n} \left| \text{row}(H) \subset \text{col}(Z^{\star})^\perp, \text{col}(H) \subset \text{col}(Z^{\star})^\perp \right.\right\}.
\end{align*}
Also define $\Delta_{G,\bmM}$ the projection of $\Delta_G\in\mathbb R^{n\times n}$ into the subspace $\bmM$, in which the projection is defined in terms of the matrix trace inner product $\langle A,B\rangle:=\text{trace}(A^\top B) $, and define the projections $G^{\star}_{\mM}, G^{\star}_{\mM^{\perp}},\Delta_{G,\bmM^{\perp}} $ similarly. 
Then
we can decompose $G^\star$ and $\Delta_G$ as 
\begin{align}
G^{\star} =  G^{\star}_{\mM} + G^{\star}_{\mM^{\perp}}, \quad \Delta_G = \Delta_{G,\bmM} + \Delta_{G,\bmM^{\perp}},\notag
\end{align}
so that
\begin{align}
\|\hat G\|_* - \|G^{\star}\|_* &= \|G^{\star} + \Delta_G \|_*  - \|G^{\star}\|_*\notag  \\
&\geqslant \|G^{\star}_{\mM} + \Delta_{G,\bmMp} \|_* -  \|G^{\star}_{\mM^\perp}\|_* - \|\Delta_{G,\bmM} \|_* - \|G^{\star}_{\mM}\|_* - \|G^{\star}_{\mM^\perp}\|_*\notag \\
&=  \|G^{\star}_{\mM}\|_* + \|\Delta_{G,\bmMp} \|_* -  \|G^{\star}_{\mM^\perp}\|_* - \|\Delta_{G,\bmM} \|_* - \|G^{\star}_{\mM}\|_* - \|G^{\star}_{\mM^\perp}\|_*\notag \\
&=  \|\Delta_{G,\bmMp} \|_* -  2\|G^{\star}_{\mM^\perp}\|_* - \|\Delta_{G,\bmM} \|_*. \label{G_decomp_lowrank}
\end{align}
When choosing $\lambda_{n,T} \geqslant 2 \|I_a\|_{\operatorname{op}} $, and substituting \eqref{G_decomp_lowrank} into \eqref{basic_ineq3A_original}, we get
\begin{align*}
\lambda_{n,T} \big(\|\Delta_{G,\bmMp} \|_* -  2\|G^{\star}_{\mM^\perp}\|_* - \|\Delta_{G,\bmM} \|_*\big) 
&\leqslant \frac{\lambda_{n,T}}2 \big(\|\Delta_{G,\bmMp} \|_* + \|\Delta_{G,\bmM} \|_* \big) 
+ \|I_b\|_{\mathrm{F}} \| \Delta_G\|_{\mathrm{F}}  \\
\|\Delta_{G,\bmMp} \|_* - 3\|\Delta_{G,\bmM} \|_* - 4\|G^{\star}_{\mM^\perp}\|_* &\leqslant \frac{2}{\lambda_{n,T}} \|I_b\|_{\mathrm{F}} \| \Delta_G\|_{\mathrm{F}}.
\end{align*}
Note that as $G^\star= Z^\star (Z^\star)^\top $, $G^{\star}_{\mM^\perp} = 0 $. Thus,
\begin{align*}
\|\Delta_G\|_* &\leqslant \|\Delta_{G,\bmMp}\|_* + \|\Delta_{G,\bmM}\|_* \\
&\leqslant 4\|\Delta_{G,\bmM}\|_* + 4\|G^{\star}_{\mM^\perp}\|_*
+\frac{2}{\lambda_{n,T}} \|I_b\|_{\mathrm{F}} \| \Delta_G\|_{\mathrm{F}} \\
&\leqslant 4\sqrt{2k} \|\Delta_{G,\bmM}\|_{\mathrm{F}} +\frac{2}{\lambda_{n,T}} \|I_b\|_{\mathrm{F}} \| \Delta_G\|_{\mathrm{F}}
\end{align*}
where the third line is because any matrix in $\bmM$ has rank at most $2k$. Using $\|\Delta_{G,\bmM}\|_{\mathrm{F}}\leqslant \|\Delta_{G}\|_{\mathrm{F}} $ gives the stated
inequality of \eqref{basic_ineq3A}.

\subsection{Proof of Lemmas  \ref{lem_bdd}--\ref{lem_Ki1t} in Section \ref{sec:lem_of_thm}} \label{subsec:lems_lem}

\subsubsection{Proof of Lemma \ref{lem_bdd}}
\label{sec:pfD4}
By the Bernstein inequality for Poisson random variables in Lemma \ref{lem:bern_pois}, we have for $d_{2t,i\cdot} := 2A_{t,ii}+ \sum_{j\neq i} A_{t,ij} $,
\begin{align}
\label{eq:bern_d}
\Pr\left(\Big|d_{2t,i\cdot} - \Big(2\mu_{t,ii}^{\star}+ \sum_{j\neq i} \mu_{t,ij}^{\star}\Big)\Big| \geqslant x \right) \leqslant 2\exp\left( -\frac{x^2}{ 2(4\mu_{t,ii}^{\star}+ \sum_{j\neq i} \mu_{t,ij}^{\star} + 10x ) } \right).
\end{align}
Setting $x = 2(s+1)\sqrt{(n+3)\exp(-M_{\Theta,1}) \log(nT) }$ for some constant $s>0$ and by a union bound, we derive that with  probability at least $1 - O((nT)^{-s_1} ) $
 there exist two 
 constants
 $0<c_2<c_1$ s.t.
\begin{align}\label{dit_bd}
c_2(n+1)\leqslant  d_{2t,i\cdot}\leqslant c_1(n+1)
\end{align}
holds for any $i=1,\ldots,n;t=1,\ldots,T$.

In the proof of this lemma, we use $\hat\alpha$ to denote $\hat\alpha(G) $.
By the definition of $\hat\alpha$, it satisfies the score equations given by
\begin{align}
\label{hatalphascore}
\frac{\partial l(G,\hat\alpha) }{\partial \alpha_{it}} = 2A_{t,ii} - 2\exp(2\hat\alpha_{it}+G_{ii}) +\sum_{j\neq i} [A_{t,ij} -\exp(\hat\alpha_{it}+\hat\alpha_{jt}+G_{ij})]=0, \forall 1\leqslant i\leqslant n.
\end{align}

Now we fix at a $t$ and show that the event \eqref{dit_bd} implies boundedness of $\hat\alpha(G)$ for $G$ bounded in the sense that there exists $M_{Z,1} > 0 $ s.t. $|G_{ij}| \leqslant M_{Z,1} $. For simplicity, in the following of the proof, we fix at a $t$, omit the $t$, and write $d_i = d_{2t,i\cdot}, \alpha_i = \alpha_{it} $.

\paragraph*{Step 1.} Prove when $\hat{\alpha}_{\max}$ is unbounded there are at least $\tilde c\times n$ coefficients $\hat{\alpha}_j$ $\geqslant \hat{\alpha}_{\max}/2$ for some $\tilde c\in (0,1) $.

\paragraph*{Step 1.1} First prove that there exists $i$ s.t. $\exp(\hat{\alpha}_i)\leqslant c\times \exp(-\hat{\alpha}_{\max})$ for some $c$. 
Let $m$ be one index that achieves the maximum $\hat{\alpha}_m=\hat{\alpha}_{\max}$. 
Suppose if all $\exp(\hat{\alpha}_j) > c\exp(-\hat{\alpha}_{\max})$ for $j\neq m$, and then from the score equation of $\hat\alpha(G) $ \eqref{hatalphascore}
\begin{align}
  d_m &=\sum_{j}(1+1_{(j=m)}) \exp(\hat{\alpha}_{\max} + \hat{\alpha}_j + G_{mj}) \notag\\
  &>   \sum_{j}(1+1_{(j=m)})c \times \exp(\hat{\alpha}_{\max} - \hat{\alpha}_{\max} + G_{mj})\notag\\
   &= c\sum_{j}(1+1_{(j=m)})\exp(G_{mj}) \geqslant c_1 (n+1)
   \label{eq:bound1}
\end{align}
where we choose
\begin{align}
  c \geqslant c_1\exp(M_{Z,1}). \label{cc1}
\end{align}
Then \eqref{eq:bound1} contradicts with $d_i\leqslant c_1(n+1)$. 

\paragraph*{Step 1.2} Assume there exist $l$ number of coefficients $\hat{\alpha}_j$ s.t. $\hat{\alpha}_j< \hat{\alpha}_{\max}/2$. We next prove $l\leqslant (1-\tilde{c})\times n$ for some constant $\tilde{c}\in(0,1)$. Choose the $i$ obtained in Step 1.1, we have $\exp(\hat{\alpha}_{i})<c \exp(-\hat{\alpha}_{\max})$, and hence
\begin{align*}
d_i=&~\sum_{j}(1+1_{(j=i)})\exp(\hat{\alpha}_{i} + \hat{\alpha}_j + G_{ij}) < c\times \sum_{j}(1+1_{(j=i)})\exp(-\hat{\alpha}_{\max}+ \hat{\alpha}_j + G_{ij}) \notag\\
< &~ c\times\exp(-\hat{\alpha}_{\max})\big\{ l\times \exp(\hat{\alpha}_{\max}/2) + (n+1-l)\exp(\hat{\alpha}_{\max}) \big\} \times \exp(M_{Z,1}) \notag\\
\leqslant &~ c \times \exp(M_{Z,1})\times \big\{l \times \exp(-\hat{\alpha}_{\max}/2)  + (n+1-l) \big\}.
\end{align*} 
By $d_i\geqslant c_2 n$, we have
\begin{align*}
  \big\{1-\exp(-\hat{\alpha}_{\max}/2)\big\}\times l\leqslant n+1-\frac{c_2n}{c\times \exp(M_{Z,1})}
\end{align*}
If $\hat{\alpha}_{\max}$ is unbounded, when it  is sufficiently large, $1-\exp(-\hat{\alpha}_{\max}/2)>0$, and then
\begin{align*}
  l\leqslant \frac{1}{1-\exp(-\hat{\alpha}_{\max}/2)}\times \biggr(\frac{n+1}{n}-\frac{c_2}{c\times \exp(M_{Z,1})}\biggr)n :\leqslant (1-\tilde{c})\times  n
\end{align*}
where from \eqref{cc1} we know that $c\exp(M_{Z,1})\geqslant c\exp(-M_{Z,1})\geqslant c_1>c_2 $.
This shows that if $\hat{\alpha}_{\max}$ is unbounded, (can take sufficiently large number{, e.g., take $c = c_1 \exp(M_{Z,1}), \tilde c =c_2 / (2c_1\exp(2M_{Z,1})) $ and $\hat\alpha_{\max} \geqslant 2 \log(2(1-\tilde c)/\tilde{c} ) $}), number of coefficients $\hat{\alpha}_i\geqslant\hat{\alpha}_{\max}/2$ is at least $\tilde{c}\times n$ for a 
constant
$\tilde c \in (0,1) $.

\paragraph*{Step 2.} By conclusion in Step 1, we next prove that the unboundedness of $\hat{\alpha}_{\max}$ will lead to a contradiction. 
Take an index set $\mathcal{N}$ that contains indexes $j$ s.t. $\hat{\alpha}_j\geqslant\hat{\alpha}_{\max}/2$,
where by Step 1, $|\mathcal{N}|\geqslant \tilde{c}\times n$ for some $\tilde{c}<1$. 
\begin{align*}
\tilde{c}\times n \times \exp(3\hat{\alpha}_{\max}/2)\times \exp(-M_{Z,1}) \leqslant &~|\mathcal{N}|\exp(3\hat{\alpha}_{\max}/2)\times \exp(-M_{Z,1})\notag\\
  \leqslant &~ \sum_{j\in \mathcal{N}}\exp(\hat{\alpha}_j+\hat{\alpha}_{m}+G_{mj})\notag\\
\leqslant &~ \sum_{j} \exp(\hat{\alpha}_j+\hat{\alpha}_{m}+G_{mj}) = d_{m}\leqslant c_1 \times (n+1).
\end{align*}
Comparing the first row and last row in the above inequalities, we note that they are contradictory  
when $\hat{\alpha}_{\max}$ is taken sufficiently large (depending only on $c_1,\tilde c$ and $M_{Z,1}$).

\paragraph*{Step 3.} Show that $\hat\alpha_{\min} $ is also bounded. Consider index $m'$ such that $\hat\alpha_{m'} = \hat\alpha_{\min} $. Then
\begin{align*}
c_2 (n+1 )< d_{m'} \leqslant (n+1)\exp(\hat\alpha_{\max}+\hat\alpha_{\min}+M_{Z,1}).
\end{align*}
Hence 
\begin{align*}
\hat\alpha_{\min} > \log(c_2) - M_{Z,1} - \hat\alpha_{\max},
\end{align*}
i.e. $\hat\alpha_{\min}$ is bounded from below (depending only on $c_1,c_2,\tilde c$ and $M_{Z,1}$). 

\subsubsection{Proof of Lemma \ref{lem_hillar_union}} \label{sec:pfD5}
The result of part (a) is already derived in \eqref{eq:myAt_star_inv_inftynorm}:
Lemma \ref{lem:thm1.1_Hillar} shows that as  $\mu_{t,ij}^{\star} \geqslant \exp(-M_{\Theta,2}) $, we have 
$$\max_{i,t} \sum_j |\omega^{\star}_{t,ij}| \leqslant 2(n\exp(-M_{\Theta,2}))^{-1} =2\exp(M_{\Theta,2}) / n. $$

Now we prove part (b). We can rewrite the Estimating Equation \eqref{eq::EstimEq_alpha} with respect to $\alpha$ given $G^\star$ into 
\begin{align*}
0&= \left.\frac{\partial l_t}{\partial \alpha_{it}} \right|_{(G^\star,\hat\alpha_t(G^\star))} =\left(2A_{t,ii} - 2 e^{G_{ii}^{\star} + 2\hat\alpha^{\star}_{it} }\right) + \sum_{j\neq i}\left(A_{t,ij}- e^{G_{ij}^{\star}+\hat{\alpha}^{\star}_{it}+\hat{\alpha}^{\star}_{jt}} \right) \\
&= M_{2t,i\cdot} - \Big[ 2\Big( e^{G_{ii}^{\star} + 2\hat\alpha^{\star}_{it} } -  e^{G_{ii}^{\star} + 2\alpha^{\star}_{it} }\Big) + \sum_{j\neq i} \Big( e^{G_{ij}^{\star} + \hat\alpha^{\star}_{it} + \hat\alpha^{\star}_{jt} } -  e^{G_{ij}^{\star} + \alpha^{\star}_{it}  + \alpha^{\star}_{jt}}\Big) \Big]
.
\end{align*}
Recall we defined $\Delta\alpha_{it} = \hat\alpha_{it}^\star -  \alpha_{it}^\star $.
By applying a mean value theorem to
the above equation between the vectors 
$(\hat\alpha_{1t}^\star,\ldots,\hat\alpha_{nt}^\star )^\top$ and $(\alpha_{1t}^\star,\ldots,\alpha_{nt}^\star )^\top$,
we have
\begin{align*} 
 M_{2t,i\cdot} =   \Big( 4e^{G_{ii}^{\star} + 2{\tilde\alpha_{it}}^{(it)} } + \sum_{j\neq i} e^{G_{ij}^{\star} + {\tilde\alpha_{it}}^{(it)}  + {\tilde\alpha_{jt}}^{(it)} } \Big)   \Delta\alpha_{it} 
 + \sum_{j\neq i} \Big( e^{G_{ij}^{\star} + \tilde\alpha^{(it)}_{it}  + \tilde\alpha^{(it)}_{jt}} \Big)  \Delta\alpha_{jt} 
\end{align*}
where for $j=1,\ldots,n$, $\tilde\alpha_{jt}^{(it)} $ is a midpoint between $\hat\alpha_{jt}^\star$ and $\alpha_{jt}^\star $ (since this midpoint is used in the expansion of the score equation with respect to $\alpha_{it} $, the midpoint depends on $i,t$, and hence the superscript $^{(it)}$).
By writing out the above equation for each $i=1,\ldots,n$ into a system and rearranging terms, we have
\begin{equation}
\label{eq::Yt_Deltaalpha_M2}
\tilde{\mathcal{Y} }_t
\begin{pmatrix} 
\Delta \alpha_{1t}\\
\vdots\\
\Delta \alpha_{nt}
\end{pmatrix}=
\begin{pmatrix}
M_{2t,1\cdot} \\
\vdots \\
M_{2t,n\cdot}  
\end{pmatrix}
\end{equation}
where $\tilde{\mathcal{Y} }_t $ is an $n\times n$ matrix whose $(i,i)$-th entry and $(i,j)$-th entry for $j\neq i$ are  $$ \tilde y_{t,ii}= \Big( 4e^{G_{ii}^{\star} + 2{\tilde\alpha_{it}}^{(it)} } + \sum_{j\neq i} e^{G_{ij}^{\star} + {\tilde\alpha_{it}}^{(it)}  + {\tilde\alpha_{jt}}^{(it)} } \Big) \quad\text{ and }\quad \tilde y_{t,ij}=\Big( e^{G_{ij}^{\star} + \tilde\alpha^{(it)}_{it}  + \tilde\alpha^{(it)}_{jt}} \Big),$$ respectively.
By Lemma \ref{lem_bdd}, we have $|\hat\alpha_{it}^\star|=|\hat\alpha_{it}(G^\star)| < M_{\alpha}' $ with probability $1-O((nT)^{-s}) $. Since $\tilde\alpha_{jt}^{(it)}$ is a midpoint between $\hat\alpha_{jt}^\star $ and  $\alpha_{jt}^\star $, we have $|\tilde\alpha_{jt}^{(it)}| < M_{\alpha}' $, and hence $\tilde y_{t,ij} > \exp(-M_{\Theta,2}') $ with probability $1-O((nT)^{-s}) $ for any $1\leqslant i,j\leqslant n$. Denote $\tilde{\mathcal{W}}_t = (\tilde{\mathcal{Y}}_t)^{-1} \in \mathbb R^{n\times n} $ with $(i,j)$-th entry $\tilde w_{t,ij} $.
Applying Lemma \ref{lem:thm1.1_Hillar} on $\tilde{\mathcal{Y}}_t$, we have $$\max_{i,t} \sum_{j=1}^n |\tilde w_{t,ij} | \leqslant 2\exp(M_{\Theta,2}') / n $$ 
with probability $1-O((nT)^{-s}) $.
From \eqref{eq::Yt_Deltaalpha_M2} we have $\Delta\alpha_{it} = \sum_{j=1}^n \tilde w_{t,ij} M_{2t,j\cdot} $;
combining with the Bernstein bound 
$$\max_{i,t} |M_{2t,i\cdot}| \leqslant 2 (s+1)\sqrt{(n+3)e^{-M_{\Theta,1}} \log(nT) }$$
 with probability at least $1 - O((nT)^{-s} )$, which is derived similarly to \eqref{eq:bern_d}, we obtain 
 \begin{align} \max_{i,t}
      |\Delta\alpha_{it}| =&~\max_{i,t} |\sum_{j=1}^n \tilde w_{t,ij} M_{2t,j\cdot}| \leqslant \left(\max_{l,t}|M_{2t,l\cdot}| \right)\left( \max_{l,t} \sum_{j=1}^n |\tilde w_{t,lj}| \right) \notag \\
      \leqslant&~ C_{s}e^{M_{\Theta,2}'-\frac12 M_{\Theta,1}}\sqrt{\log(nT)/ n} \label{eq::abs_Deltaalpha}
 \end{align}
 with probability at least $1 - O((nT)^{-s} )$ for some constant $C_s$ ($C_s$ denotes a constant whose value only depends on $s$ but could change from line to line). By the inequality $|e^x - 1| \leqslant |x| + x^2 $ for $x\in(-\infty,1] $, under the condition that $n/\log(T) $ is sufficiently large, we have that \eqref{eq::abs_Deltaalpha} also hold for $\max_{i,t} |e^{\Delta\alpha_{it}}-1| $ for some constant $C_s$, which proves part (b) of Lemma \ref{lem_hillar_union}.

\subsubsection{Proof of Lemma \ref{lem_Ki1t}}
\label{sec:pfD6}
We first consider the first term $\sum_t \mu_{t,ij'}^\star K_{i1t} $.
Note that $\{M_{2t,j\cdot}\}_{j=1,\ldots,n;t=1,\ldots,T}$ are independent centered Poisson 
variables. 
Denote $\mu_{2t,j\cdot} = 4\mu_{t,jj} + \sum_{l\neq j} \mu_{t,jl} =$ Var$(M_{2t,j\cdot})$.
By the Bernstein inequality for Poisson variables in 
Lemma \ref{lem:bern_pois}, we have
\begin{align*}
&\Pr\left( \Big|\sum_{t=1}^T \mu_{t,ij'}^{\star} \sum_{j=1}^n \omega^{\star}_{t,ij} M_{2t,j\cdot}\Big|  \geqslant x \right) \\
\leqslant& \exp\left( -\frac{x^2}{ 2(\sum_t \sum_j (\mu_{t,ij'}^{\star})^2 (\omega^{\star}_{t,ij})^2 \mu_{2t,j\cdot} + 20x\exp(M_{\Theta,2})/n )  } \right)\\
 \leqslant& \exp\left( -\frac{x^2}{ 2(T (2\exp(M_{\Theta,2})/n)^2 (n+3)\exp(-3M_{\Theta,1}) + 20x\exp(M_{\Theta,2})/n )  } \right),
\end{align*}
i.e., with probability at least $ 1 - (nT)^{-s}$, 
$$\max_{i,j'} \Big|\sum_t \mu_{t,ij'}^{\star} K_{i1t1}\Big|\leqslant C_{s} \sqrt{Te^{2M_{\Theta,2}-3M_{\Theta,1}} \log(nT)/n}.$$
Here $C_{s}$ denotes a constant that depends only on $s$ and its value could change from line to line in the proof.

Next we bound $K_{i2t1} $.
Using the bounds in Lemma \ref{lem_hillar_union}, we have
with probability $1-O((nT)^{-s-1} ) $,
\begin{align*}
\max_{t,i} |K_{i2t1}|
\leqslant&~\sum_{l=1}^n|\omega^{\star}_{t,il}|\sum_{j=1}^n|2\mu_{t,lj}^{\star}|\sum_{k_1=1}^n |\omega^{\star}_{t,jk_1}| (\max_{j',t} |M_{2t,j'\cdot}| ) \sum_{k_3=1}^n |\omega^{\star}_{t,lk_3}| (\max_{j',t} |M_{2t,j'\cdot}| )   \notag\\
\leqslant&~ C_{s} \frac{e^{M_{\Theta,2}}}{n} \cdot ne^{-M_{\Theta,1}} \cdot  \frac{e^{M_{\Theta,2}}}{n} \sqrt{n e^{-M_{\Theta,1}}\log(nT)} \cdot \frac{e^{M_{\Theta,2}}}{n} \sqrt{n e^{-M_{\Theta,1}}\log(nT)}  \\
  \leqslant&~ C_{s} \exp(3M_{\Theta,2}-2 M_{\Theta,1} ) \times \log(nT)/n.
\end{align*}
Thus, we have 
$$\max_{i,j'} \Big|\sum_t \mu_{t,ij'}^{\star} K_{i2t1}  \Big| \leqslant C_{s}  \exp(3M_{\Theta,2}-3 M_{\Theta,1} ) \times T \log(nT)/n$$
with probability $1-O((nT)^{-s}) $. 

Next we bound $K_{i2t2}$ and $K_{i2t3}$.
Using the bounds in Lemma \ref{lem_hillar_union}, we have
with probability $1-O((nT)^{-s-1} ) $,
\begin{align*}
\max_{t,i} |K_{i2t2}|
\leqslant&~\sum_{l=1}^n|\omega^{\star}_{t,il}|\sum_{j=1}^n|2\mu_{t,lj}^{\star}|\sum_{k_1=1}^n |\omega^{\star}_{t,jk_1}| (\max_{j',t} |M_{2t,j'\cdot}| ) \\
&\quad \sum_{k_3=1}^n |\omega^{\star}_{t,lk_3}|\sum_{k_4=1}^n|2\mu_{t,k_3k_4}^{\star}|  (\max_{j',t} |\exp(\Delta_{\alpha_{j't}})-1|)^2 \notag\\
\leqslant&~ C_{s} \frac{e^{M_{\Theta,2}}}{n} \cdot ne^{-M_{\Theta,1}} \cdot  \frac{e^{M_{\Theta,2}}}{n} \sqrt{n e^{-M_{\Theta,1}}\log(nT)} \\
&\quad\cdot \frac{e^{M_{\Theta,2}}}{n} \cdot ne^{-M_{\Theta,1}} \cdot \left(e^{M_{\Theta,2}' - \frac12 M_{\Theta,1} }\sqrt{\log(nT)/n}\right)^2 \\
  \leqslant&~ C_{s} \exp(2M_{\Theta,2}'+3M_{\Theta,2}-3.5 M_{\Theta,1} ) \times n^{-3/2}\log^{3/2}(nT).
\end{align*}
Thus, 
$$\max_{i,j'} \Big|\sum_t \mu_{t,ij'}^{\star} K_{i2t2}  \Big|
\leqslant C_{s}  \exp(2M_{\Theta,2}'+3M_{\Theta,2}-4.5 M_{\Theta,1} ) \times T n^{-3/2}\log^{3/2}(nT)$$ with probability $1-O((nT)^{-s}) $. 
Same conclusion holds for $\left|\sum_t \mu_{t,ij'}^{\star} K_{i2t3}  \right|$.

Finally, $K_{i2t4}$ can be bounded in a similar fashion by
\begin{align*}
\max_{t,i}
|K_{i2t4}| \leqslant C_{s} \exp(4M_{\Theta,2}'+3M_{\Theta,2}-5M_{\Theta,1}) \times n^{-2}\log^2(nT)
\end{align*}
 with probability $1-O((nT)^{-s-1}) $,
and hence $$\max_{i,j'} \Big|\sum_t \mu_{t,ij'}^{\star} K_{i2t4} \Big|
\leqslant C_s \exp(4M_{\Theta,2}'+3M_{\Theta,2}-6M_{\Theta,1}) T n^{-2} \log^2(nT)$$  with probability $1-O((nT)^{-s}) $.

\newpage
\section{Technical Results for Section \ref{subsec:manifold} on Geometric Interpretation}

\subsection{Background} \label{sec:backg}
In this section, we give a brief
review on how a Newton step on a smooth quotient manifold is specified in the literature. Recall that for an ordinary smooth multivariate function $g: \mathbb{R}^m\to \mathbb{R}$, the Newton update at a point $x\in \mathbb{R}^m$ is given by
$ x+ \nu$, 
where $\nu \in \mathbb{R}^m$ is solved from $ \operatorname{Hess} {g} (x) [\nu] = -\operatorname{Grad} {g}(x)$, and  $\operatorname{Grad} g(x) \in \mathbb{R}^{m}$ and $\operatorname{Hess} g(x) \in \mathbb{R}^{m\times m}$ denote the ordinary gradient and Hessian of $g$ at $x\in \mathbb{R}^m$. Such an idea has been generalized to smooth functions defined on manifolds 
\citep{absil2009optimization}.

Consider an abstract Riemannian manifold $\mathcal{M}$  equipped with an affine connection and a retraction $R$ \citep{boumal2020introduction}. For a smooth function $f: \mathcal{M}\to \mathbb{R}$, 
the Newton update at a point $q \in \mathcal{M}$  is given by $R_q(\nu)$,
where $\nu$ is a tangent vector specified by
\begin{align}\label{eq:netwonhessgrad}
    \operatorname{Hess}f(q) [\nu] = -\operatorname{Grad}f(q), 
\end{align} 
and 
$ \operatorname{Hess}f(q)$ and $\operatorname{Grad}f(q) $ denote the Hessian and gradient of $f$ at $q\in \mathcal{M}$.
At the point $q$, the retraction $R_q$ is a mapping from the tangent space to $\mathcal{M}$ at $q$, denoted as $T_q\mathcal{M}$, to the manifold $\mathcal{M}$.  
Conceptually,  for any $\xi \in T_q\mathcal{M}$,  $R_q(\xi)$ generalizes the notion of moving in the direction of $\xi$, while staying on the manifold. 

However, tangent vectors are defined abstractly, and it is difficult to provide an explicit expression for the Newton update $R_q(\nu)$. To solve this problem, a concrete horizontal space is proposed in the literature to represent the abstract tangent space in the quotient manifold scenario. For generality, we consider a quotient manifold $\mathcal{M}=\bar{\mathcal{M}} / \sim$ for any unspecified Euclidean space $\bar{\mathcal{M}}$ and equivalence relation $\sim$. Let $\pi :\bar{\mathcal{M}}\to {\mathcal{M}} $ denote the canonical projection and let $\bar{q} \in \pi^{-1}(q)$ be a representation of $q$ in $\bar{\mathcal{M}}$. Given any $\bar{q} \in \bar{\mathcal{M}}$, let $\mathcal{H}_{\bar{q}}$ denote a subspace of $T_{\bar{q}} \bar{\mathcal{M}}$ satisfying 
$$
\mathcal{H}_{\bar{q}} \oplus \mathcal{V}_{\bar{q}}= T_{\bar{q}}\bar{\mathcal{M}} \quad \text{where} \quad
\mathcal{V}_{\bar{q}}=T_{\bar{q}}\left(\pi^{-1}(q)\right) \subset T_{\bar{q}}\bar{\mathcal{M}}.
$$
The constructed $\mathcal{H}_{\bar{q}}$ satisfies that each tangent vector $\xi \in T_q \mathcal{M}$ has a unique counterpart  $\bar{\xi}_{\bar{q}} \in \mathcal{H}_{\bar{q}}$ with the same differential properties as $\xi$. The Euclidean subspace $\mathcal{H}_{\bar{q}}$ is called the horizontal space at ${\bar{q}}$, and the unique counterpart is called the horizontal lift of $\xi$ at
$\bar{q}$.

Furthermore, we can construct a natural retraction based on the previously mentioned horizontal lift. Suppose that for all $q \in \mathcal{M}$ and $\xi \in T_q \mathcal{M}$, 
$$
\pi\left({\bar{q}_1} + \bar{\xi}_{\bar{q}_1}\right)=\pi\left({\bar{q}_2} + \bar{\xi}_{\bar{q}_2}\right)
$$
for all ${\bar{q}_1}, {\bar{q}_2}\in \pi^{-1}(q)$, then $R_q(\xi) := \pi({\bar{q}} + \bar{\xi}_{\bar{q}})$ defines a retraction on $\mathcal{M}$. Equipped with such natural retraction, we can obtain the explicit form of $R_q(\nu)$ by simply specifying $\bar{\nu}_{\bar{q}}$.
Fortunately, when $\mathcal{M}$ is endowed with the Riemannian metric naturally induced by $\bar{\mathcal{M}}$ and the corresponding Riemannian connection, $\bar{\nu}_{\bar{q}}$ can be solved from
\begin{align*} 
    \mathcal{P}_{\mathcal{H}_{\bar{q}}} \operatorname{Hess} \bar{f}(\bar{q})\left[\bar{\nu}_{\bar{q}}\right] =-\mathcal{P}_{\mathcal{H}_{\bar{q}}}\operatorname{Grad} \bar{f}(\bar{q}),
\end{align*}
where $\bar{f} = f \circ \pi$ denotes the smooth Euclidean function induced by $f$, $\operatorname{Grad} \bar{f}(\bar{q})$ and $\operatorname{Hess} \bar{f}(\bar{q})$ denote the ordinary gradient and Hessian of $\bar{f}$ at $\bar{q}$, and $\mathcal{P}_{\mathcal{H}_{\bar{q}}}$ denotes the Euclidean projection operator to $\mathcal{H}_{\bar{q}}$.

\subsection{Proof of Proposition \ref{prop:newtonstepform}}
\label{sec:pfProp}
We first show that the  intrinsic parameter space $\mathbb{R}_0^{n\times k}/\sim$ has a natural smooth manifold structure, and the induced log profile likelihood $pL_Q$ is a smooth function on $\mathbb{R}_0^{n\times k}/\sim$.

\begin{lemma}\label{lemma:Quo_mani} 
There exists a differential structure such that 
% ...
the quotient set $\mathbb{R}_0^{n\times k}/\sim$, equipped with such differential structure, is a smooth quotient manifold of dimension $nk-k(k+1)/2$.  The induced function $pL_Q:\mathbb{R}_0^{n\times k}/\sim \ \to \mathbb R$ defined by $pL_Q([Z]) = pL(Z) $ is a smooth function on $\mathbb{R}_0^{n\times k}/\sim$. 
\end{lemma}
\begin{proof}
See Section \ref{sec:pfQuo_mani} on Page \pageref{sec:pfQuo_mani}. 
\end{proof}

\medskip

For the simplicity of notation, we let $\mathcal{M}$ denote $\mathbb{R}_0^{n\times k}/\sim$ in the following. Let \begin{align} \label{eq:hzdef_prop3}
    \mathcal{H}_Z := \operatorname{col}(\mathcal{U}_Z) =  \operatorname{col}(\Ieff(Z,\alpha)) \subset \mathbb{R}^{nk}
\end{align}
where the second equation follows by Lemma \ref{lm:presentuzmain}-(i).   
Lemma \ref{lem_HspaceDiff} below proves that each tangent vector $\xi\in T_{[Z]}\mathcal{M}$ has a unique counterpart $\bar{\xi}_Z$ in $\mathcal{H}_Z$
that preserve the differential properties. 

\begin{lemma}
\label{lem_HspaceDiff}
Let $\bar{\mathcal{F}}:= \{ \bar{f} \in C^{\infty}(\mathbb{R}^{nk}) :  \bar f(Z_v) =\bar f((ZQ)_v), \forall Z\in \mathbb{R}^{n\times k}, Q\in \mathcal{O}(k)\} $, 
where $C^{\infty}(\mathbb{R}^{nk})$ represents the class of smooth functions defined on the domain $\mathbb R^{nk}$.
Each $\bar f\in\bar{\mathcal{F}}$ naturally induces a smooth function $f:\mathcal{M} \to \mathbb{R} $ defined by $f([Z]) := \bar f(Z_v)$.
Then $\mathcal{H}_Z$
satisfies the following property:

Given any tangent vector $\xi \in T_{[Z]} \mathcal{M}$, there exists one and only one vector $\bar \xi_Z \in \mathcal{H}_Z$ such that for any $ \bar{f}\in \bar{\mathcal{F}}$,
$$
\xi f = \bar \xi_Z \bar f,
$$
where $\xi f$ denotes the directional derivative of $f$ at $[Z]$ along the direction $\xi$, and $\bar\xi_Z\bar f$ denotes the directional derivative of $\bar f$ at $Z_v$ along  $\bar \xi_Z$.
\end{lemma}
\begin{proof}
See Section \ref{sec:pf_Hspace} on Page \pageref{sec:pf_Hspace}.
\end{proof}

\medskip
In Lemma \ref{lem:expmap} below, we prove $R_{[Z]}(\xi) = [\operatorname{mat}(Z_v + \bar{\xi}_Z)]$ for any $\xi \in T_{[Z]}\mathcal{M}$. 

\begin{lemma}\label{lem:expmap}
Define $\tilde{R}_{[Z]}: T_{[Z]}\mathcal{M}\to \mathcal{M}$ with $\tilde{R}_{[Z]}(\xi)  = [\operatorname{mat}(Z_v + \bar{\xi}_Z)]$. 
Let 
    $$\tilde{R} : T\mathcal{M} \to \mathcal{M} \quad \text{ with }\quad \tilde{R}(\xi)= \tilde{R}_{[Z]_{\xi}}(\xi),$$ 
     where $T\mathcal{M} = \cup_{[Z] \in \mathcal{M}} T_{[Z]}\mathcal{M}$ is the tangent bundle of $\mathcal{M}$, and $[Z]_{\xi}\in \mathcal{M}$ is the foot of $\xi \in T\mathcal{M}$.  Then $\tilde{R}$ is the unique exponential retraction on $\mathcal{M}$.
\end{lemma}

\begin{proof}
See Section \ref{sec:pf_expmap} on Page \pageref{sec:pf_expmap}. 
\end{proof}


\medskip

Recall that the Newton-Raphson update for $pL_Q$ 
is given by the  exponential retraction $R$ in Section \ref{subsec:manifold}. 
Lemma \ref{lem:expmap}  suggests that $\tilde{R}=R$.
Therefore, $R_{[\check Z]}(\nu) = [\operatorname{mat}(\check{Z}_v + \bar{\nu}_{\check Z})]$. 
To finish the proof of Proposition \ref{prop:newtonstepform}, it remains to obtain the  closed-form formula for $\bar{\nu}_{\check Z}$, which is derived based on Lemmas \ref{lem:spec_nu} and \ref{lm_cal_grad_hess_proj} below. 


\begin{lemma} \label{lem:spec_nu}
Define $\overline{pL}: \mathbb{R}^{nk}\to \mathbb{R}$ 
with $\overline{pL}(Z_v) = pL(Z)$ for any $Z \in \mathbb R_0^{n \times k}$. 
Then  $\bar{\nu}_{\check Z}$ satisfies 
$$
\mathcal{P}_{\mathcal{H}_{\check Z}}\operatorname{Hess} \overline{pL}({{\check Z_v}})[\bar{\nu}_{\check Z}] = -\mathcal{P}_{\mathcal{H}_{\check Z}}\operatorname{Grad} \overline{pL}({{\check Z_v}}).
$$
\end{lemma}

\begin{proof}
See Section \ref{sec:pf_nu} on Page \pageref{sec:pf_nu}. 
\end{proof}

\begin{lemma}\label{lm_cal_grad_hess_proj}
For any $Z \in\mathbb{R}^{n\times k}_0$, we have
$$\operatorname{Grad} {\overline{pL}}(Z_v)=  \Seff({Z}, \hat{\alpha}(Z)), \quad \operatorname{Hess} {\overline{pL}}( Z_v) =  \Heff( Z, \hat{\alpha}(Z)),$$
where $\Seff$ and $\Heff$ refer to the efficient score and the efficient Hessian matrix defined in Section \ref{sec:newton}. 
\end{lemma}
\begin{proof}
See Section \ref{sec:pf_cal} on Page \pageref{sec:pf_cal}.
\end{proof}

\medskip
Combining Lemmas \ref{lem:spec_nu} and  \ref{lm_cal_grad_hess_proj}, $\bar{\nu}_{\check Z}$ satisfies 
 \begin{align} \label{mani-2}
\mathcal{P}_{\mathcal{H}_{\check Z}} \Heff(\check Z, \hat{\alpha}(\check Z)) [\bar{\nu}_{\check Z}] = -\mathcal{P}_{\mathcal{H}_{\check Z}} \Seff({\check Z}, \hat{\alpha}(\check Z)).
\end{align}
By the definition of $\mathcal{H}_Z$ in \eqref{eq:hzdef_prop3},  we have
$$
\mathcal{P}_{\mathcal{H}_{\check Z}} v = \check{\mathcal{U}} \check{\mathcal{U}}^{\top} v,\quad  \forall v \in \mathbb R ^{nk},
$$
where $\check{\mathcal{U}}$ is defined same as that in \eqref{newtonsolUZ}. 
 Note that $\bar{\nu}_{\check Z} \in \mathcal{H}_{\check Z}$ by the definition in Lemma \ref{lem_HspaceDiff}. 
 Then $$\bar{\nu}_{\check Z}=\mathcal{P}_{\mathcal{H}_{\check Z}} \bar{\nu}_{\check Z}= \check{\mathcal{U}} \check{\mathcal{U}}^{\top} \bar{\nu}_{\check Z} .$$ 
Plugging the above equations into \eqref{mani-2},  we obtain
\begin{align*}
\check{\mathcal{U}} \check{\mathcal{U}}^{\top} \Heff(\check Z, \hat{\alpha}(\check Z))  \check{\mathcal{U}} \check{\mathcal{U}}^{\top} \bar{\nu}_{\check Z} = -\ \check{\mathcal{U}} \check{\mathcal{U}}^{\top} \Seff({\check Z}, \hat{\alpha}(\check Z)).
\end{align*}
Multiplying $\check{\mathcal{U}}^{\top}$ on both sides of the above equation shows 
\begin{align*} 
\check{\mathcal{U}}^{\top}\bar{\nu}_{\check Z} = - \left[\check{\mathcal{U}}^{\top} \Heff(\check Z, \hat{\alpha}(\check Z)) \check{\mathcal{U}}\right]^{-1}\check{\mathcal{U}}^{\top} \Seff({\check Z}, \hat{\alpha}(\check Z)).
\end{align*} 
As $\bar{\nu}_{\check Z}\in \mathcal{H}_{\check Z}=\mathrm{col}(\check{\mathcal{U}})$, 
\begin{align*}
    \bar{\nu}_{\check Z} = - \ \check{\mathcal{U}}\left[\check{\mathcal{U}}^{\top} \Heff(\check Z, \hat{\alpha}(\check Z)) \check{\mathcal{U}} \right]^{-1}\check{\mathcal{U}}^{\top} \Seff({\check Z}, \hat{\alpha}(\check Z)). 
\end{align*}

\subsection{Proof of Lemmas \ref{lemma:Quo_mani}--\ref{lm_cal_grad_hess_proj} in Section \ref{sec:pfProp}}

\subsubsection{Proof of Lemma \ref{lemma:Quo_mani}}\label{sec:pfQuo_mani}
Since $\mathbb R^{n\times k}_0$ is an open subset of
$\mathbb S^{n\times k}_0:= \{Z\in \mathbb R^{n\times k} : 1_n^\top Z = 0 \} $, it is naturally a smooth submanifold of the linear manifold $\mathbb S^{n\times k}_0$.
Furthermore, it can be shown that the following map
$$
\psi: \mathcal{O}(k) \times \mathbb{R}_0^{n \times k} \rightarrow \mathbb{R}_0^{n \times k}:(Q, Z) \mapsto ZQ
$$
is a smooth, free, and proper Lie group action on the smooth manifold $\mathbb R^{n\times k}_0$. Then by the quotient manifold theorem \citep[][Thm. 21.10]{lee2013smooth}, $ \mathbb{R}_0^{n \times k}/ \sim$ is a topological
manifold of dimension $nk - k(k+1)/2$, and has a unique differentiable structure such that the quotient map $\pi : \mathbb{R}_0^{n \times k} \to \mathbb{R}_0^{n \times k}/ \sim$ is a smooth submersion. Endowed with such a differentiable structure, $\mathbb{R}_0^{n \times k}/\sim$ is turned into a smooth manifold. Note that $pL$ is a smooth function on $\mathbb{R}_0^{n \times k}$, then $pL_Q$ is smooth on $\mathbb{R}_0^{n \times k}/\sim$ according to the characteristic property of surjective smooth submersions \citep[][Thm. 4.29]{lee2013smooth}.


\subsubsection{Proof of Lemma \ref{lem_HspaceDiff}} 
\label{sec:pf_Hspace}
Each $\bar f \in \bar{\mathcal{F}}$ also induces a rotation invariant smooth function $\tilde f : \mathbb R^{n\times k}_0 \to \mathbb R$ by $\tilde f(Z) = \bar f(Z_v)$, differing from $f$ only in their domains.
For any $Z\in \mathbb R^{n\times k}_0$, note that
$$T_Z \mathbb R^{n\times k}_0 = \{Z\in \mathbb R^{n\times k} : 1_n^\top Z = 0 \}=\mathbb S^{n\times k}_0.$$ 
Following the classical Riemannian quotient manifold theory \citep[][Section 3.6.2]{absil2009optimization}, we have the orthogonal decomposition 
$$\mathbb S^{n\times k}_0 = \widetilde{\mathcal V}_Z  
\oplus \widetilde{\mathcal H}_Z,$$
where $\widetilde{\mathcal V}_Z = \{ZA : A\in \mathbb R^{k\times k}, A^\top = -A \}$ is the vertical space, and $\widetilde{\mathcal{H}}_Z$ denotes the horizontal space.
Moreover, the horizontal space  $\widetilde{\mathcal{H}}_Z$ satisfies that for any tangent vector $\xi \in T_{[Z]}\mathcal{M}$, there exists one and only one element $\tilde \xi_Z \in \widetilde{\mathcal{H}}_Z$ such that for any smooth $f$ on $\mathcal{M}$,
\begin{align}
    \xi f = \tilde \xi_Z \tilde f.\label{eq:lemf2_1}
\end{align}

We next show $\operatorname{mat}(\mathcal{H}_Z) = \widetilde{\mathcal{H}}_Z$, where $\mathrm{mat}(\mathcal{H}_Z)=\left\{\operatorname{mat}(v): v \in \mathcal{H}_Z\right\}$. From the definition of $\UZ$ (Definition \ref{def::UZA}) we have 
\begin{align*}
\mathbb R^{nk} = \text{col}(\UZ) \oplus \text{col}(\myDtri) \oplus \text{null}(\DZ),
\end{align*}
thus equivalently, 
$$
\mathbb R^{n \times k} = \operatorname{mat}(\text{col}(\UZ)) \oplus \operatorname{mat}(\text{col}(\myDtri)) \oplus \operatorname{mat}(\text{null}(\DZ)).
$$
By the definition of $\myDtri$, $\operatorname{mat}((\myDtri)_j) = 1_n (e_j^{(k)})^\top $, where $(\myDtri)_j$ denotes the $j$th column of $\myDtri$, and $e_j^{(k)}\in\mathbb R^k$ is the vector whose $j$th element is $1$ and others $0$. Thus,  
$$
\operatorname{mat}(\text{col}(\myDtri)) = \operatorname{span} \{\operatorname{mat}((\myDtri)_1), \ldots, \operatorname{mat}((\myDtri)_k)\} = \{1_ny_k^\top : y_k\in \mathbb R^k\},
$$
which indicates
$$ \operatorname{mat}(\text{col}(\UZ)) \oplus \operatorname{mat}(\text{null}(\DZ)) = \{Z\in \mathbb R^{n\times k}: 1_n^\top Z = 0 \} = \mathbb S^{n\times k}_0. $$
It remains to check that 
$\operatorname{mat}(\text{null}(\DZ)) = \widetilde{\mathcal V}_Z$.
We pick a linear basis of $\widetilde{\mathcal V}_Z$ given by $\{ UDB_{ij}V: 1\leqslant i<j\leqslant k \} $, where $Z=UDV$ is the singular value decomposition of $Z$ as before, and $B_{ij}$ denotes the $k\times k$ matrix whose $(i,j)$ element is $1$, $(j,i)$ element is $-1$, and others $0$. By direct calculation, we have $\operatorname{mat}(\lambda_i \beta_{ij} - \lambda_j \beta_{ji}) = UDB_{ij}V$ as follow,
\begin{align*}
UDB_{ij}V &= \begin{pmatrix}
& {{(\text{col }i)}} & & {{(\text{col }j)}} &
\\
0 \cdots 0 & -\lambda_j U_j & 0\cdots 0 & \lambda_i U_i & 0\cdots 0
\end{pmatrix}
\begin{pmatrix}
V_1^\top \\ \vdots \\ V_k^\top
\end{pmatrix} \\
&= \lambda_i U_iV_j^\top - \lambda_j U_j V_i^\top = \operatorname{mat}(\lambda_i \beta_{ij} - \lambda_j \beta_{ji}).
\end{align*}

The above argument indicates that the elements in $\mathcal{H}_Z$ and $\widetilde{\mathcal{H}}_Z$ differ only by vectorization, and so do $\bar f$ and $\tilde f$. Consequently, for any horizontal lift $ \tilde\xi_Z \in \widetilde{\mathcal{H}}_Z$, there exists one and only one vector $\bar \xi_Z = \operatorname{mat}^{-1}(\tilde\xi_Z) \in {\mathcal{H}}_Z$ such that for any $ \bar{f}\in \bar{\mathcal{F}}$, 
{\begin{align}
    \bar \xi_Z \bar f=\tilde{\xi}_Z \tilde f . \label{eq:lemf2_2}
\end{align}
Combining \eqref{eq:lemf2_1} and \eqref{eq:lemf2_2} yields Lemma \ref{lem_HspaceDiff}. }

\subsubsection{Proof of Lemma \ref{lem:expmap}}
\label{sec:pf_expmap}
We first show the constructed mapping $\tilde R_{[Z]}$ is well-defined. For any $\xi \in T_{[Z]}\mathcal{M} $, $Q \in \mathcal{O}(k)$, and $k \times k$ antisymmetric matrix $A$, 
$$
\langle \tilde{\xi}_Z Q , ZQA \rangle = 
\operatorname{tr}(Q^{\top}(\tilde{\xi}_Z)^{\top}ZQA) = \operatorname{tr}((\tilde{\xi}_Z)^{\top}ZQAQ^{\top}) = \langle \tilde{\xi}_Z  , ZQA Q^{\top}\rangle = 0,
$$
where the last equality follows from the definition of $\tilde{\xi}_Z$; see Section \ref{sec:pf_Hspace}. The above equation indicates $\tilde \xi_Z Q \in \widetilde{\mathcal{H}}_{ZQ}$, which further shows $\tilde\xi_Z Q = \tilde\xi_{ZQ}$ since $$(\tilde\xi_Z Q) \tilde{f} = \tilde\xi_Z \tilde{f} = \xi f = \tilde\xi_{ZQ}  \tilde{f}$$ for any $\tilde{f}$. As a result,
$$
[\operatorname{mat}((ZQ)_v + \bar\xi_{ZQ})] = [ZQ + \tilde\xi_{ZQ}] = [ZQ + \tilde\xi_{Z} Q] = [Z+\tilde\xi_Z] = [\operatorname{mat}(Z_v + \bar\xi_{Z})] .
$$

Let $\widetilde{\mathcal{H}}^{\star}_Z$ denote the largest star-shaped set of lifted tangent vectors that  emanate from $Z$ and maintaining full rank:
    $$
    \widetilde{\mathcal{H}}^{\star}_Z := \{\tilde \xi_Z \in \widetilde{\mathcal{H}}_Z : \operatorname{rank} (Z+ x \tilde \xi_Z) = k, \ \forall x \in [0,1]\}.
    $$
    We next illustrate that for any $\xi \in T_{[Z]}\mathcal{M}$ satisfying $\tilde \xi_Z \in \widetilde{\mathcal{H}}^{\star}_Z$, the exponential map $R^{exp}$ on $\mathcal{M}$ is given by 
    $$
    R^{exp}_{[Z]}(\xi)=[Z+\tilde{\xi}_Z].
    $$
For any $\tilde \xi_Z \in \widetilde{\mathcal{H}}^{\star}_Z$, we consider a curve $\gamma : x \mapsto Z + x \tilde{\xi}_Z$. It follows that $\gamma$ is a geodesic of $\mathbb R^{n\times k}_0$ for all $x \in [0,1]$, since $\gamma$ is a straight line  in $\mathbb R^{n\times k}_0$. Note that $\dot{\gamma}(0)$ is horizontal, then the theory of geodesics and Riemannian submersions \citep[][Prop. 2.109]{gallot1990riemannian} shows $\pi \circ \gamma$ is a geodesic of $\mathcal{M}$. By checking 
$$\pi \circ \gamma(0)=[Z] \quad\text{and}\quad \frac{\mathrm{d}}{\mathrm{d} x}(\pi \circ \gamma)(0)= \xi, $$ we have $x \mapsto [Z + x \tilde{\xi}_Z]$ is a geodesic of $\mathcal{M}$ emanating from $[Z]$ with initial velocity $\xi$. Therefore, for any $\xi \in T_{[Z]} \mathcal{M}$ satisfying $\tilde \xi_Z \in \widetilde{\mathcal{H}}^{\star}_Z$, we have
$$
R_{[Z]}^{exp}(\xi)=[Z+\tilde{\xi}_Z] = [\operatorname{mat}(Z_v+\bar{\xi}_Z)] = \tilde R_{[Z]}(\xi).
$$

\subsubsection{Proof of Lemma \ref{lem:spec_nu}} \label{sec:pf_nu}

{We note that $\overline{pL}\in \mathcal{F}$. 
By Lemma \ref{lem_HspaceDiff},  each $\bar{f}\in \mathcal{F}$ naturally induces a smooth function $f:\mathcal{M}\to \mathbb{R}$.
To prove Lemma \ref{lem:spec_nu}, we next show that for any $\bar{f}\in \mathcal{F}$  and its induced $f$, 
\begin{align}
     \overline{\operatorname{Grad} f([Z])}_Z &= \mathcal{P}_{\mathcal{H}_{ Z}}\operatorname{Grad} \bar{f}({Z_v}),\label{eq:gradproj}\\
     \overline{\operatorname{Hess} f({[Z]})[\eta]}_Z &= \mathcal{P}_{\mathcal{H}_{Z}}{\operatorname{Hess}} \bar{f}({Z_v})[\bar{\eta}_Z]. \label{eq:hessproj}
\end{align}
}

{To prove \eqref{eq:gradproj},}
firstly by seeing $\mathcal{M}$ as a Riemannian quotient manifold of $\mathbb R_0^{n\times k}$ whose Riemannian metric is induced from $\mathbb S_0^{n\times k}$, then we have $$\overline{\operatorname{Grad} f([Z])}_Z= \operatorname{mat}^{-1}(\mathcal{P}_{\widetilde{\mathcal{H}}_Z} \operatorname{Grad} \tilde{f}(Z))$$ 
by the classical theory of Riemannian quotient manifold \citep[][Section 3.6.2]{absil2009optimization},
{where $\tilde{f}$ is defined same as in Section \ref{sec:pf_Hspace}. } 
Secondly, by seeing $\mathbb R_0^{n\times k}$ as an embedded submanifold of $\mathbb R^{n\times k}$, we have $\operatorname{Grad} \tilde{f}(Z) = \mathcal{P}_{\mathbb S_0^{n\times k}} \operatorname{mat}(\operatorname{Grad} \bar{f}(Z_v)) $. Combining the two steps above, we establish \eqref{eq:gradproj} for any $\bar f \in \bar{\mathcal{F}}$ and $Z\in\mathbb{R}^{n\times k}_0$. 

{To prove \eqref{eq:hessproj},}
when the Hessian of $f$ is naturally defined by the Riemannian connection of $\mathcal{M}$, for any tangent vector $\eta \in T_{[Z]}\mathcal{M}$ we have
$$
\overline{\operatorname{Hess} f({[Z]})[\eta]} _Z= \operatorname{mat}^{-1}(\mathcal{P}_{\widetilde{\mathcal{H}}_{Z}}{\operatorname{Hess}} \tilde{f}({Z})[\tilde{\eta}_Z]).
$$
Then, from the embedded submanifold perspective,
$${\operatorname{Hess}} \tilde{f}({{Z}})[\tilde{\eta}_Z] = \mathcal{P}_{\mathbb S_0^{n\times k}} \operatorname{mat}({\operatorname{Hess}} \bar{f}({Z_v})[\bar{\eta}_Z]) .$$
As a result, we obtain \eqref{eq:hessproj} for any $\bar f \in \bar{\mathcal{F}}$, $Z\in\mathbb{R}^{n\times k}_0$ and tangent vector $\eta \in T_{[Z]}\mathcal{M}$.

Plugging $\bar f = \overline{pL}$, $Z = \check Z$, and $\eta = \nu$ into \eqref{eq:gradproj} and \eqref{eq:hessproj}, we obtain 
\begin{align*} 
\overline{\operatorname{Hess} pL_Q({[\check Z]})[\nu]}_{\check Z} &= \mathcal{P}_{\mathcal{H}_{\check Z}}\operatorname{Hess} \overline{pL}({{\check Z_v}})[\bar{\nu}_{\check Z}] \\
- \overline{\operatorname{Grad} pL_Q([\check Z])}_{\check Z}
&= - \mathcal{P}_{\mathcal{H}_{\check Z}}\operatorname{Grad} \overline{pL}({{\check Z_v}}).
\end{align*}
Moreover, by the definition of Newton update in \eqref{eq:netwonhessgrad}, and Lemma \ref{lem_HspaceDiff}, we have 
\begin{align*}
    \overline{\operatorname{Hess} pL_Q({[\check Z]})[\nu]}_{\check Z} = - \overline{\operatorname{Grad} pL_Q([\check Z])}_{\check Z}. 
\end{align*}
Combining the above three equations, Lemma  \ref{lem:spec_nu} is proved.


\subsubsection{Proof of Lemma \ref{lm_cal_grad_hess_proj}}
\label{sec:pf_cal}
Note that from the definition of $\hat\alpha(Z)$ we have $\dot{L}_{\alpha}(Z,\hat\alpha(Z)) =0$ for any $Z$. Taking the derivative with respect to $Z$ gives the important formula
$$
\left.\frac{\partial{\hat\alpha(Z)} }{\partial Z_v}\right|_{Z} = - \left.\left(\ddot{L}_{\alpha,\alpha}^{-1} \ddot{L}_{\alpha,Z}\right)\right|_{(Z,\hat\alpha(Z))}.
$$
Hence, we have
\begin{align*}
\operatorname{Grad} \overline{pL}(Z_v)&= \left.\frac{\partial pL(Z)}{\partial Z_v}\right|_{Z}=
\left.\frac{\partial L(Z,\hat\alpha(Z))}{\partial Z_v}\right|_{Z}=
\left.\left\{\dot{L}_Z + \left(\frac{\partial{\hat\alpha(Z)} }{\partial Z_v}\right)^\top \dot{L}_{\alpha}\right\}\right|_{(Z,\hat\alpha(Z))}
\\&= \left.\left\{\dot{L}_Z - \left(\ddot{L}_{\alpha,\alpha}^{-1} \ddot{L}_{\alpha,Z}\right)^\top \dot{L}_{\alpha}\right\}\right|_{(Z,\hat\alpha(Z))}
= \Seff({Z}, \hat{\alpha}(Z)).
\end{align*}
Further, note that 
$$\operatorname{Grad} {\overline{pL}}(Z_v)=\left.\left\{\dot{L}_Z - \left(\ddot{L}_{\alpha,\alpha}^{-1} \ddot{L}_{\alpha,Z}\right)^\top \dot{L}_{\alpha}\right\}\right|_{(Z,\hat\alpha(Z))} = \dot{L}_Z(Z,\hat\alpha(Z))$$
as $\dot{L}_{\alpha}{(Z,\hat\alpha(Z))} =0$. 
Thus, 
\begin{align*}
    \operatorname{Hess} \overline{pL}(Z_v)&= \left.\frac{\partial \operatorname{Grad} \overline{pL}(Z_v)}{\partial Z_v}\right|_{Z} = \left.\frac{\partial \dot{L}_Z(Z,\hat\alpha(Z))}{\partial Z_v}\right|_{Z} 
    \\&= \left.\left\{\ddot{L}_{Z,Z} + \left(\frac{\partial{\hat\alpha(Z)} }{\partial Z_v}\right)^\top \ddot{L}_{\alpha,Z}\right\}\right|_{(Z,\hat\alpha(Z))}
    \\&= \left.\left\{\ddot{L}_{Z,Z} - \left(\ddot{L}_{\alpha,\alpha}^{-1} \ddot{L}_{\alpha,Z}\right)^\top \ddot{L}_{\alpha,Z}\right\}\right|_{(Z,\hat\alpha(Z))}
    =  \Heff( Z, \hat{\alpha}(Z)).
\end{align*}




\newpage

\subsection{Additional Results for Section \ref{subsec:manifold}}\label{sec:addgeores}
In this section, we establish rigorous results for the heuristic statement after Proposition \ref{prop:newtonstepform} in the main text.
In particular, 
we first present  Propositions \ref{prop::arguments_after_prop3} and \ref{prop::9_in_R0nbyk}, with their proofs given in Sections \ref{pf::prop::arguments_after_prop3} and \ref{pf::prop::9_in_R0nbyk}, respectively.

\begin{proposition} \label{prop::arguments_after_prop3}
   Consider a  pair of initial estimate $(\check{Z},\check{\alpha}) $ satisfying Condition \ref{cond_elem_init} with constants $M_b$ and $\epsilon$. Under Condition \ref{cond:truvalueregularity}, we have
    \begin{enumerate}
    \item 
    For any constant $s>0$, there exists a constant $C_s>0$ such that when $n/\log^{2\varsigma}(nT) $ is sufficiently large ($\varsigma$ is defined same as in Theorem \ref{thm_NRerr}), 
\begin{align*}
     \Pr\left\{ \underset{1\leqslant i\leqslant n,1\leqslant t\leqslant T}{\max} \,  |\hat{\alpha}_{it}(\check Z) -   \check{\alpha}_{it}| > \frac{C_s \log^{\varsigma}(nT)}{\sqrt{n} } \right\} = O((nT)^{-s}). 
\end{align*}
    \item 
    For any constant $s>0$, there exists a constant $C_s>0$ such that when $n/\log^{2\varsigma}(nT) $ is sufficiently large, 
\begin{align*}
\Pr\left\{
\big\|\check Z_v + \bar{\nu}_{\check Z} -  \eqref{newtonsolUZ2} \big\|_2^2 > 
\frac{1}{T} \times C_s r_{n,T}
\right\} = O(n^{-s}),
\end{align*}
where $r_{n,T} = \max\{1,\frac{T}{n}\}\log^{4\varsigma}(nT)  $.
\end{enumerate}
\end{proposition}

\bigskip 
Proposition \ref{prop::arguments_after_prop3} 
provides theoretical justification for the heuristic argument after Proposition \ref{prop:newtonstepform}. 
In particular, it suggests that 
up to logarithmic factors, 
``\textit{$\hat{\alpha}(\check Z)$ and $\check\alpha $ are close}'' in the sense that $\max_{i,t} |\hat{\alpha}_{it}(\check Z)-\check{\alpha}_{it}| = O_p( n^{-1/2})$,
and ``\textit{$\check Z_v + \bar{\nu}_{\check Z} $ and \eqref{newtonsolUZ2} are close}'' in the sense that $\|\check Z_v + \bar{\nu}_{\check Z} -  \eqref{newtonsolUZ2}\|_2 =O_p\left(T^{-1/2} + n^{-1/2} \right) $, both of which decay  as $n$ and $T \to \infty$.   

\vspace{1em}
\begin{proposition} \label{prop::9_in_R0nbyk}
Let $\hat Z \in \mathbb{R}^{n\times k}$ denote the matrix version of \eqref{newtonsolUZ2}. Under Conditions \ref{cond:truvalueregularity} and \ref{cond_elem_init},  we have 
    \begin{itemize}
         \item[(i)] $1_n^\top \hat Z = 0 $; 
    \item[(ii)] 
    $\Pr\{\operatorname{det}(\hat Z^\top \hat Z )= 0\}=O(n^{-1}) $ when $n / \log^{4\varsigma}(nT) $ is sufficiently large.
    \end{itemize}
\end{proposition}

\vspace{1em}
Proposition \ref{prop::9_in_R0nbyk} suggests that  $\hat Z\in \mathbb R_0^{n\times k} $ with high probability.   
Thus, the resulting equivalence class $[\hat Z]$ belongs to the quotient set $\mathbb{R}_0^{n\times k}/\sim $.  

\vspace{1em}

\subsubsection{Proof of Proposition \ref{prop::arguments_after_prop3}}
\label{pf::prop::arguments_after_prop3}
\noindent
\textit{Proof of Statement 1:}
Denote $\check{G}_{ij} = \check z_i^\top \check z_j $; recall $L(Z,\alpha) $ is the likelihood function defined in \eqref{eq:LZalpha}, and $G^\star_{ij} = (z_i^\star)^\top z_j^\star $.
By the definition of $\hat{\alpha}$, we have that $(\check{Z}, \hat{\alpha}(\check{Z}) ) $ satisfies the estimating equations
\begin{align*}
0&= \left.\frac{\partial L(Z,\alpha)}{\partial \alpha_{it}} \right|_{(\check{Z},\hat\alpha_t(\check Z))} =\left(2A_{t,ii} - 2 e^{\check G_{ii} + 2\hat\alpha_{it}(\check{Z}) }\right) + \sum_{j\neq i}\left(A_{t,ij}- e^{\check G_{ij}+\hat{\alpha}_{it}(\check{Z})+\hat{\alpha}_{jt}(\check{Z})} \right) \\
&= 2\check M_{t,ii} + \sum_{j\neq i}\check M_{t,ij}\\
&\quad+ \left( 2 e^{ \check G_{ii} + 2\alpha_{it}^\star } - 2 e^{ \check G_{ii} + 2\hat\alpha_{it}(\check{Z}) }\right) + \sum_{j\neq i}\left(  e^{ \check G_{ij}+{\alpha}_{it}^\star+{\alpha}_{jt}^\star} -   e^{ \check G_{ij}+\hat{\alpha}_{it}(\check{Z}) +\hat{\alpha}_{jt}(\check{Z})} \right) 
\end{align*}
for $i=1,\ldots,n$ and $t=1,\ldots,T$, where
\begin{align*}
    \check M_{t,ij} := A_{t,ij} - e^{\check{G}_{ij} + \alpha_{it}^\star + \alpha_{jt}^\star  } = \left(A_{t,ij} - e^{{G}_{ij}^\star + \alpha_{it}^\star + \alpha_{jt}^\star  } \right) + \left( e^{{G}_{ij}^\star + \alpha_{it}^\star + \alpha_{jt}^\star  } - e^{\check{G}_{ij} + \alpha_{it}^\star + \alpha_{jt}^\star  } \right).
\end{align*}
Denote $\check{M}_{2,t,i\cdot} :=  2\check M_{t,ii} + \sum_{j\neq i}\check M_{t,ij}$, and $\check\Delta\alpha_{it} := \hat\alpha_{it}(\check Z) - \alpha_{it}^\star$. By applying a mean value theorem to the above equation between the vectors 
$(\hat\alpha_{1t}(\check{Z}),\ldots,\hat\alpha_{nt}(\check{Z}) )^\top$ and $(\alpha_{1t}^\star,\ldots,\alpha_{nt}^\star )^\top$, we have
\begin{equation}
\label{eq::checkM_1stTaylor}
       \check{M}_{2,t,i\cdot} = \Big( 4e^{\check G_{ii} + 2{\bar\alpha_{it}}^{(it)} } + \sum_{j\neq i} e^{\check G_{ij} + {\bar\alpha_{it}}^{(it)}  + {\bar\alpha_{jt}}^{(it)} } \Big)   \check\Delta\alpha_{it} 
 + \sum_{j\neq i} \Big( e^{\check G_{ij} + \bar\alpha^{(it)}_{it}  + \bar\alpha^{(it)}_{jt}} \Big)  \check\Delta\alpha_{jt} 
\end{equation}
where for $j=1,\ldots,n$, $\bar\alpha_{jt}^{(it)} $ is a midpoint between $\hat\alpha_{jt}(\check{Z}) $ and $\alpha_{jt}^\star $ (since this midpoint is used in the expansion of the score equation with respect to $\alpha_{it} $, the midpoint depends on $i,t$, and hence the superscript $^{(it)}$).

Since $\|z_i^\star\|_2^2 \leqslant M_{Z,1} $ and $\operatorname{dist}_i(\check z_i, z_i^\star) \leqslant M_b n^{-1/2}\log^{\epsilon}(nT) $, we have $\|\check{z}_i\|_2^2 \leqslant M_{Z,1} + 1 $ as $n/\log^{2\varsigma}(nT) $ is sufficiently large where $\varsigma = \max\{\epsilon,1/2\}$. By the same argument as in Lemma \ref{lem_bdd}, there exists a constant $M_{\Theta,2}'' $ depending only on $M_{\Theta,1}, M_{Z,1},M_{Z,2} $ and $s$, such that with probability $1-O((nT)^{-s}) $, we have $|\check{G}_{ij} +  \hat\alpha_{it}(\check{Z}) +  \hat\alpha_{jt}(\check{Z})|< M_{\Theta,2}'' $ for any $i,j,t$; without loss of generality, we can take $M_{\Theta,2}'' > M_{Z,1} + 2 M_{\alpha} +1 $ so that $|\check G_{ij} + \alpha_{it}^\star + \alpha_{jt}^\star| < M_{\Theta,2}''  $ for any $i,j,t$. Since $\bar\alpha_{jt}^{(it)} $ is a midpoint between $\hat\alpha_{jt}(\check{Z}) $ and $\alpha_{jt}^\star $, with probability $1-O((nT)^{-s}) $, we have $|\check G_{ij} + \bar\alpha_{it}^{(it)} + \bar\alpha_{jt}^{(it)}| < M_{\Theta,2}'' $ for any $i,j,t$. 
Let $\bar{\mathcal{Y} }_t $ be an $n\times n$ matrix whose $(i,i)$-th entry and $(i,j)$-th entry for $j\neq i$ are  $$ \bar y_{t,ii}= \Big( 4e^{\check G_{ii} + 2{\bar\alpha_{it}}^{(it)} } + \sum_{j\neq i} e^{\check G_{ij} + {\bar\alpha_{it}}^{(it)}  + {\bar\alpha_{jt}}^{(it)} } \Big)\quad\text{ and }\quad \bar y_{t,ij}=\Big( e^{\check G_{ij} + \bar\alpha^{(it)}_{it}  + \bar\alpha^{(it)}_{jt}} \Big),$$ respectively.
Let $\bar{\mathcal{W}}_t = \bar{\mathcal{Y} }_t^{-1} $ and denote by $\bar w_{t,ij} $ the $(i,j)$-th entry of $\bar{\mathcal{W}}_t$.
Applying Lemma \ref{lem:thm1.1_Hillar} on $\bar{\mathcal{Y}}_t$, we have 
\begin{equation} \label{eq::linfty_bd_barW}
    \max_{i,t} \sum_{j=1}^n |\bar w_{t,ij} | \leqslant 2\exp(M_{\Theta,2}'') / n  
\end{equation}
with probability $1-O((nT)^{-s}) $.

On the other hand, by applying the mean value theorem, we have 
$$e^{{G}_{ij}^\star + \alpha_{it}^\star + \alpha_{jt}^\star  } - e^{\check{G}_{ij} + \alpha_{it}^\star + \alpha_{jt}^\star  } = e^{\bar{G}_{ij} + \alpha_{it}^\star + \alpha_{jt}^\star  } (G_{ij}^\star - \check G_{ij} ) $$ 
where $\bar{G}_{ij}$ is a midpoint between $G_{ij}^\star$ and $\check G_{ij}$. Same as the above arguments for $\bar\alpha_{it}^{(it)} $, we also have $|\bar G_{ij}+ \alpha_{it}^\star + \alpha_{jt}^\star| < M_{\Theta,2}''  $ for any $i,j,t$. Since $\|z_i^\star\|_2^2 \leqslant M_{Z,1} $ and $\operatorname{dist}_i(\check z_i, z_i^\star) \leqslant M_b n^{-1/2}\log^{\epsilon}(nT) $, we have 
\begin{align*}
   |\check G_{ij} - G_{ij}^\star| \leqslant 2\sqrt{M_{Z,1}}M_b n^{-1/2}\log^{\epsilon}(nT) + \{M_b n^{-1/2}\log^{\epsilon}(nT)\}^2 <  3\sqrt{M_{Z,1}}M_b n^{-1/2}\log^{\epsilon}(nT)
\end{align*}
as $n/\log^{2\varsigma}(nT) $ is sufficiently large where $\varsigma = \max\{\epsilon,1/2\}$. Therefore we have 
\begin{equation}
    \label{eq::checkM_2ndterm}
    \begin{aligned}
            \max_{i,j,t} | e^{{G}_{ij}^\star + \alpha_{it}^\star + \alpha_{jt}^\star  } - e^{\check{G}_{ij} + \alpha_{it}^\star + \alpha_{jt}^\star  }| =&~  \max_{i,j,t}|e^{\bar{G}_{ij} + \alpha_{it}^\star + \alpha_{jt}^\star  } (G_{ij}^\star - \check G_{ij} ) | \\
            <&~ 3\sqrt{M_{Z,1}}M_b e^{-M_{\Theta,2}'' } n^{-1/2}\log^{\epsilon}(nT).
    \end{aligned}
\end{equation}
Recall from \eqref{eq::def_Mtij}, \eqref{eq::def_M2ticdot} and Lemma \ref{lem_hillar_union} that with probability at least $1 - O((nT)^{-s} )$,
\begin{equation}
   \label{eq::checkM_1stterm}
\max_{i,t} |M_{2t,i\cdot}| \leqslant 2 (s+1)\sqrt{(n+3)e^{-M_{\Theta,1}} \log(nT) } 
\end{equation}
where $M_{t,ij} = A_{t,ij} - e^{{G}_{ij}^\star + \alpha_{it}^\star + \alpha_{jt}^\star  } $ and $M_{2t,i\cdot} = 2M_{t,ii} + \sum_{j\neq i} M_{t,ij} $.
Combining \eqref{eq::checkM_1stterm} and \eqref{eq::checkM_2ndterm}, with probability at least $1 - O((nT)^{-s} )$, we have
\begin{equation}
    \label{eq::checkM_cdot_sup_bd}
    \max_{i,t} |\check{M}_{2,t,i\cdot}| \leqslant \max_{i,t} |\check{M}_{2t,i\cdot}| + (n+1) \max_{i,j,t} | e^{{G}_{ij}^\star + \alpha_{it}^\star + \alpha_{jt}^\star  } - e^{\check{G}_{ij} + \alpha_{it}^\star + \alpha_{jt}^\star  }| \leqslant C_s \sqrt{n} \log^{\varsigma}(nT).
\end{equation}

By writing out equation \eqref{eq::checkM_1stTaylor} for each $i=1,\ldots,n$ into a system and rearranging terms, we have
\begin{align*}
\bar{\mathcal{Y} }_t
\begin{pmatrix} 
\check\Delta \alpha_{1t}\\
\vdots\\
\check\Delta \alpha_{nt}
\end{pmatrix}=
\begin{pmatrix}
\check M_{2,t,1\cdot} \\
\vdots \\
\check M_{2,t,n\cdot}  
\end{pmatrix}, \quad \text{i.e., } \
\begin{pmatrix} 
\check\Delta \alpha_{1t}\\
\vdots\\
\check\Delta \alpha_{nt}
\end{pmatrix}=\bar{\mathcal{W}}_t
\begin{pmatrix}
\check M_{2,t,1\cdot} \\
\vdots \\
\check M_{2,t,n\cdot}  
\end{pmatrix}.
\end{align*}
By applying \eqref{eq::linfty_bd_barW} and \eqref{eq::checkM_cdot_sup_bd} in the above equation, we have with probability at least $1 - O((nT)^{-s} )$ that
\begin{align*}
&~\max_{i,t} |\hat\alpha_{it}(\check Z) - \alpha_{it}^\star|=
    \max_{i,t} |\check{\Delta}\alpha_{it}| \\ =&~\max_{i,t} | \sum_{j=1}^n \bar w_{t,ij} \check{M}_{2,t,j\cdot} |  \leqslant \left(\max_{l,t}|\check M_{2,t,l\cdot}| \right)\left( \max_{l,t} \sum_{j=1}^n |\bar w_{t,lj}| \right) \notag 
      \leqslant C_{s} n^{-1/2}\log^{\varsigma}(nT)
\end{align*}
where $C_s $ is a constant depending on $ M_{\Theta,1}, M_{Z,1}, M_{Z,2}, M_b $ and $s$.
Combining this result and   $\max_{i,t} |\check\alpha_{it} - \alpha^\star_{it}| \leqslant M_b n^{-1/2} \log^{\epsilon}(nT) $ from Condition \ref{cond_elem_init} with a triangle inequality, we come to the conclusion of Statement (a): 
for any  $s>0$, there exists a constant $C_s>0$ such that
\begin{align*}
     \Pr\left\{ \underset{1\leqslant i\leqslant n,1\leqslant t\leqslant T}{\max} \,  |\hat{\alpha}_{it}(\check Z) -   \check{\alpha}_{it}| > \frac{C_s \log^{\varsigma}(nT)}{\sqrt{n} } \right\} = O((nT)^{-s}). 
\end{align*}
\bigskip

\noindent
\textit{Proof of Statement 2:} In Section \ref{sec:pfRMK}, we have shown that, given a pair of initial estimators $(\check{Z},\check{\alpha}) $ satisfying Condition \ref{cond_elem_init} with constants $M_b$ and $\epsilon$, \eqref{newtonsolUZ2} defined by
$$
\eqref{newtonsolUZ2} = \check Z_v - \check{\mathcal{U}}\left\{\check{\mathcal{U}}^\top H_{eff}(\check Z,\check\alpha) \ \check{\mathcal{U}}  \right\}^{-1} \check{\mathcal{U}}^\top S_{eff}(\check Z,\check\alpha)
$$
satisfies
\begin{align}
\label{eq::9-Zstarv}
\big\| \eqref{newtonsolUZ2} -  (\ZQs)_v \big\|_2 \leqslant 
C_{s} \left(\frac{\log^{2\varsigma}(nT)\sqrt k}{\sqrt n}  + \frac{\log^{\frac12}(nT) \sqrt{k}}{\sqrt{T}}\right)
\end{align}
with probability $1-O(n^{-s}) $. (Recall that $C_s$ denotes a constant depending on $s$ and not on $n,T$; the value of $C_s$ may change from line to line.)
From Statement 1, we have that the pair $(\check{Z}, \hat{\alpha}(\check{Z}) ) $ satisfies the following two conditions:
\begin{itemize}
    \item[(i)] There exist constants $C_s >0$ and  $\varsigma $ such that with probability $1-O((nT)^{-s}) $, $|\hat\alpha_{it}(\check{Z}) - \alpha_{it}^\star| + \operatorname{dist}_i(\check{z}_i, z_i^{\star})\leqslant C_s n^{-{1}/{2}}\log^{\varsigma}(nT)$ for  $ 1\leqslant i \leqslant n$ and $1\leqslant t \leqslant T$; 
    \item[(ii)] $1_n^{\top} \check{Z}= 0$. 
\end{itemize}
These two conditions are similar to Condition \ref{cond_elem_init}, with the constant $M_b$ replaced by $C_s$, $\epsilon$ replaced by $\varsigma=\max\{\epsilon,1/2\} $, and with high probability. 
Note that 
$$\check Z_v + \bar{\nu}_{\check Z} = \check Z_v - \check{\mathcal{U}}\left\{\check{\mathcal{U}}^\top H_{eff}(\check Z,\hat\alpha(\check{Z})) \ \check{\mathcal{U}}  \right\}^{-1} \check{\mathcal{U}}^\top S_{eff}(\check Z,\hat\alpha(\check{Z})),  $$
i.e., $\check Z_v + \bar{\nu}_{\check Z}$ has the same formula as \eqref{newtonsolUZ2}  except that $\check{\alpha}$ is replaced with $\hat{\alpha}(\check{Z})$; in particular, they share the same $\check Z$, and hence the same $\check{\mathcal{U}} $ and $\check{Q}$.
Thus, following analogous arguments to Supplementary Material \ref{sec:pfRMK}, we have 
\begin{align}
\label{eq::Prop3-Zstarv}
\big\|\check Z_v + \bar{\nu}_{\check Z} -  (\ZQs)_v \big\|_2 \leqslant 
C_{s} \left(\frac{\log^{2\varsigma}(nT)\sqrt k}{\sqrt n}  + \frac{\log^{\frac12}(nT) \sqrt{k}}{\sqrt{T}}\right)
\end{align}
with probability $1-O(n^{-s}) $.
Combining \eqref{eq::9-Zstarv} and \eqref{eq::Prop3-Zstarv} with a triangle inequality, we come to the conclusion of Statement (b):
for any $s>0$, there exists a constant $C_s>0$ such that 
\begin{align*}
\Pr\left\{
\big\|\check Z_v + \bar{\nu}_{\check Z} -  \eqref{newtonsolUZ2} \big\|_2 > 
C_{s} \left(\frac{\log^{2\varsigma}(nT)\sqrt k}{\sqrt n}  + \frac{\log^{\frac12}(nT) \sqrt{k}}{\sqrt{T}}\right) \right\} = O(n^{-s}).
\end{align*}

\bigskip

\subsubsection{Proof of Proposition \ref{prop::9_in_R0nbyk}}
\label{pf::prop::9_in_R0nbyk}
First we show $1_n^\top \hat Z = 0 $. By \eqref{newtonsolUZ2} the vectorized version of $\hat Z$ is defined by
$$
\hat Z_v = \check Z_v - \check{\mathcal{U}}\left\{\check{\mathcal{U}}^\top H_{eff}(\check Z,\check\alpha) \ \check{\mathcal{U}}  \right\}^{-1} \check{\mathcal{U}}^\top S_{eff}(\check Z,\check\alpha).
$$
By Condition \ref{cond_elem_init}, $1_n^\top \check{Z} = 0$. By the definition of $\check{\mathcal{U}} = \mathcal{U}_Z$ in definition \ref{def::UZA} in Section  \ref{sec:eigenspace}, we have $(1_n \otimes \mathrm I_{k\times k})^\top \check{\mathcal{U}} = 0 $, where $\mathrm I_{k\times k}$ is the $k\times k$ identity matrix. This is equivalent to $1_n^\top \operatorname{mat}(\check{\mathcal{U}} v) = 0 $ for any $v\in \mathbb R^{nk - k(k+1)/2} $, where the operator $\operatorname{mat}(\cdot)$ is defined in Proposition \ref{prop:newtonstepform}. In particular, we have $$1_n^\top \operatorname{mat}\left[\check{\mathcal{U}} \left\{\check{\mathcal{U}}^\top H_{eff}(\check Z,\check\alpha) \ \check{\mathcal{U}}  \right\}^{-1} \check{\mathcal{U}}^\top S_{eff}(\check Z,\check\alpha)\right] = 0.$$ Therefore, we have $1_n^\top\hat Z=0$. 

Next, we show that $\operatorname{det}(\hat Z^\top \hat Z) \neq 0$ asymptotically. Denote by $\sigma_l(G) $ the $l$-th largest eigenvalue of (a symmetric) matrix $G$; we prove the following statement: ``For any constant $s>0$, there exists a constant $C_s$ such that when $n / \log^{4\varsigma}(nT) $ is sufficiently large, $\sigma_k(\hat Z^\top \hat Z^\top) \geqslant C_s n $ with probability $1-O(n^{-s})$.''

By Remark \ref{RMK:Heff_replace_Ieff} (and the corresponding proof in Supplementary Material \ref{sec:pfRMK}), we have 
$$\operatorname{dist}^2(\hat Z,Z^\star) = \min_{Q\in\mathcal{O}(k)} \|\hat Z - Z^\star Q\|_{\mathrm{F}}^2 \leqslant C_s \times\frac{1}{T} \times \max\big\{1,\frac{T}{n}\big\}\log^{4\varsigma}(nT) $$ 
with probability $1-O(n^{-s}) $ for the $\hat Z$ defined by the matrix version of \eqref{newtonsolUZ2}. 
By Lemma \ref{Lemma-29}, since $$\operatorname{dist}(\hat Z,Z^\star) \leqslant \sqrt{M_{Z,2}n}\leqslant \|Z^\star\|_{\mathrm{op}} $$ when $n/\log^{2\varsigma}(nT) $ is sufficiently large, we have 
$$
\|\hat Z \hat Z^\top - Z^\star (Z^\star)^\top\|_{\mathrm{F}}  \leqslant 3 \|Z^\star\|_{\mathrm{op}} \operatorname{dist}(\hat Z, Z^\star) \leqslant 3\sqrt{M_{Z,1}n} \operatorname{dist}(\hat Z, Z^\star).
$$
 By Weyl's inequality, we have with probability $1-O(n^{-s}) $ that
\begin{align*}
  \sigma_k(\hat Z\hat Z^\top) \geqslant&~ \sigma_k\{Z^\star (Z^\star)^\top \} - \sigma_1\{Z^\star (Z^\star)^\top -\hat Z\hat Z^\top \}  \\
   \geqslant&~ \sigma_k\{Z^\star (Z^\star)^\top \} - \|\hat Z \hat Z^\top - Z^\star (Z^\star)^\top\|_{\mathrm{F}} \\
   \geqslant&~ M_{Z,2} n - 3C_s \sqrt{M_{Z,1}n\times \max\big\{\frac{1}{T},\frac{1}{n}\big\}\log^{4\varsigma}(nT)} \\
   \geqslant&~ C_s n 
\end{align*}
when $n / \log^{4\varsigma}(nT) $ is sufficiently large. In particular, this implies  $\Pr(\operatorname{det}(\hat Z^\top \hat Z )= 0)=O(n^{-1}) $ when $n / \log^{4\varsigma}(nT) $ is sufficiently large.

Combining $1_n^\top \hat Z = 0 $ and $\operatorname{det}(\hat Z^\top \hat Z) \neq 0 $ asymptotically, we come to the conclusion that $\hat Z \in$ $\mathbb R_0^{n\times k} $. 

\newpage

\section{Supplementary Technical Lemmas}

\begin{lemma}[
\citealp{hillar2012inverses}, Theorem 1.1]
\label{lem:thm1.1_Hillar}
An $n\times n$ real matrix $J$ is called ``diagonally dominant'' if $|J_{ii}| - \sum_{j\neq i} |J_{ij}| \geqslant 0$ for $i=1, \ldots,n$.
Let $n\geqslant 3$. For any symmetric diagonally dominant matrix $J$ with $J_{ij} \geqslant l > 0$, we have
$$
\|J^{-1}\|_\infty \leqslant 
\frac{3n-4}{2l(n-2)(n-1)},
$$
where $\|A\|_\infty := \sup_{x\neq 0} \|Ax\|_\infty / \|x\|_\infty = \max_{1\leqslant i\leqslant m} \sum_{j=1}^n |a_{ij}| $ for any $A\in \mathbb R^{m\times n} $.
(In our proof we use a neater bound $\|J^{-1}\|_\infty \leqslant2/{(ln)} $ which holds when $n\geqslant 7$.)
\end{lemma}

\begin{lemma}[\citealp{nesterov2003introductory}, Theorem 2.1.12] \label{Lemma-30}
 For a continuously differentiable function $f$, if it is $\mu$-strongly convex and $L$-smooth on a convex domain $\mathcal{D}$, say for any $x,y\in \mathcal{D}$,
\begin{align*}
    \frac{\mu}{2}\|x-y\|_{\mathrm{2}}^2\leqslant
    f(y)-f(x)-\left\langle \nabla f(x),y-x\right\rangle
    \leqslant \frac{L}{2}\|x-y\|_{\mathrm{2}}^2,
\end{align*}
then
\begin{align*}
    \left\langle \nabla f(x)-\nabla f(y),x-y\right\rangle
    \geqslant \frac{\mu L}{\mu+L}\|x-y\|_{\mathrm{2}}^2
    +\frac{1}{\mu+L}\|\nabla f(x)-\nabla f(y)\|_{\mathrm{2}}^2.
\end{align*}
\end{lemma}

\begin{lemma}[\citealp{tu2016low}, Lemma 5.4]
\label{Lemma-28}
For any $Z_1,Z_2\in \mathbb{R}^{n\times k}$, let   $\operatorname{dist}(Z_1,Z_2)=\min_{Q \in \mathcal{O}(k)}\left\|Z_1-Z_2 Q\right\|_{\mathrm{F}}$, then 
\begin{align*}
\|Z_1Z_1^\top-Z_2Z_2^\top\|_{\mathrm{F}}^2\geqslant{2(\sqrt{2}-1)\sigma_k^2(Z_1)}\operatorname{dist}^2(Z_1,Z_2).
\end{align*}
\end{lemma}

\begin{lemma}[\citealp{tu2016low}, Lemma 5.3]
\label{Lemma-29}
For any $Z_1,Z_2\in \mathbb{R}^{n\times k}$, let   $\operatorname{dist}(Z_1,Z_2)=\min_{Q \in \mathcal{O}(k)}\left\|Z_1-Z_2 Q\right\|_{\mathrm{F}}$. If $\operatorname{dist}(Z_1,Z_2) \leqslant \|Z_1\|_{\mathrm{op}}$, then
\begin{align*}
\|Z_1Z_1^\top-Z_2Z_2^\top\|_{\mathrm{F}}\leqslant3\|Z_1\|_{\mathrm{op}}\operatorname{dist}(Z_1,Z_2).
\end{align*}
\end{lemma}

\begin{lemma}[Bernstein inequality.  \citealp{wellner2005empirical}, page 103-104]\label{lem:bernin}
	Suppose $X_1,\ldots,X_n$ are independent random variables with $\mathbb{E}X_i=0$ and $\mathbb{E}|X_i|^k\leqslant\frac{1}{2}\mathbb{E}X_i^2L^{k-2}k!$ for $k\geqslant2$. For $M\geqslant\sum_{1 \leqslant i\leqslant n}\mathbb{E}X_i^2$ and $x\geqslant0$,
	$$\Pr\left(\sum_{i = 1}^{n} X_i\geqslant x\right)\leqslant\exp\left(-\frac{x^{2}}{2(M+xL)}\right).$$
\end{lemma}

\begin{lemma}[\citealp{huang2018pairwise}, Lemma A.2.2]
\label{Lemma-A2}
Assume $A\sim \operatorname{Poisson}(\lambda)$ and let $X=A-\lambda$, then for any 
$0<\lambda<1/2$ and for any integer $\ell>2$,
\begin{align*}
\EXPT|X^{\ell}|\leqslant \frac{\ell!}{2}\left(\frac{e^2+1}{1-\lambda}\right)^{\ell -2}\lambda.
\end{align*}
\end{lemma}

\begin{lemma}[Bernstein inequality for Poisson random variables] \label{lem:bern_pois}
    Suppose $X_1,\ldots,X_n$ are independent Poisson random variables with $\mathbb{E}X_i=\lambda_i$. For $M\geqslant\sum_{1 \leqslant i\leqslant n}\lambda_i$ and $x\geqslant0$,
	$$\Pr\left(\sum_{i = n}^{n} (X_i-\lambda_i)\geqslant x\right)\leqslant\exp\left(-\frac{x^{2}}{2(M+10x)}\right).$$
\end{lemma}
\begin{proof}
    For $X_i\sim\text{Poisson}(\lambda_i)$, we can construct $Y_j \overset{i.i.d.}{\sim}\text{Poisson}(\lambda_i / N), j=1,\ldots,N$ for a large enough $N$ such that $\lambda_i / N<1/2$. Thus, a combination of Lemma \ref{lem:bernin} and Lemma \ref{Lemma-A2} gives the result of Lemma \ref{lem:bern_pois}.
\end{proof}

\begin{lemma}[\citealp{huang2018pairwise}, Theorem 15]
\label{Lemma-concen}
Let $A$ be a $n \times n$ symmetric random matrix with independent Poisson components. Set $\EXPT A = P = (p_{ij})_{i,j = 1, \ldots, n}$, then for any $s >0$ there exists a constant $C_s$ such that
$$
\Pr\left(\|A-P\|_{\operatorname{op}} \geqslant C_s \sqrt{n \max_{ij}p_{ij} + \log n}\right) \leqslant n^{-s}.
$$
\end{lemma}

\begin{lemma}[\citealp{tropp2012user}, Theorem 6.2]
\label{Lemma-MB}
Let $X_1,\ldots,X_m$ are independent, zero mean, $n\times n$ and symmetry
matrices, such that for $k = 1, \ldots, m$
$$\EXPT X_k^{\ell}\preceq \frac{\ell!}{2}B^{\ell-2}A_k^2 \quad\text{ for } \ell = 2, 3, \ldots.$$
Then for all $x\geqslant0$, we have 
\begin{align*}
    \Pr\left(\left\|\sum_{k=1}^m X_k\right\|_{\mathrm{op}}\geqslant x\right)
    \leqslant 2n\cdot\exp{\left(\frac{-x^2/2}{\sigma^2+B x}\right)}, 
\end{align*}
where $\sigma^2=\|\sum_{k=1}^m A_k^2\|_{\mathrm{op}}$.
\end{lemma}

\begin{lemma}[\citealp{chatterjee2015matrix}, Lemma 3.5]
\label{Lemma-3.5}
Suppose $A$ and $B$ are two $n \times m$ matrices with $n \geqslant m$. Let $A=\sum_{i=1}^m\sigma_i u_i v_i^\top$ be the SVD of $A$. 
For any fixed $\delta> 0$, define
\begin{align*}
\widehat{B}=\sum_{i:\sigma_i>(1+\delta)\|A-B\|_{\mathrm{op}}}\sigma_i u_i v_i^\top. 
\end{align*}
Then there exists a constant $C_{\delta} > 0$  depending only on $\delta$ such that
\begin{align*}
    \big\|\widehat{B}-B\big\|_{\mathrm{F}}^2
    \leqslant C_{\delta} \big\|A-B\big\|_{\mathrm{op}}  \big\|B\big\|_{\mathrm{*}}.
\end{align*}
\end{lemma}

\newpage
\section{Supplementary Numerical Results}


\subsection{Empirical Errors for Estimating \texorpdfstring{$G^\star$}{} by the Penalized MLE} \label{sec:suppsimultaion}
\begin{figure}[h]
\centering
\begin{subfigure}{0.32\textwidth}
		\centering
		\includegraphics[width=1\linewidth]{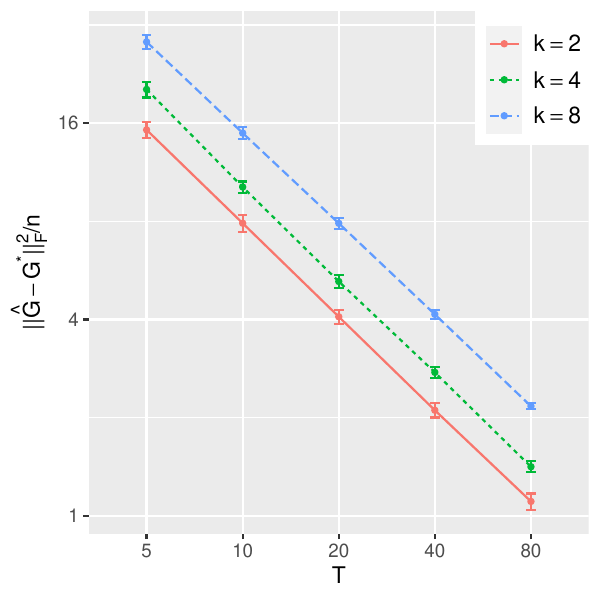}
    \caption{$\|\hat{G}-G^{\star} \|_{\mathrm{F}}^{ 2}/n$ versus $T$}
\end{subfigure}
 \begin{subfigure}{0.32\textwidth}
		\centering
		\includegraphics[width=1\linewidth]{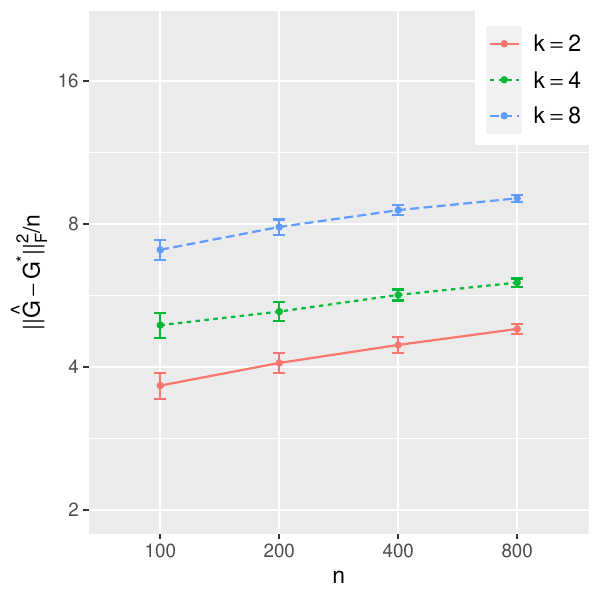}
    \caption{$\|\hat{G}-G^{\star} \|_{\mathrm{F}}^{ 2}/n$  versus $n$}
	\end{subfigure}
\begin{subfigure}{0.32\textwidth}
		\centering
		\includegraphics[width=1\linewidth]{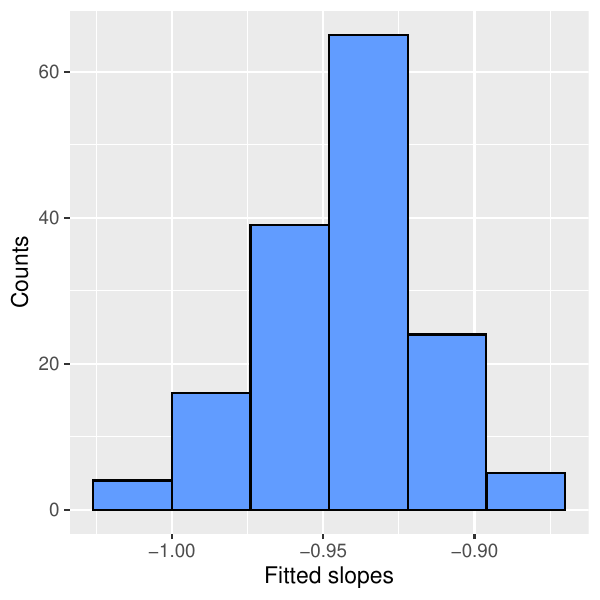}
    \caption{Histogram of fitted slopes}
	\end{subfigure}
\caption{Case (I): Empirical estimation errors of the penalized MLE. Panel (a) presents    $\|\hat{G}-G^{\star} \|_{\mathrm{F}}^{ 2}/n$ (averaged over 50 repetitions) versus $T$ in the scenario (a).  
 Panel (b) presents  $\|\hat{G}-G^{\star} \|_{\mathrm{F}}^{ 2}/n$ (averaged over 50  repetitions)  versus $n$ in the scenario (b). 
In (a) and (b), axes are in the \textit{log scale}, three lines correspond to results under $k\in \{2,4,8\}$, respectively, and error bars are obtained by $\pm$ the {standard deviation} from 50 repetitions.
 Panel (c) presents the slopes from regressing $\log (\|\hat{G}-G^{\star} \|_{\mathrm{F}}^{ 2}/n)$ on $\log T$ with fixed $(n,k)\in \{200\}\times \{2,4,8\}$ in the 50 repetitions under the scenario  (a). 
} \label{fig:resultscasea1gerror}
\end{figure}

\begin{figure}[!htbp]
\centering
\begin{subfigure}{0.32\textwidth}
		\centering
		\includegraphics[width=1\linewidth]{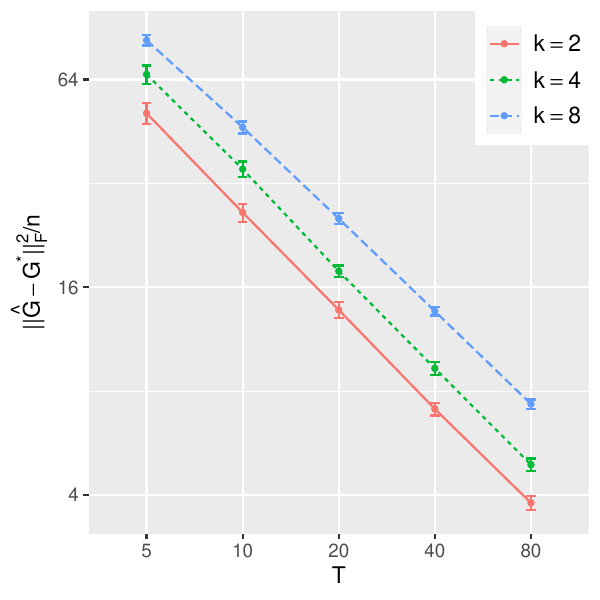}
    \caption{$\|\hat{G}-G^{\star} \|_{\mathrm{F}}^{ 2}/n$ versus $T$}
\end{subfigure}
 \begin{subfigure}{0.32\textwidth}
		\centering
		\includegraphics[width=1\linewidth]{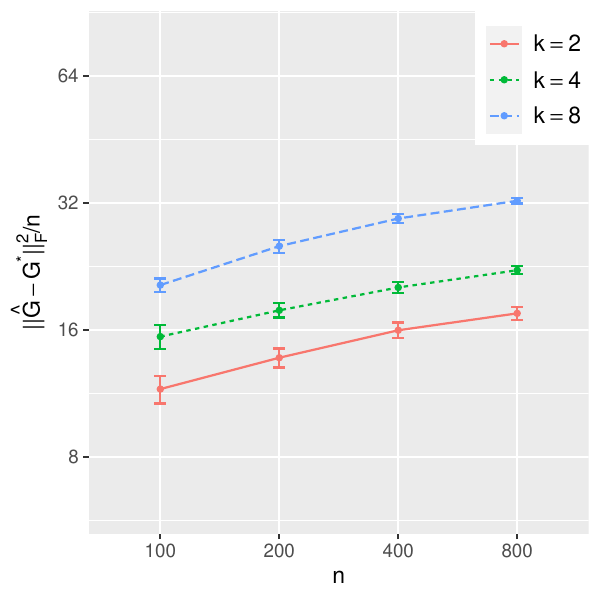}
     \caption{$\|\hat{G}-G^{\star} \|_{\mathrm{F}}^{ 2}/n$ versus $n$}
	\end{subfigure}
 \begin{subfigure}{0.32\textwidth}
		\centering
		\includegraphics[width=1\linewidth]{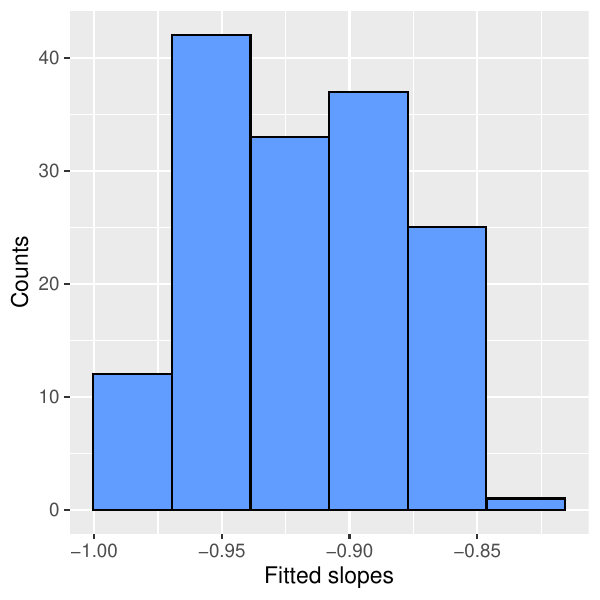}
\caption{Histogram of fitted slopes}
	\end{subfigure}
  \caption{Case (II): Empirical estimation errors of $\hat{G}$ by the penalized MLE. Panels (a)--(c) are presented similarly to Figure \ref{fig:resultscasea1gerror}.}  \label{fig:resultscase2convexgerror}
\end{figure}

Recall cases (I) and (II) introduced in Section \ref{sec::simulations} in the main text. 
We next present the estimation error $\|\hat{G}-G^{\star} \|_{\mathrm{F}}^{\blue 2}/n$ obtained by the penalized MLE under cases (I) and (II) in Figures \ref{fig:resultscasea1gerror} and \ref{fig:resultscase2convexgerror}, respectively. 
Similarly to Remark \ref{rm:choiceoflambda}, we choose $\lambda_{n,T}=c_{\lambda}\sqrt{nT\hat{\mu}}$.
When estimating $G^{\star}$, we do not employ the truncation step in Remark \ref{rm:fromgtoz}, which is utilized for estimating $Z^{\star}$.  
Consequently, in this section, we choose a larger constant $c_{\lambda}=1.5$ to achieve a higher level of shrinkage, compensating for the absence of the truncation step.


In each figure,
panel (a) suggests that $\|\hat{G}-G^{\star} \|_{\mathrm{F}}^{ 2}/n$ is inverse proportional to $T$ when $n$ is fixed. 
Particularly, panel (c) shows that all the fitted slopes are close to $-1$. 
In addition, 
panel (b) shows that $\|\hat{G}-G^{\star} \|_{\mathrm{F}}^{ 2}/n$ increases slowly as $n$ increases. 
Particularly, the slope of each curve in the panel (b) is bounded by $\log^2(n)$. 
{We emphasize that in panels (a) and (b), axes are in the log scale for clear visualization.}   
Overall, the numerical results support the theoretical error bound in Theorem \ref{MLE_mainthm} in the main text. 
It shows that up to logarithmic factors, the penalized MLE can achieve the oracle rate  $O(1/T)$ for the estimation error of $G^{\star}$.  

\newpage
\subsection{Simulation Results with \texorpdfstring{$T>n$}{}}\label{sec:suppsimultaiontn}
\begin{figure}[!htbp]
	\centering
\begin{subfigure}{0.32\textwidth}
		\centering
		\includegraphics[width=1\linewidth]{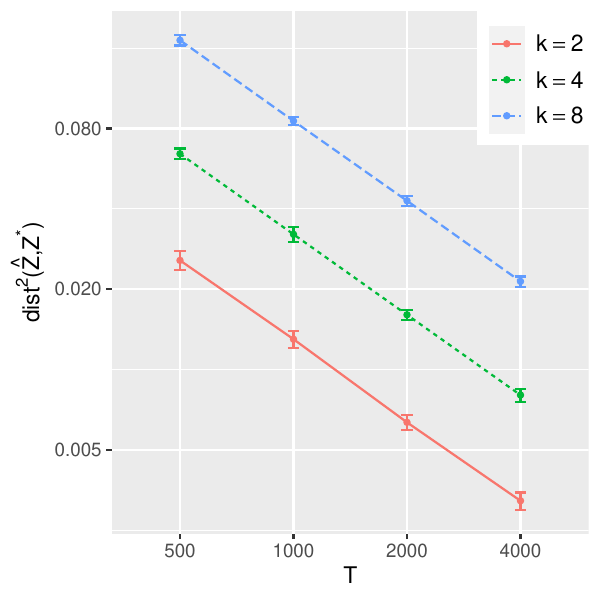}
  \caption{$\operatorname{dist}^2(\hat Z, Z^\star)$ versus $T$}
\end{subfigure}
 \begin{subfigure}{0.32\textwidth}
		\centering
		\includegraphics[width=1\linewidth]{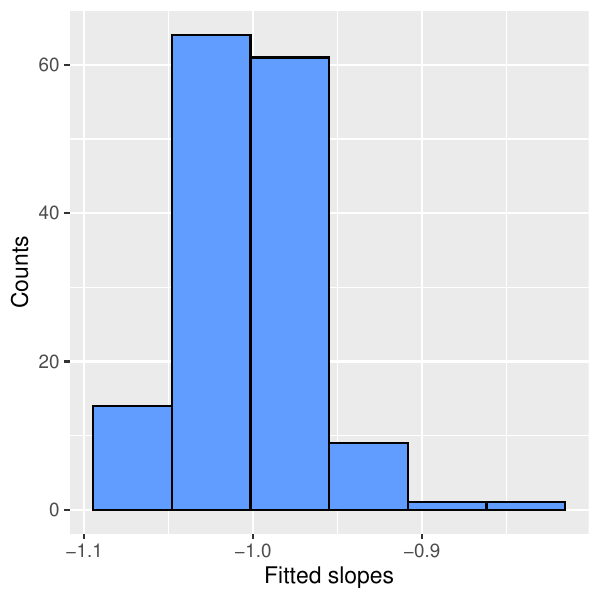}
    \caption{Histogram of fitted slopes} 
\end{subfigure}
 \caption{Case (I): Empirical estimation errors of the one-step estimator. Panel (a) presents    $\operatorname{dist}^2(\hat Z, Z^\star)$ (averaged over 50 repetitions) versus $T$ in the scenario (a). Axes are in the log scale, three lines correspond to results under $k\in \{2,4,8\}$, respectively, and error bars are obtained by $\pm$ the {standard deviation} from 50 repetitions. Panel (b) presents the slopes from regressing $\log \operatorname{dist}^2(\hat Z, Z^\star)$ on $\log T$ with fixed $(n,k)\in \{200\}\times \{2,4,8\}$ in the 50 repetitions under the scenario  (a).}  \label{fig:tlarge1onestep_supp}
\end{figure}

\begin{figure}[!htbp]
	\centering
\begin{subfigure}{0.32\textwidth}
		\centering
		\includegraphics[width=1\linewidth]{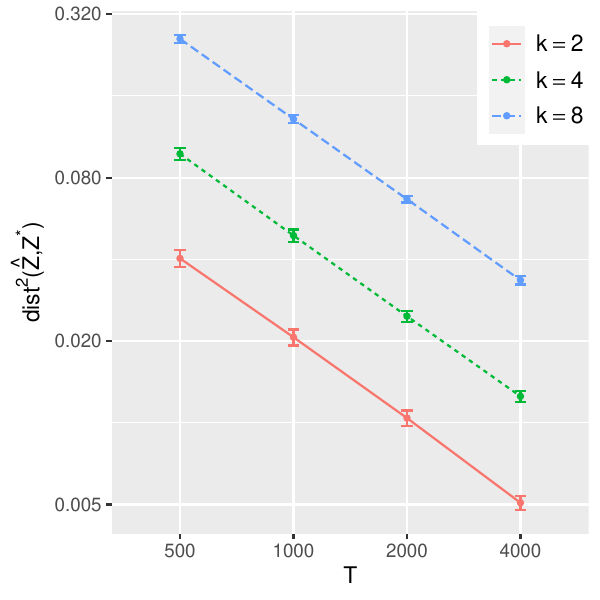}
  \caption{$\operatorname{dist}^2(\hat Z, Z^\star)$ versus $T$}
\end{subfigure}
 \begin{subfigure}{0.32\textwidth}
		\centering
		\includegraphics[width=1\linewidth]{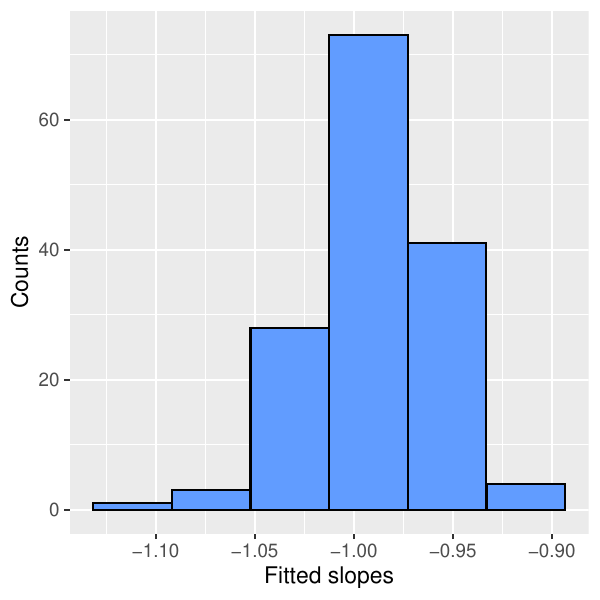}
    \caption{Histogram of fitted slopes} 
\end{subfigure}
 \caption{Case (I): Empirical estimation errors of the penalized MLE. Panels (a)--(b) are presented similarly to Figure \ref{fig:tlarge1onestep_supp}.}  \label{fig:tlarge1penalized_supp}
\end{figure}

\begin{figure}[!htbp]
	\centering
\begin{subfigure}{0.32\textwidth}
		\centering
		\includegraphics[width=1\linewidth]{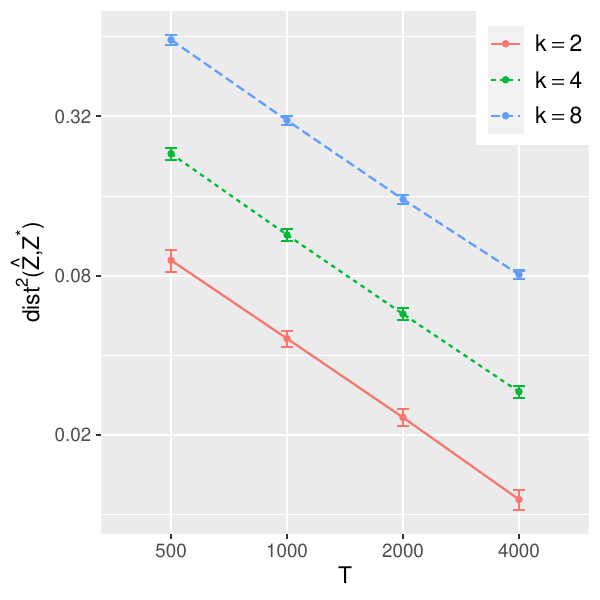}
  \caption{$\operatorname{dist}^2(\hat Z, Z^\star)$ versus $T$} 
\end{subfigure}
 \begin{subfigure}{0.32\textwidth}
		\centering
		\includegraphics[width=1\linewidth]{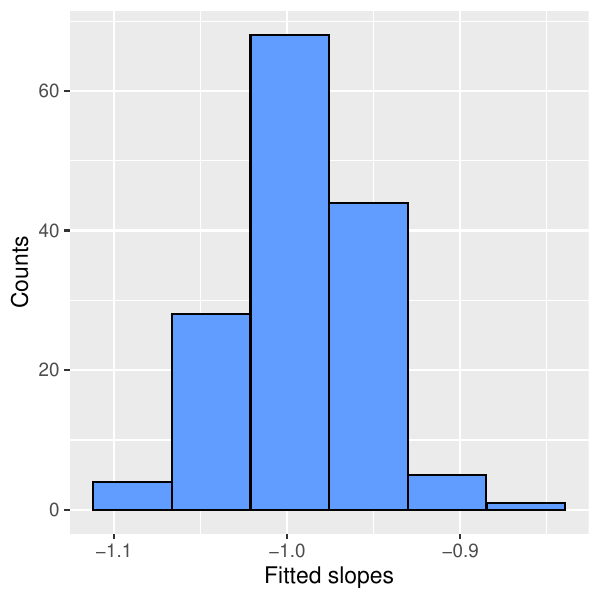}
    \caption{Histogram of fitted slopes} 
\end{subfigure}
 \caption{Case (II): Empirical estimation errors of the one-step estimator. Panels (a)--(b) are presented similarly to Figure \ref{fig:tlarge1onestep_supp}.}   \label{fig:tlarge2onestep_supp}
\end{figure}

\begin{figure}[!htbp]
	\centering
\begin{subfigure}{0.32\textwidth}
		\centering
		\includegraphics[width=1\linewidth]{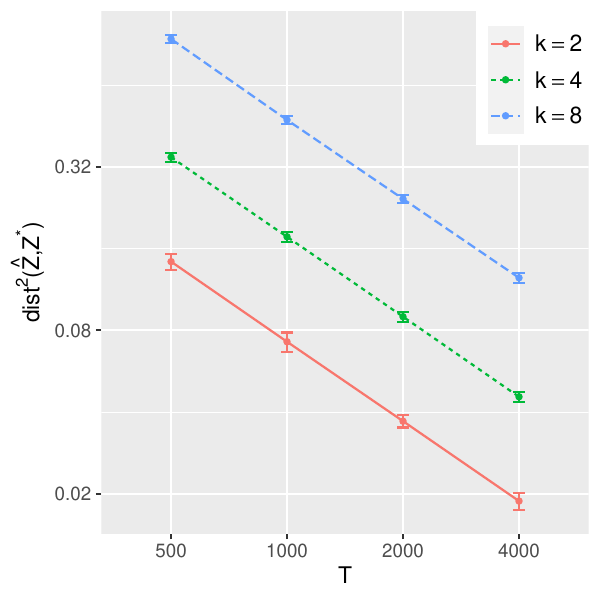}
  \caption{$\operatorname{dist}^2(\hat Z, Z^\star)$ versus $T$}
\end{subfigure}
 \begin{subfigure}{0.32\textwidth}
		\centering
		\includegraphics[width=1\linewidth]{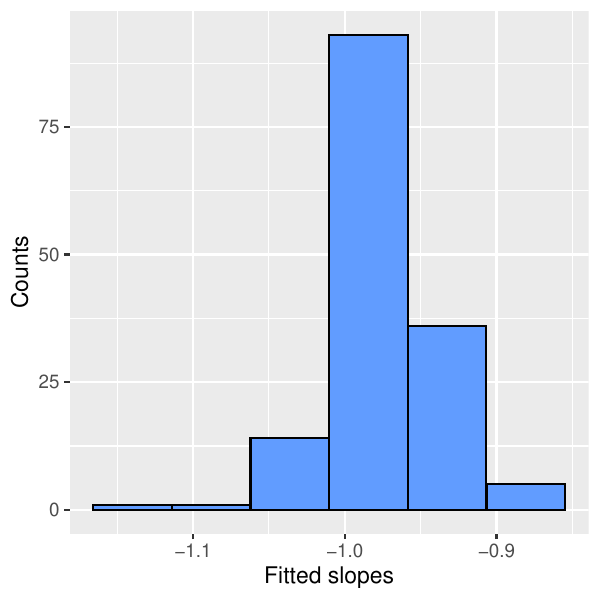}
    \caption{Histogram of fitted slopes} 
\end{subfigure}
 \caption{Case (II): Empirical estimation errors of the penalized MLE. Panels (a)--(b) are presented similarly to Figure \ref{fig:tlarge1onestep_supp}.} \label{fig:tlarge2penalized_supp}
\end{figure}

In this section, 
we present numerical results  to demonstrate that the ``inverse proportional'' pattern in Section \ref{sec::simulations} holds similarly when $T>n$.  
  In particular,  we set $n=200$ and $T \in \{500,1000,2000, 4000\}$. Then, we generate data and apply the methods in the same way as that in Section \ref{sec::simulations}. 
  The numerical results are presented in Figures \ref{fig:tlarge1onestep_supp}--\ref{fig:tlarge2penalized_supp}, respectively, showing that the estimation errors are inversely proportional to $T$ when $n$ is fixed and $T>n$.  

\newpage
\subsection{Analysis of Citi Bike Data by the Penalized MLE}
\begin{figure}[!htbp]
\begin{subfigure}{0.48\textwidth}
		\centering
		\includegraphics[width=0.5\linewidth]{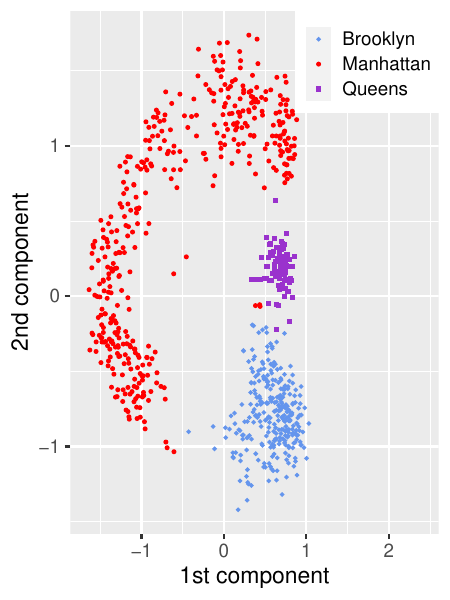}
    \caption{Estimated Latent Positions by 24-Hour Data}
\end{subfigure}
\begin{subfigure}{0.5\textwidth}
		\centering
		\includegraphics[width=0.64\linewidth]{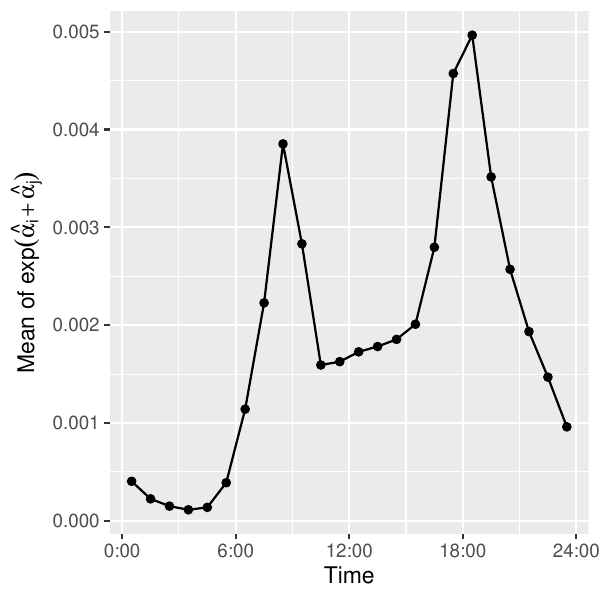}
    \caption{Estimated Baseline Levels}
	\end{subfigure}
  \caption{Analysis Results of Citi Bike Data by the Penalized MLE: Panel (a) shows the estimated latent positions $\hat{z}_i\in \mathbb{R}^2$ using all the data across the 24-hour period. Each point is colored based on which borough the corresponding bike station is located, and we present the results after a rotation for better visualization. Panel (b) presents the mean of estimated baseline levels $\sum_{i,j=1}^n\exp(\hat{\alpha}_{it}+\hat{\alpha}_{jt})/n^2$ for $t=1,\ldots, 24$, i.e., across the 24-hour period.} \label{fig:estfigpenalizedmle}
\end{figure}

In this section, we analyze the New York Citi Bike Dataset by the penalized MLE.  
Following Remark \ref{rm:fromgtoz}, we obtain estimated latent space vectors with $k=2$ and set $\lambda_{n,T}=0.5\sqrt{nT\hat{\mu}}$ as in Remark \ref{rm:choiceoflambda}. 
We visualized the estimates  $\hat{z}_i$'s  in Figure \ref{fig:estfigpenalizedmle} (a).     
Comparing Figure \ref{fig:estfigpenalizedmle} with Figure \ref{fig:citidata} in the main text,
we see that 
the estimated latent vectors display three clusters of points that are consistent with the three boroughs of New York City.
Moreover, we presented the estimated averaged baseline levels in Figure \ref{fig:estfigpenalizedmle} (b). 
We can see that the estimated averaged baseline levels across 24 hours align with the total number of ride counts over 24 hours.  
The results suggest that the penalized MLE returned reasonable and interpretable estimated parameters under the model \eqref{model_DLSM}. 
The analysis reflects both the static latent space information and the time heterogeneity simultaneously.

\newpage
\subsection{Simulation Results for Remark \ref{rm:choiceoflambda}
} \label{sec:lambdatuning}

\begin{figure}[!htbp]
	\centering
\begin{subfigure}{0.32\textwidth}
		\centering
		  \caption{$c_{\lambda} = 1.0$}
	\includegraphics[width=1\linewidth]{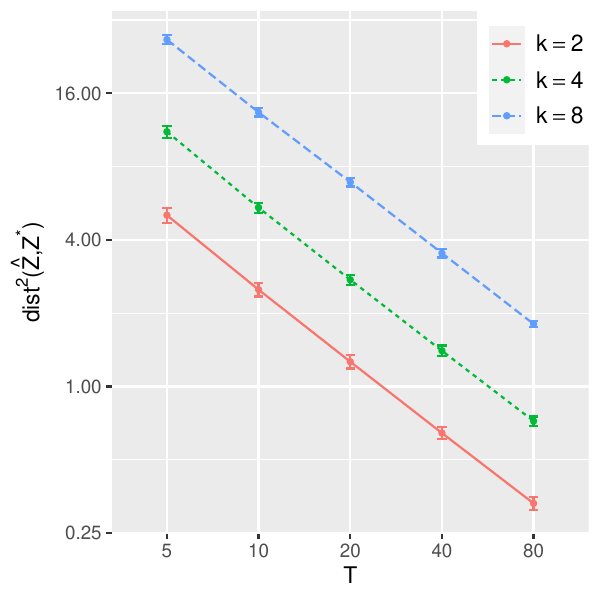}
\end{subfigure}
 \begin{subfigure}{0.32\textwidth}
		\centering
  \caption{$c_{\lambda} = 1.5$}
		\includegraphics[width=1\linewidth]{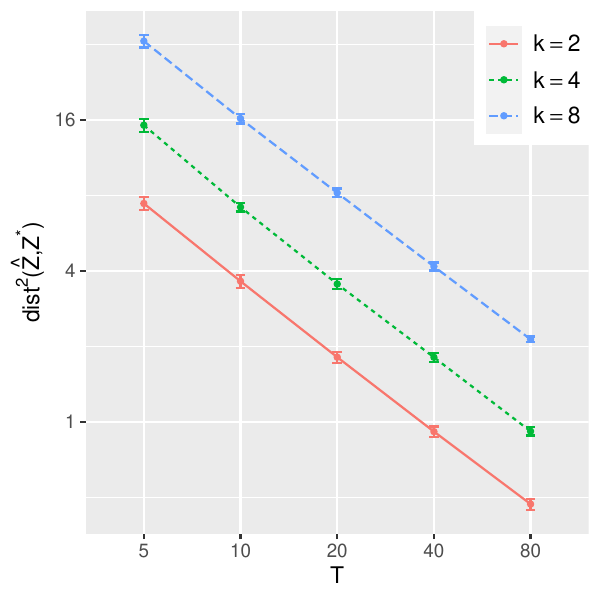}
\end{subfigure}
\begin{subfigure}{0.32\textwidth}
  \caption{$c_{\lambda} = 2.0$}
		\centering
		\includegraphics[width=1\linewidth]{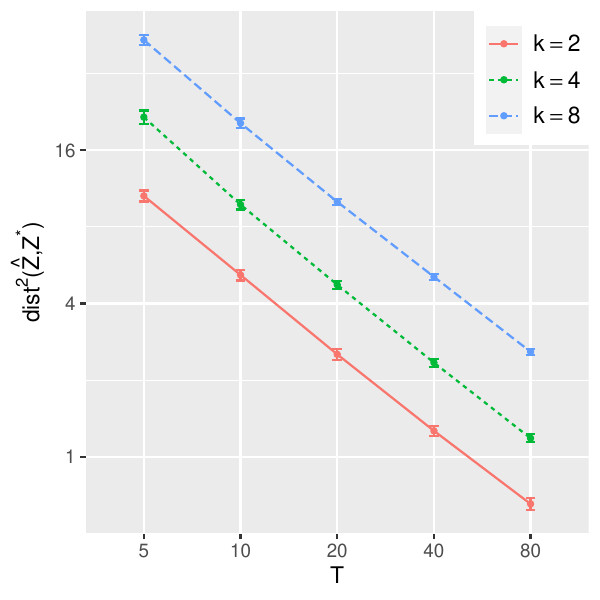}
\end{subfigure}
 \caption{Empirical estimation errors of the penalized MLE versus $T$ under Case (I) with $n=200$ and different $c_{\lambda}\in \{1.0, 1.5, 2.0\}$.}  \label{fig:errvsTlambda_supp}
\end{figure}

In this section, 
we present simulation studies to demonstrate Remark \ref{rm:choiceoflambda}. 
Specifically,
we apply the penalized MLE with  $\lambda_{n,T}= c_{\lambda}\sqrt{nT\hat{\mu}}$ and varying constants $c_{\lambda} \in \{1.0, 1.5, 2.0\}$.  
As an illustration, we generate data under Case (I) in Section \ref{sec::simulations} with $n=200$,  $T\in \{5, 10, 20, 40, 80\}$, and $k\in \{2,4,8\} $. 
Figure \ref{fig:errvsTlambda_supp} presents the estimated error $\operatorname{dist}^2(\hat Z, Z^\star)$  over 50 Monte Carlo simulations. 
It shows that the 
``inverse proportional''
pattern of the estimation error versus $T$ remains the same under different $c_{\lambda}$ values. 

\vspace{1em}

\newpage 
\section{Additional Results for Remark \ref{rm:fromgtoz}}\label{sec:estimatekdetails}

In this section, we provide  details for estimating $k$ and obtaining $\hat{Z}_k$ from $\hat{G}$, as discussed in Remark \ref{rm:fromgtoz}.  
Specifically, 
 Section \ref{sec:estimatektheory}  presents detailed 
 procedures 
and theoretical results, and  
Section \ref{sec:simuestk}  showcases simulation studies.  

\subsection{Theoretical Results}\label{sec:estimatektheory}

As mentioned in Remark \ref{rm:fromgtoz},  
we estimate $k$ as the number of significant eigenvalues of an estimator $\Ginit$ that approximates $G^{\star}$ sufficiently well.
We next provide a  general spectral strategy with the theoretical guarantee given in Proposition \ref{prop_choosek2supp}.


\paragraph*{Spectral method   for estimating $k$}
Consider an estimator $\Ginit\in\mathbb{S}_+^n$ that approximates   $G^{\star}$ sufficiently well in the sense that $\|\Ginit-G^{\star}\|_{\mathrm{F}}/n = O_p(\varepsilon_{n,T})$ with  $\varepsilon_{n,T}\to 0$ as $n,T\to \infty$. 
Then, we construct  an estimator for $k$ as 
\begin{equation}
    \label{eq::khat_defsupp}
    \hat k = \min\left\{1\leqslant l \leqslant n-1: \ \frac{\sigma_{l+1}(\Ginit) }{ \sigma_{l}(\Ginit)} \leqslant \tau_{n,T,l} \right\}, 
\end{equation}  
where $\sigma_l(G)$ denotes the $l$-th largest eigenvalue of a symmetric matrix $G$, and
$\tau_{n,T,l}$ is a threshold set to satisfy 
\begin{align}\label{eq:tauthres}
	\tau_{n,T,l}\to 0 \hspace{1em} \text{and} \hspace{1em} \tau_{n,T,k}\gg \varepsilon_{n,T}
\end{align}
as $n$ and $T\to \infty$. 

\begin{proposition}
\label{prop_choosek2supp} 
Under Condition \ref{cond:truvalueregularity}, 
the aforementioned estimator $\hat k$ satisfies that $\Pr(\hat k = k) \to 1 $ as $n,T\to\infty$, i.e., $\hat k$ 
can consistently estimate $k$. 
\end{proposition}

\medskip
We next show that   $\Ginit$ can be set as  the ``double-SVD'' estimator $\mathring{G}$ in Section \ref{sec:initest} 
or the penalized MLE $\hat{G}$ defined in 
\eqref{eq:penalizedMLEproblem}. 
Remarkably, Proposition \ref{prop::choosek_examplessupp} below proves that the general conclusion in Proposition \ref{prop_choosek2supp} is applicable  for these two estimators 
using thresholds that are independent of  unknown parameters.  
 

\smallskip

\begin{proposition} \label{prop::choosek_examplessupp} \quad

\begin{itemize}
	\item[(a)] 
	Let $\Ginit = \mathring G$. 
	Under the conditions of Theorem \ref{thm:initialerror}, 
	we have $\|\Ginit-G^\star \|_{\mathrm{F}}/n  = O_p(\varepsilon_{n,T})$ with $	\varepsilon_{n,T} = n^{-\frac{1}{2(k+3)}}\log(nT) $.  
In this case, the threshold 
\begin{align*}
	\tau_{n,T,l}= n^{-\frac{1}{2l+7}}
\end{align*}
 satisfies \eqref{eq:tauthres} 
when $\log(T)\ll n^{\frac{1}{2(k+3)(2k+7)}}$. 

\item[(b)]   Let $\Ginit=\hat{G}$. 
Under the conditions of Theorem \ref{MLE_mainthm}, we have $\|\Ginit-G^\star \|_{\mathrm{F}}/n  = O_p(\varepsilon_{n,T})$ with $\varepsilon_{n,T}= n^{-1/2}\max\{n^{-1/2}, T^{-1/2}\} \log(nT)$. In this case, the threshold 
\begin{align*}
	\tau_{n,T,l}= n^{-\frac{1}{l+2}}
\end{align*}
 satisfies \eqref{eq:tauthres} 
when $\log(T) \ll n^{\frac{k}{2(k+2)}}$.
\end{itemize}
\end{proposition}
 


\vspace{1em}
Given an estimator  for $G^\star$, 
we can construct an estimator for $Z^{\star}$ via truncation as discussed in Remark \ref{rm:fromgtoz}. Next, we will develop a theoretical guarantee for this strategy. 
We point out that the following proposition holds for a general estimator $\Ginit$, which includes the penalized MLE $\hat{G}$ or the ``double-SVD'' estimator $\mathring G$ in Remark \ref{rm:fromgtoz} as special cases.  


\begin{proposition}
    \label{prop::truncated_PMLE}
    Given an estimator $\Ginit$ of $G^\star$, let  
$\hat{Z}_k=\hat{U}_k\hat{D}_k^{1/2}$, 
where $\hat{U}_k\hat{D}_k\hat{U}_k^{\top}$ represents the top-$k$ eigenvalue components of  $\Ginit$.
Under Condition \ref{cond:truvalueregularity}, 
we have
\begin{align*}
 \|\hat{Z}_{k}\hat{Z}_{k}^{\top}- G^\star\|_{\mathrm{F}} \leqslant 2\|\Ginit- G^\star\|_{\mathrm{F}}  \hspace{1em} \text{ and }\hspace{1em}\operatorname{dist}^2(\hat{Z}_k,Z^\star) \leqslant \frac{2 \|\Ginit- G^\star\|_{\mathrm{F}}^2 }{(\sqrt{2}-1)n M_{Z,2} }, 
\end{align*}
where $M_{Z,2}$ is specified same as in Condition \ref{cond:truvalueregularity}. 
\end{proposition}  
\begin{proof}
By Eckart-Young-Mirsky theorem, we have 
$\|\hat{Z}_{k}\hat{Z}_{k}^{\top}- \Ginit \|_{\mathrm{F}} \leqslant \|\Ginit- G^\star\|_{\mathrm{F}}$, which gives  
$\|\hat{Z}_{k}\hat{Z}_{k}^{\top}- G^\star\|_{\mathrm{F}} \leqslant 2\|\Ginit- G^\star\|_{\mathrm{F}}$. 
A direct application of Lemma \ref{Lemma-28} gives the bound on $\operatorname{dist}^2(\hat{Z}_k,Z^\star)$.
\end{proof}


\subsubsection{Proof of Proposition \ref{prop_choosek2supp}}
Since $ \varepsilon_{n,T}\ll\tau_{n,T,k}\ll 1  $ as $n,T\to\infty$, we have $\|\Ginit- G^\star\|_{\mathrm{F}} / n = o_p(\tau_{n,T,k} ) $ with $\tau_{n,T,k}  \to 0 $ as $n,T\to\infty$.
This implies that for any constant $a_0$, we have $\Pr( \|\Ginit- G^\star\|_{\mathrm{F}} / n \leqslant a_0 \tau_{n,T,k}  ) \to 1 $.
In what follows we condition on the event $\|\Ginit- G^\star\|_{\mathrm{F}} / n \leqslant a_0 \tau_{n,T,k}  $, where the choice of the constant $a_0$ is specified later.  By Weyl's inequality, we have  
\begin{align*}
    &\sigma_{1}(\Ginit) \leqslant \sigma_{1}(G^\star) + \sigma_1(\Ginit-G^\star ) \leqslant \|G^\star\|_{\mathrm{F}} + \|G^\star - \Ginit\|_{\mathrm{F}} \leqslant M_{Z,1}n + a_0 n\tau_{n,T,k} ;\\
    &\sigma_k(\Ginit) \geqslant \sigma_k(G^\star) - \sigma_1(G^\star - \Ginit) \geqslant \sigma_k(G^\star) - \|G^\star - \Ginit\|_{\mathrm{F}} \geqslant M_{Z,2} n - a_0 n\tau_{n,T,k} ; \\
    &\sigma_{k+1}(\Ginit) \leqslant \sigma_{k+1}(G^\star) + \sigma_1(\Ginit-G^\star ) \leqslant 0 + \|G^\star - \Ginit\|_{\mathrm{F}} \leqslant a_0 n\tau_{n,T,k} .
\end{align*} 
Since $\tau_{n,T,k}\ll 1$, we have $\sigma_1(\Ginit) \leqslant (3M_{Z,1}/2) n $ and $\sigma_k(\Ginit) \geqslant (M_{Z,2}/2) n$. 
For any $1\leqslant l < k $, we have $ \sigma_{l+1}(\Ginit) /  \sigma_{l}(\Ginit) \geqslant \sigma_k(\Ginit) /  \sigma_{1}(\Ginit) \geqslant \frac{M_{Z,2}}{3M_{Z,1}} 
\gg \tau_{n,T,l} $. On the other hand, for $l=k$, we have $ \sigma_{l+1}(\Ginit) /  \sigma_{l}(\Ginit) =  \sigma_{k+1}(\Ginit) /  \sigma_{k}(\Ginit) \leqslant   \frac{2a_0}{M_{Z,2}}\tau_{n,T,k}   \leqslant   \tau_{n,T,k}  = \tau_{n,T,l} $; here, we have chosen $a_0=\frac{M_{Z,2}}{2} $ since $a_0$ can be any constant in the above arguments. Hence, by the definition of $\hat k$, 
we have $\hat k = k$ given that $\|\Ginit- G^\star\|_{\mathrm{F}} / n \leqslant \frac{M_{Z,2}}{2}  \tau_{n,T,k}  $.
In summary, we have shown that 
$$\Pr(\hat k = k) \geqslant \Pr( \|\Ginit- G^\star\|_{\mathrm{F}} / n \leqslant \frac{M_{Z,2}}{2}  \tau_{n,T,k}  ) \to 1.$$

\subsubsection{Proof of Proposition \ref{prop::choosek_examplessupp}}
  
\textit{(a).}
Combining \eqref{th2}, \eqref{eq:gringbd1} and Lemma \ref{lem:concen}, 
we obtain that 
under the conditions of Theorem \ref{thm:initialerror}, 
there exists a constant $C>0$ such that for any $s>0$
\begin{align*}
	\Pr\left\{ \|\mathring{G}-G^\star \|_{\mathrm{F}}/n > C(s+1)n^{-\frac{1}{2(k+3)}}\log(nT) \right\}\leqslant (nT)^{-s}, 
\end{align*} 
i.e., $\|\mathring{G}-G^\star \|_{\mathrm{F}}/n = O_p(\varepsilon_{n,T} ) $ with $\varepsilon_{n,T} = n^{-\frac{1}{2(k+3)}}\log(nT)$.
	  Therefore, the choice of $\tau_{n,T,l}  =  n^{-\frac{1}{2l+7}} $ satisfies $\tau_{n,T,l} \to 0  $ and $\tau_{n,T,k} \gg \varepsilon_{n,T} $ when $n \gg \log^{(2k+6)(2k+7)}(T) $. 

\bigskip
   
\noindent  \textit{(b).}  From Theorem \ref{MLE_mainthm}, we have that for any constant $s>0$, there exists a constant $C_s>0$ such that
$$
\Pr\left\{ \|\hat G - G^\star \|_{\mathrm{F}} / n > C_s \left(\frac{1}{\sqrt{nT}} + \frac{1}{n}\right)\log(nT) \right\}= O(n^{-s}).
$$
	i.e., $\|\hat{G}-G^\star \|_{\mathrm{F}}/n = O_p(\varepsilon_{n,T} ) $ with $\varepsilon_{n,T} = n^{-\frac{1}{2}}\log(nT)$.
	  Therefore, the choice of $\tau_{n,T,l}  =  n^{-\frac{1}{l+2}} $ satisfies $\tau_{n,T,l} \to 0  $ and $\tau_{n,T,k} \gg \varepsilon_{n,T} $ when $n \gg \log^{\frac{2(k+2)}{k}}(T) $. 

\subsection{Simulation Results} \label{sec:simuestk} 
In this section, 
we conduct simulation studies showing that $k$ can be correctly estimated with high probability.  
To this end, we consider  $n\in \{200, 400\}$, $T\in \{5,10,20\}$, and   $k\in \{2,4,8\} $. 
Given each $(n,T,k)$, 
we generate data 
under Case (I) outlined in Section \ref{sec::simulations}. 
Based on the construction in Section \ref{sec:estimatektheory}, 
we define 
\begin{align*}
    \hat k = \left\{\begin{aligned}
        &\min\{1\leqslant l \leqslant n-1:\  \sigma_{l+1}(\mathring G) / \sigma_{l}(\mathring G) \leqslant n^{-\frac{1}{2l+7}} \} \quad & \text{ for  }\Ginit= \mathring G,\\
        &\min\{1\leqslant l \leqslant n-1:\  \sigma_{l+1}(\hat G) / \sigma_{l}(\hat G) \leqslant n^{-\frac{1}{l+2}} \} \quad & \text{ for }\Ginit= \hat G.
    \end{aligned}\right.
\end{align*}
We estimate the probability of $\hat{k}=k$ over 100 simulations,
and present results for $n=200$ and $n=400$ in Tables \ref{table::estim_k_n200} and \ref{table::estim_k_n400}, respectively. 
Note that the penalized MLE $\hat{G}$ depends on the tuning parameter $\lambda_{n,T}= c_{\lambda}\sqrt{nT}$ as pointed out in Remark \ref{rm:choiceoflambda}. 
We let $\hat{G}(c_{\lambda})$ represent the estimator obtained with a given $c_{\lambda}$ for varying  $c_{\lambda}\in \{0.4,0.8,1.6\} $. 
It is noteworthy that results are robust to the choice of $c_{\lambda}$ as $n$ becomes larger.

\vspace{1em}
\setlength{\tabcolsep}{12pt} 
\renewcommand{\arraystretch}{1}
\begin{table}[H]
    \centering
    \begin{tabular}{M{10pt} M{10pt}|cccc}
    \hline
    $k$ & $T$ & 
  $\Ginit =\mathring G$ & $\Ginit =\hat{G}(0.4)$ & $\Ginit =\hat{G}(0.8)$  & $\Ginit =\hat{G}(1.6)$ \rule[-1ex]{0pt}{3.7ex} \\
    \hline
         & $5$  & 1.00 & 1.00 &1.00 & 1.00 \\
      $2$   & $10$  & 1.00 & 1.00  & 1.00 & 1.00 \\
         & $20$  & 1.00 & 1.00 & 1.00 & 1.00 \\
        \hline
            & $5$  & 1.00 & 1.00 &1.00 & 1.00 \\
      $4$   & $10$ & 1.00 & 1.00 &1.00 & 1.00\\
         & $20$  & 1.00 & 1.00 &1.00 & 1.00 \\
        \hline
                  & $5$  & 1.00 & 1.00 &1.00 & 0.97 \\
      $8$  & $10$  & 0.98 & 1.00 & 1.00 & 0.99 \\
         & $20$  & 0.95 & 1.00 & 1.00 &1.00  \\
        \hline
    \end{tabular}
    \caption{Empirical proportions of simulations satisfying $\hat k = k$ over 100 repetitions under $n=200$.} 
    \label{table::estim_k_n200} 
\end{table}

\setlength{\tabcolsep}{12pt} 
\renewcommand{\arraystretch}{1}
\begin{table}[H]
    \centering
    \begin{tabular}{M{10pt} M{10pt}|cccc}
    \hline
    $k$ & $T$ & 
  $\Ginit =\mathring G$ & $\Ginit =\hat{G}(0.4)$ & $\Ginit =\hat{G}(0.8)$  & $\Ginit =\hat{G}(1.6)$ \rule[-1ex]{0pt}{3.7ex} \\
    \hline
         & $5$  & 1.00 & 1.00 &1.00 & 1.00 \\
      $2$   & $10$  & 1.00 & 1.00  & 1.00 & 1.00 \\
         & $20$  & 1.00 & 1.00 & 1.00 & 1.00 \\
             \hline
         & $5$  & 1.00 & 1.00 &1.00 & 1.00 \\
      $4$   & $10$  & 1.00 & 1.00  & 1.00 & 1.00 \\
         & $20$  & 1.00 & 1.00 & 1.00 & 1.00 \\
             \hline
         & $5$  & 1.00 & 1.00 &1.00 & 1.00 \\
      $8$   & $10$  & 1.00 & 1.00  & 1.00 & 1.00 \\
         & $20$  & 1.00 & 1.00 & 1.00 & 1.00 \\
        \hline
    \end{tabular}
    \caption{Empirical proportions of simulations satisfying $\hat k = k$ over 100 repetitions under $n=400$.} 
    \label{table::estim_k_n400} 
\end{table}

\newpage

\section{Generalizability under Heterophilic Products}\label{sec:genhetero}

In this paper, we focus on the model \eqref{model_DLSM} with Euclidean inner product kernel between latent vectors $\langle z_i, z_j\rangle$. 
As discussed in Section \ref{sec:detailedcomareliterature}, 
heterophilic and time-varying products have also been widely adopted to model multiple networks, e.g., \cite{zhang2020flexible} considered $\langle z_i, z_j\rangle_{\Lambda_t}$ with symmetric kernel matrices $\Lambda_t \in \mathbb{R}^{k\times k}$. 
Such a model could allow more complex relationships between latent variables, including heterophilic structure. 
Following the idea, the model \eqref{model_DLSM} could be generalized 
by considering 
\begin{align}\label{eq:m2}
    \EXPT(A_{t,ij})=\exp(\alpha_{i t}+\alpha_{j t}+z_i^{\top} \Lambda_t z_j),
\end{align}
and when  $\Lambda_t=\mathrm{I}_k$, \eqref{eq:m2} reduces to \eqref{model_DLSM}. 
It is feasible to extend the analysis  under \eqref{model_DLSM} to \eqref{eq:m2}. 
To illustrate the idea, we next discuss the connections and differences between the two models in terms of identifiability, estimation,  theory, and data analysis in Sections \ref{sec:modelhetero}--\ref{sec:datahetero}, respectively. 

\subsection{{Identifiability}}\label{sec:modelhetero}
 To analyze the  latent vectors $Z$, the similarity under \eqref{model_DLSM} and \eqref{eq:m2} is that $Z=[z_1,\ldots, z_n]^{\top}$ is only identifiable up to an orthogonal transformation.  However, the constraints for identifiability under the two models are different, leading to distinct spaces to consider in the analysis.

\quad In particular, for \eqref{model_DLSM}, we consider the parameter space     
\begin{itemize}[leftmargin=3.4em]
\setlength{\itemsep}{0pt}
	\item[(A1)]  $\alpha_t\in \mathbb{R}^n$ for $t=1,\ldots, T$, 
	\item[(A2)]  $Z\in \mathbb{R}_0^{n \times k}$, where $ \mathbb{R}_0^{n \times k}=\left\{Z \in \mathbb{R}^{n \times k}: 1_n^{\top} Z=0,\,  \operatorname{det}\left(Z^{\top} Z\right) \neq 0\right\}$. 
\end{itemize}
When two sets of parameters $(\alpha, Z)$ and $(\tilde{\alpha}, \tilde{Z})$ satisfying (A1) and (A2)  generate the same model, i.e., $\exp \left(\alpha_{i t}+\alpha_{j t}+\left\langle z_i, z_j\right\rangle\right) = \exp \left(\tilde{\alpha}_{i t}+\tilde{\alpha}_{j t}+\left\langle \tilde{z}_i, \tilde{z}_j\right\rangle\right),$
we have $\alpha_t=\tilde{\alpha}_t,$ and there exists $Q\in \mathcal{O}(k) $   such that  $Z=\tilde{Z}Q.$
 For \eqref{eq:m2}, following  \cite{zhang2020flexible},  the model parameter space can be specified as 
\begin{itemize}[leftmargin=3.4em]
\setlength{\itemsep}{0pt}
	\item[(B1)] $\alpha_t\in \mathbb{R}^n$ for $t=1,\ldots, T$,
	\item[(B2)]  $Z\in  \mathbb{R}_1^{n \times k}$, 
where $	 \mathbb{R}_1^{n \times k}=\left\{Z \in \mathbb{R}^{n \times k}: 1_n^{\top} Z=0,\, {Z^{\top}Z =n\mathrm{I}_k} \right\} 
$, 
\item[(B3)] $\Lambda_t$  for $t=1,\ldots, T$ are  symmetric matrices, 
and at least one of $\Lambda_t$'s is full rank. 
\end{itemize}
When two sets of parameters $(\alpha,\Lambda, Z)$ and $(\tilde{\alpha},\tilde{\Lambda}, \tilde{Z})$ satisfying (B1)--(B3) 
generate the same model, i.e., $\alpha_{i t}+\alpha_{j t}+z_i^{\top} \Lambda_t z_j  = \tilde{\alpha}_{i t}+\tilde{\alpha}_{j t}+\tilde{z}_i^{\top} \tilde{\Lambda}_t \tilde{z}_j$,  Proposition 1 in  \cite{zhang2020flexible} showed that  $\alpha_t=\tilde{\alpha}_t,$ and there exists $Q\in \mathcal{O}(k) $   such that  $Z=\tilde{Z}Q$ and $\Lambda_t=Q^{\top}\tilde{\Lambda}_tQ$.

  Comparing (A2) and (B2) above, we can see that the  extra constraint $Z^{\top}Z=n\mathrm{I}_k$ needs to be taken into account when analyzing $Z$ under \eqref{eq:m2}, necessitating adjustments to both estimation and theory.

\subsection{{Estimation}}\label{sec:esthetero}
Under both models \eqref{model_DLSM} and \eqref{eq:m2},  
 similar estimation strategies can be applied, but adjustments are needed to accommodate the parameter space. Recall that for \eqref{model_DLSM}, we proposed to apply Algorithm \ref{algor:init} (based on projected gradient descent) and semiparametric one-step estimator \eqref{eq:ieffdef} in the main text.  
 For \eqref{eq:m2}, 
 we can also extend Algorithm 1 in 
 \cite{zhang2020flexible} (also based on projected gradient descent) under the Poisson distribution model. 
 Notably, after each iteration of gradient descent, \cite{zhang2020flexible} further transforms 
 the updated $Z$ estimates to be in $\mathbb{R}_1^{n\times k}$, whereas our proposed algorithm only requires $Z$ estimates to be in $\mathbb{R}_0^{n\times k}$.  Moreover, to establish the oracle estimation error rate of  $Z$ in theory, we can  
  generalize our semiparametric one-step update for $Z$ by treating  $\{\alpha_t,\Lambda_t: t=1,\ldots, T\}$ as nuisance parameters under \eqref{eq:m2}.  We next provide the details on the construction.  

\subsubsection{Semiparametric one-step estimator}
Under \eqref{eq:m2}, 
let $N$ denote all the nuisance parameters $\{\alpha_t,\Lambda_t: t=1,\ldots, T\}$, and then 
the log-likelihood function of $(Z,N)$ 
is $L(Z,N)=\sum_{t=1}^T L_t(Z,N)$, where 
\begin{align*}
    L_t(Z,N) =& ~\sum_{1\leqslant i\leqslant j\leqslant n} \left[A_{t,ij} (\alpha_{it}+\alpha_{jt}+z_i^\top \Lambda_t z_j) - \exp(\alpha_{it}+\alpha_{jt}+z_i^\top \Lambda_t z_j) \right].
\end{align*}
Similarly to Section  \ref{subsec:formula_efficient}, 
we define the vectorizations    $Z_v$ and $\alpha_v$ same as in \eqref{eq:vecdef}, and further define the vectorization of all nuisance parameters $N$ under \eqref{eq:m2} as 
\begin{align*}
  N_v  = \begin{pmatrix}
      \alpha_1\\
      \Lambda_{1,v}\\
      \vdots\\
      \alpha_T\\
      \Lambda_{T,v}
  \end{pmatrix}\in \mathbb R^{(n + k(k+1)/2)T \times 1}, \quad \quad \text{where}\quad \quad \Lambda_{t,v} = \begin{pmatrix}
      \Lambda_{t,11}/\sqrt{2}\\
      \vdots\\
      \Lambda_{t,kk}/\sqrt{2}\\
      \Lambda_{t,12}\\
      \vdots\\
      \Lambda_{t,k-1,k}
  \end{pmatrix}\in \mathbb{R}^{\frac{k(k+1)}{2}\times 1}, 
\end{align*}
where we point out that the diagonal entries in $\Lambda_t$ are rescaled by $1/\sqrt{2}$ for the convenience of derivation,  similarly to Remark \ref{rmk::vectorization_sqrt2}. 
With slight abuse of notation, we let $\dot L_Z(Z,N)$ and $\dot L_{N}(Z, N)$ denote the partial derivatives of $L(Z,N)$ with respect to vectors $Z_v$ and $N_v$, respectively. 
Then the efficient score and information matrix of $Z$ under \eqref{eq:m2} are 
\begin{align}
       \Seff(Z,N) =&~\dot{L}_Z - \mathbb{E}\big(\dot{L}_{Z}\, \dot{L}_{N}^{\top}\big) \big\{\mathbb{E}\big(\dot{L}_{N}\, \dot{L}_{N}^{\top}\big)\big\}^{-1} \dot{L}_{N} \ \in \  \mathbb{R}^{nk\times 1},\label{eq:seffdef_hetero}\\
    \Ieff(Z, N)=&~\mathbb{E}\big(\Seff(Z,N)\, \Seff^{\top}(Z,N) \big) \ \in \ \mathbb{R}^{nk\times nk},\label{eq:ieffdef_hetero}
\end{align} 
respectively, 
where for notational simplicity, $(Z,N)$ is omitted in $\dot{L}_Z$ and $\dot{L}_{N}$, and  $\mathbb{E}(\cdot)$ refers to the expectation  taken when data follows  \eqref{eq:m2} with parameters $(Z,N)$.
Similarly to Section \ref{subsec:formula_efficient}, 
we can 
express $ \Seff(Z,N)$ and $\Ieff(Z,N)$ as close-form functions of $(Z,N)$ and the observed data; see Sections \ref{sec:effscore_hetero} and  \ref{sec:effinfo_hetero} below. 
Given initial estimators $\check{Z}$ and $\check{N}$, e.g., similarly to Algorithm 1 in \cite{zhang2020flexible},  we can construct semiparametric one-step estimator as
\begin{align} \label{eq:newtonsolpseudo_hetero} 
 \hat{Z}_v = \check{Z}_v  + \big\{ \Ieff( \check{Z},  \check{N})\big\}^{+}\Seff( \check{Z},\check{N}). 
\end{align}

\subsubsection{Formula of the efficient score function} \label{sec:effscore_hetero} 
To derive the analytic expression of $\Seff(Z,N)$ as defined in \eqref{eq:seffdef_hetero}, we first need to determine the first-order derivatives $\dot L_Z$ and $\dot{L}_{N}$ as well as the expectations $\mathbb{E}\big(\dot{L}_{Z}\dot{L}_{N}^\top\big)$ and $ \mathbb{E}\big(\dot{L}_{N}\dot{L}_{N}^\top\big)$.

By the chain rule of the first-order derivatives, we have 
\begin{align*}
   \dot{L}_Z =  ~\frac{\partial L}{\partial Z_v} =\sum_{t=1}^T \frac{\partial  (\Theta_t)_v^{\top}}{\partial Z_v}\frac{\partial L}{\partial (\Theta_t)_v} \in \mathbb{R}^{nk\times 1}\quad \text{ and } \quad
\dot{L}_{N}=~
\begin{pmatrix}
	\dot{L}_{\alpha_1}\\
 \dot{L}_{\Lambda_{1}}\\
	\vdots\\
	\dot{L}_{\alpha_T}\\
 \dot{L}_{\Lambda_{T}}
\end{pmatrix}\in \mathbb{R}^{(n + k(k+1)/2)T\times 1}
\end{align*}
with
\begin{align*}
\dot{L}_{\alpha_t} =\frac{\partial L}{\partial \alpha_t} = \frac{\partial  (\Theta_t)_v^{\top}}{\partial \alpha_t}\frac{\partial L}{\partial (\Theta_t)_v}\in \mathbb{R}^{n\times 1} 
\  \text{ and }\ \dot{L}_{\Lambda_{t}} =\frac{\partial L}{\partial \Lambda_{t,v}} = \frac{\partial  (\Theta_t)_v^{\top}}{\partial \Lambda_{t,v}}\frac{\partial L}{\partial (\Theta_t)_v}\in \mathbb{R}^{k(k+1)/2\times 1}.
\end{align*}
To obtain formulae of $\dot{L}_Z$ and $\dot{L}_N$, it remains to derive  partial derivatives in the above expressions. 
Under \eqref{eq:m2}, the partial derivative of $(\Theta_t)_v$ with respect to $Z_v$ can be redefined as
\begin{align*}
\mathcal{D}_{Z \Theta_t} :=& \frac{\partial (\Theta_t)_v^{\top}}{\partial Z_v}=
		\begin{pmatrix}
	\frac{\partial (\Theta_{t,11}/\sqrt{2})}{\partial z_1} &  \cdots  & \frac{\partial (\Theta_{t,nn}/\sqrt{2})}{\partial z_1}  & \frac{\partial \Theta_{t,12}}{\partial z_1} & \cdots & \frac{\partial \Theta_{t,n-1,n}}{\partial z_1} \\
		\vdots &  \cdots  & \vdots  & \vdots & \cdots & \vdots \\
	\frac{\partial (\Theta_{t,11}/\sqrt{2})}{\partial z_n} & \cdots  & \frac{\partial (\Theta_{t,nn}/\sqrt{2})}{\partial z_n}  & \frac{\partial \Theta_{t,12}}{\partial z_n} & \cdots & \frac{\partial \Theta_{t,n-1,n}}{\partial z_n} 
	\end{pmatrix} \in \mathbb{R}^{nk\times \frac{n(n+1)}{2}} \notag\\[5pt]
	=&
\begin{array}{c@{}c}
\left(
  \begin{BMAT}[2pt]{c:c:c}{c}
    \begin{BMAT}[3pt]{ccc}{cccccc}
   \sqrt{2} \Lambda_t z_1   &\cdots  & 0\\
     0 &\cdots & 0\\
    0   &\cdots & 0 \\
    \vdots   & & \vdots  \\
    0  & & 0\\
     0 &\cdots &\sqrt{2} \Lambda_t z_n
    \end{BMAT}
    & 
       \begin{BMAT}[3.1pt]{cccc}{cccccc}
  \Lambda_t z_2   & \Lambda_t z_3 &\cdots  & \Lambda_t z_n\\
  \Lambda_t z_1 &0 & \cdots & 0\\
    0   & \Lambda_t z_1 & \cdots & 0 \\
    \vdots   & & &  \vdots  \\
    0  & 0 & \cdots &  0\\
     0 &0 & \cdots & \Lambda_t z_1
    \end{BMAT}    
    &
  \begin{BMAT}[3pt]{ccc}{cccccc}
\ 0   & \cdots  & 0\\
\ \Lambda_t z_3 & \cdots & 0\\
\  \Lambda_t z_2 &  \cdots & 0 \\
    \vdots   & & \vdots  \\
  \  0  & & \Lambda_t z_n\\
  \   0 &\cdots &\Lambda_t z_{n-1}
    \end{BMAT}
  \end{BMAT} 
\right)  \\[-1ex]
\ \hexbrace{3.7cm}{n}\hexbrace{3.7cm}{n-1}\hexbrace{3.3cm}{\frac{n(n+1)}{2}-(2n-1)}
\end{array} ,
\end{align*}
while the partial derivative of $(\Theta_t)_v$ with respect to $\alpha_t$ remains as in \eqref{eq:dthetaalphadefexpre} and the partial derivative of $L$ with respect to $(\Theta_t)_v$ remains as in \eqref{eq:ldot_theta_t}. Moreover, the partial derivative of $(\Theta_t)_v$ with respect to $\Lambda_{t,v}$ is
\begin{align*}
    \mathcal{D}_{\Lambda \Theta} :=& \frac{\partial (\Theta_t)_v^{\top}}{\partial \Lambda_{t,v}}=
		\begin{pmatrix}
	\frac{\partial (\Theta_{t,11}/\sqrt{2})}{\partial (\Lambda_{t,11}/\sqrt{2})} &  \cdots  & \frac{\partial (\Theta_{t,nn}/\sqrt{2})}{\partial (\Lambda_{t,11}/\sqrt{2})}  & \frac{\partial \Theta_{t,12}}{\partial (\Lambda_{t,11}/\sqrt{2})} & \cdots & \frac{\partial \Theta_{t,n-1,n}}{\partial (\Lambda_{t,11}/\sqrt{2})} \\
		\vdots &  \cdots  & \vdots  & \vdots & \cdots & \vdots \\
	\frac{\partial (\Theta_{t,11}/\sqrt{2})}{\partial (\Lambda_{t,nn}/\sqrt{2})} & \cdots  & \frac{\partial (\Theta_{t,nn}/\sqrt{2})}{\partial (\Lambda_{t,nn}/\sqrt{2})}  & \frac{\partial \Theta_{t,12}}{\partial (\Lambda_{t,nn}/\sqrt{2})} & \cdots & \frac{\partial \Theta_{t,n-1,n}}{\partial(\Lambda_{t,nn}/\sqrt{2})} \\[5pt]
    \frac{\partial (\Theta_{t,11}/\sqrt{2})}{\partial \Lambda_{t,12}} &  \cdots  & \frac{\partial (\Theta_{t,nn}/\sqrt{2})}{\partial \Lambda_{t,12}}  & \frac{\partial \Theta_{t,12}}{\partial \Lambda_{t,12}} & \cdots & \frac{\partial \Theta_{t,n-1,n}}{\partial \Lambda_{t,12}} \\
		\vdots &  \cdots  & \vdots  & \vdots & \cdots & \vdots \\
	\frac{\partial (\Theta_{t,11}/\sqrt{2})}{\partial \Lambda_{t,k-1,k}} & \cdots  & \frac{\partial (\Theta_{t,nn}/\sqrt{2})}{\partial \Lambda_{t,k-1,k}}  & \frac{\partial \Theta_{t,12}}{\partial \Lambda_{t,k-1,k}} & \cdots & \frac{\partial \Theta_{t,n-1,n}}{\partial \Lambda_{t,k-1,k}}
	\end{pmatrix} \in \mathbb R^{\frac{k(k+1)}{2} \times \frac{n(n+1)}{2}}\notag\\[5pt] 
 =& \begin{pmatrix}
     z_{11}^2  & \ldots & z_{n1}^2  &  \sqrt{2} z_{11}z_{21} & \ldots &  \sqrt{2}z_{n-1,1}z_{n1} \\
    \vdots & \ddots & \vdots &  \vdots & \ddots & \vdots\\
    z_{1k}^2 & \ldots & z_{nk}^2  & \sqrt{2} z_{1k}z_{2k} & \ldots &  \sqrt{2}z_{n-1,k}z_{nk} \\[3pt]
    \sqrt{2}z_{11}z_{12} & \ldots & \sqrt{2}z_{n1}z_{n2} & z_{11}z_{22} + z_{12} z_{21} & \ldots &  z_{n-1,1}z_{n2} + z_{n-1,2} z_{n1}\\
    \vdots & \ddots & \vdots &  \vdots & \ddots & \vdots\\
    \sqrt{2}z_{1,k-1}z_{1k} & \ldots & \sqrt{2}z_{n,k-1}z_{nk} & z_{1,k-1}z_{2k} + z_{1k} z_{2,k-1} & \ldots &  z_{n-1,k-1}z_{nk} + z_{n-1,k} z_{n,k-1}
 \end{pmatrix}.
\end{align*}
We note that $\mathcal{D}_{\Lambda \Theta}$ depends only on $Z$, and its formulae does not change across $t = 1, \ldots, T$. In contrast, $\mathcal{D}_{Z \Theta_t}$ is a function of both $Z$ and $\Lambda_t$ and may vary with respect to $t$. To facilitate the presentation, we further define $\mathcal{D}_{\Theta \Lambda} = \mathcal{D}_{\Lambda \Theta}^\mytrans$, $\mathcal{D}_{\Theta_t Z} = \mathcal{D}_{Z \Theta_t}^\mytrans$, $\mathcal{D}_{\Theta N} = \begin{pmatrix}
    \mathcal{D}_{\Theta \alpha}\,\mathcal{D}_{\Theta \Lambda}
\end{pmatrix}$, and $\mathcal{D}_{N \Theta } = \mathcal{D}_{\Theta N}^\mytrans$. In summary, we obtain 
\begin{align}
     \dot{L}_Z(Z,N) =& ~ \sum_{t=1}^T
\mathcal{D}_{Z\Theta_t}(Z,N) \dot{L}_{\Theta_t}(Z,N)\  \in \  \mathbb{R}^{nk\times 1}, \label{eq:seff_1}\\
\dot{L}_{N}(Z,N) =&~ 
\begin{pmatrix}
    \mathcal{D}_{N \Theta}(Z) \dot{L}_{\Theta_1}(Z,N) \\
    \vdots\\
    \mathcal{D}_{N \Theta}(Z) \dot{L}_{\Theta_T}(Z,N)
\end{pmatrix} \in \ \mathbb{R}^{(n + k(k+1)/2)T\times 1}. \label{eq:seff_2}
\end{align}

Following a similar analysis as in \eqref{eq:EdotL_alphadotL_Ztop}--\eqref{eq:EdotL_alphadotL_alphatop}, we have
\begin{align}
\mathbb{E}\big(\dot{L}_{Z}\dot{L}_{N}^\top\big) =&~ 
\begin{pmatrix}
 \mathcal{D}_{Z \Theta_1} \mathcal{D}_{\mu_1} \mathcal{D}_{\Theta N} & 
\cdots & 
  \mathcal{D}_{Z \Theta_T} \mathcal{D}_{\mu_T} \mathcal{D}_{\Theta N}
\end{pmatrix} 
\ \in \ \mathbb R^{nk \times (n + k(k+1)/2)T}, \label{eq:seff_3}\\[3pt]
\mathbb{E}\big(\dot{L}_{N}\dot{L}_{N}^\top\big) =&~ \text{diag}\big\{(\mathcal{D}_{N \Theta} \Dmut \mathcal{D}_{ \Theta N} )_{t=1, \ldots,T}\big\} 
 \ \in \ \mathbb R^{(n + k(k+1)/2)T \times (n + k(k+1)/2)T}. \label{eq:seff_4}
\end{align}
Then combining \eqref{eq:seffdef_hetero} with \eqref{eq:seff_1}--\eqref{eq:seff_4}, we have
\begin{align*}
     \Seff = \sum_{t=1}^T \left\{ \mathcal{D}_{Z \Theta_t} \dot{L}_{\Theta_t} -  \mathcal{D}_{Z \Theta_t} \Dmut \mathcal{D}_{\Theta N}(\mathcal{D}_{ N \Theta} \Dmut \mathcal{D}_{\Theta N})^{-1} \mathcal{D}_{ N \Theta} \dot{L}_{\Theta_t}\right\}.
\end{align*}

\subsubsection{Formula of the efficient information matrix}\label{sec:effinfo_hetero} 

By \eqref{eq:ieffdef_hetero}, we have
\begin{align*}
    \Ieff =&~\mathbb{E} \left[\dot{L}_Z - \mathbb{E}\big(\dot{L}_{Z}\, \dot{L}_{N}^{\top}\big) \big\{\mathbb{E}\big(\dot{L}_{N}\, \dot{L}_{N}^{\top}\big)\big\}^{-1} \dot{L}_{N}\right]\left[\dot{L}_Z - \mathbb{E}\big(\dot{L}_{Z}\, \dot{L}_{N}^{\top}\big) \big\{\mathbb{E}\big(\dot{L}_{N}\, \dot{L}_{N}^{\top}\big)\big\}^{-1} \dot{L}_{N}\right]^\top \\
    =&~ \mathbb{E}\big(\dot{L}_Z\dot{L}_{Z}^\top\big)
    - \mathbb{E}\big(\dot{L}_Z\dot{L}_{N}^\top\big)\big\{\mathbb{E}\big(\dot{L}_{N}\, \dot{L}_{N}^{\top}\big)\big\}^{-1}\mathbb{E}\big(\dot{L}_N\dot{L}_{Z}^\top\big).
\end{align*}
Combining with \eqref{eq:seff_3}, \eqref{eq:seff_4}, and
\begin{align*}
\mathbb{E}\big(\dot{L}_Z\dot{L}_{Z}^\top\big) =&~ \sum_{1 \leqslant t_1,t_2 \leqslant T}\mathbb E\bigg( \frac{\partial \Theta_{t_1}}{\partial Z}^\top \dot{L}_{\Theta_{t_1}}\dot{L}_{\Theta_{t_2}}^\top \frac{\partial \Theta_{t_2}}{\partial Z} \bigg) = \sum_{t = 1}^T \mathcal{D}_{Z \Theta_t}  \Dmut \mathcal{D}_{\Theta_t Z} \in\mathbb R^{nk \times nk},
\end{align*}
we can write
\begin{align*}
    \Ieff = \sum_{t=1}^T \Big\{ \mathcal{D}_{Z \Theta_t}\Dmut \mathcal{D}_{\Theta_t Z} - \mathcal{D}_{Z \Theta_t}\Dmut\mathcal{D}_{\Theta N} (\mathcal{D}_{N \Theta} \Dmut \mathcal{D}_{\Theta N})^{-1} \mathcal{D}_{N \Theta} \Dmut \mathcal{D}_{\Theta_t Z}   \Big\}.
\end{align*}

\subsection{{Theory}}\label{sec:theoryhetero}

In this section, we show the feasibility of generalizing  Theorem \ref{thm_NRerr} for \eqref{eq:newtonsolpseudo_hetero}  under the model \eqref{eq:m2}. 
We will next focus on discussing the key steps in Section \ref{sec:pfthmnerr},  whereas a rigorous proof requires a significant amount of work and is left for future study; see more discussions in Remark \ref{rm:future_hetero} below.

\medskip 

 \textit{Preliminary:} Section \ref{sec:pfthmnerr} utilizes 
a restricted eigenspace 
 in Definition \ref{def::UZA}. 
 We next adopt the same definition and establish  Lemma \ref{lem_presentingUZ_hetero} in Section \ref{sec:lem_presentingUZ_hetero} on properties of the efficient information matrix similar to Lemma \ref{lem_presentingUZmain}.  
Moreover, we assume that the initial estimator $(\check{Z},\check{N})$ used in \eqref{eq:newtonsolpseudo_hetero}  can achieve error rates similar to those in Condition \ref{cond_elem_init}, which is expected to be mild as discussed in Section \ref{sec:formulaofupdateandtheory} of the main text.
 




\smallskip 
 \textit{Step 1:} 
 Let $\{\check{\mathcal{B}}_{C,j}: j=1,2,3\}$ denote basis vectors constructed following Definition \ref{def::UZA} at initial estimator $\check{Z}$. 
 Following Step 1 in Section \ref{sec:pfthmnerr}, we have $\hat Z_v -  (\ZQs)_v  \in \mathrm{span}(\cup_{j=1}^3\check{\mathcal{B}}_{C,j})$, as the proof only depends on the definition of $\check{Q}$, construction of $\check{\mathcal{B}}_{C,j}$, and Lemma \ref{lem_presentingUZ_hetero}-(iii).  
Therefore, similarly to \eqref{pfthm1-3A}, we have 
\begin{align*}
    \| \hat Z_v -  (\ZQs)_v \|_2^2= \sum_{b\, \in \,  \cup_{j=1}^3\check{\mathcal{B}}_{C,j}}| b^{\top} ( \hat Z_v -  (\ZQs)_v) |_2^2. 
\end{align*}

 \textit{Step 2:} 
Define $\bar{\mathcal{U}}_{\check{Z}}$ as a matrix whose columns are formed by vectors in $\check{\mathcal{B}}_{C,3}$. 
Define $\check{I}_{e,\bar{U}}=\bar{\mathcal{U}}_{\check{Z}}^\top \Ieff(\check{Z}, \check{N})\bar{\mathcal{U}}_{\check{Z}} /(nT)$ as the  $\bar{\mathcal{U}}_{\check{Z}}$-restricted efficient information matrix. 
By Lemma \ref{lem_presentingUZ_hetero}, we have
\begin{align*}
    \sigma_{\min}(\check{I}_{e,\bar{U}}) \geqslant \check{M}_{\Lambda,2} \exp(-\check{M}_{\Theta,2}),
\end{align*} 
where we define $\check{M}_{\Lambda,2} = \min_{1 \leqslant t \leqslant T} \sigma_{\operatorname{min}}(\check \Lambda_t^2)$ and $\check{M}_{\Theta,2} = \max_{1 \leqslant i,j \leqslant n, 1 \leqslant t \leqslant T} \{| \check \alpha_{it} + \check \alpha_{jt} + \check z_i^\mytrans \check{\Lambda}_t \check{z}_j |\}$. 
Similarly to  Step 2 in Section \ref{sec:pfthmnerr}, we \textit{expect} that $ \sigma_{\min}(\check{I}_{e,\bar{U}})$ has a constant upper bound when the true parameters $(Z^{\star},N^{\star})$ satisfy regularity conditions similar to Condition \ref{cond:truvalueregularity}, and initial estimators $(\check{Z},\check{N})$ are consistent for estimating $(Z^{\star},N^{\star})$. 

 \textit{Step 3:}
By Lemma \ref{lem_presentingUZ_hetero},  
\begin{align*}
    \sum_{b\in \check{\mathcal{B}}_{C,3}}| b^{\top} ( \hat Z_v -  (\ZQs)_v) |_2^2 = &~ \big\| \check{I}_{e,\bar{U}}^{-1}\ \check{I}_{e,\bar{U}} \ \bar{\mathcal{U}}_{\check{Z}}^{\top} \, (\hat Z_v -  (\ZQs)_v)\big\|_2^2\\
    \leqslant &~ \frac{1}{\sigma_{\min}^2(\check{I}_{e,\bar{U}})}\big\|\check{I}_{e,\bar{U}} \ \bar{\mathcal{U}}_{\check{Z}}^{\top} \, (\hat Z_v -  (\ZQs)_v)\big\|_2^2.
\end{align*} 
For the numerator above, by plugging the definition of $\check{I}_{e,\bar{U}}$ and the construction of $\hat Z_v$ in \eqref{eq:newtonsolpseudo_hetero}, we obtain 
\begin{align*}
    \big\|\check{I}_{e,\bar{U}} \ \bar{\mathcal{U}}_{\check{Z}}^{\top} \, (\hat Z_v -  (\ZQs)_v)\big\|_2 = &~\big\| \bar{\mathcal{U}}_{\check{Z}}^{\top}\big\{\Ieff(\check Z, \check N) (\check Z_v-  (\ZQs)_v ) +   \Seff(\check{Z},\check{N})\big\}\big\|_2 /(nT)\notag\\   
    \leqslant &~\big\| \Ieff(\check Z, \check N) (\check Z_v-  (\ZQs)_v ) +   \Seff(\check{Z},\check{N})\big\|_2/(nT), 
\end{align*}
which takes a form similar to \eqref{eq:multipliediff2A}. We \textit{expect} that the analysis in \eqref{eq:thm1s1termA}--\eqref{pfthm1-12A} can be similarly extended, which would give that with high probability, the numerator $\big\|\check{I}_{e,U}\,  \bar{\mathcal{U}}_{\check{Z}}^{\top} ( \hat Z_v -  (\ZQs)_v)\big\|_2$ is of the order of 
$\max\{T^{-1/2}, n^{-1/2}\}$ up to logarithmic factors.  


 \textit{Step 4:}
In Section \ref{sec:pfbc12errorbd}, we prove that there exists a constant $C>0$ such that 
\begin{align}\label{eq:bc12errorbd}
     \sum_{b\in \check{\mathcal{B}}_{C,1}\cup \check{\mathcal{B}}_{C,2}}| b^{\top} ( \hat Z_v -  (\ZQs)_v) |_2^2 \leqslant \frac{C}{n}  \| \check Z_v -  (\ZQs)_v \|_2^4. 
\end{align}
When the initial estimator $\check{Z}$ satisfies Condition \ref{cond_elem_init}, we know that $\eqref{eq:bc12errorbd} $ is of the order of $n^{-1}$ up to logarithmic factors. 

 \textit{Step 5:} Combining Steps 1--4, we could obtain  that with high probability, 
 $  \| \hat Z_v -  (\ZQs)_v \|_2^2$ is of the order of  
$\max\{T^{-1}, n^{-1}\}$ up to logarithmic factors, which gives the same conclusion as in Theorem \ref{thm_NRerr}.  


\begin{remark}[Connections and differences with Section \ref{sec:pfthmnerr}]
\label{rm:conndisc}
Intuitively, the analysis for \eqref{eq:newtonsolpseudo_hetero}  can be similar to  Section \ref{sec:pfthmnerr} because  $Z$ is unidentifiable up to an orthogonal transformation under both \eqref{model_DLSM} and \eqref{eq:m2} by Section \ref{sec:modelhetero}.
The unidentifiability of $Z$ causes the singularity of the efficient matrix $I_{eff}(Z,N)$ in \eqref{eq:ieffdef_hetero} similarly to Remark \ref{lb:singularity}. 
Thus the theoretical techniques under \eqref{model_DLSM} would still be useful under \eqref{eq:m2}. 
However, given the constraint $Z^{\top}Z=n\mathrm{I}_k$ under \eqref{eq:m2}, 
the number of free parameters in $Z$ is  reduced compared to that under \eqref{model_DLSM}, and the column space of  $I_{eff}(Z,N)$ under \eqref{eq:m2} becomes smaller. As a result,  theoretical analysis would not be exactly the same as that under \eqref{model_DLSM}.

Technically, the reduced dimension of the column space of $I_{eff}$ is substantiated by comparing Lemma \ref{lem_presentingUZmain}-(i) and Lemma \ref{lem_presentingUZ_hetero}-(i), showing that the column space of $I_{eff}$ under \eqref{eq:m2} is a subset of that under \eqref{model_DLSM}. 
To address this new challenge, 
we propose to consider projection of $ \hat Z_v -  (\ZQs)_v$ onto only $\check{\mathcal{B}}_{C,3}$ in Step 3 so that an upper bound similar to \eqref{eq:multipliediff2A} can be obtained.  
Meanwhile, we develop new Step 4 to control the projection of  of $ \hat Z_v -  (\ZQs)_v$ onto $\check{\mathcal{B}}_{C,1}\cup \check{\mathcal{B}}_{C,2}$.  

\end{remark}

\begin{remark}[Future work]\label{rm:future_hetero}
The discussion above aims 
to demonstrate that the developed analysis for the model \eqref{model_DLSM} can also pave the way for analyzing other more complex models such as  \eqref{eq:m2}. 
We acknowledge that
a substantial amount of work is still needed to obtain rigorous results. 
Specifically, it remains to establish sufficiently good rates for   
 the initial estimator $(\check{Z},\check{N})$,  
and the error rates anticipated in Step 3. 
We expect that those conclusions could be  developed similarly to Sections \ref{sec:pfthm1} and \ref{sec:techinitial} and leave the details for future work. 
\end{remark}


\subsubsection{Lemma \ref{lem_presentingUZ_hetero} and its proof} \label{sec:lem_presentingUZ_hetero}
Given $Z\in \mathbb{R}^{n\times k}$, we define $\beta_{ij}$,  $\{\mathcal{B}_{C,j}: j=1,2,3\}$, and $\{\mathcal{B}_{N,j}: j=1,2\}$ in the same way as in \eqref{eq:betaijdefA}, Definition \ref{def::UZA}, and \eqref{eq:bn1bn2setsupp}, respectively. Let $\bar{\mathcal{U}}_{Z}$ be the matrix whose columns are formed by vectors in ${\mathcal{B}}_{C,3}$. 

\begin{lemma}\label{lem_presentingUZ_hetero}
   For any parameter $(\alpha, Z, \Lambda)$ satisfying the identifiability constraints (B1)--(B3), define $\Ieff(Z,N)$ as in \eqref{eq:ieffdef_hetero}. Then we have 
   \begin{enumerate}
       \item[(i)] $\operatorname{col}(\Ieff(Z,N)) = \operatorname{col}(\bar{\mathcal{U}}_Z) = \operatorname{span}(\mathcal{B}_{C,3})$, $\operatorname{rank}(\Ieff(Z,N)) = nk - k(k+1)$, and 
       \begin{align*}
\sigma_{\operatorname{min}}\big(\bar{\mathcal{U}}_Z^\mytrans \Ieff(Z,N) \bar{\mathcal{U}}_Z\big)  \geqslant nT \min_{1 \leqslant t \leqslant T}\Big\{\sigma_{\operatorname{min}}(\Lambda_t^2) \min_{1 \leqslant i , j \leqslant n}\exp(\alpha_{it} + \alpha_{jt} + z_i^\mytrans \Lambda_t z_j) \Big\}.
       \end{align*}
        \item[(ii)] $\operatorname{null}(\Ieff(Z,N)) = \operatorname{col}(\bar{\mathcal{U}}_Z)^{\bot} = \operatorname{span}(\mathcal{B}_{C,1} \cup \mathcal{B}_{C,2} \cup \mathcal{B}_{N,1} \cup \mathcal{B}_{N,2})$,  where $\operatorname{col}(\bar{\mathcal{U}}_Z)^{\bot}$ denotes the orthogonal complement of $\operatorname{col}(\bar{\mathcal{U}}_Z)$.
        \item[(iii)] $\bar{\mathcal{U}}_Z$ satisfies $\bar{\mathcal{U}}_Z\bar{\mathcal{U}}_Z^\top \mathcal{D}_{Z \Theta_t} \myXt = \mathcal{D}_{Z \Theta_t} \myXt$ with $\myXt=  \Dmut- \Dmut \mathcal{D}_{\Theta N} (\mathcal{D}_{N \Theta} \Dmut \mathcal{D}_{\Theta N})^{-1} \mathcal{D}_{N \Theta } \Dmut$,  
     and thus  
    $$\bar{\mathcal{U}}_Z\bar{\mathcal{U}}_Z^\top \Seff(Z,N) =  \Seff(Z,N), \quad \bar{\mathcal{U}}_Z \bar{\mathcal{U}}_Z^\top \Ieff(Z,N) =  \Ieff(Z,N),$$i.e., $\Seff$ lies in the eigenspace corresponding to the non-zero eigenvalues of $\Ieff$.
   \end{enumerate} 
   
\end{lemma}

\smallskip
\begin{proof}
To prove (i) and (ii), 
it suffices to show that 
(a) for any $v\in \text{col}(\bar{\mathcal{U}}_Z)^\perp$, $\Ieff v =0$, 
and (b) for any $v \in$ col$(\bar{\mathcal{U}}_Z)$, $v^{\top}\Ieff v$ has a positive lower bound given as in (i).

Part (a): By Definition \ref{def::UZA} in the main text and Lemma \ref{lem_DZtDz_eigenstructure}, we have
\begin{align*}
    \text{col}(\bar{\mathcal{U}}_Z)^\perp = \text{span}(\mathcal{B}_{N,1}) \oplus \text{span}(\mathcal{B}_{C,1} \cup \mathcal{B}_{C,2} \cup \mathcal{B}_{N,2}) = \text{col}(\mathcal{D}_{I}) \oplus \text{col}(\mathcal{D}_{ZI}),
\end{align*} where $\mathcal{D}_I = 1_n \otimes \mathrm I_k \in \mathbb R^{nk \times k}$ and $\mathcal{D}_{ZI} = Z \otimes \mathrm{I}_k \in \mathbb R^{nk \times k^2}$. To prove $\text{col}(\bar{\mathcal{U}}_Z)^\perp \subseteq \text{null}(\Ieff)$, 
it suffices to prove 
$\Ieff \myDtri = 0$ and $\Ieff \mathcal{D}_{ZI} = 0$. 
In particular, we have
\begin{align*}
    \mathcal{D}_{\Theta_t Z} \mathcal{D}_{I} = \begin{pmatrix}
        \sqrt{2}  z_1^\mytrans \Lambda_t \\
        \vdots \\
        \sqrt{2}  z_n^\mytrans \Lambda_t \\
         z_2^\mytrans \Lambda_t + z_1^\mytrans \Lambda_t \\
         \vdots \\
         z_n^\mytrans \Lambda_t + z_{n-1}^\mytrans \Lambda_t \\
    \end{pmatrix} \hspace{2.8cm} \text{and} \hspace{3cm}
\end{align*} 
\begin{align*}
    \mathcal{D}_{\Theta_t Z} \mathcal{D}_{ZI} = \begin{pmatrix}
        \sqrt{2} z_{11} z_1^\mytrans \Lambda_t & \ldots & \sqrt{2} z_{1k} z_1^\mytrans \Lambda_t \\
        \vdots & \ddots & \vdots \\
        \sqrt{2} z_{n1} z_n^\mytrans \Lambda_t & \ldots & \sqrt{2} z_{nk} z_n^\mytrans \Lambda_t \\
        z_{11}z_2^\mytrans\Lambda_t + z_{21} z_1^\mytrans \Lambda_t & \ldots & z_{1k}z_2^\mytrans \Lambda_t + z_{2k} z_1^\mytrans \Lambda_t \\
        \vdots & \ddots & \vdots \\
        z_{n-1,1}z_n^\mytrans\Lambda_t + z_{n1} z_{n-1}^\mytrans \Lambda_t & \ldots & z_{n-1,1}z_n^\mytrans\Lambda_t + z_{n1} z_{n-1}^\mytrans \Lambda_t \\
    \end{pmatrix},
\end{align*}
which imply $\operatorname{col}(\mathcal{D}_{\Theta_t Z} \mathcal{D}_{I}) \subseteq \operatorname{col}(\mathcal{D}_{\Theta \alpha}) \subseteq \operatorname{col}(\mathcal{D}_{\Theta N}) $ and $\operatorname{col}(\mathcal{D}_{\Theta_t Z} \mathcal{D}_{ZI}) \subseteq \operatorname{col}(\mathcal{D}_{\Theta \Lambda}) \subseteq \operatorname{col}(\mathcal{D}_{\Theta N})$ by the formulae in Section \ref{sec:effscore_hetero}. Combining with $\Ieff = \sum_{t=1}^T \mathcal{D}_{Z \Theta_t} \myXt \mathcal{D}_{\Theta_t Z}$ and $\operatorname{null}(\myXt) = \operatorname{col}(\mathcal{D}_{\Theta N})$ by an argument similar to Lemma \ref{lem:property_myXt}-2, we obtain $\Ieff \myDtri = 0$ and $\Ieff \mathcal{D}_{ZI} = 0$.

Part (b): 
For any $v \in$ col$(\bar{\mathcal{U}}_Z)$, 
we first show that  
\begin{align}\label{eq:ieffquadineq_heter}
     v^{\top} \Ieff v =\sum_{t=1}^T v^\top \mathcal{D}_{Z \Theta_t} \myXt \mathcal{D}_{\Theta_t Z} v
    \geqslant T\min_{t,i,j}(\mu_{t,ij}) \|\mathcal{D}_{\Theta_t Z} v\|_2^2.
\end{align}
To prove the inequality in \eqref{eq:ieffquadineq_heter}, 
it suffices to show that $\mathcal{D}_{\Theta N}^{\top}\mathcal{D}_{\Theta_t Z} v =0$ by $\operatorname{null}(\myXt) = \operatorname{col}(\mathcal{D}_{\Theta N})$. 
In particular, we have
\begin{align*}
    \mathcal{D}_{\alpha \Theta} \mathcal{D}_{\Theta_t Z} =  \begin{pmatrix}
         z_1^\mytrans \Lambda_t & \ldots &  z_1^\mytrans \Lambda_t \\
        \vdots & \ddots & \vdots \\
         z_n^\mytrans \Lambda_t & \ldots &  z_n^\mytrans \Lambda_t
    \end{pmatrix} \hspace{2.8cm} \text{and} \hspace{3cm}
\end{align*}
\begin{align*}
    \mathcal{D}_{\Lambda \Theta} \mathcal{D}_{\Theta_t Z} =  \begin{pmatrix}
        \sqrt{2} z_{11} e_1^\mytrans \Lambda_t & \ldots & \sqrt{2} z_{n1} e_1^\mytrans \Lambda_t \\
        \vdots & \ddots & \vdots \\
        \sqrt{2} z_{1k} e_k^\mytrans \Lambda_t & \ldots & \sqrt{2} z_{nk} e_k^\mytrans \Lambda_t \\[3pt]
        z_{11} e_2^\mytrans \Lambda_t + z_{12} e_1^\mytrans \Lambda_t & \ldots & z_{n1} e_2^\mytrans \Lambda_t + z_{n2} e_1^\mytrans \Lambda_t \\
        \vdots & \ddots & \vdots \\
        z_{1,k-1} e_k^\mytrans \Lambda_t + z_{1k} e_{k-1}^\mytrans \Lambda_t & \ldots & z_{n,k-1} e_k^\mytrans \Lambda_t + z_{nk} e_{k-1}^\mytrans \Lambda_t \\
    \end{pmatrix},
\end{align*}
then for any $v = (v_{1}^\mytrans, \ldots, v_{n}^\mytrans)^\mytrans \in \operatorname{col}(\bar{\mathcal{U}}_Z)$, we have $\mathcal{D}_{\alpha \Theta} \mathcal{D}_{Z \Theta_t} v = 0$ and $\mathcal{D}_{\Lambda \Theta} \mathcal{D}_{Z \Theta_t} v = 0$ since 
\begin{align*} \sum_{l = 1 }^n  v_{l} = 0 \quad\text{and}\quad
    \sum_{l = 1 }^n z_{li} v_{l} = 0 \quad\text{ for any }1 \leqslant i\leqslant k.
\end{align*}

We next derive a lower bound for $\|\mathcal{D}_{\Theta_t Z} v\|_2^2$ in \eqref{eq:ieffquadineq_heter}. 
For any $v\in \operatorname{col}(\bar{\mathcal{U}}_Z)$, 
\begin{align}\label{eq:dzvquadlowbd_heter}
    \|\mathcal{D}_{\Theta_t Z} v\|_2^2 = v^\top \mathcal{D}_{ Z \Theta_t} \mathcal{D}_{\Theta_t Z} v \geqslant  \sigma_{\min}(\Lambda_t^2) \|v\|_2^2, 
\end{align}
where the inequality follows by an argument similar to Lemma \ref{lem_DZtDz_eigenstructure}. 
 In summary, combining   \eqref{eq:ieffquadineq_heter} and  \eqref{eq:dzvquadlowbd_heter}, 
  $$\sigma_{\min}(\bar{\mathcal{U}}_Z^{\top} \Ieff(Z,N) \bar{\mathcal{U}}_Z) \geqslant nT\, \min_{1 \leqslant t \leqslant T}\Big\{\sigma_{\operatorname{min}}(\Lambda_t^2)\min_{{1\leqslant i,j\leqslant n}}\exp(z_i^\top \Lambda_t  z_j +\alpha_{it} +  \alpha_{jt})\Big\}.$$

  To prove (iii), it suffices to show that $\myXt \mathcal{D}_{\Theta_t Z} v = 0$ for any $v \in \operatorname{col}(\bar{\mathcal{U}}_Z)^\perp$. This can be demonstrated using a similar technique as in part (a).
\end{proof}

\subsubsection{Proof of \eqref{eq:bc12errorbd}} \label{sec:pfbc12errorbd}

By the construction of $\hat{Z}_v$ in \eqref{eq:newtonsolpseudo_hetero}, we know 
\begin{align}
   b^{\top} ( \hat Z_v -  (\ZQs)_v)=&~b^{\top}\Big[  \check Z_v-  (\ZQs)_v  + \big\{ \Ieff( \check{Z},  \check{N})\big\}^{+}\Seff( \check{Z},\check{N})\Big] 
   = b^{\top} ( \check Z_v -  (\ZQs)_v) \notag 
\end{align}
where the second equation follows since $\{ \Ieff( \check{Z},  \check{N})\}^{+}\Seff( \check{Z},\check{N}) \in \check{\mathcal{B}}_{C,3} $ by Lemma \ref{lem_presentingUZ_hetero}-(i)
and $ b \in \check{\mathcal{B}}_{C,1}\cup \check{\mathcal{B}}_{C,2} $. Then
\begin{align} \label{eq:bc12errorbd_hetero1}
    \sum_{b\in \check{\mathcal{B}}_{C,1}\cup \check{\mathcal{B}}_{C,2}}| b^{\top} ( \hat Z_v -  (\ZQs)_v) |_2^2  =  \sum_{b\in \check{\mathcal{B}}_{C,1}\cup \check{\mathcal{B}}_{C,2}}| b^{\top} ( \check Z_v -  (\ZQs)_v) |_2^2.  
\end{align}

As $\check{Z}^{\top}\check{Z} = n \mathrm{I}_{k}$,  
we can write the singular value decomposition  $\check{Z}=\sum_{i=1}^n \sqrt{n}\check{U}_i e_i^{\top} $, where $\check{U}_i\in \mathbb{R}^{n\times 1}$ denotes the $i$-th left singular vector of $\check{Z}$, and $e_i\in \mathbb{R}^{k\times 1}$ is an indicator vector for the $i$-th entry. 
By Definition \ref{def::UZA}, $b\in \check{\mathcal{B}}_{C,1}\cup \check{\mathcal{B}}_{C,2}$ is linear combinations of 
$\check{\beta}_{ij}:=\check{U}_i\otimes e_j$ for $1\leqslant i,j\leqslant k$. 
Therefore, there exists a constant $C>0$ such that 
\begin{align}
\eqref{eq:bc12errorbd_hetero1}  \leqslant &~ C \sum_{1\leqslant i,j\leqslant k} \big| (\check{U}_i\otimes e_j)^{\top} ( \check Z_v -  (\ZQs)_v) \big|_2^2\notag \\
  =&~  C \sum_{1\leqslant i,j\leqslant k} \big| \check{U}_i ^{\top} ( \check Z  -  Z^{\star}\check{Q}) e_j\big|_2^2   
  =  \frac{C}{n}  \big\| \check{Z}^{\top}(  \check Z  -  Z^{\star}\check{Q} )\big\|_{\mathrm{F}}^2.\label{eq:bc12errorbd_hetero2}  
\end{align}

By the property of orthogonal procrustes problem, 
\begin{align}\label{eq:qcheckq_hetero}
    \check{Q}= \arg\min_{Q\in \mathcal{O}(k)} \| \check{Z}-Z^* \check{Q} \|_F=  Z^{\star\top}\check{Z}  (\check{Z}^{\top} Z^{\star}Z^{\star\top}\check{Z})^{-1/2}. 
\end{align}
Then 
\begin{align}\label{eq:projnorm_hetero}
    \big\| \check{Z}^{\top}(  \check Z  -  Z^{\star}\check{Q} )\big\|_{\mathrm{F}}^2 = \|n\mathrm{I}_k - (\check{Z}^{\top} Z^{\star}Z^{\star\top}\check{Z})^{1/2} \|_F^2= n^2 \sum_{i=1}^k (1-\sigma_i)^2, 
\end{align}
where $\sigma_i$ represents the $i$-th singular value of $\check{Z}^{\top} Z^{\star}/n$. 
Moreover, by $\check{Z}^{\top}\check{Z} = Z^{\star \top} Z^{\star}= n \mathrm{I}_{k}$ and \eqref{eq:qcheckq_hetero}, 
\begin{align}
     \big\|    \check Z  -  Z^{\star}\check{Q}  \big\|_{\mathrm{F}}^2 =   2 n k - 2\mathrm{tr}( (\check{Z}^{\top} Z^{\star}Z^{\star\top}\check{Z})^{1/2})   =&~2n\sum_{i=1}^k (1-\sigma_i).  \label{eq:fnormdelta_hetero}
\end{align}
By  $0 \leqslant \sigma_i\leqslant 1$, $\sum_{i=1}^k (1-\sigma_i)^2 \leqslant\{\sum_{i=1}^k (1-\sigma_i)\}^2  $. 
By \eqref{eq:projnorm_hetero} and \eqref{eq:fnormdelta_hetero}, 
$  \big\| \check{Z}^{\top}(  \check Z  -  Z^{\star}\check{Q} )\big\|_{\mathrm{F}}^2 \leqslant   \big\|    \check Z  -  Z^{\star}\check{Q}  \big\|_{\mathrm{F}}^4 $. 
Combining with \eqref{eq:bc12errorbd_hetero2}, \eqref{eq:bc12errorbd} is proved.



\bigskip 

\subsection{{Data Analysis}}\label{sec:datahetero}

\begin{figure}[!htbp]
	\centering
    \includegraphics[width=0.55\linewidth]{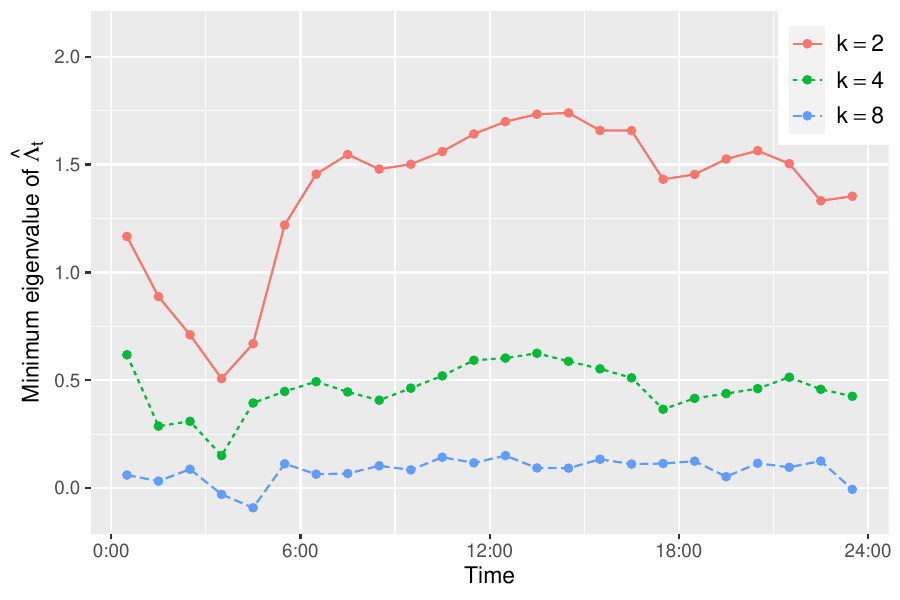}
\caption{Minimum eigenvalue of estimated $\hat{\Lambda}_t$ under \eqref{eq:m2} versus $t$.} \label{fig:lambdatvst}
\end{figure}

\begin{figure}[!htbp]
	\centering
\begin{subfigure}{0.22\textwidth}
		\centering
		  \caption{$k = 2$}
	\includegraphics[width=1\linewidth]{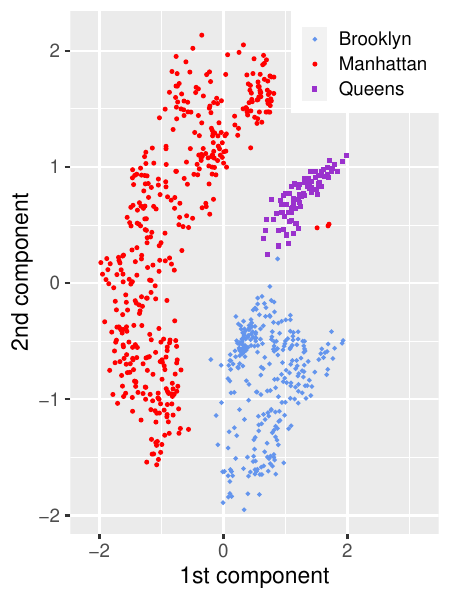}
\end{subfigure}
 \begin{subfigure}{0.22\textwidth}
		\centering
  \caption{$k = 4$}
		\includegraphics[width=1\linewidth]{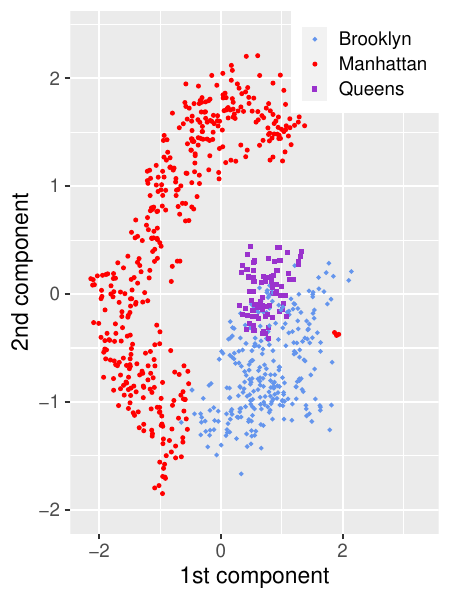}
\end{subfigure}
\begin{subfigure}{0.22\textwidth}
  \caption{$k = 8$}
		\centering
		\includegraphics[width=1\linewidth]{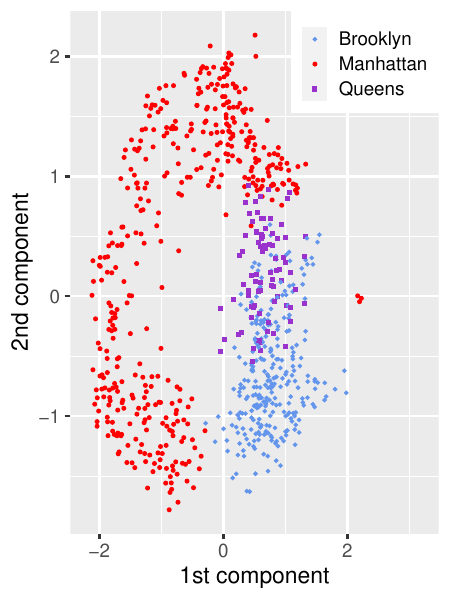}
\end{subfigure}
\begin{subfigure}{0.3\textwidth}
  \caption{True Locations}
		\centering
		\includegraphics[width=0.98\linewidth]{Bike_map.pdf}
\end{subfigure}
 \caption{Estimated top two components of latent vectors under the heterophilic model \eqref{eq:m2} when $k\in \{2,4,8\}$ and true locations of bike stations.}   \label{fig:latentzm2}
\end{figure}

\begin{figure}[!htbp]
	\centering
\begin{subfigure}{0.22\textwidth}
		\centering
		  \caption{$k = 2$}
	\includegraphics[width=1\linewidth]{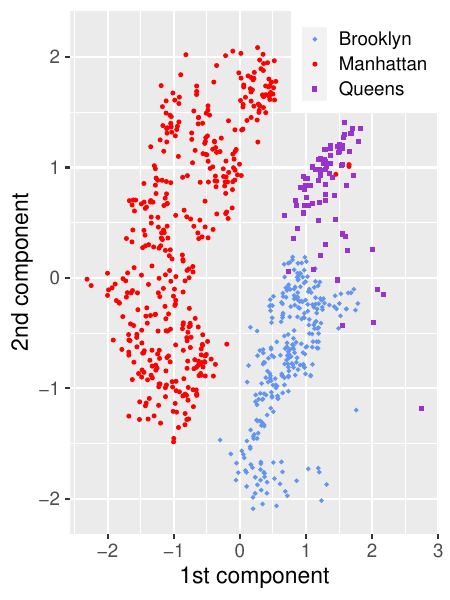}
\end{subfigure}
 \begin{subfigure}{0.22\textwidth}
		\centering
  \caption{$k = 4$}
		\includegraphics[width=1\linewidth]{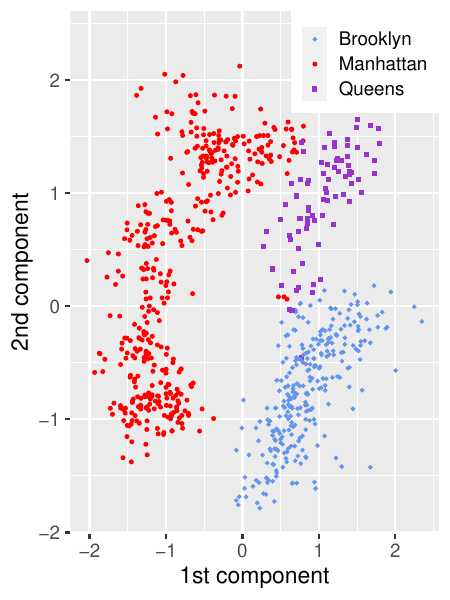}
\end{subfigure}
\begin{subfigure}{0.22\textwidth}
  \caption{$k = 8$}
		\centering
		\includegraphics[width=1\linewidth]{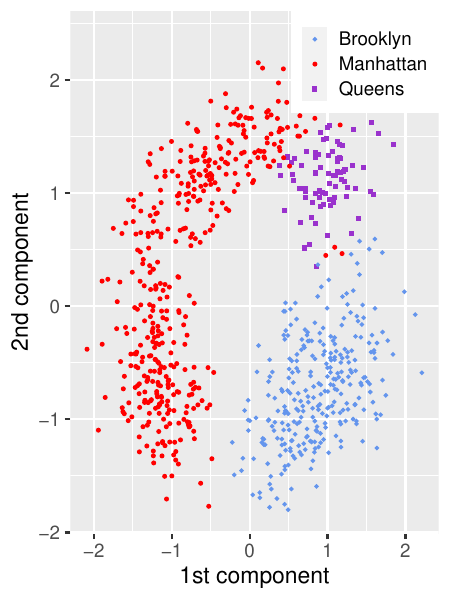}
\end{subfigure}
\begin{subfigure}{0.3\textwidth}
  \hspace{10em} 
\end{subfigure}
 \caption{Estimated top two components of  latent vectors under the homophilic model \eqref{model_DLSM} when $k\in \{2,4,8\}$.}   \label{fig:latentzm2_homo}
\end{figure}

As a comparison, we fit the heterophilic model \eqref{eq:m2} to the NYC bike data in Section \ref{sec::bikedata} and compare with the estimation results under the homophilic model \eqref{model_DLSM}.
Specifically, Figure \ref{fig:lambdatvst} presents the minimum eigenvalue of estimated $\hat{\Lambda}_t$  versus $t$ when $k \in \{2, 4, 8\}$ under \eqref{eq:m2}.  
The results show that most $\hat{\Lambda}_t$'s have  positive minimum eigenvalues, suggesting they are positive definite matrices. 
When $k=8$, only a few $\hat{\Lambda}_t$'s display negative eigenvalues, but their magnitudes are close to zero.  
In addition, we compare the estimated latent vectors under models \eqref{model_DLSM} and \eqref{eq:m2} when $k\in \{2, 4, 8\}$. 
Figure \ref{fig:latentzm2} presents the estimated latent vectors corresponding to the leading two diagonal components in $\hat{\Lambda}_t$ under the model \eqref{eq:m2}, and Figure \ref{fig:latentzm2_homo} similarly presents the leading two principal components of estimated $Z$ under \eqref{model_DLSM}. 
Although the estimated latent vectors under the two models are not identical, they reveal  
similar underlying patterns, particularly in forming the clustering of stations according to  three boroughs of New York City.


In summary, the studied data does not exhibit significant evidence for heterophilic structure in \eqref{eq:m2} when compared to \eqref{model_DLSM}. 
In future studies, it would be of interest to 
develop rigorous quantitative measures such as goodness-of-fit test to  assess the model fit or develop results under other more complex models.

\clearpage



\end{document}